**The life and work of Major Percy Alexander MacMahon**
**PhD Thesis by Dr Paul Garcia**

History of Mathematics Thesis submitted to the Open University in fulfilment of the Degree of Doctor of Philosophy

# The life and work of Major Percy Alexander MacMahon

Paul Garcia
B.A. (Hons), MSc (Mathematics), MSc (Information and Knowledge Technology)

May 2006

# The life and work of Major Percy Alexander MacMahon
# PhD Thesis by Dr Paul Garcia

**Abstract**

This thesis describes the life and work of the mathematician Major Percy Alexander MacMahon (1854 - 1929). His early life as a soldier in the Royal Artillery and events which led to him embarking on a career in mathematical research and teaching are dealt with in the first two chapters. Succeeding chapters explain the work in invariant theory and partition theory which brought him to the attention of the British mathematical community and eventually resulted in a Fellowship of the Royal Society, the presidency of the London Mathematical Society, and the award of three prestigious mathematical medals and four honorary doctorates. The development and importance of his recreational mathematical work is traced and discussed. MacMahon's career in the Civil Service as Deputy Warden of the Standards at the Board of Trade is also described. Throughout the thesis, his involvement with the British Association for the Advancement of Science and other scientific organisations is highlighted. The thesis also examines possible reasons why MacMahon's work, held in very high regard at the time, did not lead to the lasting fame accorded to some of his contemporaries. Details of his personal and social life are included to give a picture of MacMahon as a real person working hard to succeed in a difficult context.

**Published work from this thesis**.

Some of the biographical details discovered during the work for this thesis were used to create an entry for MacMahon in the *Dictionary of Nineteenth Century British Scientists,* published by Thoemmes in 2004.

Biographical details and a summary of MacMahon's recreational work were used as the introduction to a reprint of MacMahon's 1921 book *New Mathematical Pastimes*, published by Tarquin in 2004.

Some of the material on recreational mathematics was used in two articles published in the Mathematical Association's *Mathematics in School* magazine: *The mathematical pastimes of Major Percy Alexander MacMahon: Part 1 - slab stacking*, Mathematics in School **34**(2) (2005) 23 - 25, and *Part 2 - triangles and beyond*, Mathematics in School **34**(4) (2005) 20 - 22; and an article in the magazine *Infinity* published by Tarquin: *From puzzles to tiling patterns*, Infinity, **2005**(2) (2005) 28 - 32.

# The life and work of Major Percy Alexander MacMahon
## PhD Thesis by Dr Paul Garcia

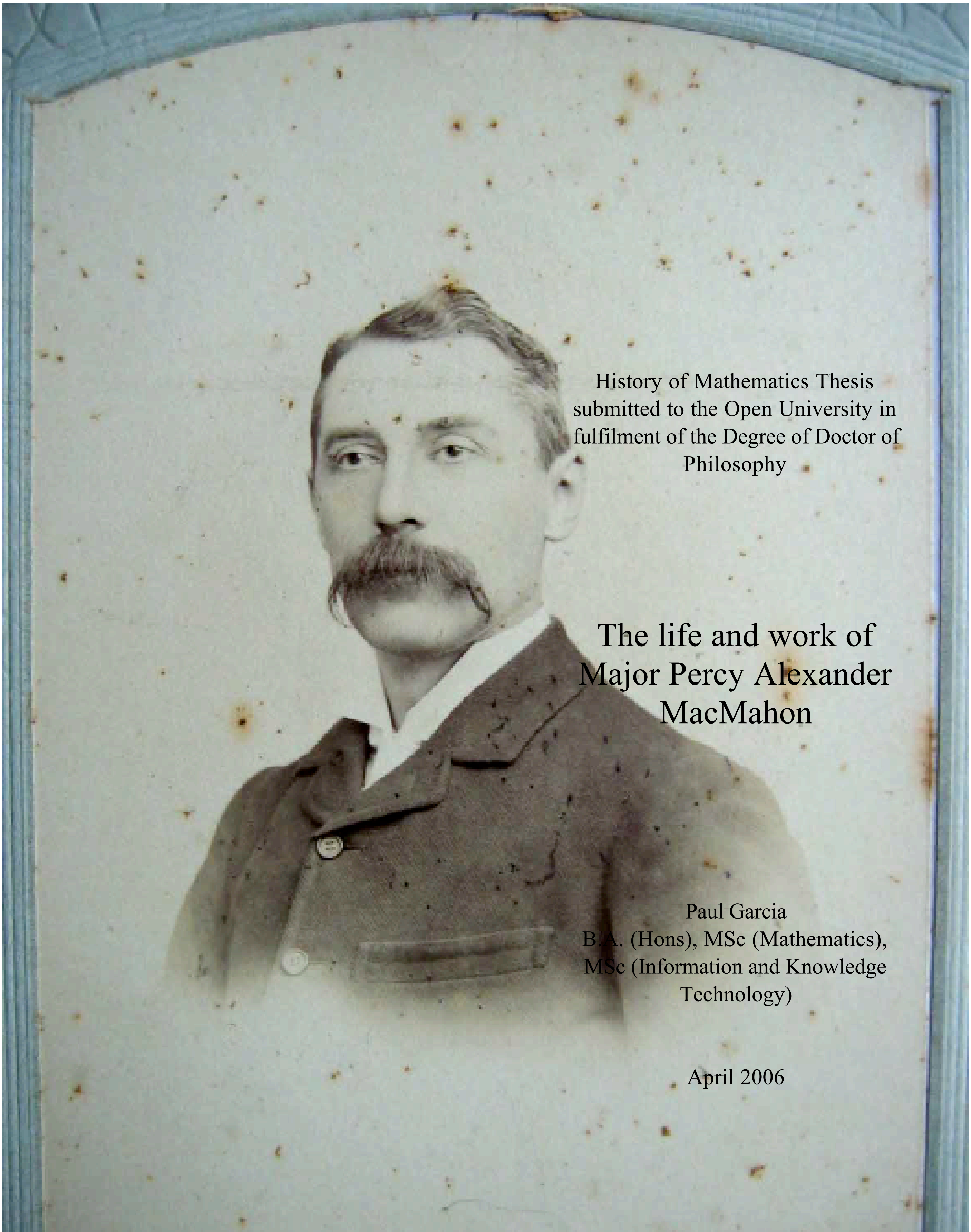

History of Mathematics Thesis submitted to the Open University in fulfilment of the Degree of Doctor of Philosophy

# The life and work of Major Percy Alexander MacMahon

Paul Garcia
B.A. (Hons), MSc (Mathematics), MSc (Information and Knowledge Technology)

April 2006

## Acknowledgements

"Remember, it's supposed to be fun"

Richard Feynman, quoted in *Some time with Feynman*, L. Mlodinow, 2003

I have pieced together this story of MacMahon's life and work from the writings of his three obituarists (H. F. Baker, H. H. Turner and A. R. Forsyth), notices of his death published in *The Times*, and information from various other sources, most notably St John's College, Cambridge, Winchester College, The Royal Military Academy (Sandhurst Collection), the Royal Irish Academy, Aberdeen University, St. Andrew's University, Trinity College, Dublin, the Royal Society, the Royal Astronomical Society, and the London Mathematical Society. The librarians at Cambridge University Library, the Bodleian Library, the British Library, the Public Record Office at Kew, the Family Record Office in Islington, Claire Harrington of the Glasgow City Council Library Service, and Maria Martinez, a genealogist in Monte Vista, Colorado, have all provided invaluable help and pointed me at material I might otherwise have taken decades to discover.

I would also like to thank all the people who have helped and shown great kindness and patience over the 7 years I have been working on MacMahon. In particular, George Andrews (the editor of MacMahon's *Collected Papers*); Anthony Unwin; the Humphreys family (Charles & Rosie Humphreys, Mick Humphreys, and Clare Middleton, all descendants of MacMahon's brothers); the late John Fauvel for his encouragement and advice in my first year of study; my supervisors, Robin Wilson for all the time and advice he has lavished on me, and June Barrow-Green for many interesting and helpful conversations.

Special thanks are due to Professor Margaret Farrell for details of the *Mayblox* puzzle, and Professor Robert N. Andersen for donating his aunt's *Mayblox* puzzle to my collection of MacMahon memorabilia , and to Kev Stenning of Rapido UK who helped me make a copy of the bronze bust of MacMahon.

Finally, I would like to thank my wife Mary, for putting up with this obsession for so many years.

# The life and work of Major Percy Alexander MacMahon
# PhD Thesis by Dr Paul Garcia

## Contents

**Figures**

**Tables**

# The life and work of Major Percy Alexander MacMahon
# PhD Thesis by Dr Paul Garcia

## Chapter 1 - Introduction

> .. there to go at a spot of combinatorial analysis, that favourite pastime of retired army officers, with a rattling intense devotion.
>
> Thomas Pynchon, *Gravity's Rainbow*

> ... historians have also to examine the way mathematicians as people have functioned within their society and their mathematical tradition.
>
> John Fauvel, *MA290 Unit 17*, Open University

### Scope of the thesis

The conventional view of Victorian combinatorics in Britain is that it was the province of amateurs[1] - barristers, clergymen, military men - and it might be argued that MacMahon was the among the most important of these. However, this thesis will argue that MacMahon was indeed a most important Victorian combinatorialist but was by no means an amateur.

Percy Alexander MacMahon was an unusual mathematician. He did not come from a family of mathematicians - his father and brothers were soldiers - and he did not have the university education usually associated with a prominent mathematician. Further, he entered the world of mathematics at the relatively late age of 26, and did not become well known until he was 30. Before he joined the advanced class in mathematics at the Royal Military Academy in 1880 he had shown no hint of any particular interest in, or talent for, doing mathematics.

The speed with which MacMahon rose to prominence in the mathematical community of the late nineteenth century is remarkable. From the start of his career he was held in high regard by several famous and well-established mathematicians; James Joseph Sylvester[2], for example, specifically mentioned MacMahon three times in his inaugural speech as Savilian Professor of Geometry at

[1] [Biggs, 1979] and [Biggs, et al., 1995]

[2] James Joseph Sylvester, 1814 - 1897, was second Wrangler in the Cambridge tripos in 1837, Professor of Mathematics at Woolwich from 1855 - 1870 and Savilian Professor of Geometry in Oxford from 1883 until his death in 1897.

Oxford[3]. A review[4] of MacMahon's two-volume work, *Combinatory Analysis*[5], compared him with such luminaries as Fermat, Pascal, Euler, Lucas and Sylvester, at least in the field of combinatory analysis, and described his mathematical style as 'impeccable'.

MacMahon did pioneering work in invariant theory, symmetric function theory, and partition theory. He brought all these strands together to bring coherence to the discipline we now call combinatorial analysis[6]. According to Peter J. Cameron[7], combinatorial analysis was not well thought of in the first half of the twentieth century. MacMahon was ahead of his time, since in the second half of the twentieth century, attitudes to combinatorial analysis changed dramatically, as illustrated by this quotation from Professor Steven Pinker[8], referred to by Cameron:

> It may not be a coincidence that the two systems in the universe that most impress us with their open ended complex design - life and mind[9] - are based on discrete combinatorial systems.

MacMahon also carried out research in recreational mathematics, creating visual puzzles, but again in a pioneering and non-traditional way. This thesis will describe only as much of the mathematics as is necessary to make the context and extent of MacMahon's work clear; it is not intended that this will be a detailed technical exposition of the topics he worked on.

How did MacMahon come from a non-mathematical background to a position of such prominence so quickly ? What was special about his work ? Why did he fade from prominence ? What was his legacy to mathematics ? To answer these questions, his personal life, which was also unconventional, and his working life will be described. MacMahon worked as a soldier, a teacher and a government official, whilst simultaneously pursuing his mathematical researches.

---

[3] Sylvester gave his inaugural speech for the Savilian professorship on 12 December 1885 [Sylvester, 1886]. MacMahon's presence in the audience was noted by Sylvester, and there were two references to the process developed by MacMahon in his 1883 paper on seminvariants - see chapter 6.
[4] [Anonymous, 1916].
[5] [MacMahon, 1984]. The date of this reference is to the 1984 Chelsea reprint.
[6] Described in chapters 6, 7 and 8.
[7] [Cameron, 2001]. Peter Cameron is Professor of Mathematics at Queen Mary College, London.
[8] [Pinker, 1994], p. 85. Steven Pinker (1954 - ), is the Johnstone Family Professor of Psychology at Harvard University [Pinker, 1994].
[9] The two systems Pinker is referring to are the DNA code and language, respectively.

During his career he published over 125 papers and four books. *Combinatory Analysis* is still in print[10] and *New Mathematical Pastimes* was reprinted in 2004 as part of the work for this thesis. A list of all MacMahon's mathematical papers, in chronological order, may be found in *Percy Alexander MacMahon Collected Papers*[11]. This list is reproduced in the bibliographies.

**Content of the thesis**

MacMahon's life and career are described in succeeding chapters in more or less chronological order. The breaks between the chapters divide his life into six natural periods, so the chapters are somewhat unequal in length.

Chapter 2 outlines the type of mathematical education that MacMahon received as a Gentleman Cadet at the Royal Military Academy in Woolwich, the training institution for the Royal Artillery and Royal Engineers, and his subsequent military service in India. MacMahon's education was not particularly mathematical, and there is no evidence that his gifts were recognised by his teachers.

The beginnings of MacMahon's mathematical career are covered in Chapter 3. Whilst working as an instructor at his *alma mater*, the Royal Military Academy at Woolwich, MacMahon made a significant discovery in the *algebra of forms* (also known as *invariant theory*), which brought him to the attention of the wider mathematical community in London. He also began to study *partition theory*, in which field he was to become the most important researcher of the late nineteenth and early twentieth centuries.

During the 1890s, MacMahon's reputation as a mathematician grew. Chapter 4 chronicles his career in the Royal Artillery, his continuing work on partitions, and his election to the Royal Society and to the

[10] [MacMahon, 1984]. *Combinatory Analysis* was originally published as two separate volumes in 1915 and 1916; it was republished by Chelsea as single volume in 1960, 1964 and 1984.
[11] [Andrews, 1978]

Presidency of the London Mathematical Society. He also became involved with the Royal Astronomical Society, received his first honorary degree and experienced a disappointment at the hands of the electors of the Savilian Professor of Geometry in Oxford.

The material in Chapter 5 covers the period between MacMahon's retirement from the Royal Artillery and his appointment as Deputy Warden of the Standards at the Board of Trade. It shows the breadth of MacMahon's involvement with the mathematical and scientific community, and the respect he commanded from the distinguished organisations he had joined whilst teaching at Woolwich.

Chapter 6 describes one of the most eventful periods of MacMahon's life. He was by now an established figure enjoying a distinguished reputation, with the Presidency of the London Mathematical Society and two honorary degrees already to his credit. A further two honorary degrees were awarded in 1911. He also wrote four books in this period, consolidating his work in partitions and symmetric functions into *Combinatory Analysis*, published in two volumes in 1915 and 1916[12], and very favourably reviewed. An introductory volume to *Combinatory Analysis* followed in 1920[13], and in 1921 MacMahon wrote *New Mathematical Pastimes*, developing the puzzle work he had carried out in the early 1890s. MacMahon also wrote more than 30 papers between 1906 and 1922, in which he developed the ideas he had expounded earlier in symmetric function theory and partition theory. In 1906 he took up the post of Deputy Warden of the Standards with the Board of Trade.

Finally, Chapter 7 records MacMahon's work at St John's College, Cambridge, his influence on some other mathematicians, and his continuing mathematical work. Ill-health forced him to move to the south coast in 1928, where he died in 1929.

[12] [MacMahon, 1984].
[13] [MacMahon, 1980].

## Sources

George Andrews' *Collected Works of Percy Alexander MacMahon*, published in two volumes in 1977 and 1986, was the main source for the mathematical papers. References to these papers include the numbering system used by Andrews, which is given as [Chronological paper number; chapter in the *Collected Papers*], so that the interested reader may readily locate items. MacMahon's military papers do not include such a reference, as Andrews did not include them in the *Collected Papers*. Obituaries of MacMahon were published by A. R. Forsyth[14], H. H. Turner[15] and H. F. Baker[16].

Private information provided by descendants of MacMahon's brothers has also been invaluable.

The archive of the British Association for the Advancement of Science held at the Bodleian Library, the Public Record Office[17] at Kew, London, the Family Records Centre in Islington, London, and the India Institute Library in Oxford provided much background material, some of which has not found its way into the final thesis. A full list of the archival sources consulted is given in the bibliography.

[14] [Forsyth, February 1930]. Andrew Russell Forsyth, 1858 - 1942, was senior Wrangler in 1881, elected FRS in 1886, and was Sadleirian Professor of Pure Mathematics at Cambridge from 1895 - 1910.

[15] [Andrews, 1986, pp. 940 - 945]. Herbert Hall Turner, 1861 - 1930, was educated at Cambridge and worked at the Royal Greenwich Observatory. He was appointed Savilian Professor of Astronomy in Oxford in 1893 and was Bruce Medal winner in 1927.

[16] [Baker, June & October 1930]. Henry Frederick Baker, 1866 - 1956, was senior Wrangler (bracketed) in 1887, won the Smith's prize in 1889 and was elected FRS in 1898. He won the Sylvester medal in 1910, and the De Morgan Medal of the London Mathematical Society in 1905. He held the Cayley lectureship in Cambridge from 1903 - 1914, was Lowndean Professor of Astronomy and Geometry in Cambridge from 1914 - 1936, and edited Sylvester's *Mathematical Papers* 1904 - 1912.

[17] During the course of the work on this thesis, the Public Record Office was renamed The National Archive.

## Chapter 2 Early Life 1854 - 1881

Percy Alexander MacMahon was born in Sliema, Malta, on 26 September 1854, the second son of Brigadier-General Patrick William MacMahon and Ellen Curtis, daughter of George Savage Curtis of Teignmouth, Devon. There are no official records of his birth, and apart from an anecdote about noticing the way cannon balls were stacked that MacMahon himself told in later life, nothing is known of his early childhood.

### Early education

MacMahon attended the Proprietary School in Cheltenham. At the age of 14 he won a Junior Scholarship to Cheltenham College, which he attended as a day boy from 10 February 1868[18] until December 1870. MacMahon was always ranked in the top five in his classes[19], but never came first. It is likely that one of his mathematics teachers was Edward Walker[20] from Trinity College, Cambridge, who had been 8th Wrangler in 1844 and was a master at the school from 1861 until 1873.

At the age of 16 MacMahon was admitted to the Royal Military Academy at Woolwich, just a year after J. J. Sylvester had been forced to retire from the post of Professor on reaching retirement age. (MacMahon later came to work in the same field as Sylvester, and there are some parallels between the careers of both men. He was eventually to write the obituary of Sylvester published in *Nature* in 1897 and in the *Proceedings of the Royal Society* in 1898[21]).

[18] MacMahon's nomination form was signed by Edward Thomas Wilson of 6 Montpellier Terrace, Cheltenham, on 8 February 1868. As MacMahon was a day boy rather than a boarder, and his father was then stationed in Peshawar, it seems likely that Wilson was acting as his guardian.
[19] Class sizes varied between 8 and 16 pupils.
[20] Edward Walker joined Trinity College in 1839 and studied under George Peacock. He received his B.A. in 1844 and worked as an Assistant Tutor from 1846 - 1847.
[21] [MacMahon, 1897, [50;19]].

**Royal Military Academy, Woolwich**

The Royal Military Academy (RMA) at Woolwich was founded in 1741 to train gentleman cadets for the Royal Artillery and the Royal Engineers. The only other similar institutions were the Royal Military College (RMC) in Camberley, founded in 1800 as a school for staff officers, and the East India Company's Military Seminary, founded in 1809 at Addiscombe House in Croydon, which was taken over by the RMC in 1860. Sandhurst replaced both institutions in 1939.

***Mathematics at Woolwich***

The syllabus that MacMahon would have followed was the subject of revision and controversy during the quarter century prior to his enrolment. In the late 1840s, Samuel Hunter Christie[22] rewrote the course, which was regarded as inadequate[23]. According to *The Records of the Royal Military Academy*[24], Christie recommended these texts:

1. *Algebra* - Young's "Elementary Treatise." Young's "Elements of the Theory of Equations."
2. *Geometry* - Simson's "Euclid" (used at present).
3. *Trigonometry* - Hind's
4. *Tables of Logarithms, &c.* - Hutton's (used at present).
5. *Conic Sections.* - Hamilton's "Analytical System."
6. *Differential and Integral Calculus* - La Croix's (in French).
7. *Mechanics* - Whewell's "Elementary Treatise."

However, Professor Christie's efforts were not well received:

> Sir Thomas Hasting's Committee decided [in 1848] that it [Christie's course] was not suitable for the Cadets. They recommended that the Royal Engineer Officers at the Academy should compile a Mathematical Course under the superintendence of Captain Harness, R.E., Professor of Fortification, which was done in spite of much remonstrance from Mr. Christie. Dr. Hutton's[25] "Course of Mathematics" was then discontinued, after having been the textbook for Mathematics at the Royal Military Academy for nearly 50 years[26] .

[22] Samuel Hunter Christie, 1784 - 1865, was second Wrangler in 1805, elected FRS in 1826, and Professor of Mathematics at the Royal Military Academy from 1838 - 1854.
[23] [Rice, 1996], page 400.
[24] [Buchanan-Dunlop, 1895]. This book is a reprint of a volume originally published in 1851, which is all that was left of the RMA records after a fire in 1873. Later records taking it up to 1892 were added at the time of the 1895 reprint.
[25] Charles Hutton, 1737 - 1823, was professor of mathematics at the RMA from 1773 - 1807, and editor of the *Lady's Diary* from 1773 - 1818.
[26] [Buchanan-Dunlop, 1895].

According to Rice[27], the replacement for Christie's course was actually written by Stephen Fenwick, William Rutherford[28] and Thomas Stephen Davies and published between 1850 and 1852. This course was also soon subject to proposals for revision:

> The Rev. Canon Moseley, M.A., F.R.S., a member of the Council of Military Education, was then invited to suggest a Mathematical Course for the Gentlemen Cadets. This he did in a long paper, which was submitted to the Mathematical Professor and his Assistants for their remarks in October, 1859. These Gentlemen pointed out that Canon Moseley's course was chiefly copied from that of the [École] Polytechnic and other Continental Schools in which, often, Officers of the Civil Administration were trained[29].

Moseley's suggestions were also deemed too difficult for the Cadets, and his proposals were not followed up.

Sylvester was appointed Professor of Mathematics at Woolwich in 1855, a post he had unsuccessfully applied for in 1854. He was pleased to accept it since he was ineligible for appointments at Oxford or Cambridge because he was a Jew. It was not an easy appointment, for Sylvester did not like being under military command, and found that his duties were far too time consuming. He was forced to retire in 1870, and had to fight to get his pension[30].

By the time MacMahon started at Woolwich, the syllabus covered arithmetic, including proportion and compound interest; algebra, including Horner's method for the solution of equations, and permutations and combinations; geometry and conic sections, including the application of algebra to geometrical problems; plane trigonometry and mensuration, including surveying; analytical plane geometry and spherical trigonometry, analytical geometry in two and three dimensions; differential and integral calculus, statics and dynamics. Throughout the course, there was a strong emphasis on practical applications useful for military purposes; for example, the statics included material on terracing, piers and arches and the dynamics covered ballistics.

---

[27] [Rice,1996, p. 401].
[28] William Rutherford had revised Hutton's book himself only a few years previously. The 1832 edition of Hutton's three volume work still taught Newtonian calculus in the form of fluxions and fluents. Not until the 1840 edition, revised by Rutherford and still in print in 1846, was the Continental calculus using Leibniz notation taught.
[29] [Buchanan-Dunlop, 1895].
[30] See [Parshall, 1998, p. 136], for more details.

*Entering the Academy*

By the age of 16 MacMahon would already have been extremely well-educated, able to speak French and German, with a good knowledge of Latin and Greek, and a grounding in pure and applied mathematics.

To be admitted to the RMA[31], MacMahon had to pass a stiff entrance examination in Mathematics, including the first six books of Euclid, English language, Classics, French and German languages and history, Experimental sciences, Natural sciences and Drawing.

To illustrate the standard required at entry, below are some questions from an entrance examination session. These examples are from the 1858 Pure Mathematics session, held at King's College, London[32]. They show that MacMahon would have been familiar with the concepts of symmetric functions and the geometry of circles from an early age.

> Section II Q3: Find the sum of the fourth powers of the roots of the equation $x^4 + x^3 - 7x^2 - x + 6 = 0$; and prove that any symmetrical function of the roots of an equation may be expressed in terms of the coefficients. (Set by the Rev. C. Graves)

> Section III Q3: Prove that the product of the radii of the four circles each of which touches the same three intersecting right lines is equal to the square of the area of the triangle included between these right lines. (Set by J. J. Sylvester)

A comparison with the mathematical entrance requirements for the army, navy and civil service[33] shows that the Royal Military Academy was by far the most demanding. The entrance examination for Sandhurst covered only algebra and geometry, and the navy and civil service asked for a knowledge of Euclid and some mechanics. Woolwich demanded a more sophisticated grasp of algebra and the theory of equations, as well as trigonometry, statics and dynamics. Many of the questions for the Woolwich examinations were placed in a military context; for example, a trigonometrical question on angles of declination referred to a reconnaissance observer in a balloon, and a question about dynamics was placed in the context of a pistol shot fired at a moving train.

[31] As a Gentleman Cadet, MacMahon's tuition fees, board and the cost of books and equipment, including his uniform, would have been paid by his parents. However, because he was the son of a serving officer, the fees would have been subject to a discount.
[32] Taken from the *Report on the examination for admission to the Royal Military Academy as Woolwich and to The Military Class, University of Dublin*, HMSO, 1858, pp. 57 and 59.
[33] See [Austin, 1880], [Austin, 1882] and [Arnett 1874].

The standard required in some topics compared very favourably with that required for the award of a university degree. The questions in applied mathematics posed in degree examinations at Cambridge and London[34], for example, could have been tackled by a Woolwich candidate at entry. In particular, the questions in [Kimber, 1880] and [Gantillon, 1852] are remarkably similar to the questions on the Woolwich papers.

***MacMahon at Woolwich***

On 7 February 1871 MacMahon entered the Royal Military Academy in Woolwich as a Gentleman Cadet.

The period of instruction at Woolwich lasted a maximum of three years, after which time a Cadet would either have to have passed or would be required to leave. MacMahon passed after two years[35]; indeed, in 1882, when he was an Instructor at the Academy, the maximum length of the course was reduced to two years, as a result of the increased demand for officers resulting from various wars across the Empire. MacMahon started in February 1871 in the 5th Class, and worked through the 4th, 3rd, and 2nd Classes to the 1st, and final, Class, in November 1872. He was passed on 20 February 1873.

A weekly study timetable would have looked something like the one in Appendix 1, which is actually from 1892, some twenty years after MacMahon had been a Cadet. A fire in 1873[36] destroyed most of the Academy's records, so few examples of such ephemera have survived. An 1872 daily timetable is preserved in Guggisberg's *The Shop: the story of the Royal Military Academy*[37]. It reads:

> 6:15 am Extra drill
> 7 am Breakfast (9am Sundays)
> 8 - 9 Drill. First & second class riding
> 9:30 - 11:30 Study
> 11:30 Luncheon (none on Sundays) - Bread, biscuits, butter *ad lib*, beer 1 pint each
> 12 - 2 Study
> 2:15 Dinner (Beer, 1 pint per head)
> 5pm Afternoon lunch
> 6 - 8 Study
> 8pm Supper in own rooms
> 10 Rounds
> 10:30 Lights out

[34] [Kimber, 1880], [Gantillon, 1852].

[35] The length of time taken by a Cadet to complete the course varied according to the personal circumstances of the Cadet and the requirements of the Royal Artillery for officers.

The Professor of Mathematics and Mechanics at the time was Morgan W. Crofton BA FRS[38], assisted by three Instructors, Major W. H. Wordell, John MacLeod and Lieutenant E. Kensington.

It was customary for the Governor of the Academy[39] to write a report to the Commander-in-Chief, His Royal Highness the Duke of Cambridge, KG, on the success and conduct of each group of finishing cadets, certifying that they were ready to take commissions. Normally, such reports commended the good behaviour and example of this most senior class, and a sword was awarded to one student for good conduct. MacMahon's classmates, however, although well behaved, suffered from a "want of punctuality" and "general slackness in their duties"[40]. Seven had their marks reduced at the end of the two years (thus adversely affecting their choice of regiments) and the Governor was unable to recommend anyone for the good conduct sword. This is mentioned to illustrate that MacMahon was not especially distinguished as a Cadet, and had spent his time at Woolwich in the company of boisterous and not particularly academic young men. So at this stage, there was no hint of the path he was going to follow.

[36] Described by a witness in [Guggisberg, 1900, pp. 145 - 147].
[37] [Guggisberg, 1900, p. 138]. This book also contains interesting contemporary anecdotes concerning the daily life of Cadets (pp. 127 - 150).
[38] Morgan W. Crofton, 1826 - 1915, succeeded J. J. Sylvester on 1 August 1870. He was elected FRS in 1868; the election citation read: "Author of original researches and discoveries in the Geometry of Curves, the Theory of Probability, and the Integral Calculus, contained in papers in the Annales de Mathématiques, the Comptes Rendus de l'Académie des Sciences, the Proceedings of the London Mathematical Society, the Oxford Cambridge and Dublin Messenger of Mathematics, and in a Memoir recently presented for publication to the Royal Society of London: distinguished for his knowledge and ability as a Mathematician (having been First Senior Moderator, and First Gold Medallist at the Degree Examination of the University of Dublin in the year 1847)." He retired in 1884.
[39] At this time the Governor was Major-General J. Lintorn A. Simmons. This was the rank he used when signing the reports, although the Army List shows him as a Lieutenant-General. He died in 1875.
[40] These remarks are quoted from the letter written by the Governor to the Commander-in-Chief, dated 20 February 1873, to be found in the Sandhurst Collection.

MacMahon's father died[41] a few months after MacMahon had started at Woolwich. The lack of personal correspondence or diaries means that there is no record of how MacMahon may have been affected by this event.

**India**

On 12 March 1873, MacMahon was posted to Madras, India, with the 1st Battery 5th Brigade, with the temporary rank of Lieutenant. The Army List showed that in October 1873 he was posted to the 8th Brigade in Lucknow. A map showing the places where MacMahon served in India can be found in Appendix 2.

The Royal Artillery record notes that in 1874 he was moved to C/11 in Dinapore[42] where he stayed until 1876. In March 1874, MacMahon's temporary rank of Lieutenant was made permanent and backdated to 12th September 1872. The Royal Artillery records show him spending the two years from 1876 with the 8th Brigade, first in Multan[43] and then in Meerut.

MacMahon's final posting was to the No. 1 Mountain Battery with the Punjab Frontier Force at Kohat on the North West Frontier[44]. He was appointed Second Subaltern on 26 January and joined the Battery on 25 February 1877.

[41] On 14 October 1871 at 10 Devonshire Place, Brighton.
[42] Dinajpur or Dinagepour is now in Bangladesh.
[43] Multan is now in Pakistan.
[44] Kohat is in Pakistan, close to the Khyber Pass on the border with Afghanistan.

MacMahon is reported by Turner to have taken part in a punitive raid on a tribe known as the Jowaki Afridis in August 1877, but he was not a part of this expeditionary force. In the *Historical Record of the No. 1 (Kohat) Mountain Battery, Punjab Frontier Force*[45] it is recorded that he was sent on sick leave to Muree (or Maree), a town north of Kohat on the banks of the Indus river, on 9 August 1877. He thus missed the beginning of the raid by three weeks. On 22 December 1877 he started 18 months leave on a medical certificate granted under GGO[46] number 1144.

The nature of his illness is unknown[47]. Officers did not receive discharge papers in the same way as ordinary soldiers[48], whose documents contain a wealth of interesting information. This period of sick leave was one of the most significant occurrences in MacMahon's life. Had he remained in India he would undoubtedly have been caught up in Roberts's War against the Afghans, a bloody adventure which lasted two years and achieved nothing, neither in a military nor a political sense[49].

In early 1878 MacMahon returned to England and the sequence of events began which led to him becoming a mathematician rather than a soldier. The Army List records a transfer to the 3rd Brigade in Newbridge at the beginning of 1878, and then shows MacMahon as 'supernumerary' from May 1878 until March 1879. The Royal Artillery records show him spending a little time in Malta in early 1878.

In January 1879 MacMahon was posted to the 9th Brigade in Dover, moving to Sheerness in 1880. In the same year he enrolled in the Advanced Class for Artillery Officers at Woolwich. This was a two year course covering technical subjects and a foreign language. Successful completion of the course resulted in the award of the letters "p.a.c" (passed advanced class) after MacMahon's name in the

[45] Published by the Punjab Government Press in Lahore in 1886 (page 12).
[46] Governor's General Order.
[47] Malaria or dysentery are the most likely candidates. The records which might have contained the information were destroyed by the Ministry of Defence in 1969.
[48] This information was given in an e-mail from the Public Record Office, dated 18 September 2001.
[49] [Farwell, 1999, p. 216].

Army List. It was not a popular course. Commanding officers disliked it because it meant the loss of an officer for two years, and it was unpopular amongst the officers themselves for many reasons, six of which were identified in the *Report of a Committee on the Advanced Class of Artillery Officers,* published in 1884[50]:

- unattractive appointments, both in terms of pay and position, upon completion of the class;
- uncertainty about getting an appointment at all;
- the prospect of being retired early;
- uncertainty over the benefits of the course;
- being removed from ordinary regimental work for two years;
- unwillingness of "officers of a certain age and standing" to submit themselves to a regime they felt was 'too scholastic'.

It is likely that the nature of MacMahon's illness was such that a return to combat duty was impossible. His work at Dover and Sheerness as Inspector of Stores may have been dull. The Advanced Class would have been the only avenue open to MacMahon to enjoy any prospect of advancement at all.

The military training that MacMahon received at Woolwich and in India would prove useful in providing the self-discipline necessary to carry out the complex calculations demanded by combinatory analysis. The social status conferred by military rank also enabled MacMahon to overcome his lack of a traditional university education. It allowed him to meet other mathematicians on an equal footing.

[50] [PRO, 1884]

# Chapter 3 Woolwich 1881 -1890

This chapter describes the events leading from MacMahon's first entry into the British mathematical community to his election to the Royal Society in 1890. In particular, his first significant discovery in invariant theory will be described, a result that also enabled him to collaborate with Cayley for many years in the calculation and tabulation of seminvariants and perpetuants and the syzygies, or linear relations, between them. These mathematical terms will be explained later.

### Instructor at RMA Woolwich

After he passed the Advanced Course and had been promoted to the rank of Captain on 29 October 1881, MacMahon took up a post as Instructor at the Royal Military Academy on 23 March 1882. Here he met Alfred George Greenhill[51], Professor of Mathematics at the Royal Artillery College. Joseph Larmor[52], in a letter to *The Times*[53] published after MacMahon's death, wrote, 'The young Captain threw himself with indomitable zeal and insight into the great problems of the rising edifice of algebraic forms, as was being developed by Cayley, Sylvester and Salmon[54].'

MacMahon had mathematical work to do for the military, and in 1881 he wrote a paper on the trajectories of different types of artillery guns[55], published in *The Proceedings of the Royal Artillery Institution.*

The Royal Artillery needed to ensure that guns were aimed accurately and that as few rounds as possible were wasted, since munitions were expensive to manufacture and deliver to the far-flung battlefields of the Empire. A heavy gun had only a limited life before it required reboring, although it was usually cheaper to replace the barrel. In India, for example, artillery had to be transported overland by elephant for many hundreds of miles; once engaged in battle, it was vital that every shell hit the intended target, and so a great deal of work went into the accurate determination of trajectories. The weight and speed of the shells meant that air resistance was not negligible, and experimental work was carried out to determine the velocity of a projectile during its flight.

[51] Alfred George Greenhill, 1847 - 1927, was second Wrangler and joint winner of the Smith's prize in 1868, and Professor of Mathematics at the Royal Military Academy from 1876 - 1908. He worked mainly in ballistics and elliptic functions, but also in the theory of elasticity and hydrostatics. He was LMS President from 1890 - 1892, was awarded the De Morgan Medal in 1902 and the Royal Medal in 1906, and was knighted in 1908.
[52] Joseph Larmor, 1857 - 1942, was senior Wrangler 1880, Lucasian Professor of Mathematics in Cambridge from 1903 - 1932, and was knighted in 1909.
[53] This letter from Larmor was published in the edition dated 31 December 1929.
[54] The Irish mathematician George Salmon, 1819 - 1904, was elected FRS in 1863, won the Royal medal in 1868 and was Provost of Trinity College, Dublin, from 1888 - 1902.
[55] [MacMahon, 1881].

An important experimenter in the mid-19th century was the Reverend Francis Bashforth[56]. He carried out two series of experiments (during 1865 - 1870 and 1878 - 1880), to determine the motion of shells fired from heavy artillery pieces. Bashforth used an electrically operated timing device of his own invention, known as *Bashforth's Chronograph*, which timed the passage of a shell as it passed through a series of screens at 150 ft intervals[57].

Bashforth's device was a great improvement on the ballistic pendulum, a machine invented in the mid-1700s by Benjamin Robins, author of *New Principles of Gunnery* [1742], which had been used by Charles Hutton at the Royal Military Academy at the end of the 18th century.

Bashforth himself carried out extensive analysis of his results, described in several articles and books, using actuarial techniques developed by W. S. B. Woolhouse[58] to smooth the tabulated results.

Bashforth's data was used by Greenhill and MacMahon to develop a mathematical model of the shells' trajectories. The model was used to prepare firing and range tables for artillery officers. Greenhill noted that the production of the tables was much improved by the purchase by the Royal Artillery Institution of an 'Arithmometer', a mechanical calculating device invented in the 1820s and popular with banks and other commercial institutions.

[56] Francis Bashforth, 1819 – 1912, was educated at Doncaster Grammar School and St. John's College, Cambridge. Second Wrangler 1843, he worked for several years as a civil engineer in the railway industry, surveying potential routes for new lines, where he learned the careful measurement skills needed for his work in ballistics.
[57] This device now resides in the Science Museum in London, and has an interesting history.
[58] Wesley Stoker Barker Woolhouse, 1809 - 1893, was Deputy Superintendent of the Nautical Almanac 1830 - 1837, editor of *The Lady's & Gentleman's Diary* from 1844 - 1865, and a cofounder of the Institute of Actuaries in 1848.

In 1884 MacMahon calculated the tables for Greenhill's article on the motion of a particle in a resisting medium[59], as well as writing three articles on trajectories of his own, published in the *Proceedings of the Royal Artillery Institution*[60]. It is in the second of these articles, *On the motion of a projectile*, that the link with the study of algebraic forms is first made clear. In it, MacMahon had to solve a symmetrical cubic equation describing the trajectory of a projectile. An outline of the work which led from ballistics to symmetric functions is given in Appendix 4.

The complexity of using symmetric functions and the associated tables necessitated quick methods of calculating the range of a gun in the field. In 1885 MacMahon published a note in the *Proceedings of the Royal Artillery Institution*[61] entitled 'A practical rule for range finding', which read thus:

> The following practical rule for range-finding with the sextant, when the vertical angle is less than $5^\circ$, and the triangle nearly isosceles, may or may not be well known. It is *very* accurate, and only requires a quarter of a minute on the back of an envelope:-
> Multiply the base in yards by 57.3 and divide by the vertical angle in degrees. The result is the range in yards.
> E.g. Base 50 yds, Vertical angle $1^\circ$, range = 57.3 x 50/1 = 2865 yds

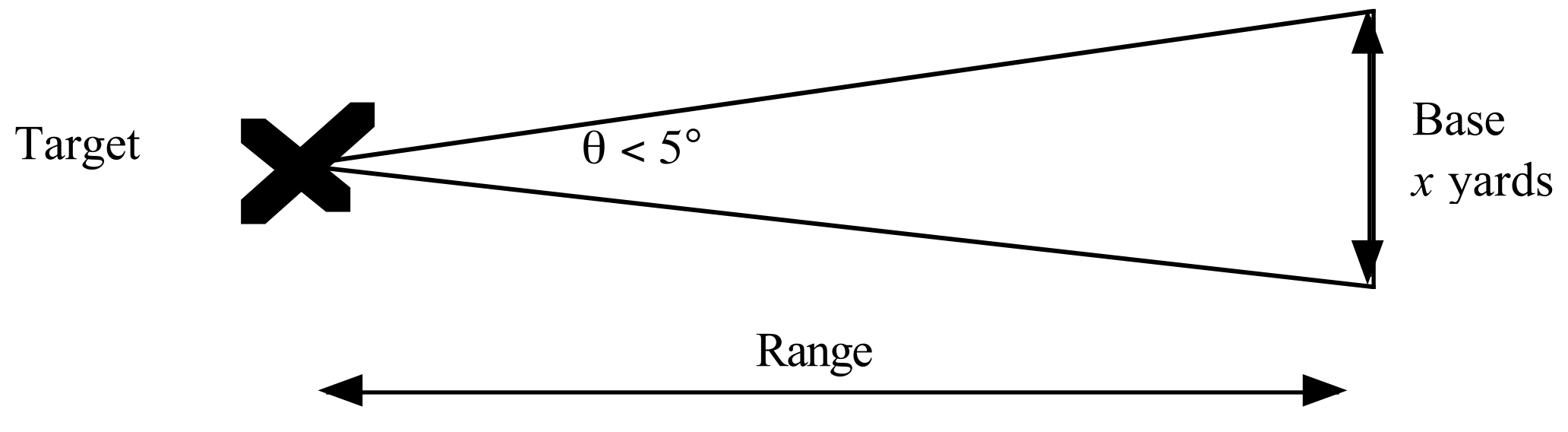

Figure 0: A rule for range finding

[59] [Greenhill, 1884].
[60] [MacMahon, 1884].
[61] [MacMahon, 1885].

This works because $\theta \approx \tan\theta$ for small angles, where $\theta$ is measured in radians, and 1 radian is about 57.3°. So if two sightings of the target $x$ yards apart (the *base*) give an angle of less than 5°, then the range in yards is the base divided by the angle in radians. To convert an angle in degrees to an angle in radians, divide by 57.3. This gives the formula $Range \approx \frac{base}{\frac{\theta}{57.3}} = \frac{base \times 57.3}{\theta}$, as described. This was the young Captain's last contribution to the *Proceedings of the Royal Artillery Institution* until the mid-1890s, by which time he had become a Major and an Instructor at the Royal Artillery College.

During this early period as an Instructor, MacMahon began to get involved with the wider mathematical community in England. He joined the British Association for the Advancement of Science (BAAS) in 1883. The minutes of the meeting held on 21 September 1883 in Southport are the first evidence of his attendance at a BAAS meeting, and record his joining the Section A (Mathematical and Physical Science) committee. On 22 September, MacMahon read a paper entitled *On symmetric functions, and in particular on certain inverse operations in connexion therewith.* On 14 December 1883, he was officially admitted to the London Mathematical Society (LMS) and 'signed the book' at the 19th session[62].

It is likely that Greenhill had a hand in suggesting that MacMahon should join both these organisations, as a way of encouraging his young protegé to meet with other mathematicians and scientists. MacMahon's colleagues did not join the LMS, suggesting that this was not common practice for Woolwich instructors at that time, although other Woolwich staff before and since have been members.

## Teaching

The professor of mathematics at Woolwich at this time was Harry Hart[63], a civilian who had been appointed as an Instructor in 1873. There were three other instructors[64]. The cadets were taught in whole-year groups, some 40 to 50 at a time in the classroom, with the professor and three of the four instructors present. There was little formal input, since cadets were expected to interpret the syllabus and work from the set texts on their own, only asking for help from the staff in the room as necessary. There was no check that progress was being made until the end of term examinations, unless a candidate was obviously not working, which would be reported to the senior officer as idleness. When MacMahon had been a Cadet it was the practice for a Cadet corporal to be present in the room to see to discipline, but this was not the system in force when he became an Instructor. The professor and the instructors did give lectures on various topics according to a schedule given to the cadets at the beginning of each term, but attendance by Cadets was entirely voluntary[65].

[62] Amongst the many signatures in black or blue ink, MacMahon's violet ink stands out.
[63] Harry Hart was a member of the LMS; he joined in 1874.

There were two syllabuses - the *Obligatory*, which every cadet was required to pass, and the *Voluntary*, which students could elect to take as an additional subject. The Obligatory syllabus contained much applied mathematics useful in a military context:

I Geometry: Plane Geometry, Elements of Conic Sections, Elements of Solid Geometry;
II Mensuration of Planes and Solids;
III Mechanics and Theoretical Mechanics; Statics and Dynamics;
IV Hydrostatics;
V Applied Mechanics: Stability of Structures and Strength of Materials;
VI Elements of Mechanism; Theory of Work.

The Voluntary syllabus contained the same topics, but at a more advanced level, and also the areas of Probability, Theory of Equations, Co-ordinates, and the Differential and Integral Calculus.

At an Academic Board meeting on 30 September 1887, the only such meeting at which MacMahon's attendance is recorded[66], the Governor wanted to change the system of teaching so that all four instructors and the Professor were present for all classes, and that the Professor was the only staff member allowed to lecture. Professor Hart protested vigorously that this would place an unduly heavy workload on the staff and thus make it impossible to recruit Instructors of a suitably high calibre. There is no indication that MacMahon contributed to this discussion.

[64] Captain F. W. Boteler, Mr E. F. S. Tylecote and Mr W. Foord-Kelsey.
[65] These arrangements are described by Harry Hart in the transcript of an inspection interview held in December 1885 at the War Office [PRO, 1886].
[66] This document is in the Sandhurst Collection.

This method of instruction implies a lack of any serious preparation to do, and thus provided MacMahon with plenty of opportunity to pursue his own mathematical interests, which were far more advanced than the syllabus he was teaching. He quickly became an important contributor to the field of algebraic forms, largely by the invention of new methods. Principally, he introduced the use of *Hammond operators* to the study of *invariants* and *covariants* of *symmetric functions* (these technical matters will be discussed later). So deep was his absorption in this work that his military colleagues described him as a "good soldier spoiled"[67]. During his time as an Instructor, he wrote over 25 non-military papers, mostly on symmetric functions and their properties. Symmetric functions are discussed in the next section.

MacMahon remained an Instructor[68] at the Royal Military Academy until 1888[69]. He was "Inspector Warlike Stores" until 1891, when he took up a post as Instructor in Electricity at the Royal Artillery College. This period of his life is described in the next chapter.

## MacMahon's mathematical work

### *Symmetric functions*

A symmetric function is "an algebraic function of a number of numerical magnitudes ... unaltered when *any* two of the magnitudes are interchanged[70]". For example, the sum c is symmetric, as is

[67] Quoted in Larmor's letter to *The Times*, 31 December 1929.
[68] There are a few references to MacMahon's time as a teacher at Woolwich amongst the obituaries of Fellows of the Royal Society. Col. C. F. Arden-Close, 1865 - 1952, Director General of the Ordance Survey from 1911 - 1922, remembered being taught by MacMahon [de Graaff-Hunter, 1953], as did Colonel G. F. Lenox-Conyngham, 1866 - 1956, Superintendent of the Trigonometrical Survey of India from 1912 - 1920, and founder of the School of Geodesy at Trinity College, Cambridge, in the 1920s.
[69] MacMahon's successor, Mr E. Brooksmith BA LLM, was appointed on 1 March 1889.
[70] This is MacMahon's own definition from *An Introduction to Combinatory Analysis* [MacMahon, 1980].

$\alpha^i\beta^j + \alpha^j\beta^i + \ldots + \mu^i\nu^j + \mu^j\nu^i$. If all the exponents are 1, then the function is called an *elementary* symmetric function.

MacMahon's first paper on the subject was published in 1883[71]; in it he considered a result due to C. H. Prior, where, for $a_1 + a_2 + a_3 = 0$:

$$\left(\frac{a_1}{a_2 - a_3} - \frac{a_2}{a_1 - a_3} + \frac{a_3}{a_1 - a_2}\right) \times \left(\frac{a_2 - a_3}{a_1} - \frac{a_1 - a_3}{a_2} + \frac{a_1 - a_2}{a_3}\right) = 3^2$$

MacMahon generalised the result to:

$$\left[\frac{a_1}{\Delta(a_2, a_3, \ldots a_n)} - \frac{a_2}{\Delta(a_1, a_3, \ldots a_n)} + \ldots\right] \times \left[\frac{\Delta(a_2, a_3, \ldots a_n)}{a_1} - \frac{\Delta(a_1, a_3, \ldots a_n)}{a_2} + \ldots\right] = n^2$$

where[72] $\sum_{i=1}^{n} a_i^j = 0$, for $1 \le j \le n-2$, and $\Delta(y_1, \ldots, y_{n-1}) = \prod_{1 \le i < j \le n-1} (y_i - y_j)$.

In subsequent papers, he developed the use of special differential operators invented by Hammond to expand symmetric functions in terms of elementary symmetric functions. The use of such operators is discussed later.

### *Invariant theory*

Aged only 30, MacMahon made a significant contribution to *invariant theory*, also known as the *algebra of forms*. In England this had been largely the province of two giants of Victorian mathematics, Cayley and Sylvester. The interest of these two important mathematicians in MacMahon's result caused MacMahon to emerge from the obscurity of the Royal Military Academy. As already described in Chapter 1, Sylvester mentioned MacMahon and his work three times during his inaugural speech as Savilian Professor of Geometry[73]. With such an endorsement, MacMahon's acceptance into the wider British mathematical community was assured.

[71] [MacMahon, 1883, [7;1]].
[72] This formulation of the result is from [Andrews, 1977, p. 13].
[73] [Sylvester, 1886]. Sylvester also proposed MacMahon for a Royal Medal in 1885 for this work, but was thwarted by G. G. Stokes. See [Crilly, 2006, p. 550] for details.

Invariant theory began with the discovery by Boole in 1841 of an *invariant* of the *binary quantic*[74] of degree 2 (i.e. a form involving only two variables, each with highest power 2, such as $ax^2 + bxy + cy^2$). This prompted an inquiry by Cayley into the number of invariants that exist for a quantic of a given degree, and methods for their construction. An invariant is a function of the coefficients of the quantic which remains unaltered except by a power of the determinant under a linear transformation of the variables.

For example, if the binary *p*-ic: $(a_0, a_1, a_2, \cdots, a_p)(x, y)^p = a_0 x^p y^0 + a_1 x^{p-1} y^1 + a_2 x^{p-2} y^2 + \cdots + a_p x^0 y^p$ is transformed by:

$$x = lX + mY$$
$$y = l'X + m'Y$$

to become $(A_0, A_1, A_2, \cdots, A_p)(X, Y)^p = A_0 X^p Y^0 + A_1 X^{p-1} Y^1 + A_2 X^{p-2} Y^2 + \cdots + A_p X^0 Y^p$ then an invariant is a function *f* such that $f(A_0, A_1 \cdots A_p) = M^s f(a_0, a_1 \cdots a_p)$, where $M = \begin{vmatrix} l & m \\ l' & m' \end{vmatrix}$.

A covariant is a similarly defined function involving both the coefficients and the variables of the quantic. An important feature of the theory is that, for a particular degree of quantic, some of the invariants are linearly independent whilst others are linear functions, or *syzygies*, of a small *basis set* of invarants.

At the time that MacMahon became interested in algebraic forms, some German mathematicians were engaged in seeking an extension of Gordan's[75] proof of the existence of a finite basis for binary quantics to quantics of higher order. At the same time, the British were busy calculating extensive tables of invariants, and also searching for an alternative proof of Gordan's finite basis theorem using direct algebraic methods, rather than the 'symbolic' method[76]. Crilly[77] has noted the period 1863 - 1895 as one of decline in the development of Cayley's invariant theory, as a result of the superiority of the German methods in proving general results.

[74] This is the name given to the algebraic form we now call an homogeneous polynomial.
[75] P. A. Gordan, 1837 - 1912, studied under Kummer in Berlin and taught at Giessen and Erlangen.

MacMahon's first forays into the world of invariant theory were published in the *Quarterly Journal of Mathematics* ('On Professor Cayley's canonical form'[78]) and the *American Journal of Mathematics* ('Seminvariants and symmetric functions'[79]). In this latter paper, MacMahon drew attention to a one-to-one correspondence[80] between *seminvariants*[81] and the elementary non-unitary symmetric functions related to the roots of the polynomial $x^n - a_1 x^{n-1} + a_2 x^{n-2} - \cdots = 0$. This was an unexpected and important breakthrough which hinted at the possibility of a British triumph over the German efforts.

This paper came to the attention of Sylvester, who wrote to Cayley about it in March 1884[82]. He made some criticisms about the detail of MacMahon's work, but nevertheless acknowledged it as a true demonstration of a result he had hinted at in 1882[83]. Indeed, Sylvester published a short note in the *Messenger of Mathematics* in 1884 entitled *Note on Capt MacMahon's transformation of the Theory of Invariants*, which although staking a prior claim to the idea, nonetheless allowed that MacMahon had made a fundamental discovery in the field. Parshall has pointed out that Sylvester was keen to claim priority because the theorem gave renewed hope that Gordan's Theorem could be proved using his computational methods, rather than by the German symbolic method, which involved replacing a quantic by a single symbol, and investigating the effect of transformations on the form of the symbol, and expressing the in- and covariants in a symbolic form.

---

[76] See Appendix 4 for more detail.
[77] [Crilly, 1988].
[78] [MacMahon, 1883, [8;18]].
[79] [MacMahon, 1884, [13;18]].
[80] This is known as the *Correspondence Theorem*.
[81] A *seminvariant* is the leading coefficient in the expression of a covariant; from it, the entire covariant can be deduced by the use of Cayley's Theorem. See, for example, [Crilly, 1988].
[82] [Parshall, 1998, p. 247].
[83] [Parshall, 1998, p. 249].

To understand the significance of MacMahon's contribution, it is necessary to see how covariants can be calculated from seminvariants. The following explanation shows how the *Hessian*[84] covariant of the binary quartic may be calculated from just one of its coefficients.

The full covariant expression is

$$H \equiv \left(ac-b^2\right)x^4 + 2\left(ad-bc\right)x^3y + \left(ae+2bd-3c^2\right)x^2y^2 + 2\left(be-cd\right)xy^3 + \left(ce-d^2\right)y^4.$$

Generally, the leading coefficient $C_o$ of a covariant is annihilated by a differential operator[85]

$$\Omega \equiv a_0\frac{\partial}{\partial a_1} + 2a_1\frac{\partial}{\partial a_2} + 3a_2\frac{\partial}{\partial a_3} + \cdots + pa_{p-1}\frac{\partial}{\partial a_p}.$$

Given the last coefficient $C_\varpi$ of a covariant, all the coefficients may be calculated by successive applications of $\Omega$; in this example for a binary quartic, the form of the covariant is:

$$\frac{\Omega^4}{4!}C_4x^4 + \frac{\Omega^3}{3!}C_4x^3y + \frac{\Omega^2}{2!}C_4x^2y^2 + \frac{\Omega^1}{1!}C_4xy^3 + C_4y^4.$$

Starting with $C_4 = ce - d^2$ we get successively[86]:

$$C_3 = \frac{\Omega}{1!}C_4 = \frac{\Omega}{1!}\left(ce-d^2\right) = 2\left(be-cd\right)$$

$$C_2 = \frac{\Omega^2}{2!}C_4 = \frac{\Omega}{2}C_3 = \left(ae+2bd-3c^2\right)$$

$$C_1 = \frac{\Omega^3}{3!}C_4 = \frac{\Omega}{3}C_2 = 2\left(ad-bc\right)$$

$$C_0 = \frac{\Omega^4}{4!}C_4 = \frac{\Omega}{4}C_1 = \left(ac-b^2\right)$$

A further application of the operator $\Omega$ would annihilate the expression.

A similar process can be employed by starting with $C_0$ and successively operating with another differential operator[87] $O$ (the final coefficient $C_\varpi$ is annihilated and is sometimes known therefore as the *anti-seminvariant*), or by starting with any coefficient and using both operators. Thus, the problem of generating covariants can be reduced to that of finding seminvariants. To illustrate this point in the context of the binary quartic, we note that there are three possible seminvariants of degree 2 for *weight*[88] 4 or less:

$$a_0^2 \equiv a^2$$
$$a_0a_2 - a_1^2 \equiv ac - b^2$$
$$a_0a_4 - 4a_1a_3 + 3a_2^2 \equiv ae - 4bd + 3c^2$$

[84] It is not necessary to know what a Hessian is to follow this explanation, but more detail is given in Appendix 4.
[85] See Appendix 4 for more detail of the operator $\Omega$. The notation is due to Cayley.
[86] See [Elliott 1895, pp. 128 - 140] for a more detailed general argument.
[87] See Appendix 4 for more detail of the operator $O$.

Taking each of these as the leading coefficient $C_0$ of a covariant, we can generate the covariants using the operator $\mathrm{O}$. Applied to each of the possible seminvariants of degree 2, the first is not annihilated at all and the third is annihilated immediately (and so is an invariant); only the second, $ac - b^2$, generates the Hessian covariant:

$$\mathrm{O}\left(ac - b^2\right) = 2\left(ad - bc\right)$$
$$\mathrm{O}^2\left(ac - b^2\right) = ae + 2bd - 3c^2$$
$$\mathrm{O}^3\left(ac - b^2\right) = 2\left(be - cd\right)$$
$$\mathrm{O}^4\left(ac - b^2\right) = ce - d^2$$
$$\mathrm{O}^5\left(ac - b^2\right) = 0$$

To explain MacMahon's discovery, we need to see the connection between the elementary symmetric functions and the roots of a binary quantic. Again, the quartic is used as an illustration.

If we take the roots of the quartic to be $\alpha, \beta, \gamma, \delta$, and expand the equation $\left(x - \alpha y\right)\left(x - \beta y\right)\left(x - \gamma y\right)\left(x - \delta y\right) = 0$ then the result is:

$$x^4 - \left(\alpha + \beta + \gamma + \delta\right)x^3y + \left(\alpha\beta + \alpha\gamma + \alpha\delta + \beta\gamma + \beta\delta + \gamma\delta\right)x^2y^2$$
$$- \left(\alpha\beta\gamma + \alpha\beta\delta + \alpha\gamma\delta + \beta\gamma\delta\right)xy^3 + \left(\alpha\beta\gamma\delta\right)y^4 = 0$$

From this we can see that the coefficients of the quartic $\left(a_0, a_1, a_2, a_3, a_4\right)\left(x, y\right)^4$ can be expressed as elementary *unitary* symmetric functions; that is, symmetric functions where no symbol has a power greater than 1. The notation $(1^n)$ denotes the symmetric function $\sum \alpha_1\alpha_2\ldots\alpha_n$.

[88] The *weight* is the sum of the coefficient suffices in each term.

$$\begin{aligned}
a_0 &= 1\\
a_1 &= -(\alpha+\beta+\gamma+\delta) = (1)\\
a_2 &= (\alpha\beta+\alpha\gamma+\alpha\delta+\beta\gamma+\beta\delta+\gamma\delta) = (1^2)\\
a_3 &= -(\alpha\beta\gamma+\alpha\beta\delta+\alpha\gamma\delta+\beta\gamma\delta) = (1^3)\\
a_4 &= (\alpha\beta\gamma\delta) = (1^4)
\end{aligned}$$

Hammond[89] had described the solutions to the differential equation $\Omega z = 0$ (one of the conditions for invariance)[90] in terms of the following *sources* (or *protomorphs*)

$$\begin{aligned}
U &= a_0\\
H &= a_0a_2 - a_1^2\\
C_3 &= a_0^2a_3 - 3a_0a_1a_2 + 2a_1^3\\
Q_4 &= a_0a_4 - 4a_1a_3 + 3a_2^2\\
&\vdots
\end{aligned}$$

Hammond was working with quantics with binomial coefficients, but MacMahon chose to work with the binary quantic $a_0x^n + na_1x^{n-1}y + n(n-1)a_2x^{n-2}y^2 + \cdots$, where the coefficients are derived, so his list looks slightly different (although he provided a simple method for converting between the forms):

$$\begin{aligned}
q_2 &= 2a_2 - a_1^2\\
c_3 &= 3a_3 - 3a_1a_2 + a_1^3\\
q_4 &= 2a_4 - 2a_1a_3 + a_2^2\\
&\vdots\\
q_{2m} &= 2a_{2m} - 2a_1a_{2m-1} + 2a_2a_{2m-2} - \cdots + (-)^{\frac{1}{2}n}a_m^2\\
c_{2m+1} &= (2m+1)a_{2m+1} - (2m+1)a_1a_{2m} + (2m-3)a_2a_{2m-1} - (2m-5)a_3a_{2m-2} + \cdots\\
&\quad + a_1\left\{2a_1a_{2m-1} - 2a_2a_{2m-2} + 2a_3a_{2m-3} - \cdots - (-)^{\frac{1}{2}n}a_m^2\right\}
\end{aligned}$$

Any solution of the differential equation is a function of these sources (also called *ground sources*). All the sources with an even suffix are of degree 2, and are known as *quadrinvariants* (hence the symbol *Q*). These are not the only possible ground sources, but they are of lower degree than any other set[91].

[89] James Hammond, 1850 - 1930, was 35th Wrangler in 1874, a teacher in Londonderry and Bath, a council member of the LMS from 1883 - 1892, and private secretary to James Joseph Sylvester from 1883 - 1894.
[90] [Hammond, 1882].

So the source of a seminvariant is a function of the ground sources. Suppose we want a seminvariant of weight 6. The non-unitary partitions of 6 are (6), (2, 4), (3, 3) and (2, 2, 2) which correspond to the ground sources $Q_6$, $HQ_4$, $C_3^2$ and $H^3$. The seminvariant source $S$ is a function of these ground sources. Since they have degrees 2, 4, 6 and 6, we multiply as necessary by powers of $a_0$ (denoted by Hammond as $U$) so that the form of $S$ is

$$SU^{6-\mu} = \alpha_1 H^3 + \alpha_2 C_3^2 + \alpha_3 U^2 HQ_4 + \alpha_4 Q_6$$

where $\mu$ is the degree of $S$, the terms of the equation all have the same degree and weight and the $\alpha_i$ are arbitrary constants. Once a source seminvariant is known, an entire covariant can be calculated as described above.

A *syzygant* is a seminvariant formed from the ground sources which is itself a ground source, not in the base set, and results from appropriate choices of the arbitrary constants. The simplest example of a syzygant is of weight 6. For the cubic this is $U^2\Delta = 4H^3 + C_3^2$ where $\Delta = \left(a_0a_3 - a_1a_2\right)^2 - 4\left(a_0a_2 - a_1^2\right)\left(a_1a_3 - a_2^2\right)$, the discriminant of the cubic. This gives the ground source seminvariant $-3a_0a_1^2a_2^2 + 4a_0^3a_2^3 + 4a_0^2a_1^3a_3 - 6a_0^3a_1a_2a_3 + a_0^4a_3^2$ which by repeated application of the operator O gives a covariant of degree-order (6.5).

Another example, from which the cubic syzygant above is deduced, may be determined as follows. We choose $\mu = 3$, so we are seeking a set of constants $\alpha_i$ such that $S$ in $SU^3 = \alpha_1 H^3 + \alpha_2 C_3^2 + \alpha_3 U^2 HQ_4 + \alpha_4 Q_6$ is itself a ground source invariant distinct from the basis set above.

[91] [Hammond, 1882, p. 220]

Since *S* has to be distinct from any of the ground sources, we choose $\alpha_4 = 0$. Multiplying out the right-hand side and dividing through by $a_0^3$ gives three conditions to be satisfied: $\alpha_1 - 4\alpha_2 = 0,\ \ \alpha_2 + \alpha_3 = 0,\ \ \alpha_1 - 3\alpha_2 + \alpha_3 = 0$. We can choose $\alpha_1 = -4, \alpha_2 = -1, \alpha_3 = 1$, so that $S = -4H^3 - C_3^2 + U^2 HQ_4 = -a_1^2 a_4 + 2a_1 a_2 a_3 - a_2^3 + a_0 a_2 a_4 - a_0 a_3^2$ which is equivalent to the invariant *J* of the quartic described in Appendix 4.

Hammond's paper concluded with a table of all the ground sources up to weight 10, including the syzygants calculated as above, showing the degree-order of the covariants which each ground source would generate for quantics of any order. Part of the table was copied from tables published by Sylvester[92]. By explicitly calculating the ground source and the syzygant, Hammond showed that a postulate of Sylvester's, that a quantic of order 7 (*septimic*) has no ground source nor syzygant of degree-order (5.13), is in fact false.

MacMahon pointed out that each of the ground sources corresponds to a non-unitary partition function. For example, $q_2$ for the quartic becomes:

$$\begin{aligned} 2a_2 - a_1^2 &= 2(\alpha\beta + \alpha\gamma + \alpha\delta + \beta\gamma + \beta\delta + \gamma\delta) - \left(-(\alpha + \beta + \gamma + \delta)\right)^2 \\ &= -(\alpha^2 + \beta^2 + \gamma^2 + \delta^2) = -(2) \end{aligned}$$

- that is, the sum of the squares of the roots, which is a non-unitary symmetric function of the roots. He deduced that every seminvariant corresponds to a non-unitary-partition function of the roots of $(1, a_1, a_2, \ldots)(x, 1)^n = 0$, using derived coefficients, rather than binomial, to simplify some of the calculations. MacMahon then deduced the important result which interested Sylvester, that the number of linearly independent seminvariants of degree-weight (*j,w*) is "equal to the number of non-unitary partitions of *w* which contain *j* and no part greater than *j*".

[92] Hammond remarked on p. 227 that all of Sylvester's results would need a table that extends up to weight 85.

To summarise, the key to generating covariants and syzygies was finding seminvariants, and linearly independent seminvariants are in one-to-one correspondence with non-unitary partitions.

The German mathematicians R. F. A Clebsch[93] and P. Gordan had developed a method for generating covariants (of which invariants are a special case), but were unable to locate the syzygies between them, and so overestimated the size of the basis set. Sylvester had hoped that his 1876 work would get around the problem by constructing the set algorithmically, but his method underestimated the size of the basis set. MacMahon's result gave renewed hope to Sylvester by providing another avenue to explore.

***Partition theory***[94]

The consideration of symmetric functions and their relationship to partitions led MacMahon to study partitions, the field in which he did his most original and creative work. Whilst at Woolwich he published two major papers and one minor paper related to partitions.

If $N$ is a sum of positive integers, then that sum is called a *partition* of $N$. For example, if $N = 4$, then the possible partitions are 4, 3 + 1, 2 + 2, 2 + 1 + 1 and 1 + 1 + 1 + 1. The order of the parts is not important (so that 3 + 1 is the same partition as 1 + 3), and it is customary to write the parts in non-ascending order, in brackets, thus: (4), (31), (22) or ($2^2$), (211) or ($21^2$), and (1111) or ($1^4$). The basic problem is to enumerate the partitions of a given integer $N$. Related problems include the enumeration of partitions subject to certain restrictions - for instance, where no part is to be greater than some given integer, or the number of parts is restricted, or only certain types of integer may be used (even numbers, primes, multiples of 3, etc.). By the time MacMahon came to the study of partitions in the mid-1880s, it was already a well-established field. The question of calculating a value for $p(n)$ for any $n$ had first been raised by Leibniz in a letter to Johann Bernoulli in 1669 [95] Euler, prompted by Naudé, worked on various aspects of the enumeration problem during the 1740s. I do not intend to give a complete history of the topic, although a description of the early English work is given in Appendix 5, to show the tradition from which MacMahon came.

[93] R. F. A. Clebsch, 1833 - 1872, studied at Königsberg, and taught at Giessen and Göttingen.
[94] A complete list of MacMahon's papers dealing with partition theory is given in Appendix 5.
[95] [Dickson,1992].

Amongst those working on partition problems in the middle of the 1880s, the most eminent in England was J. J. Sylvester. In 1886 he published in *Mathematical Questions with their solutions from the Educational Times*[96] a list of unsolved questions, including five on partitions, one of which is reproduced below to illustrate the style.

3535. 1. Euler has shown that the number of modes of composing $n$ with $i$ distinct numbers is equal to the denumerant (that is, the number of solutions in positive integers, zeros included) of the equation

$$x + 2y + 3z + \cdots + i\omega = n - \tfrac{1}{2}\left(i^2 + i\right)$$

Show more generally that the number of modes f composing n with i numbers, alll distinct except the largest, which is always to be taken $j$ times, is the denumerant of the equation

$$jx + (j+1)y + \cdots + i\omega = n - \tfrac{1}{2}(i - j + 1)(i + j)$$

2. Show also that, if all the partitions of $n$ into $i$ parts are distinct except the *least*, which is to be taken $j$ times, then the number of such partitions is the denumerant of the equation

$$x + 2y + \cdots + (i - j)\phi + \omega = n - \tfrac{1}{2}(i - j)^2 + 3i - 1$$

The questions show that, to a large extent, research in the subject involved ever more complicated conditions on the type of partition to be enumerated. Many authors found formulas for particular sets of conditions, often by creating recurrence formulas to carry out the enumeration. MacMahon's unique contribution was to generalise the idea of a *partition* to that of a *distribution*, so that he could use the calculus of symmetric functions to solve the enumeration problem. Partitions were then a special case of this more general method. MacMahon was to become the most prolific writer on the subject in the period from 1886 to 1918. Dickson[97] cites the names of 66 mathematicians who wrote papers on partitions. Of those, the majority are mentioned only once; 8 are mentioned twice, three are mentioned three times and one, R. D. Sterneck, is mentioned six times. MacMahon merits 19 citations. Until around 1902, the writers of papers on partitions were principally from France, Germany or Great Britain, but after that date, some Italian mathematicians also published work on partitions.

[96] [Sylvester, 1886].
[97] [Dickson,1992].

In his first paper[98], *Certain Special Partitions of Numbers,* MacMahon introduced the novel concept of a *perfect partition.* This is a partition which contains exactly one partition of every lower number. For example, the partition (421) of 7 contains the partitions (42) of 6, (41) of 5, (4) of 4, (21) of 3, (2) of 2 and (1) of 1.

These partitions were introduced in the context of producing sets of weights so that any integral weight from 1 lb. to $u$ lb. can be made in a unique way. Perfect partitions of $u$ are needed to achieve this if weights may only be used on one side of a balance. For example, if $u$ = 7lb., then a set of weights of 4lb., 2lb. and 1lb. is needed, corresponding to the perfect partition explained in the previous paragraph. If weights may be placed on both sides of a balance[99], then a *subperfect* partition suffices; this is a partition in which the parts are allowed to be negative. This context presages his later work for the Board of Trade (see Chapter 6). This theme is picked up again in a later paper, *Weighing by a series of weights*[100].

MacMahon developed methods for counting and producing perfect partitions of numbers. First he defined $\phi_{p.q} = 1 + x^q + x^{2q} + \ldots + x^{pq}$ and noted that for $q$ = 1, $\phi_{p.1} = 1 + x + x^2 + \ldots + x^p = \sum_{k=0}^{k=p} x^p = \frac{1-x^{p+1}}{1-x}$.

If $p$ + 1 is not prime, then $\phi_{p.1}$ may be factorised as $\phi_{p.1} = \phi_{l.\lambda} + \phi_{m.\mu} + \phi_{n.\nu} + \ldots$, and $\left(\lambda^l \mu^m \nu^n \cdots\right)$ is a perfect partition of $p$. The Greek letters represent the parts, and Roman letters the repetitions of those parts in the partition of $p$.

MacMahon's ability to invent new ideas to generalise existing concepts and theorems characterises his work in partitions over the following four decades. Perfect partitions were also the subject of a paper in 1891, *The theory of perfect partitions and the compositions of multipartite numbers*[101].

[98] [MacMahon,1886, [20;6]]

[99] It is interesting to note that in this paper MacMahon used the phrase "both scale pans", but in the equivalent chapter of *Combinatory Analysis*, written nearly 30 years later and after he had been working as Deputy Warden of the Standards for almost a decade, he used the phrase "weights pan and goods pan", reflecting a professional familiarity with the technical terms of the scale trade.

[100] [MacMahon,1890, [37;19]].

As described above, a syzygy is a linear relationship between invariants, and an important part of the search for complete lists of invariants was the need to determine the basis set of independent invariants from which the entire set could be created. In his second partition theory paper, *The expressions of syzygies amongst perpetuants by means of partitions*[102], MacMahon brought together non-unitary partitions, monomial symmetric functions and a process in invariant theory known as 'transvection' (at the time still known by its German name of *Überschiebung*, meaning 'slide over', because two algebraic forms are combined to form a third) to make explicit the connection between invariant theory and partition theory.

A minor paper, *On play à outrance*[103] was only peripherally concerned with partition theory; in this work, MacMahon considered the probabilities of the outcomes of a game for two players with a total of $n$ counters between them. However, it does show that MacMahon had an early interest in games and recreational mathematics, a topic which is discussed in later chapters, as is MacMahon's later work on partition theory.

**Personal events**

The following outline of MacMahon's first marriage and its subsequent break-up is included because it shows that there was more to his life than the mathematics which has come down to us. It is easy to forget that the creators of mathematics were human beings with many facets to their character. It is not clear how the events described might have affected MacMahon's work, nor whether they were a consequence of that work.

101 [MacMahon, 1891,[38;6]].
102 [MacMahon, 1888, [27;18]].
103 [MacMahon,1889, [33;10]].

On 19 April 1881, in St Mark's Church, Lewisham, MacMahon married a 21 year old American girl, Aimee Rose Leese, daughter of Robert Henry Leese[104], a 70 year old retired US Consul (although he had been born in Portsmouth, England) who lived with his family in Kent. The day after the wedding, Robert died. On 11 October 1882, Aimee Rose gave birth to a daughter, Florence Aimee[105]. On 18 January 1886 MacMahon filed for a divorce from Aimee. Aimee[106] and Florence had moved to the United States, to a town called Lariat (known as Henry between 1884 and 1886 and is now called Monte Vista) in Rio Grande County, Colorado, where her elder brother, Percy Henry Leese, was living with his family. Aimee claimed that she had obtained a divorce from MacMahon under Colorado State law on 8 November 1884[107]. She later married Edward Bevan in September 1885. MacMahon had not recognised the Colorado divorce and his 1886 petition cited bigamy and adultery as grounds for a UK divorce. He also claimed £5000 damages and asked for custody of baby Florence. After two years of transatlantic legal wrangling[108], MacMahon was granted the divorce he sought. He received no compensation, and Florence stayed with her mother and stepfather in the United States[109] .

In 1893, Aimee and her mother were in Chicago[110], on their way to California[111]. A letter written to Percy Henry Leese suggests they were going to visit him in New Mexico before completing the journey to California to settle. There is no mention of Florence in the letter, so she may already have gone to California to start school. Certainly Florence was in California by 1901[112]. It is not certain whether she married or not[113]. She had no children, so there are no direct descendants of MacMahon today.

---

[104] Robert Henry Leese, 1811 - 1881.
[105] An announcement appeared in *The Times* on 14 October 1882.
[106] It is likely that Aimee's mother, Catherine, also moved to Colorado at the same time. Indeed, it is possible that Aimee's original intention was simply to acccompany her mother on the journey and not to stay away permanently.
[107] The state of Colorado did not start keeping records of divorces in a systematic way until 1890, so it is not possible to verify this claim independently, nor to discover the grounds used by Aimee for the divorce.
[108] Bevan was not happy about being cited as co-respondent when he was legally married under Colorado law, and certainly saw no reason to pay damages.
[109] She is listed on the school census for Monte Vista in 1889.
[110] Possibly to visit the 1893 Columbian exhibition held in Chicago between 1 May 1893 and 31 October 1893.
[111] Private communication from Linda Cates, Percy Henry Leese's great granddaughter.
[112] She died on 7 May 1961 in Los Angeles. Her death certificate lists her profession as a dentist.
[113] In the US census of 1930, Florence was living in Oakland, California, and she described herself as a widow. She gave her occupation as stewardess on a steamship! It seems that her life may have been quite eventful, but documentary evidence is very sparse. Enquiries made at Florence's last known address in the 1970s and again in 2002 produced no information.

# Chapter 4 Royal Artillery College 1890 - 1898

This chapter discusses the period from MacMahon's election to the Royal Society to his retirement from the Royal Artillery, during which he consolidated his reputation as a first class mathematician and laid the foundation for his most important work.

### Royal Artillery College

In 1891 MacMahon took up a new post as Military Instructor in Electricity[114] at the Royal Artillery College, Woolwich[115]. This institution arose out of the Advanced Class for Artillery Officers in Mathematics[116] that MacMahon had attended just ten years previously. No surviving records shed light on the work he was required to do in this post, or on any students he may have taught.

MacMahon continued to work on military problems and to take part in Royal Artillery Institution meetings, but to a lesser extent, since his mind was now occupied with invariant theory and partitions. In 1894 he published a short paper entitled 'Note on the correction of artillery fire' in the *Proceedings of the Royal Artillery Institution*, but made no further military contributions to the journal before he retired from the Army, although on 7 February 1895 he gave a lecture on *Terrestrial refraction and mirages* to the Royal Artillery Institution[117]. MacMahon's part in the discussion following a talk on eclipses by Captain Hills[118] in 1897 is his last recorded contribution to the Royal Artillery Institution.

---

[114] Some sources (e.g. the three obituarists) have said that this post was 'Professor of Physics', but this is not correct, as Greenhill held that post until his retirement in 1908.

[115] MacMahon's friend and collaborator Major J. R. J. Jocelyn had previously held the post.

[116] This class was introduced in 1864 under the Director of Artillery Studies at Woolwich. Until 1885, instruction took place at the premises of the Royal Artillery Institution, when it was moved to the 'Red Barracks' at Woolwich. At the same time, the class was thrown open to the Army and Royal Marines with the new name "Senior Class Artillery College". The College was renamed the "Ordnance College" in 1899.

[117] [MacMahon, 1895]

[118] Captain, later Major, Edmund Herbert Hills, the treasurer of RAS from 1904, was educated at Winchester and the RMA. He entered the army in 1884, was appointed Instructor at the School of Military Engineering in 1898; and contested Portsmouth at the 1905 general election for the Conservative party [source: *Irish Times*, 2 September 1908].

**The Royal Society**[119]

In 1890 MacMahon was elected a Fellow of the Royal Society. The proposers were: "from general knowledge - A. J. Ellis, J. Henry Lefroy and J. T. Walker" and "from personal knowledge - A Cayley, J. J. Sylvester, Percival Frost, A. G. Greenhill, J. W. L. Glaisher[120], A. R. Forsyth, James Cockle[121], J. J. Walker and M. W. Crofton." The citation in support of the application for Fellowship of the Society was first made on 26 February 1889 and then 'resuspended'[122] in 1890, and read:

> As author of numerous papers in the Quarterly Journal of Mathematics, vols XIX - XXI, Proceedings of the London Mathematical Society vols XV - XIX, American Journal of Mathematics vols VI - XI on various subjects in Pure Mathematics, connected with Invariants, Seminvariants, Perpetuants, Reciprocants, Partitions, Distributions and Symmetric Functions. Associate Member of the Ordnance Committee. Instructor in Mathematics at the Royal Military Academy, Woolwich, 1882 - 1888.

The election was confirmed on 5 June 1890. This recognition of his contribution to mathematics firmly established him in the British mathematical community.

As further evidence of his close association with the mathematical community (and to some extent the entire scientific community), we note that MacMahon personally signed the proposal forms for 47 applicants for membership of the Royal Society during his lifetime, including William Burnside, Karl Pearson, G. H. Hardy, Arthur Eddington and S. Ramanujan (a complete list is given in Appendix 8). This does not prove that he knew all of them personally, but indicates that he was at least familiar with their work.

MacMahon published five papers in the *Proceedings of the Royal Society*, and twelve in the *Philosophical Transactions*, including, in 1894, the paper that described his 'Master Theorem'

[119] The Royal Society has some original letters from MacMahon which are not mathematical in nature, as well as various referees' reports on some of his papers. These are listed in Appendix 8. MacMahon also refereed a Royal Society paper by R. F. Gwyther on *Differential Coefficients of Twisted Curves* (20 October 1894), two papers by Karl Pearson, and one by A. R. Forsyth.

[120] J. W. L. Glaisher, 1848 - 1928, was second Wrangler in 1871, elected FRS in 1875, and was President of the LMS from 1884 - 1886.

[121] James Cockle, 1819 - 1895, was elected FRS in 1865, knighted in 1869, and was President of the LMS from 1886 - 1888.

[122] Resuspension refers to the procedure for resubmitting a proposal for election.

(described later).

### The London Mathematical Society

Further recognition of his status as a mathematician came in 1894 - 1896, when MacMahon became President of the London Mathematical Society (LMS). Two years was, and still is, the standard period for a person to hold the Presidency of the Society.

MacMahon had joined the London Mathematical Society on 14 December 1883, whilst still an Instructor at Woolwich, and had been a Council member since 1890. E. W. Hobson[123] and A. E. H. Love[124] proposed him for the presidency on 11 October 1894. The first President after the incorporation[125] of the LMS on 23 October 1894, and successor to A. B. Kempe, at 40 years of age, MacMahon was the third youngest person ever to hold the post of President. Only Lord Rayleigh[126] (1876 - 78) and J. W. L. Glaisher (1884 - 1886) had been younger, at ages 34 and 36 respectively[127] .

When his term of office came to an end, he delivered his valedictory President's address on 12 November 1896, before an audience of 24 members and one visitor. Entitled *Combinatory Analysis: A review of the present state of knowledge*[128], it was in this speech that the foundation for his greatest work, the two volume *Combinatory Analysis*, was set.

#### *1896 Presidential address to the London Mathematical Society*

MacMahon's address was in three parts: an introductory section with some history of Cayley's work on invariant theory, a middle section concerning the place of combinatory analysis in classification schemes, and a final part on the state of partition theory in 1896. There are some splendid phrases dotted about the text (some of which are quoted in the following paragraphs), showing MacMahon's skill as an orator and communicator. This outline of the main points of the speech also demonstrates MacMahon's understanding of the importance of the role of history in understanding the development of a mathematical subject.

[123] E. W. Hobson, 1856 - 1933, was senior Wrangler in 1878, and was Sadleirian professor at Cambridge from 1910. He was president of the LMS 1900 - 1902, and was awarded the De Morgan Medal in 1920.
[124] Augustus Edward Hough Love, 1863 - 1940, held the Sedleian Chair of Natural Philosophy, Oxford, from 1899.
[125] "Incorporation" is a legal process whereby an organisation becomes a legal entity and can own property in its own right.
[126] John William Strutt, 3rd Baron Rayleigh, 1842 - 1919, was senior Wrangler 1865, and elected FRS in 1873.
[127] No President since MacMahon has been younger than 41 at the time of election, the median age being 51.
[128] [MacMahon, 1897, [49;19]].

The address began with a tribute to Cayley, "in consequence of his exceeding kindness to me at the time, some fourteen years ago, when I first placed my foot upon the thorny path of mathematical research." In his description of Cayley's work in invariant theory, MacMahon referred to his own contribution in 1884[129], although he was too modest to claim the credit for it. He described Cayley's dissatisfaction with the purely symbolic work of the German mathematicians Gordan and Hilbert[130], but noted that "the symbolism employed is essentially that initiated by Cayley, who had little fancy for the method, which, dropping from his hands, was carried on with great results by the mathematicians of Germany."[131]

The second part of the address is arguably the most important. Here, MacMahon asserted that the areas of mathematics ascribed to combinatory analysis were incorrect. He used the classification published in Paris in 1893 by *La commission permanente du répertoire* as a typical example[132]. In it, combinatory analysis was shown as including permutations and combinations, but partition theory was classed as part of number theory, and symmetric functions as part of the theory of equations. This latter classification MacMahon claimed was totally without foundation, and described it as "the thing itself confounded with one of its applications." MacMahon was concerned not just with proving results, but also with providing the correct framework for mathematical research.

[129] *Seminvariants and symmetric functions* [MacMahon, 1884, [13;18]] .
[130] David Hilbert, 1862 - 1943, taught at Göttingen from 1895.
[131] MacMahon also referred to several letters from Cayley which he had in his possession; these are now lost.
[132] The text he quoted was:
Classe J1 Analyse combinatoire, with the subdivision
- a. Groupes où l'on tient de l'ordre:
  - a. Permutations et arrangements de toutes sortes, nombre et loi de formation;
  - b. Structure, sequences, &c.,
  - c. Manière de ramener les permutations les unes aux autres par permutation circulaire de groups d'elements, permutations, semblables, &c.
- b. Groupes ou l'ordre est indifferent:
  - a. Combinaisons simples at avec repetition;
  - b. Combinaisons regulieres, completes.
- c. Fonctions qui presentent dans l'analyse combinatoire, theorie des chémins.
- d. Applications:
  - a. A la recherche des termes généraux des series;
  - b. Des tables;
  - c. Des coefficients de certaine developpements.

MacMahon then drew attention to some of his own work in the field of symmetric functions, including his Master Theorem[133] and his extension of symmetric functions[134] to those with negative exponents, that is, in the general symmetric function $\sum \alpha^p \beta^q \gamma^r \ldots$, the exponents $p$, $q$ and $r$ may be positive, zero or negative, and still be interpreted sensibly in terms of partitions. This is accomplished by using a generalisation of the concept of distribution, an idea central to MacMahon's work on partition theory. Briefly put, a partition $(p_1 p_2 p_3 \ldots p_m)$ of $n$ may be considered as a distribution of $n$ similar objects into $m$ differently sized parcels, or alternatively as $m$ different objects. For example, the partition (421) of 7 may be thought of as seven 1s in three parcels of size 4, 2 and 1 respectively, or as three different numbers[135], 4, 2 and 1. The generalisation that MacMahon made was to consider distributions into two groups ('upper' and 'lower'), where one is subtracted from the other. The generalised subject of distributions was, he said, "virgin soil laden with germinative seeds which only await the application of the husbandman."

The third section of the address dealt in considerable depth with partition theory, and began with a detailed history of the subject - necessary, MacMahon claimed, in order fully to appreciate the advances that had been made. The true story of partitions began in 1750 with Euler's use of generating functions, MacMahon asserted, but little progress was made on the Continent by Paoli, Lacroix,

[133] This result in partition theory is described later in this chapter.
[134] [MacMahon, 1889, [31:2]] and [MacMahon, 1890, [35:2]].
[135] That is, three different kinds of object, 4 of one kind, 2 of a second kind, and 1 of a third kind.

Legendre and Jacobi[136]. In Britain, interest in partitions was initiated by De Morgan[137] in a letter to Henry Warburton; this event and its consequences are described in Appendix 5.

Sylvester and Cayley became involved in partition theory in 1855, and MacMahon described Sylvester's outline of his researches, published in the *Quarterly Journal of Mathematics* in 1855, as "incomparably the finest contribution that has ever been made to combinatory analysis."[138]

After a summary of Sylvester's analytical work on partition theory, MacMahon turned to the diagrammatic approach introduced by Sylvester, which MacMahon himself found both attractive and fruitful. He described his own extensions of the technique, "for the most part arithmetical and intuitive, and not algebraical." This led to a discussion of the possibility of constructing an analytical theory to enumerate the partitions of multipartite numbers, which Euler had formulated in the bipartite case as concerned with the fraction $\left\{\prod\left(1-x^{\alpha}y^{\beta}\right)\right\}^{-1}$, traditionally known as the "Problem of the Virgins"[139]. Sylvester had claimed to have a complete solution to the problem, but MacMahon expressed some astonishment that this solution had never been published, perhaps due to some "deplorable mischance." From a perusal of the printed outlines of Sylvester's 1859 series of lectures on partition theory[140] and Cayley's paper[141] on the subject of the special problem of double partitions, MacMahon hoped to reconstruct Sylvester's line of thought, since he believed that the rediscovery of these lost theorems of Sylvester was a matter of great importance to combinatory analysis.

The talk concluded with the observation that "the subject of combinatory analysis has for forty years been much neglected", and with the following quotation from Sylvester, whom he described as the "greatest living English master of pure mathematics."

[136] Pietro Paoli (1759 - 1839) in *Opuscula analytica*, 1780 - reference in L. E. Dickson's *History of the theory of numbers*, Volume 2 p. 107. Paoli was Professor of Mathematics at the University of Pisa; S. F. Lacroix, 1765 - 1843; A. M. Legendre, 1752 - 1833; C. G. J. Jacobi, 1804 - 1851.
[137] Augustus De Morgan, 1806 - 1871.
[138] A full account of this outline by Sylvester, which included results but no detailed explanations, was not published until 1882 [Sylvester, 1882]. The work is very technical and a full description is beyond the scope of this thesis, although some details are given in Appendix 5.
[139] Sylvester [Sylvester, 1858, pp 113 - 117] confessed that he had no idea why the solution in integers of two simultaneous equations with an indefinite number of variables should be called 'the problem of the virgins.'
[140] [Sylvester, 1973, pp 119 - 175].
[141] [Cayley, 1860].

> Partitions constitute the sphere in which analysis lives, moves and has its being; and no power of language can exaggerate or paint too forcibly the importance of this till recently almost neglected, but vast, subtle and universally permeating, element of algebraical thought and expression.

MacMahon's presidential address shows his skill as a speaker in drawing together different threads to form a coherent whole, his mastery of the subject matter, his ability to summarise difficult topics in a succinct way, and his passion for mathematics. He could deal with both the fine technical detail of the area of research, as well as having a clear view of its place within, and value to, mathematics as a whole.

There is no evidence that the address provoked any reaction on the part of the audience or on subsequent readers. It was at the time atypical in dealing with a specific mathematical topic: outgoing presidents usually chose to address general questions about the nature of mathematics or the state of mathematics education. For example, the immediate past president, A. B. Kempe, had spoken on the question 'What is mathematical knowledge?'[142], and before that, Greenhill had talked about collaboration in mathematics[143]. The following president, E. B. Elliott, spoke about the importance of problem making and solving in continuing mathematical education[144] . But the next holders of the office, Lord Kelvin and E. W. Hobson, both spoke about specific mathematical topics (the transmission of force through a solid[145] and the role of the infinite and the infinitesmal in mathematical analysis[146]). So MacMahon may have set the trend for presidential addresses..

MacMahon continued to serve on the Council of the LMS for many years, and was Vice-President in 1896-98, 1910-11, 1915-16 and 1919-20. During this last term of office, MacMahon had a disagreement with the Council over the involvement of the LMS in the formation of the International

[142] [Kempe, 1894].
[143] [Greenhill, 1892].
[144] [Elliott, 1898].
[145] This talk was never published.
[146] [Hobson, 1902].

Mathematical Union, which appears to have led to his resignation from the Council. The minutes of the LMS Council record the events thus:

> Oct. 26 1919…
> The letter from the International Research Council, read at the meeting of 12 June 1919, was further considered. It was proposed by Prof. Macdonald, seconded by Mr. Littlewood, and carried nemine contradicente, that
> ‘The Council of the London Mathematical Society is not prepared to make definite proposals as to the formation of an International Union of Mathematical Societies.’
>
> Dec 11 1919…
> A letter from the Secretary of the International Research Council, asking the Society to send 4 additional delegates to a meeting of the Conjoint Board on Jan 8, 1920, and other letters, were read.
> Mr. Hardy (seconded by Prof. Watson) moved the following resolution :–
> ‘The Council of the L.M.S. does not consider the organization of an International Union of Mathematical Societies to be of urgent importance in the interests of mathematical research; and is of the opinion that the formation of such a Union should be postponed until such time as it can be constituted on a fully international basis.’
> Prof. Love (seconded by Major MacMahon) moved as an amendment :–
> To delete all words after ‘research’, and to substitute: ‘but is of the opinion that if such a Union is formed the Society should join it’.
> After considerable discussion the amendment was lost by 4 votes to 7.
> The President (seconded by Prof. Watson) proposed as an amendment :–
> To omit all words after ‘postponed’ and this amendment was accepted by the mover of the original resolution. The resolution as amended was then put to the meeting and carried by 9 votes to 3. Major MacMahon, Prof. Love and Prof. Filon voted against the resolution. It was agreed that additional delegates should be sent; and the President, the Treasurer, Prof. Eddington & Prof. Love were appointed.
>
> 15 Jan 1920…
> Major MacMahon resigned his membership of the Council as from 11 Dec 1919. It was agreed unanimously that ‘The Council of the L.M.S. learns with extreme regret that Major MacMahon proposes to resign his membership of the Council and hopes that he will reconsider the matter’.

MacMahon did not reconsider and is the only Council member in the history of the LMS to have resigned in mid-term. It seems an extreme action to have taken in the circumstances. Hardy had clearly felt that if a truly international union, that is, one that included the Axis powers, could not be formed immediately, then the matter should wait until Germany and her allies had been restored to favour. MacMahon seems to havce adopted a more relaxed approach, and was happy to accept a less rigorous interpretation of the meaning of 'international'. Perhaps, in his 65th year, MacMahon felt that the much younger Hardy was exercising too much influence, and not according sufficient respect to the views of a former LMS president and 'grand old man' of mathematics. It is worth noting that a mere two years after this episode, MacMahon abandoned London altogether for the quieter life of a Cambridge academic.

### The Royal Astronomical Society[147]

The President of the Royal Astronomical Society, A. A. Common[148], proposed on 13 November 1896 that MacMahon be elected a Fellow of that Society. The Election Certificate was also signed by J. W. L. Glaisher and Sir Robert Stawell Ball. The election was ratified by the Council on 8 January 1897 and MacMahon signed his acceptance and subscription undertaking on 16 January 1897.

MacMahon's involvement with observational astronomy was limited. As a young soldier in India he had witnessed a total eclipse in 1874, and he accompanied Common on a trip to Norway to see a total eclipse on 9 August 1898, but was thwarted by the weather a few minutes before totality[149]. MacMahon probably joined the RAS because it was another club where Victorian gentlemen of science could meet socially, rather than because he wanted or needed to meet with astronomers. MacMahon

[147] The Royal Astronomical Society (RAS) was formed in 1820 by Sir James South, 1785 -1867, and was granted a Royal Charter by William IV during the first year of his reign. South was President of the RAS in 1829 and was knighted in 1830.
[148] A. A. Common, 1841 - 1903, was awarded the RAS Gold Medal in 1885, elected FRS in 1885, and served as RAS President from 1895 - 1897.
[149] An account of the trip can be found in the *Journal of British Astronomy* 101, 3 (1991).

did write one paper on an astronomical theme, in 1909 (see Chapter 6).

In 1917 he was elected to the Presidency of the Society, a post he held for the customary two years until 1919 when he was succeeded by Alfred Fowler, FRS. During his term of office he gave two Presidential addresses, at the Annual General Meetings in 1918 and 1919. In the first he described the work of Mr. John Evershed, who had been awarded the Society's Gold Medal for his investigations into sunspots, and the Rev. T. E. R. Phillips, an amateur astronomer who had been awarded the Hannah Jackson (née Gwilt) Gift for his observations of planets, double stars and variable stars. In 1919 M. Guillaume Bigourdan was awarded the Gold Medal for his observations of nebulae, which MacMahon described in considerable detail in his address.

**Honorary degree**

The first of MacMahon's four honorary degrees was awarded by Trinity College, Dublin, in 1897. The citation was not officially published for the D.Sc. *honoris causa*, but a copy exists in the Trinity college archives, although it does not say who read it out, nor why MacMahon was chosen for this honour:

> *Comitiia Aestiva, die 11 Julii 1897, habita Sc.D. PERCY ALEXANDER MACMAHON*
> *Postremo duco ad vos PERCY ALEXANDRUM MACMAHON, quem Marti deditum Minerva totum sibi vindicat; qui natales Hibernicos illustravit longa serie disquitionum quibus per regiones altiores mathematicae expatiatus est. Etenim cum horum libellorum titulos lego, qui plus quam quinquaginta sunt, reverentia tangor, atque optare subit munus laudandi metorophrontistem tam subtilem contigisse aut Praeposito nostro aut amico meo Gulielmo Burnside, qui pretium idoneum tantis merritus dicere potuissent: man me non solum Latini sermonis egestas sed crassa quaedam Minerva talia vetat culpa ingeni. Quare, Academici suppetias mihi veniat plausus vester multiplex*

This translation was provided by the University's present archivist:

> Summer Commencements on 11 July, 1897, Sc. D. PERCY ALEXANDER MACMAHON
> I finally introduce to you PERCY ALEXANDER MACMAHON, one dedicated to Mars but totally claimed by Minerva; he has shed lustre on his Irish descent by a long series of scholarly publications in which he has pursued an expansive course through the realms of higher mathematics. For when I read the titles of the articles, which number more than 50, I am struck

> with admiration, and the wish creeps into my mind that the task of praising such a high and subtle thinker had fallen to the lot of our Provost[150] or to my friend William [Snow]Burnside, both of whom could have rendered an appropriate tribute to such great merits: in my case not only the poverty of the Latin tongue, but also the obtuse texture of my thought forbids me from diminishing such achievements by a faulty use of my talents. And so, Members of the University, let the multiplication of your applause come to my assistance.

**Savilian professorship**

In 1897, just a year before MacMahon's retirement from the Army, J. J. Sylvester died, leaving vacant the post of Savilian Professor of Geometry at Oxford University. Since 1892, Sylvester's failing health had rendered him unable to fulfil the lecturing duties of the Savilian post[151]. Although Sylvester was well enough to lecture in the Michaelmas Term 1892[152], John Griffiths[153] continued to stand in for Sylvester during the Hilary Term of 1893.

It was then decided to make a more formal arrangement to cover Sylvester's indisposition.. On 23 January 1894 a proposal was made "that a Deputy be appointed to perform the statutory duties of the Savilian Professor of Geometry during the continuance of his inability."[154] A deputy was chosen in April 1894. This was William Esson[155], from Dundee, who had been a Fellow of Merton College, Oxford since 1860 and was a Lecturer and Tutor from 1865. He was the main architect of the intercollegiate scheme still operated by the University of Oxford.

[150] The Provost referred to was George Salmon, who may have had something to do with proposing MacMahon for the award, given their common interest in invarrant theory.
[151] A proposal granting leave of absence and dispensation from his duties until 30 June 1892 on the grounds of ill health was reported in the *Oxford University Gazette* on 9 February 1892 (carried *nemine contradicente* on 16 February 1892). The lecture list for the Hilary Term (from 8 February 1892 to 9 April 1892 shows John Griffiths of Jesus College lecturing on behalf of the Professor, although the Report of the Savilian Professor made on 29 March 1892 is credited to Sylvester.
[152] Hilary term ran from19 October 1892 to 17 December 1892. Sylvester's ill health continued in 1893: the *Oxford University Gazette* of 7 March 1893 records a proposal "that leave of absence be granted to James J Sylvester MA Hon DCL Savilian Professor of Geometry, for one year, and that he be allowed to appoint a Deputy Professor during that period, to be approved by the Visitatorial Board, and to assign to the Deputy Professor such stipend as the Board may approve". This proposal was reported in the *Gazette* for 14 March 1893 as "carried *nemine contradicente*".
[153] John Griffiths, 1837 - 1916, attended Cowbridge School until 1856, studied at Jesus College from 1856 - 60, and was Senior mathematical scholar in 1862. He gave three lectures on Sylvester's behalf on geometrical subjects, in 1891, 1892 and 1893. He was a Fellow of Jesus College from 1863 - 1916.
[154] There is a note appended to this proposal explaining that the Professor was permanently disabled and that the deputy be appointed for the entire term of his inability or his life if he should be unable to resume his duties or resign. This proposal was agreed on 30 January 1894, and electors were duly appointed (Edwin B. Elliott on 13 February 1894 and Edward H. Hayes on 6 March 1894). The post was advertised on 20th March 1894 with a stipend of £500 p.a.
[155] William Esson, 1838 - 1916, was elected to the Royal Society in 1869 for his investigations with A. V. Harcourt into the rate of chemical change.(*RAS Monthly Notices,* February 1917). As Deputy Savilian Professor he lectured on the 'Theory of Plain[sic] Curves' on Tuesdays and Thursdays, and on the 'Construction of Curves by Projective Pencils' on Saturday mornings During the Summer months he also held "Informal Instruction" in Geometry on Saturday mornings.

Sylvester died on 15 March 1897. On 13 April 1897, the appointment of H. H. Turner, then Savilian Professor of Astronomy, as an Elector for the Geometry post was recorded in the *Oxford University Gazette*, and on 4 May the Hebdomadal Council announced the appointment of Edwin B. Elliott, Fellow of Magdalen College and Waynflete Professor of Pure Mathematics, as a second Elector. The post was advertised in the *Oxford University Gazette* on 18 May 1897, as follows:

> The Savilian Professorship of Geometry being vacant by the death of the late Professor Sylvester, the Electors will proceed to the appointment of a successor in the course of the present term.
> The duties of the professor are defined in the following provisions of the Statutes:-
> "The Savilian Professor of Geometry shall lecture and give instruction in pure and analytical Geometry."
> "The Professor shall reside within the University during six months, at least, in each academical year, between the first day of September and the ensuing first day of July."
> "He shall give not less than 42 lectures in the course of the academical year, six at least of such lectures shall be given in each of the three University terms, and in two at least of the University terms he shall lecture during seven weeks not less than twice a week."
> The emoluments of the Professorship as determined by the Statue are as follows:-
> "He shall be entitled to the emoluments now assigned to the Professorship and derived from the benefaction of Sir Henry Savile, Knight, or from the University chest, and shall receive in addition the emoluments appropriated to the Professorship by the Statutes of New College"
> The total amount of these emoluments is £900 a year and cannot exceed that amount.
> Applications, together with such papers as a Candidate may desire to submit to the Electors, must be sent to the Registrar of the University, Clarendon Building, Oxford, on or before June 12th 1897.

The salary of £900 was a considerable sum for the late 1890s, which made the post attractive even without the status conferred by it. MacMahon had almost completed 25 years as an officer in the Royal Artillery, and was no doubt thinking about what to do in his retirement. His successes up to this point (particularly his presidency of the LMS), combined with his friendships with many of the leading mathematicians of the day and the fact that Sylvester had also been an 'outsider' - a non-Oxford man with a career that included teaching at the Royal Military Academy - gave him cause to think that he would be a serious candidate for the post. It would provide a secure financial future and considerable status in the mathematical community.

The University did not keep records of the applications for the post, nor of the meeting of the Electors. However, in addition to MacMahon's, there were at least two other applications: from William Esson, the incumbent deputy, and from Ernest William Hobson[156] of Cambridge.

The Electors made their decision on 8 July 1897; they chose Esson[157]. MacMahon was clearly disappointed, as can be seen from the letter he wrote to his friend Joseph Larmor on the following day, 9 July 1897[158]:

> Dear Larmor
> I am so pleased at your letter of sympathy.
> Hobson knows that had I been successful my pleasure would have been marred by his being unsuccessful - It would have been a real pleasure to me to have congratulated him. As it is we have been defeated by a man who I think never ought to have been a candidate. However, I must now settle down more seriously to my military work as long as I remain in the army & not be disturbed anymore by dreams of a great position which are clearly never to be realised.
> I think that on the whole my work has met with more recognition than I could have hoped for.
> Thanking you very much
> Yrs very sincerely
> Percy A MacMahon

Not only is his disappointment very clear, but also that he did not have a very high opinion of Esson. The fact that up until 1897 MacMahon had published 48 papers on a wide variety of topics, as against Esson's two articles on equations of straight lines and axial co-ordinates, may be a justification for this low opinion. Clearly, mathematical achievement and publication record were not the only criteria used by the Electors, but no records appear to exist, either in the papers of Turner or Elliott, which can shed any light on the decision making process. It is possible to speculate that the more conservative elements of the University, already discontented with the reforms that followed from the University

[156] E. W. Hobson, 1856 - 1933, was senior Wrangler in 1878, and was Sadleirian Professor at Cambridge from 1910 - 1933. He was President of the LMS from 1900 - 1902, and was awarded the De Morgan Medal in 1920.
[157] The *Oxford University Gazette* of 7 August 1897 recorded the decision thus: "At a meeting of the Electors to the Savilian Professorship of Geometry holden on Thursday July 8th, William Esson, MA, FRS, Fellow of Merton College, was elected Professor in place of James J. Sylvester MA Hon DCL deceased."
[158] This letter is to be found in the St. John's College, Cambridge, archive of MacMahon's correspondence.

Commission of the late 1870s, which included a call for more research and the widening of fellowships, felt that they had already done their duty by appointing the non-Oxonian Sylvester. They may have put pressure on the Electors to bring the Chair back to an Oxford man. It is unlikely that we will ever be able to learn the truth behind what appeared to MacMahon to be a bizarre appointment. The letter also shows that MacMahon was aware of his unusual position in the mathematical community, and may even have been surprised by his own success.

Esson[159] remained in the post[160] until his death in 1916, but MacMahon went on to achieve far more mathematically.

After Sylvester's death, the chemist and naturalist Raphael Meldola enlisted MacMahon's help to establish the Royal Society's Sylvester Medal. MacMahon was a founder member of the fund-raising committee, along with Lord N. M. Rothschild. According to Geoffrey Cantor[161], it was the use of a passage from MacMahon's obituary for Sylvester[162] that established Sylvester as "a mathematical genius and as solely responsible for the renaissance of British mathematics during Queen Victoria's reign." The passage in question does not appear in the obituary, but occurs in an earlier piece that MacMahon wrote about Sylvester's work in *Nature*[163], published in March 1897; this ends with a paragraph beginning, "The revival of the mathematical reputation of England, dating from the Queen's accession to the throne [in 1837], is to a large degree due to his genius"[164] Cantor has reinterpreted this to emphasise Sylvester's role.

## Mathematical work

### *Partition theory*

During this period, MacMahon continued the work he had begun at Woolwich with a series of four papers on partition theory. The *Memoir on the theory of the compositions of numbers*[165] is the one which led MacMahon to his Master Theorem, discussed later. This paper was written in 1892, and refereed for the Royal Society by Arthur Cayley[166] and J. W. L. Glaisher[167]. The referees did no more than comment that any paper by MacMahon would be worthy of printing.

[159] Esson was also a delegate at the 1912 ICM, at which he presented two short papers.
[160] It is only fair to point out that Esson was a very able administrator, and did much to shape the University's tutorial system.
[161] [Cantor, 2004, p. 83].
[162] [MacMahon, 1897, [50;19]].
[163] [MacMahon, 1897b].
[164] This quote is reported in the *Jewish Chronicle,* dated 4 June 1897, p. 8, in a letter from Meldola.

In the paper MacMahon developed interesting graphical methods to represent compositions, of *unipartite* numbers, *bipartite* numbers and higher order *multipartite* numbers[168]. For example, two compositions of 7 are (214) and (13111); by using a graphical method it is easy to see that these are *conjugate* compositions, so that given a composition of $n$ into $p$ parts, a second composition of $n$ into $n - p + 1$ parts is obtained simply by swapping circled and uncircled dots:

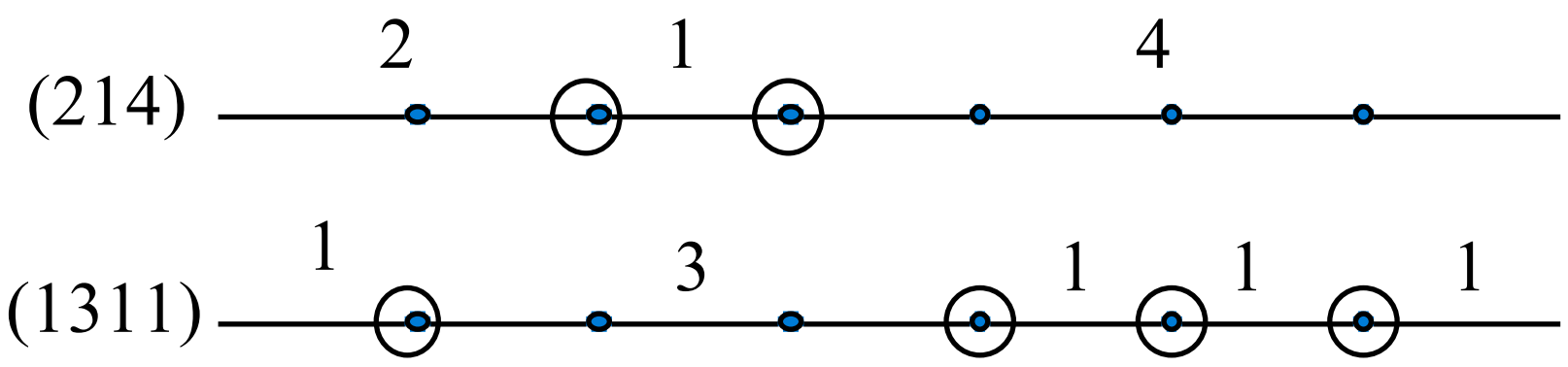

Figure 1: Conjugate partitions

Compositions of bipartite numbers are represented on a grid. The diagram below shows the composition (21)(11)(11) of the bipartite number $\overline{43}$. The paths or routes are read from bottom left to top right, so the first blue circle (*node*) represents (21).

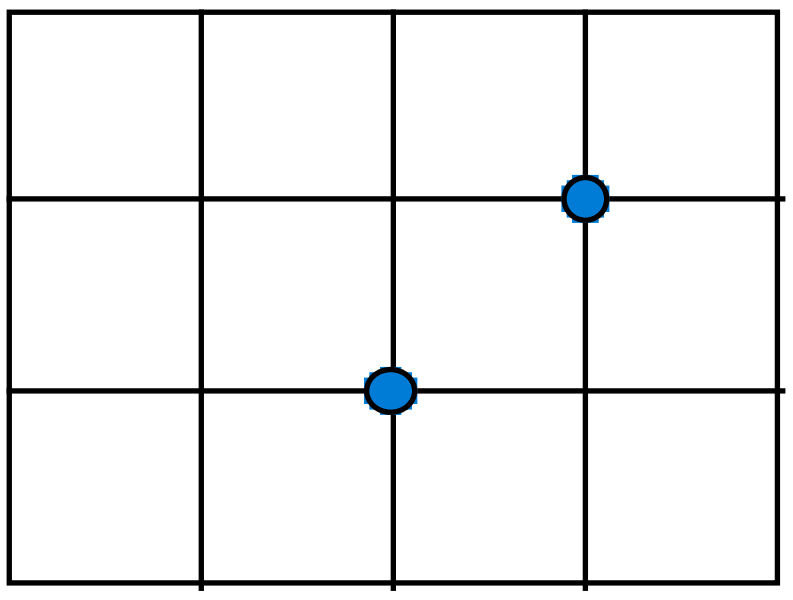

Figure 2: The bipartite number $\overline{43}$

An alternative composition is (31)(01)(11):

[165] [MacMahon, 1894, [46;5]].
[166] Royal Society archive reference RR 12.3.
[167] Royal Society archive references RR 12.4 and 12.5.
[168] An $n$-partite number is equivalent to a $n$-dimensional vector using only integers.

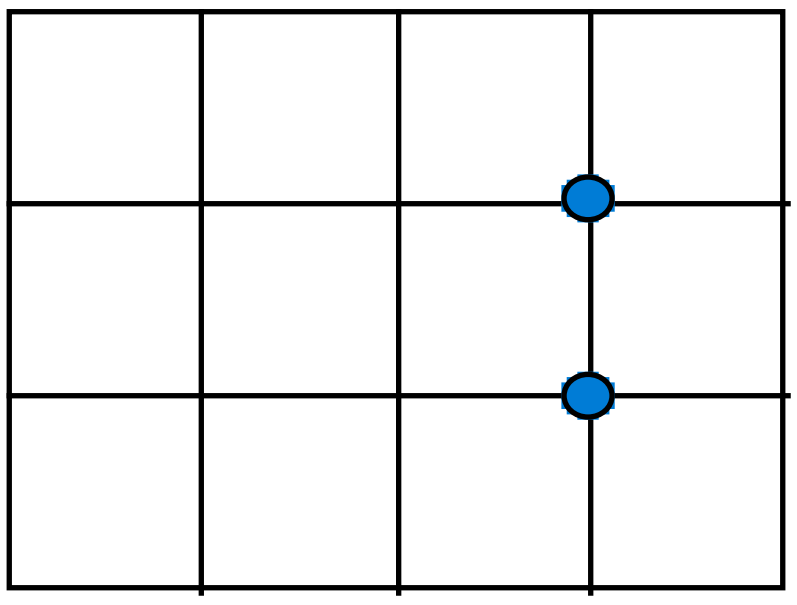

Figure 3: Alternative composition (31)(01)(11)

A node is *essential* if it occurs where the direction of the route changes from the vertical direction to the horizontal direction. The lower blue node is *non-essential* to this route; so if we delete it and put red crosses at the nodes not previously used (and the essential upper blue node) on the route we get the conjugate composition (10)(10)(10)(02)(10)(01):

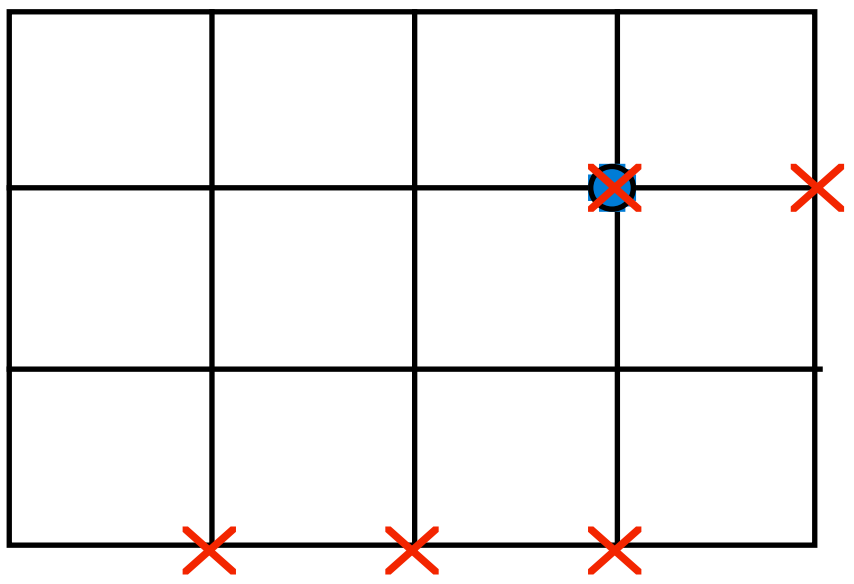

Figure 4: Conjugate composition

The paper develops some sophisticated methods for enumerating the lines of route using permutations. The number of different lines of route on the graph of the bipartite number $\overline{pq}$ is the number of permutations of $p$ symbols $\alpha$ and $q$ symbols $\beta$, where the $\alpha$ and $\beta$ represent a single unit in the horizontal and vertical directions, respectively. The compositions can be arranged in conjugate pairs, and the investigation then moves on to the enumeration of compositions represented by lines of route with exactly $s$ essential nodes. The result is $\sum^{s} \binom{p}{s}\binom{q}{s} 2^{p+q-s-1}$ which MacMahon states (but does not prove or justify) is analytically equivalent to the algebraic expression

$$\frac{1}{1-2(x+y-xy)} = \frac{1}{(1-2x)(1-2y)} + \frac{2xy}{(1-2x)^2(1-2y)^2} + \frac{2^2x^2y^2}{(1-2x)^3(1-2y)^3} + \cdots$$

The extension of this method to three and more dimensions leads him to the generating function for compositions of multipartite numbers of order $n$:

$$\frac{1}{2}\frac{1}{\{1-s_1(2\alpha_1+\alpha_2+\cdots+\alpha_n)\}\{1-s_2(2\alpha_1+2\alpha_2+\cdots+\alpha_n)\}\ldots\{1-s_n(2\alpha_1+2\alpha_2+\cdots+2\alpha_n)\}}$$

In a further step, MacMahon then generalised the concept of a composition thus: if we think of the number $p$ as succession of $p$ units arranged in a row, we have $p$ - 1 spaces between them, which we may fill with a + sign or leave blank. The $2^{p-1}$ different expressions so obtained are the compositions of $p$. This is illustrated in the table below for $p$ = 4, where $\circ$ is used to signify a blank.

| | |
|---|---|
| $1\circ1\circ1\circ1$ | $(1111)=(1^4)$ |
| $1+1\circ1\circ1$ | $(211)=(21^2)$ |
| $1\circ1+1\circ1$ | $(121)$ |
| $1\circ1\circ1+1$ | $(112)=(1^22)$ |
| $1\circ1+1+1$ | $(13)$ |
| $1+1+1\circ1$ | $(31)$ |
| $1+1\circ1+1$ | $(22)=(2^2)$ |
| $1+1+1+1$ | $(4)$ |

Table 1: Generalised composition

MacMahon then established a one-to-one correspondence between these compositions and rooted trees[169]. In the example above, each composition of 4 corresponds to one of the eight rooted trees of height 2, as shown below.

[169] Rooted trees had already been studied by Cayley, in [Cayley, 1857], [Cayley, 1875] and [Cayley, 1881].

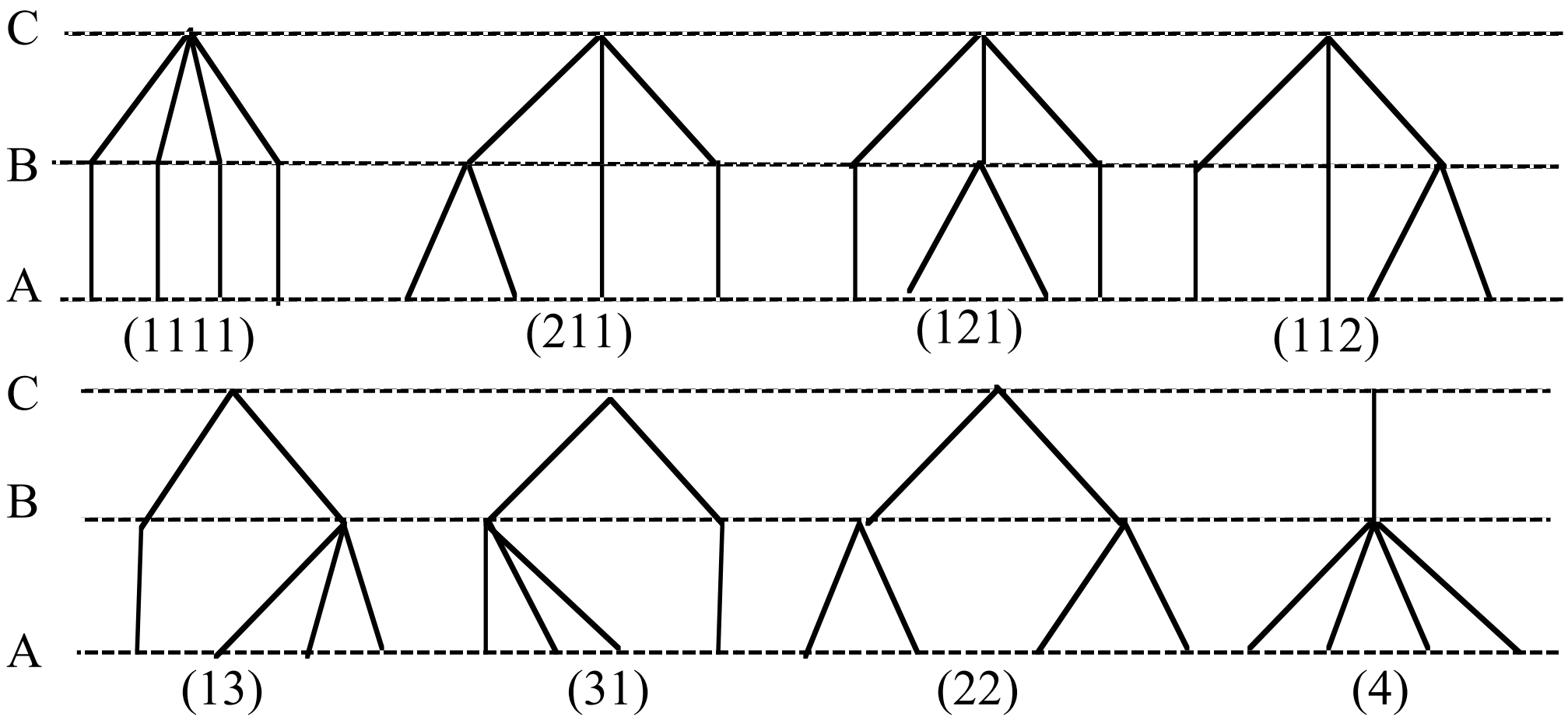

Figure 5: Composition as rooted trees

In this scheme, a node on line C corresponds to a blank, a node on line B corresponds to +, and the number of terminal nodes on line A is $p$. (Note that there is always one more 'down' branch than the number of nodes on a line, because the number of numbers in a composition must be one more than the number of symbols between them).

MacMahon then examined the case for trees of height $k$, corresponding to $k$ different symbols between the $p$ units, called *combinations* (so that compositions are the special case where $k$ = 2). This led to the generating function

$$\frac{1}{k}\bullet\frac{1}{\left\{1-s_1\left(ka_1+a_2+\cdots+a_n\right)\right\}\left\{1-s_2\left(ka_1+ka_2+\cdots+a_n\right)\right\}\cdots\left\{1-s_n\left(ka_1+ka_2+\cdots+ka_n\right)\right\}}$$

where the number of combinations of order $k$ of the multipartite number $\overline{p_1p_2\cdots p_n}$ is the coefficient of $\left(s_1a_1\right)^{p_1}\left(s_2a_2\right)^{p_2}\cdots\left(s_na_n\right)^{p_n}$ in the expansion of the function. This result follows from earlier results in the paper.

The paper ends with a short postscript, added on 25 August 1893, in which MacMahon explained that he had discovered a general result, later to be called the 'Master Theorem', in the period between reading the paper and its publication. The various enumeration problems introduced in the paper all required the evaluation of a coefficient in the expansion of a generating function of the form $\frac{1}{(1-s_1X_1)(1-s_2X_2)\cdots(1-s_nX_n)}$ where the $X_n$ are linear functions of the form $a_{11}x_1 + a_{12}x_2 + \cdots + a_{1n}x_n$. The coefficient required is that of terms which are functions of the products $s_1x_1, s_2x_2, \ldots, s_nx_n$. The Master Theorem states that the part of the expansion involving these functions is given by $V_n^{-1}$, where

$$V_n \text{ is the determinant } (-1)^n x_1x_2\ldots x_n \begin{vmatrix} a_{11}-\frac{1}{x_1} & a_{12} & \cdots & a_{1n} \\ a_{21} & a_{22}-\frac{1}{x_2} & \cdots & a_{2n} \\ \vdots & \vdots & \ddots & \vdots \\ a_{n1} & a_{n2} & \cdots & a_{nn}-\frac{1}{x_n} \end{vmatrix} \text{ and all the } s_n = 1.$$

A paper[170] the following year provided a proof of this result and examples of applications, and is summarised at the end of this section.

Following the 1894 compositions paper described above, there were no further papers on partitions for three years, until in 1897 MacMahon published the first of a series of seven papers on partition theory (the seventh and last appeared in 1917). The first *Memoir on the theory of the partition of numbers*[171] is important for two reasons: the further development of graphical methods to assist in the study of partitions, and MacMahon's efforts to extend partitions to three dimensions and beyond.

This paper was sent by the Royal Society to W. Burnside[172] for comment, but he returned the manuscript because "the subject of the memoir deals with a branch of Pure Mathematics of which I have at present no knowledge."[173] It fell to J. W. L. Glaisher, E. B. Elliott[174] and an unidentified third

170 [MacMahon, 1894, [45;3]].
171 [MacMahon, 1897,[51;19]].
172 William Burnside, 1852 - 1927, was second Wrangler in 1875, elected FRS 1893,and served as LMS President from 1904 - 1906. He was Professor of Mathematics at the Royal Naval College, Greenwich, from 1885 - 1919.
173 Royal Society archive references RR13.76 & 77.
174 Edwin Bailey Elliott, 1851 - 1937, was Waynflete Professor of Pure Mathematics in Oxford from 1892 - 1921, and served as LMS President 1896 - 1898.

referee to recommend the paper for publication.[175]

In this paper, MacMahon introduced the idea of a *separation* of a partition. This is a partition that has been broken into parts; for example, 5 + 3 + 3 + 1 + 1 is a partition of 13 which may be separated into 5 + 3 and 3 + 1 + 1. The weights of the parts are 8 and 5 respectively, so the specification of the separation is (8,5); the degree of the separation is 5 + 3 = 8, the sum of the greatest integers in each part of the separation, and the multiplicity is 2. MacMahon drew attention to the correspondence between the number of separations of a partition and the number of partitions of the associated multipartite number; that is, the number of separations of the partition $\left(p_1^{\pi_1} p_2^{\pi_2} p_3^{\pi_3} \ldots\right)$ is identical to the number of partitions of the multipartite number $\left(\overline{\pi_1 \pi_2 \pi_3 \ldots}\right)$.

From this, MacMahon was able to develop a graphical representation of partitions by combining the idea of representing a partition using dots, introduced by Sylvester and known as a *Ferrers graph*[176], and the diagram used to illustrate compositions described above (a figure known as a *reticulation*). The line representing the composition divides the diagram into two regions which can be considered as the Ferrers graphs for the partitions of two numbers, as in the following example:

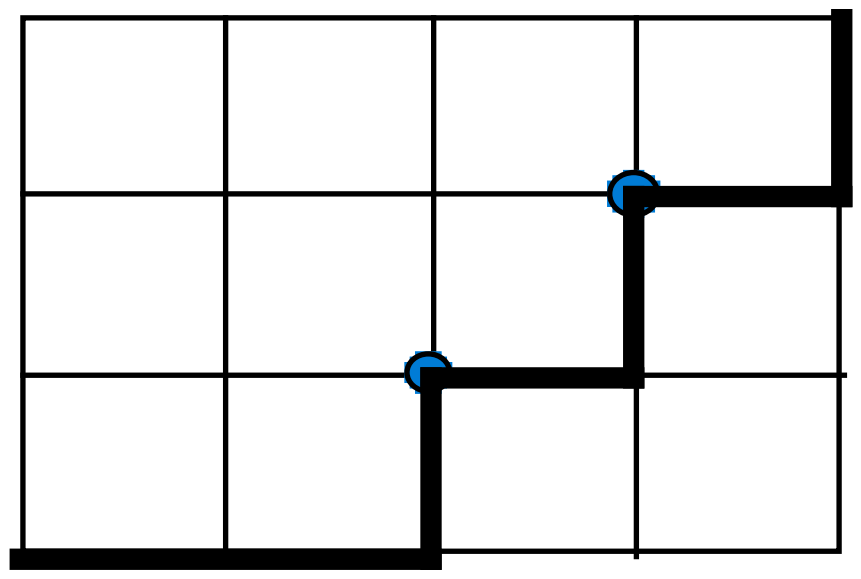

Figure 6: Line of route of a composition

The line (called a *line of route*) is the composition (21)(11)(11) of $\overline{43}$; the upper part of the diagram is the partition (432) of 9 and the lower part is the partition (21) of 3. MacMahon examined in some detail the connection thus demonstrated between compositions of bipartite numbers and partitions of unipartite numbers. From graphical considerations he was able to discover many interesting formulas concerning partitions that he claimed would have been diffficult or impossible to derive algebraically; for example, the number of partitions of a number with exactly $q$ parts, a greatest part equal to $p$ and $s$ different parts is $\binom{p-1}{s-1}\binom{q-1}{s-1}$.

[175] Royal Society archive references RR13.78 to 81. Glaisher was not very enthusiastic due to the technical nature of the paper, but nevertheless recommended publication.

[176] Sylvester attributes the first use of these diagrams to Dr. N. M. Ferrers, 1829 - 1903, in a proof of the reciprocity theorem of partitions - that the number of partitions of $n$ into $j$ parts each not greater than $i$ is the same as the number of partitions into $i$ parts not greater than $j$. This was not published by Ferrers, but by Sylvester himself [Sylvester, 1882, p. 253].

***The Master Theorem***

The Master Theorem referred to above was described more fully in a paper entitled *On a certain class of generating functions in the theory of numbers* published in 1894[177]. In it, MacMahon provided first a proof[178] and then some examples of its use in a variety of situations, including the Problème des Rencontres (a permutation in which no symbol remains in its original place) and a number of variants thereof. An example of the use of the result is given in more detail in Appendix 10. He used the Theorem in much of his later work, whenever an enumeration problem was encountered, since it enabled him to create generating functions more easily.

***Recreational mathematics***

It was during the period covered by this chapter that MacMahon's recreational work was begun, partly in collaboration with a fellow officer from the Royal Artillery, Major J. R. J. Jocelyn[179].

MacMahon's total output on recreational topics comprised seven papers and one book, listed in chronological order in Table 2 below. The first paper in the table has already been mentioned in Chapter 3, and is of minor interest. The second paper, *Weighing by a series of weights*, is essentially concerned with partition theory, albeit in the context of a puzzle. The later papers are of more importance and are discussed below and in chapter 6.

[177] MacMahon noted that, strictly speaking, the paper was concerned primarily with the theory of determinants, rather than partition theory.
[178] Andrews [Andrews, 1977, pp377 - 378] gives details of further, shorter proofs and examples of applications discovered by mathematicians later in the 20th century.
[179] J. R. J. Jocelyn,1851 - 1929, was a colleague of MacMahon's at Woolwich, and Commandant of the Ordnance College from 1899 - 1903.

| Year | Andrews reference | Title |
|---|---|---|
| 1889 | [33:10] | On play "à outrance" |
| 1890 | [37:19] | Weighing by a series of weights |
| 1893 | [44:15] | On the thirty cubes constructed with six coloured squares |
| 1902 | [60:19] | Magic squares and other problems upon a chess board |
| 1921 | [96:15] | *New Mathematical Pastimes* (CUP) [book] |
| 1922 | [98:15] | The design of repeating patterns |
| 1922 | [99:15] | The design of repeating patterns for decorative work |
| 1922 | [100:15] | Pythagoras's theorem as a repeating pattern |

Table 2: MacMahon's recreational papers

***Puzzle patents***

In addition to his papers, MacMahon was granted three patents for puzzles. The first patent for which he applied (number 21,118) was a puzzle comprising nine wooden blocks linked by flexible tapes to form a chain. The length of the links was to be such that the blocks could be formed into a stack in 4527[180] different ways. A design applied to the edges of the blocks would have to be reconstructed by finding the correct stacking. This patent was accepted on 8 October 1892, but there are no records to indicate whether the puzzle was ever manufactured commercially.

It is another example of MacMahon being ahead of his time; the idea is equivalent to the modern puzzle of counting the number of ways that a strip of postage stamps may be folded. This puzzle was described by Martin Gardner in 1983[181], and he claimed the earliest description of the puzzle dated from 1937[182]. MacMahon had created the puzzle more than forty years earlier, but by not writing about it, this contribution to recreational mathematics was lost.

180 MacMahon does not say how he calculated this figure, which is incorrect. The correct number is 4536. [Eric W. Weisstein, "Stamp folding", from *Mathworld* - a Wolfram web resource. http://mathworld.wolfram.com/StampFolding.html].
181 [Gardner, 1983].
182 [Sainte-Lague, 1937].

In 1892, MacMahon and Jocelyn jointly applied for two patents. The first of these, patent number 3927, was for *Appliances to be used in playing a new class of games*. The provisional specification was submitted on 29 February 1892, and the patent was granted on 28 January 1893. The patent described equilateral triangular pieces divided into three compartments thus:

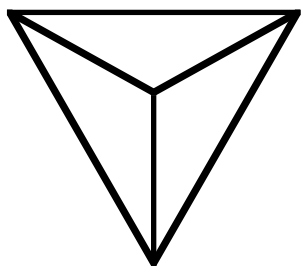

Figure 7: Three compartments of a triangle

The compartments were to be numbered or coloured in a manner analogous to dominoes[183]; for example, a set of 'treble 3s' would comprise 24 such pieces shown below, where colours are used instead of the numbers 0, 1, 2, 3 (Figure 8 is adapted from *New Mathematical Pastimes*[184], page 2, figure 2). A set of 'treble 4s' would require 45 pieces.

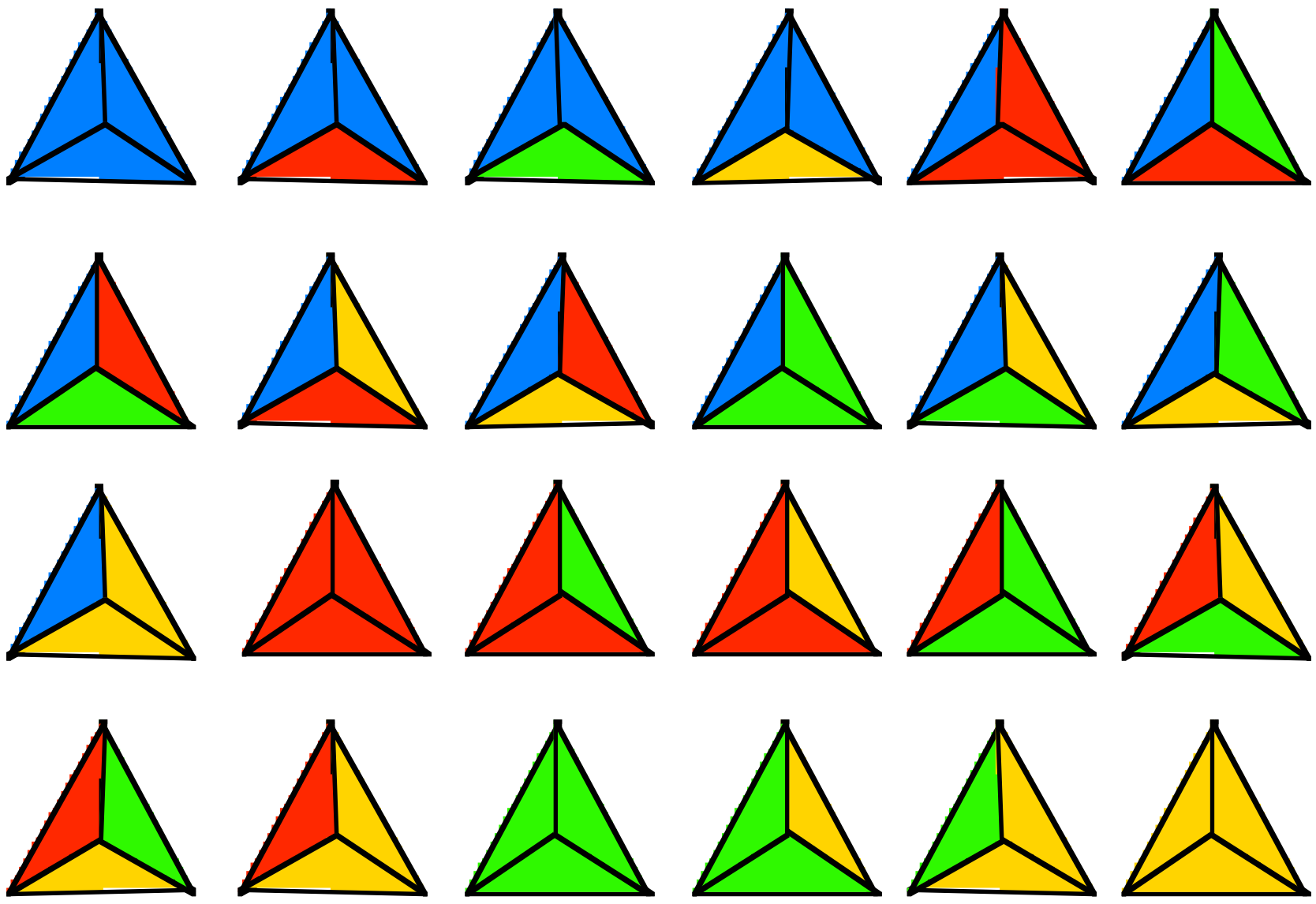

Figure 8: 24 four colour triangles

The rules of two different edge-matching games to be played by two players on specially marked boards were described in the patent in some detail, and a suggestion was made for a puzzle where the pieces were to be fitted together to form a hexagon with (say) a single colour around the perimeter, this latter idea was to be expanded upon in greater detail in *New Mathematical Pastimes* some 29 years later. It is not known whether any sets of pieces were manufactured or sold.

[183] In a standard domino set, all possible pairs of the numbers 0 to 6 (inclusive) are used, and the set is called 'double sixes', after the domino made from the two largest numbers.

[184] [MacMahon, 1921, [96;15]]. Although an Andrews reference exists for this book, it is not reproduced in the *Collected Papers*.

The second joint patent application was made on 2 May 1892. It described how a set of 27 coloured cubes could be constructed using three colours, the puzzle being to assemble the set into a larger cube with a uniform colour on each external face and matching between the internal faces. The existence of a set of 30 cubes that could be made with six colours was also mentioned, and the puzzle described was to select two 'associated' cubes, with the same pairs of colours on opposite faces, and then to locate amongst the remaining 28 cubes a set of sixteen cubes in which none of those opposites occurs. These sixteen can then be assembled into two larger copies of the two associated cubes, where internal faces must also match.

This set of cubes was described MacMahon's 1893 paper, *On the thirty cubes constructed with six differently coloured squares*[185]. This is the most important of his recreational papers, in the sense that it is the first public description of original recreational work by MacMahon. It was read at a meeting of the London Mathematical Society chaired by A. B. Kempe on 9 February 1893; the meeting and MacMahon's paper were reviewed in *Nature* on 23 February 1893.

In his paper, MacMahon first described how to calculate the number of different cubes that can be constructed with six differently coloured squares, and then explained how they may be paired up according to their symmetry group. This is done by noting that for any particular cube there is a group of 24 'substitutions' or rotations, that leave the cube unaltered[186]. The cube formed by exchanging any single pair of opposite colours is represented by the same group with the rotations in the opposite senses, so the cubes fall naturally into pairs of 'associated cubes'. From this basis, he described the essential puzzle: given any one of the thirty cubes (the *target* cube), to construct a 2 x 2 x 2 cube with the same colour arrangement from eight of the remaining cubes, such that the internally touching faces match. This is done by using the octahedral dual of the selected cube (with the colours labelled 1 to 6) and performing a sequence of substitutions on the faces of the octahedron. The method selects the correct cubes and places them in their correct orientations and positions. In *New Mathematical Pastimes*, the interested reader is referred to this paper for a detailed description of the method.

[185] [MacMahon, 1893, [44;15]].
[186] MacMahon lists the substitutions for a particular cube explicitly in [MacMahon, 1893, [44;15], p.2]

A puzzle called *Mayblox* was manufactured by R. Journet and Co. of London[187] (best known for their glass-topped dexterity puzzles) with the legend "Invented by Major P. A. MacMahon F. R. S." on the box. This version of the puzzle, described by Margaret A. Farrell in the *Journal of Recreational Mathematics* in 1969[188], asked the solver to construct a large cube from 8 smaller cubes as in MacMahon's puzzle, but without the benefit of the target cube. Farrell described how to set about solving the puzzle, and also how to create a set of suitable cubes. A footnote to the article states:

> The *Mayblox* puzzle was invented by Major P. A. MacMahon, R.R.S [sic] and was produced by a London firm, R. Journet. This puzzle was found among the effects of Maj. MacMahon's late aunt by Mr R. N. Andersen, a colleague at State University of New York at Albany

This footnote is misleading, as explained by Professor Farrell in a private communication where she said that the footnote had been altered during editing to read "that the puzzle had been found in the effects of Major M[a]cMahon's aunt. It should have said, 'found in the effects of his late aunt,' meaning the effects of R. Andersen's late aunt." The particular puzzle in question was purchased in Hamley's toy shop in Regent Street, London, sometime between 1901 and 1915, according to a private letter from Professor Andersen.

A version of a puzzle similar to *Mayblox* called 'The cubes of Mahomet' was described in the magazine called *Work* in 1916 by someone known as 'Toymaker'. *Work* was a journal for hobby woodworkers and gave instructions for the manufacture of a variety of puzzles and illusions. It is not known whether these instructions were published with MacMahon's knowledge or permission. There is no acknowledgement of the origin of the puzzle, and no hint that a solution is possible by any means other than patient 'trial and error'. However, the existence of the article demonstrates the wide appeal of MacMahon's brainchild.

[187] *Mayblox* may have been sold during the period 1895 to 1915, but records to confirm this no longer exist. The trade mark for the puzzle was registered to Major Jocelyn in 1897. From examples of the puzzle in the author's collection, there were at least two different versions made, one with cubes 1 inch on a side, the other with cubes of side $^7/_8$ inch. It is not possible to tell with absolute certainty which was the earlier set, although the better quality of the larger set suggests it may been an improved version of the puzzle.
[188] [Farrell, 1969].

As further evidence of the attractiveness of the puzzle, in 1930, Gerhard Kowalewski[189] wrote a chapter on the 'MacMahonsche' cubes in his *Alte und Neue Mathematische Spiele.*[190] This book gave instructions for constructing a set of cubes, with a detailed procedure for the colouring process. Kowalewski described the MacMahon version of the puzzle, and also claimed that a version of the puzzle was available commercially, in which eight cubes must be assembled to make a larger cube with faces each of a single colour - but without the benefit of the target cube, as in the *Mayblox* puzzle described above. Kowalewski also gave a version of his own where two sets of opposite faces have the *same* colour and the remaining two faces have a single colour each, so that two colours do not appear on the outside of the large cube.

In 1934, Ferdinand Winter, a student of Kowalewski's, published *MacMahon's Problem - Das Spiel der 30 bunten Würfeln* [191], a book reviewed, not entirely favourably, by W. L. Ayres[192] in the widely read *American Mathematical Monthly*[193], a journal published by the American Mathematical Society. The book discussed both MacMahon's puzzle and Kowalewski's variation. The level of detail was described by Ayres as tedious, which seems harsh given that the work was Winter's thesis, submitted to the Sächsische Technische Hochschule in Dresden for the degree of Doktor-Ingenieur. This is the list

[189] Gerhard Kowalewski was Professor of Pure Mathematics at the *Technische Hochschule* in Dresden.
[190] [Kowalewski, 1930].
[191] There are three surviving copies of this work: at the University of California in Berkeley, at the University of Cincinnati, and in the Niedersachsische Staats und Universitäts library at Göttingen.
[192] W. L. Ayres, 1905 - 1976, was a topologist, and a student of J . R. Kline; he worked at the University of Michigan from 1929 - 1941, and at Purdue University from1941.
[193] The review is in Vol. 42 (Issue 9, Part 1), November 1935, pp. 563 - 564.

of contents (in translation):

1. Description of the puzzle
2. A new version of the puzzle
3. MacMahon's problem
4. The problem of two MacMahon cubes
5. Kowalewski's game with 16 cubes
6. The connection between MacMahon cubes and Kowalewski cubes
7. Locating the necessary cubes
8. A new game with three cubes
9. Cube puzzle: top and bottom the same colour
10. Cube puzzle: top and bottom with four colours
11. Cube puzzle with four colours on every outer face
12. Cubes with three hidden colours
13. Games with 12 and 16 cubes
14. Games with 16 cubes
15. A new game with 30 squares
16. A new game with magic cubes

From the chapter headings, it can be seen that Winter had conducted a very detailed analysis of the possibilities afforded by the original idea, and extended the number of puzzles beyond those suggested by MacMahon. Winter claimed that the puzzle was obtainable directly from the firm S. F. Fischer in Oberseiffenbach, or in toyshops and bookstores. There are no known surviving examples, and enquiries made with S. F. Fischer in 2004 were unsuccessful in locating any records of such a puzzle.

As part of Kowalewski's 60th birthday celebrations, the *Monatshefte für Mathematik und Physik* published a series of articles, including one by Winter[194] containing a further new variant on the puzzle. This involved building two 2 x 2 x 2 cubes such that only four colours appeared on the top and bottom faces, and the middle horizontal internal faces, with the vertical internal faces containing a repeating pair of the four colours, one pair in each cube. The remaining two colours had to appear in pairs on the side faces. Clearly Winter was pushing at the limits of the puzzle, turning it from an amusing pastime into an obsession[195].

[194] [Winter, 1936].
[195] However, since MacMahon's purpose, enunciated later in *New Mathematical Pastimes*, had been to involve the reader in actively creating puzzles by providing a mathematical framework, the project must be judged a success, certainly in Winter's case.

### *Magic squares*

MacMahon's fourth recreational paper was *Magic Squares and other problems on a chessboard*[196]. This was a transcript of a lecture he gave at the Royal Institution on 14 February 1902, published in the *Proceedings of the Royal Institution*[197] and reprinted in *Nature* in March 1902. MacMahon had been invited by the Royal Institution to present the talk at one of the Friday evening discourses for the general public.

In his talk, MacMahon recounted some history of magic squares and some methods of construction (in particular the method of De la Hire, where two Latin squares are added cell by cell), and then introduced the problem of enumerating the magic squares of a given order. This led seamlessly into a more general discussion of Latin squares, Graeco-latin squares and group theory. He then introduced two famous problems:

(1) the *Problème des Rencontres* ('in how many ways can $n$ letters be placed in $n$ envelopes so that no letter is in the correct envelope?'), as a way of calculating the number of ways that the second row of a Latin square might be completed given the first row;

(2) the *Problème des Ménages* ('in how many ways can $n$ married couples be seated at a circular table so that no married pair sit together?') to calculate the number of possible third rows in the Latin square.

The problem of placing eight castles on a chessboard so that none is attacked by another was used to provide a link between algebra (the study of continuous quantities) and arithmetic (the study of discrete quantities), by showing how the process of differentiation could be used to enumerate the placing of the castles on successive rows of the chessboard. Briefly, choosing the placement of the castle on the first row is enumerated by differentiating $x^8$ to get $8x^7$, where the coefficient 8 is the

196 [MacMahon, 1902, [60;19]].

197 Volume XVII (96), February 1902, pp. 50 - 61.

number sought. A second differentiation gives $56x^6$, enumerating the ways that two castles can be placed on the first two rows of the chessboard, and so on.

The talk was illustrated with slides, and MacMahon was aware of the non-specialist nature of his audience, for he gave clear explanations of terms and processes that a mathematical audience would not need. For example, in describing a variation of the Latin square where the number of symbols to be placed is greater than the number of cells in a row, so that there is no restriction on the number of symbols that may appear in each cell but each symbol must still appear once in each row and column, he claimed that the number of arrangements is "$(4!)^7$". He gave this explanation of the notation: "4, the order of the square, must be multiplied by each lower number, and the number thus reached multiplied seven times by itself."

MacMahon's point was to show the audience that apparently widely disparate fields of study, such as the problem of the general nature of the magic square, and the infinitesimal calculus, could be linked with great advantage and reveal hidden depths to the processes involved. He finished the talk by thanking the Royal Institution managers for allowing him "to display some of the chips from a pure mathematician's workshop."

***Probability***

MacMahon was also interested in the combinatorial aspects of probability. The Karl Pearson[198] collection at University College London (UCL) contains three letters from MacMahon to Pearson[199], which are transcribed below.

*Letter 1; dated 8 May 1894, from MacMahon to Pearson*

Dear Sir,

[198] Karl Pearson, 1857 - 1936, was third Wrangler in 1879, and Gresham Professor of Geometry from 1890 -1894. He founded the journal *Biometrika* in 1901, and was the Galton Professor of Eugenics at UCL from 1911 to 1933. He used statististical methods to problems in evolution and heredity.

[199] Karl Pearson archive reference [756/7].

> I am much interested in your Monte Carlo paper in the Fortnightly Review. During my study of it it becomes necessary for me to know the following:- You give a list of runs in 4274 throws at roulette. Does this number include the zeroes thrown or is it the number after the zeroes have been thrown out ?
> Kindly excuse the trouble I give you and believe me, yrs sincerely, P A MacMahon
> P.S. I presume Black Red. Zero. Red Black would be regarded (in your paper) as involving a sequence of two reds ?

The article to which MacMahon referred appeared in the *Fortnightly Review*[200] under the title 'Science and Monte Carlo.' In this article, Pearson described his efforts to amass some experimental evidence for the 'Laws of Chance' to provide material for a series of 'popular lectures' he gave in 1893[201], for he feared that he would be unable "to lead [the] audience though the mazy and not over-sure paths of mathematical theory." Pearson described how he had spent a good portion of his vacation in tossing a shilling 25,000 times, and noted that a "machine to show within the brief period of a lecture the result of several millions of tosses of 20 coins at a time, or the like number of throws of dice, only failed owing to the views of the British carpenter on the variability of the British inch."

Pearson was thus pleased to discover a publication called *Le Monaco*, issued weekly in Paris, which recorded the results of all the roulette games at the tables of Monte Carlo. Using the following table of data (Table 3) collected by himself and others (and one can only marvel at the patience of the data collectors!), Pearson addressed the question of how far from the theoretical 50:50 split would results have to deviate "in order for us to assert that we are not dealing with a game of chance ?"

[200] [Pearson, 1894].
[201] Although Pearson did not say so, these were probably his Gresham lectures.

| Method | Success (%) | Failure (%) | Trials | Experimentor or Calculator |
|---|---|---|---|---|
| Roulette | 50.15 | 49.85 | 16141 | Pearson |
| Roulette | 50.27 | 49.73 | 16019 | De Whalley |
| Bag of balls | 50.11 | 49.89 | 10000 | Westergaard |
| Bag of balls | 50.4 | 49.6 | 4096 | Quetelet |
| Tossing | 51 | 49 | 4040 | Buffon |
| Tossing | 50.05 | 49.95 | 4092 | De Morgan's pupil |
| Tossing | 50.04 | 49.96 | 8178 | Griffith |
| Tossing | 50.16 | 49.84 | 12000 | Pearson |
| Tossing | 50.05 | 49.95 | 24000 | Pearson |
| Lottery | 50.034 | 49.966 | 7275 | Westergaard |

Table 3: Pearson's results

Having decided that, in terms of the occurrences of Red and Black numbers in roulette it appeared indeed to be a game of chance, Pearson then turned to the frequency with which the individual numbers appeared in roulette games, and asked, "What is a reasonable amount for the standard deviation of an experiment of this kind to differ from its theoretical value by?" The results from Mr. de Whalley's data surprised him by showing a divergence from the theoretical standard deviation to be expected only once in 270 years. Not wanting "to put aside as useless [this] very improbable month's returns", Pearson decided to examine the frequencies of runs of the same colour. Zeros he disregarded (hence MacMahon's postscript in the above letter).

Examining a particular fortnight's play, Pearson was led to remark, "If Monte Carlo roulette had gone on since the beginning of geological time on this earth, we should not have expected such an occurrence as this fortnight's play". The conclusion was inescapable: "Roulette as played at Monte Carlo is not a scientific game of chance." Pearson finished the article with a call for the 'hierarchy of science' to put pressure on the French government to close the casinos.

MacMahon's letter received a prompt reply from Pearson, and the second letter from MacMahon was dated the following day.

*Letter 2; dated 9 May 1894, from MacMahon to Pearson:*

> Dear Sir,
> I am much obliged for your quick reply to my letter. It is my intention to critically examine your conclusions and to tender my own for publication. I have read your article in the Fortnightly Review but will not have a copy of my own till tomorrow. I have sent to Paris for various numbers of Le Monaco and on receiving them hope to set to work at once.
> So many of us try our luck at the tables sometime during our careers that I conceive the subject to be [of] great importance and your conclusions to seriously ‘give pause’.
> Believe me, Yrs very sincerely, P A MacMahon

Pearson must have supplied MacMahon with some data, and by the time he wrote the third letter, MacMahon had had time to examine the results presented by Pearson in the *Fortnightly Review*.

*Letter 3; dated 16 May 1894, from MacMahon to Pearson:*

> Dear Sir,
> I am very much obliged to you for the information you have been kind enough to send me - I will certainly examine the (some of the) numbers in the summer and autumn of 1892. I do not expect to ‘smite’ you at all! tho’ the method I employ and the results I reach may not be wholly in agreement with yours. I am wrestling with some of your figures. So far as I can understand the “student” the number 16 in the fourth line of p.192 should read 8 - If this be so the succeeding 6 lines are open to criticism[202]. However as Monte Carlo repudiates the numbers in Le Monaco it hardly seems worthwhile undertaking the labor tho’ I believe that I will as I am greatly interested in the matter, With kind regards, Yrs very sincerely,P A MacMahon

There is no evidence that MacMahon completed the promised analysis.

## Retirement

In 1898, at the age of 44, MacMahon took paid retirement from the Royal Artillery. He had been a serving officer for a quarter of a century; the lack of any surviving personal documents means that it is not possible to be certain why he chose to retire at this particular time. His career had encompassed five years of active front line duty, two years of study and 17 years of teaching. It is possible that he was ready for a rest; he may have felt that his future lay with mathematical research, particularly since his army pension would have allowed him the time and leisure to pursue such intellectual pursuits. He was at the height of his powers, and still young, so the unencumbered bachelor life would have held great appeal.

[202] This refers to a long quote made by Pearson taken from an unnamed evening paper, chosen to illustrate the common fallacious thinking surrounding the game of roulette (namely, that a run of one colour is more likely to be followed by the other colour) which begins, ‘Some time since, at Monte Carlo, a student at the game of roulette sat for 48 days at one particular table noting the spin of the ball’ and the subsequent discussion, which depends on the average difference between observation and theory being 16, as apparently noted by the “student”.

# Chapter 5 London 1898 - 1906

During the years covered by this chapter, MacMahon lived a bachelor life in the heart of London, occupying rooms in Shaftesbury Avenue until 1901, when he moved to St Anne's Mansions in St James's (now the site of the Home Office). He continued his work on invariant theory, and although he continued to publish papers on the subject, the balance of his interests had shifted substantially towards partition theory. He was very involved with the British Association for the Advancement of Science (BAAS), and this forms a large part of this chapter.

### *Royal medal*

Two years after retiring from the Army, in 1900, he was awarded the Royal Medal of the Royal Society "for the number and range of his contributions to mathematical science." This medal was awarded annually by the Sovereign to two people upon the recommendation of the Council; the other 1900 winner was Alfred Newton[203], "for his eminent contributions to the science of ornithology and the geographical distribution of animals." Other mathematicians awarded the Royal Medal specifically for their mathematical work have included Boole in 1844, Cayley in 1859, Sylvester in 1861, Salmon in 1868, Rayleigh in 1882, Hirst in 1883, Tait in 1886, Forsyth in 1897, Burnside in 1904, Greenhill in 1906, Hobson in 1907, Larmor in 1915 and Hardy in 1920. Of the 80 medals awarded between 1890 (when MacMahon was elected to the Society) and 1929 (the year of his death), only 10 were awarded for mathematical work.

### *Winchester College*

In 1899, MacMahon was made a Governor of Winchester College. According to the Register of the College, he was nominated for the post on 16 February by "the President and Council of the Royal

[203] Alfred Newton, 1829 - 1907, was an ornithologist, elected FRS in 1870.

Society, *vice* the Rev. B. Price[204] , Master of Pembroke College deceased," and admitted on 26 April 1899. He remained a Governor until 12 March 1924[205] .

### British Association for the Advancement of Science

As mentioned in Chapter 3, MacMahon had joined the British Association for the Advancement of Science in 1883. During his long involvement with the Society, which lasted until 1923, MacMahon is known to have attended 11 annual meetings, as shown in the following table[206]:

| Year | Location |
|---|---|
| 1883 | Southport |
| 1901 | Glasgow |
| 1903 | York |
| 1904 | Cambridge |
| 1905 | South Africa |
| 1907 | Leicester |
| 1908 | Dublin |
| 1909 | Winnipeg, Canada |
| 1910 | Sheffield |
| 1911 | Portsmouth |
| 1912 | Dundee |

Table 4: BAAS meetings attended by MacMahon

He does not appear to have attended meetings on a regular basis until after his retirement from the Royal Artillery in 1898, when he was elected to the Council of the Association[207].

From 1895 to 1901, MacMahon was a member of the BAAS Mathematical Tables Committee, although no evidence that he ever calculated any tables for them has been found.[208] His experience of calculating artillery tables for George Greenhill was the most likely reason that he had been invited to join the committee. MacMahon also held various other offices (described below) until 1913, when he became one of the trustees in charge of the Society's funds.

[204] Bartholomew Price, 1818 - 1898, mathematical tutor at Pembroke College, Oxford, from 1845, was elected FRS in 1852 and was Sedleian Professor of Natural Philosophy at Oxford from 1853 - 1898.
[205] The archives at Winchester College do not contain any details of the activities that MacMahon may have undertaken in this role.
[206] The information was obtained from the BAAS Archives held in the Bodleian Library, Oxford.
[207] Council report dated 7 September 1898 in Box 12 of the BAAS Archive.
[208] [Croarken, 2000].

In 1901 MacMahon was President of the Mathematical and Physical Section A of the 71st meeting in Glasgow from 11 - 19 September[209]. This important post gave MacMahon the chance to speak publicly on a topic of his choosing, with the knowledge that his words would be fully reported in the press. In his presidential address on 12 September[210], MacMahon mourned the death of several leading mathematicians of the time (among them, Hermite, Fitzgerald, Rowland and Tait), before giving an amusing history of the Mathematical Society of Spitalfields, which had been wound up in 1845.

MacMahon then went on to bemoan the poor state of mathematics education in England. He began by drawing attention to the international cooperation between astronomers at the time of the close approach of the asteroid Eros, saying, "to my mind, this concert of the world ... is ... a grand example of an ideal scientific spirit." Such cooperation did not occur in applied physics, he claimed, partly due to the influence of commercial interests, but also as a result of the poor education of engineers, and in particular the bad teaching of mathematics. These remarks on education and their effect are described in more detail below.

There then followed a spirited defence of the scientific specialist, using the discipline of the theory of numbers to illustrate his point that specialists were not to be derided for the narrowness of their study, since the ability to make headway in a difficult field requires an understanding of general principles and an awareness of the tools which are available from other disciplines. In this respect MacMahon said, "he would be found to have realised that analogy is often the finger-post that points the way to useful advance." This is certainly the approach MacMahon himself adopted in commencing his work on plane repeating patterns, described in in Chapter 6.

[209] This meeting sent a message of sympathy to the American President McKinley following the attempt on his life; McKinley died from his wounds three days after the message was drafted.
[210] [MacMahon, 1901, [59;19], pp. 519 - 528].

MacMahon continued the theme of general principles by emphasising the importance of invariance in mathematics. The report finished with a summary of the importance of combinatorial analysis, including remarks on Latin squares, differentiation as a combinatorial process which provided a solution to the 'problème des rencontres', and Diophantine equations recast as partition problems.

Two days after the discussion on education, MacMahon gave a talk entitled *On the partition of series each term of which is a product of quartics.*[211]

***Views on education***

It is worth describing the discussion mentioned above of the state of mathematics education in a little more depth, since it shows MacMahon's interest in educational principles and purposes. He could speak with some experience in this area, having been a teacher for 17 years.

Since the 1860s there had been a debate about the suitability of Euclid's *Elements* in mathematical education, a consequence of which was the formation in 1871 of the Association for the Improvement of Geometrical Teaching[212]. The argument was essentially between those who felt that the logical structure of Euclid and its status as a part of the classical Greek tradition were very good reasons to continue using it as the basis for geometrical teaching and training of the mind, and those who saw that the use of Euclid had been reduced to mere mechanical routines and memory exercises which failed to teach anything about the nature and practical use of geometry. The emphasis on axiomatic structure was not relevant for understanding or applications such as surveying or navigation. A full discussion of the debate is to be found in [Price, 1986, pp. 13 - 41].

MacMahon's comments on the state of education were reported in the press, and were referred to by

[211] The exact content of the paper is lost, as it was not published.
[212] The Association for the Improvement of Geometrical Teaching was founded in 1871 with Thomas Archer Hirst as its first president. It later became the Mathematical Association.

Professor John Perry[213] at a special joint meeting with Section L (Education) to discuss the teaching of mathematics, held on 14 September 1901[214]. Professor Perry had proposed a new syllabus for the teaching of geometry[215]:

> *Geometry*: - Dividing lines into parts in given proportions, and other illustrations of the 6th Book of Euclid. Measurement of angles in degrees and radians. The definitions of sine, cosine and tangent of an angle; determination of their values by drawing and measurement; setting out angles by means of a protractor when they are given in degrees and radians, also when the value of the sine, cosine or tangent is given. Use of tables of sines, cosines and tangents. The solution of a right angled triangle by calculation and by drawing to scale. The construction of a triangle from given data; determination of the area of a triangle. The more important propositions fo Euclid may be illustrated by actual drawing; if the proposition is about angles, these may be measured by means of a protractor; or if it refers to the equality of lines, areas or ratios, lengths may be measured by a scale and the necessary calculations made arithmetically. This combination of drawing and arithmetical calculation may be freely used to *illustrate* the truth of a proposition.

Perry's proposal also covered the measurement of areas and volumes, as well as the use of squared paper, coordinate geometry and vectors. The scheme is very practical, and MacMahon had been very scathing in his Section A address about the thinking behind it: he said that it was an attempt "to overcome defects in training for scientific pursuits by the construction of royal roads to scientific knowledge."[216] In particular, he noted that, "Engineering students have been urged to forgo the study of Euclid, and, as a substitute, to practise drawing squares and triangles." During the discussion on 14 September MacMahon observed that no substitute for Euclid had yet been proposed, even by the Association for the Improvement of Geometrical Teaching[217]. The objects of education, MacMahon asserted, were "the habits of thought and observation, the teaching of how to think, and the cultivation of the memory."

[213] John Perry, 1850 - 1920, taught mathematics and physics at Clifton College, and was elected a Fellow of the Royal Society whilst still a teacher. He was an engineer by training, so was regarded as something of an outsider by the Cambridge mathematical fraternity. Perry was a leader in the reform of mathematical education in the early 20th century, and invented the concept of 'practical mathematics' for use in general secondary education. His 'laboratory approach', with its emphasis on experimental geometry, and the use of squared paper to draw graphs and find slopes, maxima and minima, was not popular with public schools or pure mathematicians. A detailed history may found in [Price, 1986, pp. 42 - 73], and [Price, 2003, pp. 464 - 467].

[214] The participants at this discussion had eschewed the opportunity to go on one of 13 possible excursions to, inter alia, Loch Lomond, the Trossachs, Stirling, Dunoon, or Fossil Grove, Whiteinch.

[215] This version is quoted from [Price, 2003, pp. 465 - 466].

[216] This was a reference to the famous quotation, attributed to Euclid, that 'there is no royal road to geometry.'

[217] This observation was in error, as noted later, although the suitability of the books produced was debatable. Although the Association by then was officially known as the Mathematical Association, Perry, MacMahon and most of the other participants in the discussion continued to refer to it by its old name.

Professor Perry was stung by MacMahon's remarks, and the publicity they had received. In his opening address to the joint Section A and Section L meeting on 14 September[218], he lamented, "I am sorry to think that I have had so little success in explaining my proposed reform." Perry's motives were noble[219]; he felt that existing methods of geometry teaching[220] were damaging to the majority of pupils: "10,000 Toms, Dicks and Harrys mentally destroyed for the sake of producing one man fit to be a mathematical master of a second rate public school; 10 million destroyed for the sake of producing one great mathematician."

MacMahon was not convinced by Perry's defence of his syllabus. He felt that Perry regarded "all students as boys who are to become engineers", and that before Euclid could be abandoned, some suitable substitute should be proposed. He saw Euclid as important training for the mind, and felt that it should be taught alongside geometrical drawing. In defending the use of Euclid in secondary education, MacMahon was in a minority. Most of the other speakers and correspondents, whilst taking Perry to task over points of detail, were in favour of Euclid being dropped, with or without a replacement, and were clear that major reform in mathematics education was essential to the future of the nation. Michael Price[221] has observed that MacMahon was not really in very close touch with ordinary secondary education[222], and there is no record that he was a member of the Mathematical Association.[223]

---

[218] This meeting was chaired by the Rt. Hon Sir John E. Gorst, MP, 1835 - 1916, knighted in 1885, who represented Cambridge University in Parliament 1892 - 1906. At the time of this discussion, he was Vice-president of the committee of the Privy Council on education. His presence reflects the Government's continuing concern over the state of mathematcs education in Britain, which was being compared unfavourably to Continental methods (particularly French and Prussian).

[219] The Perry Movement also crossed the Atlantic. See [Mock, 1963] for details.

[220] This approach was called by Perry the 'Alexandrian method.'

[221] Michael Price is a lecturer in education at the University of Leicester and a historian of secondary mathematics development and the Mathematical Association.

[222] The October 1902 edition of the *Mathematical Gazette* has a report of the British Association Committee on the Teaching of Mathematics, originally presented to the Education Section at a meeting in Belfast, which lists MacMahon as a member of the Committee. So being out of touch with ordinary secondary education had not in anyway disqualified him from being asked to consider the matter.

[223] Private correspondence, 6 August 2003.

Mr E. M. Langley spoke immediately after MacMahon and took him to task for apparently having failed to notice the AIGT's textbook on geometry[224]. Mr Langley was also responsible for the splendid observation that teachers should be paid more, since then "they are not liable to be turned off at forty to attempt to get a living by poultry or fruit, or by editing a socialist newspaper."

The discussions were recorded in a book edited by Perry[225], which included written comments from those who had been unable to attend the meeting in person. Some discussants were also in favour of introducing the history of mathematics into the curriculum. Examiners and inspectors were seen as a force for change by some, and as totally reactionary by others.

It is interesting to note MacMahon's conservative view of the role of Euclid in the teaching of geometry. The two greatest influences on MacMahon, Sylvester and Cayley, took diametrically opposed views on the matter; Sylvester wanted to see Euclid "honourably shelved", whilst Cayley had wanted to see a return to the original treatise[226]. MacMahon's own education at Woolwich included the traditional study of Euclid. The austerity of Euclid and the practicality of a military training are clearly embodied in the way MacMahon regarded mathematics, summarised best by his own words in describing the emergence of analytic number theory at the turn of the century. He spoke of the work of Legendre, Gauss, Eisenstein and H. J. S. Smith[227] in adapting the research tools of continuous quantity for use in number theory, and said that, "these adaptations are of so difficult and ingenious a nature that they are today, at the commencement of a new century, the wonder, and I may add, the delight of the beholders." Noting, however, that the 'beholders' were few in number, he explained thus[228] :

"To attain to the point of vantage is an arduous task demanding alike devotion and courage. I am

[224] *The elements of plane geometry*, published by Swan Sonnenschein in 1884.
[225] [Perry, 1901].
[226] [Price, 1994, pages 23 and 30].
[227] Henry J. S. Smith, 1826 - 1883, was appointed Savilian Professor of Geometry in 1860 and elected FRS in 1861. He worked mainly in number theory and algebra.
[228] [MacMahon, 1901, [59;19], p. 525].

> reminded, to take a geographical analogy, of the Hamilton Falls[229], near Hamilton Inlet, in Labrador. I have been informed that to obtain a view of this wonderful natural feature demands so much time and intrepidity, and necessitates so many collateral arrangements, that a few years ago only nine white men had feasted their eyes on falls which are finer than those of Niagara."

To appreciate the full beauty of mathematics requires similar devotion and intrepidity; Perry, of course, did not feel that the education system should be geared solely to the needs of the intrepid few who would look upon the mathematical equivalent of the Hamilton Falls. There are no records that would show how MacMahon's practice as a teacher fitted in with this philosophical stance. Perry's attitude to education was clearly influenced by his experiences as a teacher of engineering students, but MacMahon's experience of teaching soldiers had left him with a very traditional point of view. There are no surviving archives of MacMahon's individual work as a teacher - only the general practice described in Chapter 2 - and as a non-university man, he had no students who might have commented on his abilities as a teacher. It is possible that the military environment did not encourage dissension, so MacMahon may not have been aware of any difficulties encountered by his students, whereas Perry's engineers may have been more vocal about their problems.

***Other meetings***

In 1902 MacMahon was appointed a General Secretary of the British Association[230], as successor to Sir William Roberts-Austen[231].

At the 74th meeting in Cambridge (17 - 24 August 1904) MacMahon was awarded an honorary degree by the University (see below). He read a paper entitled *On the theory of Partial Differential Equations*, but again, there is no record of the content of this paper.

[229] It was renamed "Churchill Falls" in 1965, and is now the site of a 5.4 Mw hydro-electric power station.
[230] Report to the Council dated 10 September 1902 in the BAAS archive in the Bodleian Library, Oxford. There were sometimes two General Secretaries.
[231] William Roberts-Austen, 1843 - 1902, was a metallurgist, elected FRS in 1875.

In 1905 the British Association annual meeting was held in South Africa[232]. MacMahon was a guest at the Commemoration Dinner on 18 September 1905 hosted by the South African College of Cape Town[233]. There is no evidence that he read a paper on this occasion.

## Evidence of Esteem

### *Cambridge*

As part of the 1904 British Association for the Advancement of Science meeting held in Cambridge, a number of honorary degrees were awarded by the University to distinguished delegates. Joseph Larmor proposed MacMahon for one such award[234] and on 17 August the Senate of Cambridge University passed the following Grace:

> 14. That the degree of Doctor in Science, *honoris causa*, be conferred upon Major Percy Alexander MacMahon, F.R.S., late Royal Artillery, formerly Professor of Physics, Ordnance College[235] , Woolwich, under Statute A, Chapter II, Section 18, Paragraph 3[236].

The award ceremony was held on 22 August 1904 and was reported in the Cambridge University Reporter dated 23 August:

> The Vice Chancellor entertained the recipients of Honorary Degrees on the occasion of the visit of the British Association for the Advancement of Science, and a large party of guests, at luncheon in the Hall of Queen's College, at 1:30 pm.
>
> The recipients of Honorary Degrees:
>
> 1. Dr. Backlund
> 2. Professor Becquerel
> 3. Professor Bruehl
> 4. Professor Engler
> 5. Professor Von Groth
> 6. Dr. Kabbadias
> 7. Professor Kossel
> 8. Professor Osbor
> 9. Dr. Pierson
> 10. Professor Volterra

[232] MacMahon's name appears in the list of guests printed in the *Transvaal Leader,* dated 2 September 1905. He lodged at the Kimberley Club (*Diamond Fields Advertiser* 5 September 1905, which printed this biographical note: "[He] is the General Secretary of the British Association, and is an expert in mathematics and astronomy. He was for some time Professor of Physics at Woolwich").
[233] [Ritchie, 1918].
[234] There are no records to indicate why Larmor chose to make this recommendation.
[235] The Royal Artillery College, Woolwich, was so renamed in 1899. MacMahon was never Professor of Physics, but the Military Instructor in Electricity. There may always have been some confusion over the difference between a Professor and an Instructor at the College, which MacMahon may not have felt worth the trouble of explaining.
[236] [University of Cambridge, 1914]. This statute allows the University to award degrees *honoris causa* to certain classes of people.

11. Sir David Gill
12. Mr. Howitt
13. Sir Joseph Norman Lockyer
14. Major MacMahon
15. Sir William Ramsey
16. Professor Schuster
17. Sir William Turner Thiselton-Dyer

The Cambridge University Reporter carried a notice on 10 August 1904 advertising the event, admission to which was by ticket only. All members of the University were required to wear academic dress, and doctors had to wear scarlet. An oration in Latin was read for each recipient of a degree. MacMahon's was as follows:

> *"Adest deinceps militis insignis filius, miles mathematicus paesertim in studiis spectandus, qui praeter alios laudis titulos etiam scientiarum societatis Britannicae inter ministros praecipuos numerator. Studiorum suorum in caelo puro, in regione illa sublimi a Cayleio nostro feliciter peragrata, diu versatus, studiis illis caelestibus sermonis Latini Musam pedestrem, longe infra in terris relictam, nihil aliud quam numerorum theoriam quandam e longinquo contemplari patitur. Cetera omnia scientiae tam sublimis mysteria, peritis patefacta, a nobis certe palam divulgari non concessum. illud autem unum dixerim. Si, Syracusis captis, Archimedem, intentum formis, quas in pulvere illo erudito descripserat, Marcelli in exercitu miles talis aspexisset, caeli spectatorem illim unicum, tormentorum bellicorum machinatorem illim mirabilem, sine dubio nunquam interfecisset, sed velut socium et fratrem statim esset amplexus.*
> *Duco ad vos militem quem non iam Mars sed Mathesis inter cultores suos numerat, regionis non terrestris victorem, ALEXANDRUM MACMAHON."*

> "The next to come before us is the son of a distinguished soldier, himself a soldier and mathematician whose diligent studies are considered among the members of the British Association for the Advancement of Science to be particularly worthy of praise and recognition. His studies in pure mathematics over a long time, following the ground covered so well by our own well-known and distinguished Cayley, have made the words of the Latin Muses seem pedestrian, left behind by the turning of the world beneath them, have left nothing in number theory unchanged by virtue of his dogged persistence. The lofty mysteries of the whole of science are being expertly brought to light and made plain by our efforts and will not be relinquished. On this we will speak as one. Consider Archimedes, a captive in Syracuse, drawing his forms in the sand, designing for Marcellus those magnificent engines of war, observing the heavens, and who without doubt never harmed anyone but embraced everyone as an ally and even a brother. We consider you to be not a soldier of Mars but number you amongst the husbandmen of Mathematics, a victor in areas not of the earth, ALEXANDER MACMAHON[237]."

Unlike the Trinity College, Dublin, oration cited in Chapter 4, the Cambridge orator appears to have ignored MacMahon's first name. The reason for this is not obvious; it may be that the orator did not feel that there was a suitable Latin version of 'Percy'.

[237] Translation by the author.

### *St John's College, Cambridge.*

Joseph Larmor proposed that MacMahon join St. John's College in 1904, consequent upon the honorary degree awarded at the occasion of the BAAS meeting in Cambridge that year. From that date, MacMahon maintained a loose association with St. John's College, which became stronger when he moved to Cambridge in 1922. There are 11 letters in the MacMahon archive[238] at St John's College, two of which are of interest. One letter to Larmor concerning the Oxford Savilian Professorship has already been described. The other, to H. F. Baker, will be discussed later.

## Mathematical work

### *Partitions*

The second in the series of seven memoirs on the partitions of numbers, *Memoir on the theory of the partition of numbers part II*, was published in 1899[239] and refereed by E. B. Elliott and G. B. Matthews.[240] In this paper, MacMahon generalised the concept of compositions and partitions by considering $n = \alpha_1 + \alpha_2 + \alpha_3 + \cdots + \alpha_s$, where in $\alpha_1 \circ \alpha_2 \circ \alpha_3 \circ \cdots \circ \alpha_s$ the symbol represented by the small circle between each pair can be chosen from the set $\{>, =, <, \geq, \leq, \substack{<\\=\\>}; \lessgtr\}$. In a partition, for example, the symbol is always $\geq$. He developed generating functions for these generalised partitions, and extended many of his results to *two dimensional partitions,* where the parts are placed at the vertices of a square. In a normal one dimensional partition, the parts are placed in non-ascending order; for example, a partition of 8 is (4211); as a two dimensional partition this would be written:

$$\begin{matrix} 4 & 2 \\ 1 & 1 \end{matrix}$$

[238] There is other correrspondence in the Sylvester and Cayley archives, that came to the attention of the author too late for inclusion in this thesis.
[239] [MacMahon, 1899, [54;10]].
[240] Royal Society archives RR14 62 & 63. The referees merely expressed their view that the paper was worthy of publication.

The non-ascending rule now applies from top to bottom as well as from left to right.

An obvious extension is to extend the idea to three dimensions, but this proved to be much less straightforward than the two-dimensional case. This paper contained the fundamentals of what MacMahon called "Partition Analysis" and led to further work over the next six years.

In his 1899 paper, *Partitions of numbers whose graphs possess symmetry*[241], MacMahon took an idea he developed in the first *Memoir*, where a Ferrers graph representing a partition using a two-dimensional array of dots was extended into three dimensions as a representation of a multipartite number by piling the two dimensional diagrams one on top of another. For example, the multipartite number $\left(\overline{13,\ 11,\ 6}\right)$ has partitions of its three parts $\left(\overline{643}, \overline{632}, \overline{411}\right)$ and the three Ferrers graphs are:

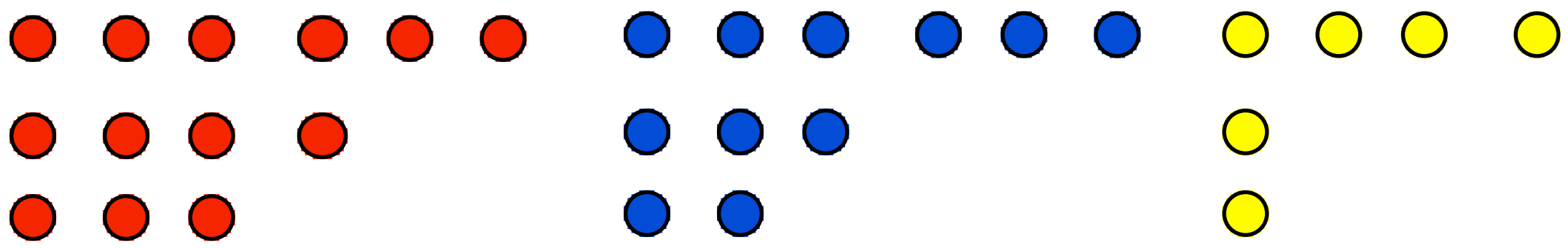

Figure 9: Three Ferrers graph

Stacked one on top of the other, they look like this (seen from above):

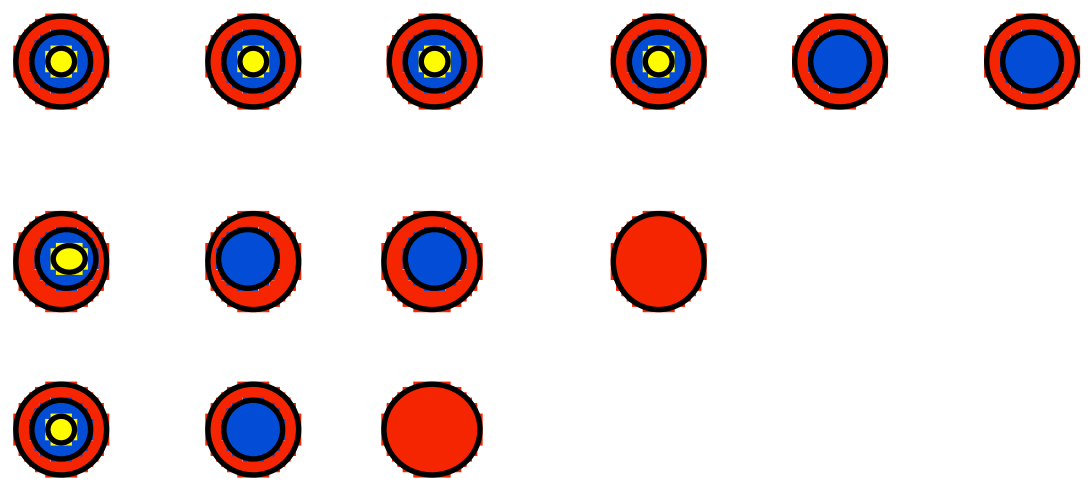

Figure 10: Stacked Ferrers graphs

The particular partition we started with can be thought of as being read in the $z$-direction, out of the page. If $x$ and $y$ are the horizontal and vertical directions in the plane of the page, then an alternative reading of the three dimensional diagram in the $y$-direction (down the page) is the multipartite number

[241] [MacMahon, 1899, [56;12]].

$\left(\overline{16, 8, 6}\right)$, partitioned as $\left(\overline{664}, \overline{431}, \overline{321}\right)$. When read in the $x$-direction from left to right, it is the partition $\left(\overline{333}, \overline{331}, \overline{321}, \overline{211}, \overline{110}, \overline{110}\right)$ of $\left(\overline{9, 7, 6, 4, 2, 2}\right)$[242].

These graphs may be symmetrical in four different ways: in the *xy*-plane, the *yz*-plane, the *xz*-plane, or in all three (*xyz*-symmetry). MacMahon began by considering graphs which are *xy*-symmetrical and made a conjecture, described by Andrews as 'striking'[243], concerning their generating function.

MacMahon used an operator he called $\Omega$, ascribed to Cayley. The effect of this operator is, when an expression of the form[244]

$$a_1 a_2 \ldots a_i x^{2i-1}\left(1 + a_1 x\right)\left(1 + a_1 a_2 x^3\right)\ldots\left(1 + a_1 a_2 \ldots a_{i-1} x^{2i-3}\right) \times \left\{\left(1 + \frac{x}{a_1}\right)\left(1 + \frac{x^3}{a_1 a_2}\right)\left(1 + \frac{x^5}{a_1 a_2 a_3}\right)\ldots\right\}$$

is multiplied out, to delete all the terms containing negative powers of the coefficients and set those remaining equal to unity. An example will illustrate the process. An *xy*-symmetrical graph with at most 4 dots along the $x$ and $y$ axes and two layers can be constructed from 'angles':

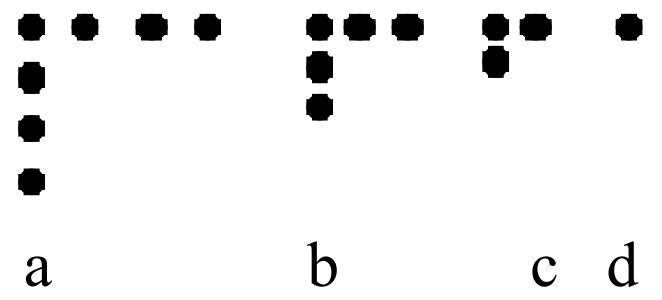

Figure 11: Angle graphs

Each layer of the graph comprises one or more angles nested together. For example, a graph with 11 dots (*weight* 11) in two layers can be made like this:

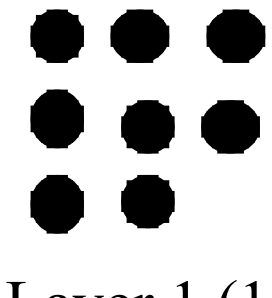

Figure 12: Layer 1 (11b + 11c)

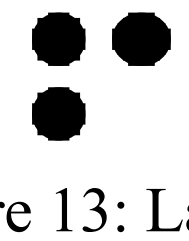

Figure 13: Layer 2 (11c)

[242] This example is taken from MacMahon's 1902 article on combinatorial analysis, written for the *Encyclopaedia Britannica*.
[243] [Andrews, 1977, p. 1380].
[244] This example is the generating function for an *xy*-symmetrical graph with two layers, but with any number of dots.

The generating function that enumerates such graphs is:

$$a_1a_2a_3a_4x^7(1+a_1x)(1+a_1a_2x^3)(1+a_1a_2a_3x^5)\times\left(1+\frac{x}{a_1}\right)\left(1+\frac{x^3}{a_1a_2}\right)\left(1+\frac{x^5}{a_1a_2a_3}\right)\left(1+\frac{x^7}{a_1a_2a_3a_4}\right)$$

The first part of the expresssion enumerates the lower layer, the second part enumerates the upper layer. The coefficients $a_i$ are needed only to pick out the unsymmetrical cases; the operator Ω is used to delete from the product any terms with a negative power of $a_i$; then the remaining $a_i$ are set to unity. In this case the result is:

$$x^7+2x^8+x^9+2x^{10}+3x^{11}+4x^{12}+4x^{13}+5x^{14}+5x^{15}+7x^{16}+5x^{17}+6x^{18}+6x^{19}$$
$$+7x^{20}+5x^{21}+5x^{22}+4x^{23}+5x^{24}+3x^{25}+3x^{26}+2x^{27}+2x^{28}+x^{29}+x^{30}+x^{31}+x^{32}$$

This shows, for example, that there are three ways to construct an *xy*-symmetrical graph of weight 11 with 2 layers. They are (using the angles as labelled in Figure 11):

(1) Layer 1 (b) +( d), Layer 2 (b);

(2) Layer 1 (a) + (c), Layer 2 (d);

(3) Layer 1 (b) + (c), Layer 2 (c).

The next step is to increase the number of layers, so that the generating function becomes:

$$\Omega(1+a_1x)(1+a_1a_2x^3)(1+a_1a_2a_3x^5)\ldots(1+a_1a_2\ldots a_ix^{2i-1})$$
$$\times\left(1+\frac{b_1}{a_1}x\right)\left(1+\frac{b_1b_2}{a_1a_2}x^5\right)\left(1+\frac{b_1b_2b_3}{a_1a_2a_3}x^5\right)\ldots\times\left(1+\frac{c_1}{b_1}x\right)\left(1+\frac{c_1c_2}{b_1b_2}x^5\right)\left(1+\frac{c_1c_2c_3}{b_1b_2b_3}x^5\right)\ldots$$

This generating function is crude, in the sense that it contains far more terms than are needed, and the operation Ω is not carried out algebraically. MacMahon's aim was then to try to work out the form of a 'reduced generating function' that involves only *x*. His strategy was to examine a large number of particular cases to try to guess the form of the required function. From his conjecture, MacMahon deduced that, for example, the number of 'at-most-two layer *xy*-symmetrical graphs of weight 2*w'* with at most *i* dots along the *x* and *y* axes is the same as the number of two dimensional graphs of weight *w* with the same restriction on the dots.

The generalised study begun in the second *Memoir* led to four further papers in 1900, 1904 and 1905[245] in which some algebraic ideas are explored, and to a fifth paper, the third *Memoir*, in 1906. In this paper, MacMahon applied his partition analysis to the study of magic squares, a topic on which he had already written in 1898 and 1900[246].

***Encyclopaedia Britannica article (1902)***

The tenth edition of the *Encyclopaedia Britannica*, published in 1902, contained an extensive article on combinatorial analysis newly written by MacMahon[247]; there is no article on the subject in the 9th edition published between 1875 and 1889.

It began with a history of the subject, which was an expanded version of the history given in the LMS Presidential address described earlier, followed by a detailed description of the theory of distributions and symmetric functions, leading to compositions, partitions and the Master Theorem. MacMahon made several references to his own work during the course of the article, as well as the work of Cayley, Sylvester, F. E. A. Lucas (1842 - 1891) and W. A Whitworth[248]. He illustrated the theory with several worked examples[249], such as the number of compositions of the multipartite number (22), calculated to be 26, but he did not list the actual compositions[250].

The article finished with a section on the construction and use of differential operators in combinatorial

[245] [MacMahon, 1900, [58:10]], [MacMahon, 1904, [62;10]], [MacMahon, 1905, [ 65;10 and 66;10]].
[246] [MacMahon, 1898, [52;7]], [MacMahon, 1900, [57;7]].
[247] [MacMahon, 1902]. Andrews was unaware of the publication of this article in the 1902 10th edition, and has it listed in the 11th edition published in 1910, reference [76;19].
[248] [Whitworth, 1867]. Whitworth, 1840 - 1905, was 16th wrangler in 1862 and vicar of All Saints, Margaret Street, London, from 1885 until his death.
[249] The example above is the one MacMahon gave of the extension of a Ferrers' diagram to three dimensions.
[250] The 26 compositions of (22) are (22), (11)(11), (21)(01), (01)(21), (12)(10), (10)(12), (20)(02), (02)(20), (20)(01)(01), (01)(20)(01), (01)(01)(20), (11)(01)(10), (11)(10)(01), (01)(11)(10), (10)(11)(01), (01)(10)(11), (10)(01)(11), (02)(10)(10), (10)(02)(10),(10)(10) (02), (01)(10)(01)(10), (01)(01)(10)(10), (01)(10)(10)(01), (10)(10)(01)(01), (10)(01)(10)(01), (10)(01)(01)(10)

analysis, illustrated by three examples:

- the differential operator $\left(\dfrac{d}{dx}\right)^n \equiv \delta_x^n$ applied to the function $x^n$ to enumerate the permutations of $n$ different objects;
- the more complex operator $D_p = \dfrac{1}{p!}\left(\delta_{a_1} + a_1\delta_{a_2} + a_2\delta_{a_3} + \ldots\right)^p$ to enumerate the ways of placing numbers in an $n \times m$ lattice so that the sums of each row and column have specified totals;
- the operator $D_{2^n-1}^n$ applied to the function $\left(\sum a_1^{2^n-1} a_2^{2^n-2} \ldots a_{n-1}^2 a_n\right)^n$ to enumerate the number of Latin squares of order $n$.

***Invariant Theory***

As MacMahon's interests moved toward partition theory, his researches in invariant theory declined, and in the period under consideration, he wrote only two papers concerned with invariant theory, both published in 1904.

The first was *On the applications of quaternions to the orthogonal transformation and invariant theory*[251], in which MacMahon derived a more symmetrical form for orthogonal transformation using quaternions, and showed that invariants can be expressed in quaternion form. If two ternary quantics $a_x^n = \left(a_1x_1 + a_2x_2 + a_3x_3\right)^n$ and $A_x^n = \left(A_1x_1 + A_2x_2 + A_3x_3\right)^n$ are connected by the orthogonal transformation

$$\begin{pmatrix} x_1 \\ x_2 \\ x_3 \end{pmatrix} = \begin{pmatrix} p_1 & q_1 & r_1 \\ p_2 & q_2 & r_2 \\ p_3 & q_3 & r_3 \end{pmatrix} \begin{pmatrix} X_1 \\ X_2 \\ X_3 \end{pmatrix}$$

then the transformation may be represented by the quaternion identity $x_i = \left(\alpha_0 + \alpha_i\right)X_i\left(\alpha_0 - \alpha_i\right)$,

[251] [MacMahon, 1904, [63;18]].

where $\alpha_0$ is the real part and $\alpha_i$ is the imaginary part of a quaternion, with $\alpha_0^2 + \alpha_1^2 + \alpha_2^2 + \alpha_3^2 = 1$. Further work shows that the problem of finding invariants for these orthogonal ternary functions may be reduced to finding functions of the quaternions which are just scalars (i.e., have zero imaginary part). The second paper, *Semivariants of systems of binary quantics, the order of each quantic being infinite*[252], extended earlier work on *perpetuants*[253] of single binary quantics to systems of quantics.

## Social life

MacMahon led a busy social life. In March 1903 he became a member of the Athenaeum Club, having been nominated by Sir William Huggins[254], proposed by Professor G. H. Darwin[255] and seconded by Professor R. B. Clifton[256]. The citation in support of the nomination read:

> Major Percy Alexander MacMahon, R.A. (Retired), D.Sc. Dublin, F.R.S. Distinguished for the number and range of his contributions to the Mathematical Sciences. The generality of his methods has opened up wide vistas in the higher algebra and the theories of combinations and partitions, and has led to the solutions of problems that have been outstanding since the time of Euler. In 1900 the Royal Society awarded him a Royal Medal "for his work, which is distinguished for great originality, research and precision."

In his obituary for MacMahon, the astronomer Herbert Hall Turner recalls MacMahon's fondness for a game of post-prandial billiards; this is borne out by the fact that MacMahon was billiards champion of the Athenaeum in 1919 and 1922.

On 1 October 1904, MacMahon chaired the morning session of the sixth meeting of the San Francisco Section of the American Mathematical Society held at the University of California, and gave a talk entitled "*Groups of linear differential operators*."[257] There is no record of how this came to pass, but since his daughter, Florence, had been living in California since 1900, it is possible that the invitation to speak gave him the opportunity to visit her.

---

[252] [MacMahon, 1904, [64;18]].

[253] A perpetuant is an irreducible seminvariant of a binary form of infinite order.

[254] William Huggins, 1824 - 1910, Astronomer, was President of the Royal Society from 1900 - 1905.

[255] G. H. Darwin, 1845 - 1912, son of Charles Darwin, was second Wrangler in 1868, and Plumian Professor of Astronomy and Experimental Philosophy at Cambridge from 1883 - 1912.

[256] R. B. Clifton, 1836 - 1921, was President of the Physical Society from 1882 - 1884, and designer of the original Clarendon Laboratory in Oxford.

[257] Sources: *Science*, 4 November 1904, Vol XX No 514, p. 604, and the *Bulletin of the American Mathematical Society*, 1904, Vol XI, pp. 68 - 69.

## Chapter 6 Board of Trade 1906 - 1922

This chapter deals with the period during which MacMahon was employed as Deputy Warden of the Standards at the Board of Trade. By this time he was a mature mathematician with an established repuation, which was such that he was consulted for advice or help by eminent people. His work for the Electoral Commission, his correspondence with D'Arcy Thompson, and the help he provided to G. H. Hardy and S. Ramanujan are described.

### Board of Trade

The Board of Trade had existed in various forms since the early 1600s. The form of the Board that MacMahon knew was created by William Pitt the Younger in 1786. Responsibility for the maintenance of the weights and measures standards was given to the Board of Trade by the Standards Act of 1866.

With effect from 6 November 1905, MacMahon was appointed as Deputy Superintendent of Weights and Measures at a salary of £300 p.a. on a one year contract[258]. There are no records in the Board of Trade archives from 1904 and 1905 which might shed any light on how or why this appointment was made. It is not known whether he applied for the post, or was invited to take it up. The Board had been seeking to appoint an assistant to the Superintendent of Weights and Measures since mid-1904. As MacMahon was about to remarry - and thus require more money than his military pension provided - his involvement with intellectual and social life in London (such as the Royal Society and the Athenaeum Club) would have placed him in an ideal position to be considered for a job requiring analytical skills and a meticulous attention to detail. There was a tradition of mathematicians' involvement in such work. For example, in 1795 Legendre was placed in charge of the Temporary

[258] Archive BT5/111 in The National Archive, Kew.

Agency of Weights and Measures by the National Convention of the French Revolution, to oversee the introduction of metric measures into France[259].

On 13 February 1906 the Superintendent, H. J. Chaney, died. The Treasury then approved a promotion for MacMahon, effective from 14 February[260], to the post of Deputy Warden of the Standards at a salary[261] of £600.

Until this appointment, MacMahon's only related work had been a recreational paper on weighing by a series of weights in 1890. In January and February 1907, after a year in post, MacMahon gave two talks on weights and measures to the Royal Institution, but these took place on a Thursday afternoon and no record of their content was kept. It is not possible to tell, therefore, how a year of professional work in the field had affected his outlook on the subject.

Much of MacMahon's work required the interpretation of law or statutory regulations and the resolution of disputes. His specialised mathematical knowledge was not explicitly called upon, but the analytical and logical skills he had acquired were very valuable. Tact and patience were also necessary, and it is clear that occasionally these were tested to the limit. Some of the more interesting duties and episodes from this period in his life are described below.

***MacMahon's duties***

The Deputy Warden of the Standards was head of the Standards Department of the Board of Trade[262]. The duties were to be responsible for the imperial, parliamentary and secondary weights and measures standards, to provide guidance to the Board on its functions under the legislation, and to verify the accuracy of weights and measures used by traders. In addition to routine verifications, the courts would often seek the advice of the Board on questions of the accuracy or efficiency of weights and measures.

[259] [Alder, 2002, p. 162].
[260] Archive T9/35 p. 832, p. 844, p. 882, The National Archive, Kew.
[261] Ibid BT5/112.
[262] The previous incumbent had been Hugh Young, father of the mathematician William Henry Young who later became the husband of Grace Chisholm.

The Deputy Warden was also required to keep abreast of the latest developments in measuring technology, which MacMahon did from the very start. In the period 1907 to 1909, MacMahon commissioned for the Standards Department a comparator that used the wavelength of light to measure distances. This device had been invented by A. E. H. Tutton[263] , the crystallographer, in 1898. In a paper read before the Royal Society in 1909[264] Tutton described the construction of the instrument for the Standards Department, and wrote, “The author desires to record his appreciation of the most kind help and encouragement invariably received from the Deputy Warden of the Standards, Major P. A. MacMahon, F.R.S., to whose initiative the whole of this advance in fine measurement at the Standards Department is due.” MacMahon had seen the potential of Tutton’s interferometer, and persuaded him to design and construct the comparator.

A regular duty of the Deputy Warden of the Standards was to supervise the “trial of the pyx” every year. The Goldsmiths’ Company has a statutory responsibility, under the Coinage Act, to undertake the annual examination of coins manufactured by the Royal Mint. The pyx[265] refers to the large chests in which samples of coins minted in ancient times were placed following minting. These chests were originally housed in the Pyx Chamber in Westminster Abbey, where the early trials were held. Since 1878 the trial has taken place at Goldsmiths’ Hall. The jury, comprising members of the Goldsmiths’ Company, was summoned by the Queen’s Remembrancer (a senior judge at the High Court of Justice) and was a proper court of law. In February officials from the Royal Mint brought to the Hall boxes containing samples of the coins, put aside during the previous year. The coins were required to conform to standards of size in diameter, weight and composition of metal, and this was verified by the

[263] Alfred Edwin Howard Tutton, 1864 - 1938, was HM Inspector of Technical Schools, Oxford, from 1895 - 1905, and was elected FRS in 1899.
[264] [Tutton, 1911].
[265] From the Latin: *pyxis* = chest or box.

London Assay Office. The verdict was delivered separately, usually in May, when the Master of the Mint (the Chancellor of the Exchequer) learned whether coins made in the previous year were up to standard. The freshly minted gold sovereign and half-sovereign coins were assayed individually for purity, after which a random number of coins were melted into an ingot, which was itself then reassayed. A typical report is this one made by MacMahon, in 1914:

> ...the millesimal fineness of the gold ingot produced by melting 127 sovereigns and 31 half sovereigns was 916.66676, the standard being 916.666… and the mean millesimal fineness of the gold coins tested separately was 916.755, the fineness of individual coins varying from 916.51 to 916.94

The dubious use of statistics in the trial of the pyx was an issue where MacMahon's mathematical expertise was useful. It was customary to calculate an overall allowable error in the pyx trials by simply summing the products of the number $n$ of each coin by the allowance $r$ for that coin, $\sum nr$; this was known as the "Master's Remedy". MacMahon pointed out[266] that the square root of the sum of the products of the number of coins and the square of the allowance, $\sqrt{\sum nr^2}$, would make more sense, since any errors would be distributed normally. It is not clear from the surviving documents whether this suggestion was adopted.

The Deputy Warden also met regularly with delegations from trade bodies, to receive complaints and petitions. For example, in 1908 a delegation from the local authority in Durham visited the Board to meet with MacMahon and Winston Churchill MP, who was President of the Board at that time. The delegation complained that the statutory charge of 6d a dozen for 'stamping' glassware[267] was too high, and that, even with the maximum allowed rebate of 3d per dozen, it was double the actual costs incurred. The surplus money retained by the stamping authority was felt to be an unnecessary tax on local manufacturers; the authority was actually granting an extra but illegal rebate to prevent the manufacturer from taking his business elsewhere.

[266] Archive BT101/734 in The National Archives, Kew.
[267] 'Stamping' refers to the process of applying an official mark to, say, a beer glass to show that it has been tested and will hold a full measure. In the case of glassware, the mark is applied by sandblasting.

he purpose of the delegation was to persuade the President of the Board of Trade to alter the legislation to allow for more local variation in the charge. Durham was able to do the job very cheaply by simply renting a space in the one factory where the glassware was made. The London stamping authority, on the other hand, maintained special premises to which manufacturers and importers were obliged to bring their glassware, which resulted in much higher costs. Upon his arrival in post, MacMahon had conducted a survey of the costs[268], and had determined that 6d per dozen would be a fair average cost. He was determined to defend this position, which eventually brought him into a polite but tetchy exchange with Churchill. It was clear that MacMahon did not appreciate the difficulty that Churchill was having grasping the issue. Churchill, with the skill of a politician, avoided making a decision by calling for more information, but intimating that he was minded to accede to the delegation's request should the new data support their position. MacMahon was charged with obtaining details of the quantity of imported foreign glassware stamped in London.

The matter rumbled on, with further meetings between the Board of Trade and representatives of the Sunderland, South Shields and Gateshead local authorities, who were very upset by Durham's disregard for the regulations, which they had applied meticulously. They were agreeable to a maximum rebate of 4d per dozen, so long as it was applied uniformly - that is, so long as Durham cooperated. At a further meeting with Durham, it was suggested that the President of the Board of Trade was minded to allow the 4d rebate, but Durham wanted a larger rebate and suggested that a conference be held in London at the Board of Trade, between themselves and the Tyneside authorities, who should supply detailed information about their costs. The Tynesiders were agreeable to this, although reluctant to provide cost data, but before the conference could be set up, Durham lost interest. A rise in alcohol tax meant that publicans in Durham had decided to sell beer by the 8 fluid ounce glass, rather than by the half pint. Selling beer in this way meant that the glasses did not need to be stamped, so the factory had stopped making half pint glasses, and all stamping activity there had ceased.

[268] Archive BT101/640 in the The National Archive, Kew.

Another regular issue was the question of metrication. Metric measures had been legal since 1897, but not compulsory, and there were recurrent calls for them to be imposed on the country. Then, as now, there were those who were vehemently opposed to any such compulsion. In 1907, the British Weights and Measures Association (BWMA) made an application to have the weight of a cubic inch of water declared a legal measure[269], as part of their campaign to defeat the Weighing and Measuring Bill, then about to have its second reading in the House of Commons. The letter from the BWMA was circulated around the senior officers of the Board, who all agreed that it was not a matter with which to bother the Prime Minister; the last note scrawled across the back of the letter read, "We must have Major MacMahon." This suggests that MacMahon's mathematical expertise was highly regarded at the Board of Trade, and his presence at the discussion of such very technical matters was felt to be essential.

An 11-page memorandum on the weight of a cubic inch of water was prepared by another official, Mr H. B. C. Darling, who was later to be the proof reader for MacMahon's *Combinatory Analysis*. The matter was raised again in 1910 by Mr Thomas Parker, who wanted to use one thousandth of an inch as the standard length, and define the weight of a cubic inch of water as 250 grains at 122°F. MacMahon could see no value in this at all[270], partly because the temperature suggested was not one normally encountered in everyday life.

A group opposed to compulsory metrication were the Cotton Traders. In 1907 MacMahon attended a meeting with the Cotton Traders Association and the President of the Board, David Lloyd George MP, where these views were expressed. This is further evidence that MacMahon was known to several of the leading politicians of the era.

[269] Archive BT101/660 in The National Archive, Kew.
[270] Archive BT101/730 in The National Archives, Kew.

A briefing paper prepared by the Board for the second reading of the Bill included this comment from MacMahon:

> In discussing the matter with retail traders it is not easy to point out to them the advantages that would accrue to them by a compulsory metric system, and no advantage at all would apparently be the lot of numerous poor; for a very long period certainly they would be subject to all kinds of fraud as a result of the change. As an instance take the case of liquid measures of capacity; under a compulsory metric system a man could not go to a public house and ask for a pint of beer; he would have to ask for either half a litre, which is 7/8 pint or for 0.6 litre which is more than a pint (about 21/20 pint); in neither case is it likely that he would get the value for his money that he does now. It is, moreover, not certain that 0.6 litre would be a legal denomination of the metric system; it is not so at the present moment. It has been calculated that the provisions of local metric standards by local authorities would cost about £60 000 and the cost to traders of metric weights and measures would be little short of £2 000 000 sterling. The cost to manufacturers and engineers of new plant is not ascertainable approximately, but is certain to be large.

The Bill was dropped. MacMahon's final report on the subject was made in 1916, when he expressed the view that the costs and confusion that would arise from compulsory metrication would outweigh any benefit to the export trade. In any event, metric measures were perfectly legal, so there was nothing to prevent exporters from using metric measures in their dealings with countries using the metric system[271].

This was not the only peculiar set of proposals that crossed MacMahon's desk. In 1912, the Weights and Measures Committee of the Central and Associated Chambers of Agriculture suggested that the pound and the 'cental' (100 lb) be made the only legal measures of weight, with a 'ton' of 2000lb or 20 centals.[272] In 1917, a petition from some colleges and public schools called for a change to decimal coinage[273]. One of the arguments put forward in support of this latter proposal was that much time was wasted teaching the system of pounds, shillings and pence. MacMahon's comment on this was: "Although a compulsory metric system would save time by simplifying arithmetic it may be questioned whether any time is <u>wholly</u> wasted by the Imperial system as the mental training is useful."

[271] Archive BT101/819 in The National Archives, Kew.
[272] Archive BT101/762 in The National Archive, Kew.
[273] Ibid - BT101/825.

In the former case, MacMahon's view was that the industry could adopt whatever system it liked for trading purposes, since the pound was already a legal measure, but that abolition of the hundredweight (1 cwt = 112 lb) and the existing ton (2240 lb) would serve no useful purpose.

The problem of short weight being given by traders, both retail and wholesale, was another issue that occupied MacMahon. In 1907, Howard Cunliffe, the local inspector from Smethwick, wrote to MacMahon with a detailed report on the weights of wrapping paper used by merchants, and the practices of paper salesmen in persuading retailers to use heavy or light paper depending on whether they sold by nett or gross weight[274]. In 1913 a report in *The Grocer* gave details of the unfair practice of using very heavy paper to wrap sausages sold by weight[275]; in the same year MacMahon met with a delegation from the National Federated Association of Fruiterers and Florists, who wanted legislation to enforce the printing of weights or counts on wrappers, so that merchants selling short weight could then be prosecuted[276]. MacMahon recommended that the legislation be restricted to the retail trade. A Select Committee was appointed, to which MacMahon gave evidence in 1914[277].

In 1909, a complaint about short weight deliveries of coal from collieries was placed before the Board by a deputation of coal merchants[278]. Many reasons for the problem were given. Collieries bulk loaded and weighed the coal, but misdeclared the *tare*[279] weight of the trucks (this was measured when the trucks were clean and dry, but not adjusted when they were wet and dirty). The weight of coal had to be checked by the merchants when it was off-loaded into smaller cart or sack loads, and the weighing errors were uncontrollable. Any loss of weight was blamed by the railway companies on the collieries, and the collieries blamed pilferage during transit. Many railway companies made handsome profits by reselling coal that had fallen off trucks during violent shunting manoeuvres, it was alleged by the merchants. The most annoying aspect was that merchants transporting coal by road could be inspected at any time and prosecuted if their load was underweight, but such regulations did not apply to transport by rail. MacMahon suggested that the merchants should fund a band of roving inspectors who could pay surprise visits to collieries to check their weighing apparatus and practices, since it was felt that legislation would be very difficult to draft - for example, how would loss of coal due to a truck striking the buffers be dealt with, or the theft of coal by children from trucks left in sidings overnight? The situation was complicated by the change in practice recently introduced by the railway companies, where they allowed only 20 cwt to the ton, instead of 21, because they felt that they were being charged too much for the coal used by their locomotives. MacMahon instituted a set of detailed studies of the practices at collieries and dockyards around the country. The inspector made a final report in late 1912, recommending that weighbridges be inspected more carefully, that checkweighers be employed to prevent fraud, the regular reassessment of the tare weight of wagons, and the covering of coal in transit, or at least careful trimming of the load before weighing, to avoid accidental loss and pilferage[280].

[274] Ibid - BT101/658.
[275] Ibid - BT101/780.
[276] Ibid - BT101/792.
[277] Archive BT101/798 in The National Archive, Kew.
[278] Archive BT101/695 in The National Archives, Kew.
[279] The tare of a container used for commercial goods sold by weight is the empty weight marked upon that container.

MacMahon's obituarists and his *Dictionary of National Biography* entry all state that he worked for the Board of Trade until 1922. However, MacMahon was of retirement age in 1920, and after the end of that year he was no longer on the circulation list for internal memoranda. In 1921, J. E. Spears was Deputy Warden of the Standards, and attended the 1921 international conference on weights and measures in that capacity. As his weights and measures swansong MacMahon attended the same conference, but as the British representative on the *Comité Internationale des Poids et des Mésures*. In his report, dated 18 October 1921[281], he described the financial position of the Bureau, and pressed for the Committee to take no action on the standardisation of electrical units before the next conference.

280 Archive BT101/778 in The National Archives, Kew.
281 Archive BT101/843 in The National Archives, Kew.

**The cube law for elections.**

Amongst analysts of elections, there is a 'rule' that in a two party system with constituencies, where each constituency is won separately[282], the ratio of the number of seats won by each party varies as the ratio of the cubes of the number of voters for each party. This is the 'cube law', attributed by the Glasgow MP James Parker Smith to MacMahon. It will be shown that MacMahon was reluctant to become involved with the problem of electoral analysis, and that the method by which he obtained this result is obscure.

In 1908, the Government set up a Royal Commission to investigate electoral reform. The problem of electoral reform was not new. In the early 1880s there had been a great deal of interest in the introduction of proportional representation to replace the system of constituencies with many members, and with voters who were allowed to vote in more than one constituency. The choice at the time was between 'redistribution', the creation of roughly equal sized constituencies with one or two members, and proportional representation using larger multi-member constituencies. Charles Lutwidge Dodgson had written many letters and pamphlets in favour of proportional representation, but redistribution had taken place as a result of the Reform Bill, in time for the General Election of November 1885.

The dramatic swings in elections, particularly the huge Liberal victory of 1906, resulted in Parliament once again wanting to consider alternative forms of election. A Royal Commission under the chairmanship of Lord Richard Cavendish was set up, "to examine the various schemes which have been adopted or proposed, in order to secure a fully representative character for popularly elected legislative: and to consider whether, and how far, they, or any of them, are capable of application in this country."

[282] This method of election is often known as 'the first past the post' system.

***James Parker Smith***

In 1909 MacMahon was contacted by the Rt. Hon. James Parker Smith[283] in connection with the Royal Commission on electoral systems. Parker Smith had been chosen to give evidence before the Commission, probably because he had published two articles in the *Law Magazine and Review* in 1883 and 1884 entitled *University Representation*, describing the history of the representation of universities in Parliament, and had also contributed to the arguments about proportional representation leading up the 1885 Reform Bill.

MacMahon and Parker Smith had probably met at London Mathematical Society meetings, and they were certainly both present at the 71st British Association for the Advancement of Science meeting in Glasgow in 1901. Parker Smith had been elected a Vice-president of the meeting[284] and MacMahon was President of Section A. On 11 January 1909[285] Parker Smith wrote to MacMahon requesting his assistance with the work of the Royal Commission. Parker Smith described the problem of the unrepresentative proportions of MPs. Parliament was interested in the possibility of introducing proportional representation to replace the 'first past the post' system, which had produced some dramatic swings in recent elections. It was clear that the majority party was getting a disproportionately large share of the seats in the Houses of Parliament, to the detriment of minority interests and the reputation of democracy. Parker Smith explained the difficulty by using the analogy of a large bin of white and red balls in the ratio 53:47. A sample selected at random is likely to have a majority of white balls. If each constituency is represented by a sample drawn at random, without replacement, then it is clear that the whites will have a larger majority than the original proportion. Parker Smith reminded MacMahon that he (Parker Smith) had worked on the problem twenty years earlier, but claimed that his "mathematical equipment is not sufficient to deal with the question adequately, and moreover is now quite rusty." The Commission chairman, Lord Cavendish, wanted a

[283] James Parker Smith, 1854 - 1929, was an MP for Glasgow. He was the eldest son of the mathematician Archibald Smith, and had been a member of the London Mathematical Society since 1881. He knew Cayley [Crilly, private note 2006], and was private secretary to Joseph Chamberlain.
[284] The election was recorded in the minutes of the meeting of the BAAS council held on 7 September 1900 in Bradford.
[285] [Parker Smith, 1909]. A full transcript of this letter is given in Appendix 7.

"mathematician of high rank" to work on the problem.

MacMahon's cautious reply was made two days later[286]. He was concerned about the time the calculations would take, and felt that he was spending quite enough time working for the Government as Deputy Warden of the Standards work at the Board of Trade. He said he would think about it, and consider who else might be interested in working on the problem.

***The evidence for the Cube Law***

There are no further letters to or from MacMahon in the Parker Smith archive, but on 19 May 1909 Parker Smith himself gave evidence to the Royal Commission[287] . In paragraph 1253, Smith used the analogy of a 'great box' filled with 11000 blue marbles and 9000 red marbles to represent the electorate. The votes of each constituency were represented by taking a shovelful of marbles and letting each side score a win for those shovelfuls with a majority of that side's colour. He continued:

> I have been going into the mathematics of it, if you are willing to take it from me, and the chances of a single marble being drawn blue is 11 to 9; but the chance of a whole shovelful having a blue majority is something very near 2 to 1. I have had the help in working this of my friend Major MacMahon, who is one of the leading mathematicians of the day, and he gives me this as the formula: that if the electors are in the ratio of A to B, then the members will be at least in the ratio $A^3$ to $B^3$.

It seems that, despite his misgivings, MacMahon had indeed worked on the problem for Parker Smith. In their 1950 paper *The law of cubic proportion in election results*, M. G. Kendall and A. Stuart[288] remark, "How MacMahon arrived at his result is a mystery."[289]. A possible derivation is given below, although the 'law' does not work very well in practice.

Further evidence for MacMahon's interest in the problem is in an undated and unsigned document in the Parker Smith Archive, attached to the letters above. It is in MacMahon's handwriting and is decribed below[290].

[286] [MacMahon, 1909]. A full transcript of this letter is given in Appendix 7.
[287] Royal Commission on Systems of Election, Minutes of Evidence, 19 May 1909, pp. 77 - 86.
[288] [Kendall, 1950].
[289] This comment is repeated in a letter to the *Observer* on 4 November 1951.

The calculation set out in MacMahon's document was intended to show that the composition of a sample of $n$ voters taken at random from the whole set $A$ of voters has a high probability of being within $r$ voters of the composition of the set $A$. It began by considering a problem known as the *Scrutin de Ballotage*, in which the probability was required of the winning candidate in an election always being ahead during the counting. MacMahon had already written about the problem in 1894[291], and had spoken about it at the BAAS meeting in Dublin in 1908[292].

Equation (1) from the document represented the probability $C_{pq}$ of getting voters in the ratio $p{:}q$ in a sample of size $n$ chosen from a population of size $A$, $b$ of one persuasion and $c$ of another. Stirling's formula was used to approximate the factorials.

$$C_{pq} = \frac{b!c!\left(A-\left(p+q\right)\right)!}{\left(b-p\right)!\left(c-q\right)!A!} \times \frac{\left(p+q\right)!}{p!q!} \qquad (1)$$

Using the assumption that the voters were equally distributed in $A$, further simplifications were made to enable the calculation of values of $C_r$, the probability of the sample containing $r$ more voters of one sort than the exact proportion. From this the sum $S_r$ of all values of $C_r$ from $+r$ to $-r$, could be calculated, giving a cumulative probability that the sample would be within $r$ of exact proportionality. The values of $C_r$ and $S_r$ were tabulated for an electorate of 10,000 voters whose sympathies were equally divided between two parties.

It is not clear whether the calculation in the document led MacMahon to the cube rule, which would amount to the claim that MPs would be in the ratio $b^3{:}c^3$. Kendall and Stuart[293] observed, "If

[290] [MacMahon, undated] A full trancript of the document is given in the Appendix 7.
[291] [MacMahon, 1894, [46;5]]. A detailed discussion of MacMahon's work on the problem of the *Scrutin de Ballotage* follows later in this chapter.
[292] There is no record of the content of the talk that MacMahon gave in Dublin in 1908.
[293] [Kendall, 1950].

MacMahon regarded the law as a phenomenon of simple sampling of the kind mentioned by Smith then it could not have been derived mathematically", and later, "It is unthinkable that a man of MacMahon's calibre should have made an error in mathematics."

However, in a paper written in 1980, R. G. Stanton of the University of Manitoba[294] showed that MacMahon could indeed have derived the result. Stanton's argument calculated the expected number of constituencies with a win for the majority party, using arguments similar to those used by MacMahon, and showed that this does indeed agree with an approximate cube law.

So although MacMahon's original argument is lost, the combined evidence of the document and Stanton's paper shows that he could easily have derived the cube law, contrary to Kendall's view. A recent paper[295] shows that a 'square root of three rule', where voters in the ratio *a:b* produce seats in the ratio $a^{\sqrt{3}}:b^{\sqrt{3}}$, would be a more accurate model of recent election results. This paper mentions the work of Kendall, but attributes the rule to Parker Smith, rather than to MacMahon. It builds on the work of R. Taagepera[296], who showed how a cube rule might arise, as follows. It is not known whether MacMahon used an argument of this nature, but it does show that the cube law can be derived from fairly simple considerations.

Taagepera assumed that for two parties *i* and *j* with $V_i$ and $V_j$ votes and $S_i$ and $S_j$ seats, the value of *n* in the relation $\frac{S_i}{S_j} = \left(\frac{V_i}{V_j}\right)^n$ can be estimated by $\frac{\ln V}{\ln S}$, where *V* and *S* are the total number of voters and seats respectively. This can be estimated by a combinatorial argument, different from Stanton's, which runs as follows. Each MP has $\frac{V}{S}$ constituents, so has $\frac{2V}{S}$ channels of two way communication with those constituents, plus $2(S-1)+\frac{(S-1)(S-2)}{2} \cong \frac{S^2}{2}$ channels with all the other MPs. The minimum value of $\frac{2V}{S}+\frac{S^2}{2}$ occurs at $2V = S^3$, so that $n = \frac{\ln V}{\ln S} \cong 3$.

[294] [Stanton, 1980]. Stanton was less accurate with his history than his mathematics; he described Parker Smith as a 'semi-anonymous civil servant' rather than a prominent Glasgow MP from a famous Glasgow family.
[295] [Mahoney et al., 2003].
[296] [Taagepera et al., 1989].

In 1910 the Royal Commission reported, and recommended that the alternative vote system should be adopted in single member constituencies, and that two member constituencies should be abolished as soon as possible.

## Honorary degrees

### *University of Aberdeen*

At a meeting of the Senate on 28 February 1911, the University of Aberdeen decided to award MacMahon an LL.D degree; the minutes do not record who proposed him for this award. The degree was conferred at a ceremony held on 5 April at Mitchell Hall, Marischal College. The event was reported the following day in the *Aberdeen Daily Journal*:

> **Major P. A. MacMahon, D.Sc., F.R.S.**
> Major MacMahon holds a position of proud eminence among the pure mathematicians of this country. His contributions to that science, remarkable in number, in importance, and, above all, in originality, have elucidated many provinces of mathematical learning and furnished new demonstrations of the far-reaching power of mathematical analysis. In particular, his brilliant memoirs on the theory of numbers and the theory of algebraic forms vindicate his title to rank in the direct line of succession to the great English mathematicians of the last generation. (Applause.) As Instructor in Mathematics in the Royal Military Academy, and, later, as Professor of Physics in the Ordnance College, Major MacMahon proved signally successful in training others in the principles and methods of science of which he himself is so great a master. Of recent years his special faculties and attainments have been utilised in a very practical sphere of public service; for, in 1906, he was appointed Deputy Warden of the Standards, in which capacity he is charged with the task of testing and verifying the precision of our instruments of weight and measure - a task of great importance alike to our science and our commerce. (Applause.) Major MacMahon has filled the president's chair of the London Mathematical Society, and holds a Royal medal from the Royal Society. To these distinctions the University of Aberdeen now rejoices to add her Doctorate of Laws. (Loud applause.)

### *St Andrews University*

On 14 September 1911 MacMahon was also awarded an LL.D. degree by St Andrew's University, at a ceremony to celebrate the 500th anniversary of the foundation of that University. A document in the possession of Professor Antony Unwin, a distant relative of MacMahon, reads:

> *Universitas Sancti Andreae*
> *Viro Pretissimo Percivallo Alexander MacMahon*
> *S. P. D.*
> *Nos, Universitatis Andreanae Cancellarius, Vice-Cancellarius, Rector, Facultatum Decani, Colegorium Praefecti, Professores, cum anno salutis MCCCCXI Academia Nostra, in hoc regno Scotorum vetustissima, a viro admodum reverendo Henrico de Wardlaw, Episcopo Andreapolitano, condita sit, mox autem Apostolica et Regia auctoritate confirmata, ut rite celebretur almae matris natalis quingentesimus, constituimus viros doctos undique advocare, ut caerimoniis nostris quingenariis intersint. Oramus igitur et summo studio a te petimus ut per festos dies quos in mensem septembrem MCMXI indiximus hospitium nostrum accipere digneris*
>
> *Dabamus Andreapoli, Mense Martio MCMXI*
> *(signed by the Chancellor, Rector and Vice-Chancellor)*

This is simply an invitation to accept the honorary degree and attend the anniversary celebrations. In this case the writer has not shirked the task of translating 'Percy' into Latin, although the use of 'Percival' as a starting point may be regarded as cheating. The award was made by the Chancellor, the Right Hon. Lord Balfour of Burleigh. MacMahon was one of 96 candidates "capped" that day; it is likely that he was recommended for the honour by D'Arcy Wentworth Thompson[297] , who was Professor of Biology at Dundee at the time, and became Professor of Natural History at St Andrew's in 1917. He was the author of a seminal work, *On Growth and Form*[298] which is regarded as having founded the discipline of biomathematics. MacMahon and Thompson had been friends for some time and corresponded regularly.

The University itself has no record of the Lauration address, and the local newspaper reported the event thus: "The graduands were now called in turn to the platform and duly 'capped' by the Chancellor, who with Principal Donaldson, Principal Stewart, and Professor Scott Lang (who introduced them) cordially shook them by the hand."

---

[297] D'Arcy Wentworth Thompson, 1860 - 1948, was elected FRS in 1916 and knighted in 1937.
[298] [Thompson, 1997].

### The Royal Irish Academy

MacMahon was recommended for Honorary membership of the Royal Irish Academy by its Science Committee at a meeting on 15 February 1917. There are no records to shed any light on the reason for this. The election was confirmed on 16 March 1917, and acknowledged by MacMahon on 21 March 1917 thus:

> Sir, I have the honour to acknowledge the receipt of your letter informing me that I have been elected an Honorary Member of the Royal Irish Academy, and also a Diploma. I beg to assure you of my deep appreciation of the great honour that has been conferred upon me. I am, Sir, Your obedient Servant, Percy A. MacMahon

There is no record of his having participated in the work of the Academy, and the only other mention of him in the minutes is on 15 March 1930, noting his death.

### British Association for the Advancement of Science

MacMahon was a fairly regular attender of BAAS meetings between 1901 and 1912[299], and remained an active member of the Association until 1923. He often gave talks, and was entrusted with supervision of the Association's finances for a decade from 1913.

#### *1907*

At the 1907 meeting in Leicester MacMahon presented the draft report of the Council but did not present any mathematical papers. A photograph of the Council of the BAAS standing on the steps of the town hall in which MacMahon is standing on the far right, was published in the local press.

#### *1908*

At the 1908 meeting in Trinity College, Dublin[300], MacMahon gave a talk entitled *On a generalisation on the question in probabilities known as Le Scrutin de Ballotage*. No record of the content of this talk exists, although the ideas are covered in his 1909 *Memoir on the theory of partitions of numbers Part*

[299] See Table 4 on p.72.
[300] On 4 September 1908.

*IV*[301].

A brief biography of MacMahon appeared in the *Irish Times* of 2 September 1908, which drew attention to MacMahon's Irish roots. A second talk was mentioned in a report by the same newspaper on 8 September, but no details of the title were given.

***1909***

The 1909 meeting was held in Winnipeg, Canada. At the inaugural meeting on the evening of 25 August, the proceedings were opened by Professor George Carey Foster[302], who called upon MacMahon to read a letter from the retiring President, Francis Darwin[303], who was unable to attend[304]. This shows the high esteem in which MacMahon was held. Later, MacMahon gave a talk entitled *On a correspondence in the theory of the partition of numbers.*

The meeting was widely covered in the local Canadian press. There are articles in the *Manitoba Free Press* (26 June - a 'word sketch' of MacMahon), the *Winnipeg Telegram* (28 August - an article with a picture of MacMahon and a quote from him; 30 August - an article describing a meeting of the Ladies' Circle at which Mrs MacMahon was a 'top table' guest) and the *Calgary Daily Herald* (7 September where 'MacMahon' is misspelled as 'McMann')[305].

***1913***

In 1913 MacMahon gave up the post of General Secretary[306] and became a Trustee of the British Association. This meant that he shared responsibility for the financial affairs of the Association. There is no record of why he chose to take this action.

***1920***

[301] [MacMahon, 1909, [73:10]].
[302] George Cary Foster,1835 - 1919, was professor of natural philospohy at Strathclyde University from 1862 - 1865, and an honorary member of IEE in 1916.
[303] Francis Darwin, 1848 - 1925, was a botanist, a son of Charles Darwin; he was knighted in 1913.
[304] Source: *Science*, Vol XXX No 767, p. 351, 10 September 1909.
[305] BAAS archive Box 421.
[306] The committee appointed to appoint a successor drew up a short list which comprised J. Larmor, A. E. H. Love, A. W. Porter and H. H. Turner. The post went to Turner.

In 1920, the committee appointed to organise the next meeting proposed that MacMahon be invited to present a paper "on the lines of his recent Royal Institution lectures." This would have been his work on the design of repeating patterns, described in detail later, but there is no record of MacMahon having presented a paper at the 1921 Edinburgh meeting.

***1923***

MacMahon was still active as a trustee in 1923, after his move to Cambridge, as evidenced by a letter[307] to the Secretary of the Society in the BAAS archive dated 11 February concerning some investments.

## Mathematical work

### *Partition theory*

MacMahon's second memoir on compositions was published in 1908[308]. In this paper, he dealt with the problem, drawn to his attention by Simon Newcomb[309], which he stated thus:

> A pack of cards of any specification is taken - say that there are *p* cards marked 1, *q* cards 2, *r* cards 3 and so on - and, being shuffled, is dealt out on a table; so long as the cards that appear have numbers that are in descending[310] order of magnitude, they are placed in one pack together - equality of number counting as descending order - but directly descending order is broken a fresh pack is commenced, and so on until all the cards have been dealt. The result of the deal will be *m* packs containing, in order, *a*, *b*, *c*, ... cards respectively, where, *n* being the number of cards in the whole pack, (*abc*...) is some composition of the number *n*, the number of parts in the composition being *m*.

Among the questions to be answered are, how many ways of arranging the cards yield a given composition and how many arrangements lead to *m* packs? These could also be regarded as recreational problems.

To solve the problem, MacMahon introduced an hybrid adaptation of the Ferrers graph and the line of route graph described above, which he called a *zig-zag graph*. The composition (3321) of 9 is

[307] BAAS archive Box 111.
[308] [MacMahon, 1908, [71;5]].
[309] Simon Newcomb, 1835 - 1909, was an American mathematician and astronomer who is reputed to have played a solitaire card game that gave rise to this problem in combinatorial analysis.
[310] In the statement of this problem in his *Combinatory Analysis*, MacMahon used 'ascending' instead of 'descending.'

represented by horizontal rows of dots thus:

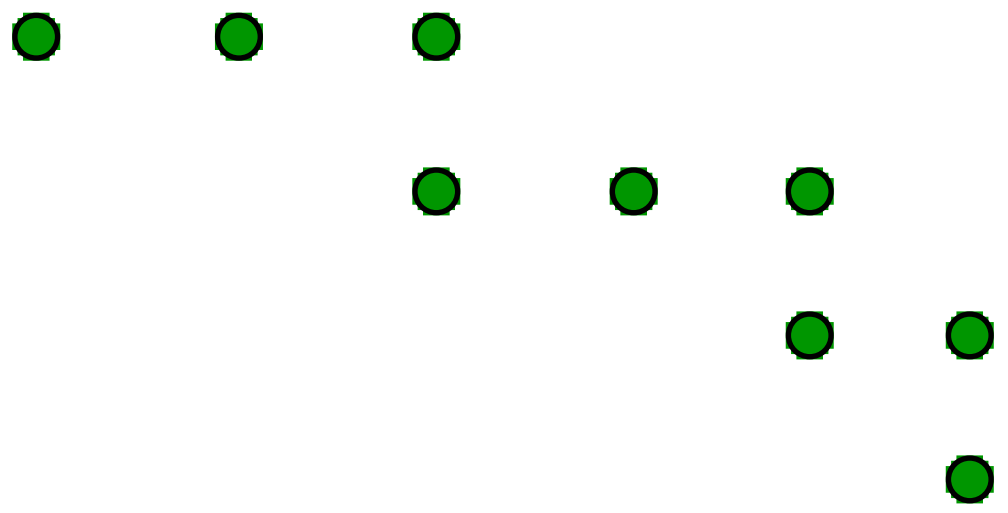

Fig. 14. Zig-zag graph

The conjugate composition (112122) can be obtained by reading the vertical columns from left to right. MacMahon was confident that this representation of compositions would become "an important instrument of research"[311].

The fourth *Memoir*[312] in 1909 was devoted to the application of the partition analysis to various ballot or voting problems. Use was made of the graphical representation of composition first described in the 1894 *Memoir on the composition of numbers* to deal with 'certain questions of probability'. The question was one posed by J. Bertrand and solved by Désiré André[313] in 1887:

> Peter and Paul are the candidates in a secret ballot; the urn contains *m* votes for Peter and *n* votes for Paul; *m* is greater than *n*, and Peter will be elected. What is the probability that, during the counting of the votes, Peter is always ahead?

The number of votes $m$ and $n$ can ber represented as an $m \times n$ rectangular lattice or *reticulation*, and the problem can be recast as counting the number of lines of route (compositions of ($mn$)) that lie below, and do not cross, a line at 45° from the bottom left corner of the reticulation. Having solved this problem[314], MacMahon extended it to any number of candidates and concluded that, for $n$ candidates who receive votes in their favour of $a_1, a_2, \ldots, a_n$, (where $a_1 \geq a_2 \geq \ldots \geq a_n$), the probability that the order of candidates at any point in the election is the same is $\prod_{s=t+1}^{s=n} \prod_{t=1}^{t=n-1} \left(1 - \frac{a_s}{a_t + s - t}\right)$.

[311] There is no evidence that this hope was realised.
[312] [MacMahon, 1909, [73;10]].
[313] "Solution directe du problème résolu par M. Bertrand" in the *Comptes rendus*, 1887, vol. 105, pp. 436 - 437. André was Belgian, born in 1829, the author of *Developpements en series des fonctions elliptiques et de leur puissances*, published in 1877, and *Des notations mathématiques: enumeration, choix et usage*, published in 1909 in Paris.
[314] The solution involves Catalan numbers.

It was in the two following *Memoirs*[315], parts V and VI, both read in 1911 and published in 1912, that MacMahon developed this work to what Andrews has described as his 'crowning achievement' in the subject[316]. This is a generating function for the number of *plane partitions*[317] of $n$.

A plane partition of $n$ is a partition whose parts are arranged in non increasing rows and columns. For example, the partitions of 4 are (4), (31), (22), (211) and (1111). The plane partitions are

| | | | | | | | | | | | | |
|---|---|---|---|---|---|---|---|---|---|---|---|---|
| 4 | 31 | 3 | 22 | 2 | 211 | 21 | 2 | 1111 | 111 | 11 | 11 | 1 |
| | | 1 | | 2 | | 1 | 1 | | 1 | 11 | 1 | 1 |
| | | | | | | | 1 | | | | 1 | 1 |
| | | | | | | | | | | | | 1 |

The aim was to find a generating function for the number GF($l$, $m$, $n$) of plane partitions with no part greater than $l$, and with $m$ rows and $n$ columns. To facilitate the task, MacMahon used a notation developed by Cayley and Sylvester, where the expression $\left(1-x^{s}\right)$ is written simply as (**s**). This has the effect of concentrating the attention on the exponent only.

Only two generating functions had been established by MacMahon at the time of writing the paper: GF($l$, 1, $n$), for partitions written on a line, and GF($\infty$, 2, $n$) a result known also to A. R. Forsyth[318]. This latter result admits of two interpretations: either two rows, $n$ columns and no restriction of the size of the parts, or unrestricted columns, $n$ rows and no part greater than 2.

The next paper on partitions, *On compound denumeration*[319], was published in 1912. Here, MacMahon brought his work on the theory of distributions and symmetric functions from 1887 and 1888 to bear on the problem. *Compound denumeration* was a term coined by Sylvester to describe the partition of multipartite numbers. No new graphical methods are described in this work, which is purely algebraic in character. The seventh and final *Memoir*[320] was of a similar character, but did not follow until 1917. There is no evidence that MacMahon was involved in any way with work related to the First World War which might have contributed to this delay, so it was probably due to the demands of the post of Deputy Warden of the Standards and the work involved in preparing *Combinatory Analysis* for publication.

[315] [MacMahon, 1912, [77;11]] and [MacMahon, 1912, [78;12]].
[316] [Andrews, 1977, p. 1315].
[317] Plane partitions are similar, but not identical, to the two dimensional partitions introduced on p. 79, as the parts are not constrained to be at the vertices of a square.
[318] Alfred Russell Forsyth,1858 - 1942, was senior Wrangler in Cambridge in 1881, and Sadlerian Professor of Pure Mathematics from 1895 - 1910.
[319] [MacMahon, 1912, [79;8]].

After a further gap of three years, MacMahon published *Divisors of numbers and their continuations in the theory of partitions*[321], in which he generalised the concept of a divisor in a novel way. MacMahon sought generating functions for power series of the form $A_k = \sum a_{n,k} q^n$ where the coefficients $a_{n,k}$ have a definite meaning in the theory of partitions. The formal definition is $a_{n,k} = \sum s_1 s_2 \ldots s_k$, where the products $s_1 s_2 \ldots s_k$ are all possible $k$-tuples satisfying $\sum_{i=1}^{k} s_i m_i = n$. If $k$ = 1, then $a_{n,1}$ is the sum of the divisors of $n$ and $A_1$ is a power series whose coefficients are the values of the divisor function, usually known as $\sigma(n)$ or $d(n)$. The generating function for this series was already well known, but MacMahon generalised it to all values of $k$, and to several other series that he denoted by $B_k$, (where each coefficient represents the excess of the sum of the odd divisors of $n$ over the sum of the even divisors), $C_k$ (where each coefficient is the sum of the divisors with odd[322] conjugates; that is, in the product $sm = n$, $m$ is an odd number), and so on. The paper concludes with some tables of the coefficients for the various series for $n \leq 16$ and $k \leq 5$.

The connection with partition theory arises as follows: if we have $sm = n$, then $m^s$ is a partition of $n$. Similarly, if $m_1 s_1 + m_2 s_2 + \cdots + m_k s_k = n$, then $m_1^{s_1} m_2^{s_2} \ldots m_k^{s_k}$ is a partition of $n$. For example, $(2 \times 2) + (1 \times 2) = 6$, and 2211, written as $2^2 1^2$, is a partition of 6.

320 [MacMahon, 1917, [87;8]].
321 [MacMahon, 1920, [93;16]].
322 MacMahon himself almost always uses 'uneven' where we would use 'odd'.

In the same year, 1920, MacMahon published a short paper, *On the partitions into unequal and uneven parts*[323] in which he also generalised Euler's well-known theorem that the number of partitions of an number into unequal parts is the same as the number of partitions into odd parts. The generalised result is that if **N** is the set of all positive integers, **O** is the set of all odd positive integers, and **P** is any set of prime numbers, and we delete from **N** and **O** all those elements that are divisible by elements of **P**, then we are left with two sets of numbers from **N** and **O** that are not divisible by any element of **P**. The number of partitions of an integer *n* into distinct parts from the first set is then the same as the number of partitions of *n* into parts from the second set, with repetitions allowed.

***Combinatory analysis***

<u>Combinatory Analysis - Volumes 1 & 2 (1915 & 1916)</u>

MacMahon collected together his work on symmetric functions and partition theory in two volumes, published in 1915 (at 15 shillings) and 1916 (at 18 shillings). These two volumes form MacMahon's most important and influential work, but before more is said about MacMahon's reasons for writing the book, its content and legacy, the contemporaneous reviews of the work written by mathematicians at home and abroad are described, in order to show that MacMahon's work was already known and held in high regard, and that the books met fully the expectations of his peers. Even more than a decade after publication, reviews of books by other authors saw fit to mention *Combinatory Analysis*; for example, in 1929[324] the anonymous reviewer of *Statistique Mathématique* by G. Darmois[325] wrote:

> The reviewer has long been of the opinion that the method of generating functions as extended by MacMahon in his classical treatise on Combinatory Analysis should be of immense value in theoretical statistics; but, unfortunately, highly trained mathematical statisticians - and only highly trained mathematicians could handle so delicate an instrument - have so far not been interested in the matter.

[323] [MacMahon, 1920, [94;9]].
[324] *Journal of the Royal Statistical Society*, Volume 92, pp. 100 - 101.
[325] Georges Darmois, 1888 - 1960, was a French applied statistician, and President of the International Statistical Institute at the time of his death.

At the time of its publication, volume 1 was reviewed anonymously in *Nature*[326]. The unknown writer considered himself privileged to have heard MacMahon speak, and hoped that volume 2 "will not be long delayed."

A longer and more detailed review was written by an author known only as "C" in *Science Progress in the Twentieth Century*[327]. After describing the contents, C remarked that "any one who wishes to measure the magnitude of the task achieved in the book is recommended to take the article on this subject published in the *Encyclopédie des Sciences Mathématiques* and reflect upon the changes which are now required to bring that portion of this great work into line with the subject of which it now treats." C was clearly an admirer of MacMahon, for the review concluded with a lament that MacMahon had never been able to direct a research school, and had not been able to inspire a generation of young mathematicians in the manner of Cayley. This sad state of affairs was attributed to the poor organisation of mathematical teaching in England.

In July 1915, G. B. Mathews[328] reviewed Volume 1 for the *Mathematical Gazette*[329]. Mathews described MacMahon as "a past-master in every kind of symmetrical algebra." He remarked particularly on the consistency of the approach to a wide variety of problems. Mathews also reviewed Volume 2[330]. He commented particularly on MacMahon's use of graphical methods, and hoped that the work would encourage younger analysts to make further discoveries.

The book was also reviewed in America and on the Continent. In the *Bulletin of the American Mathematical Society,* W. V. Lovitt[331] described the content of the book, and noted that "in general the work is well written and clear." In the *Comptes Rendus*, E. Cahen[332] remarked that combinatorial analysis was a subject only barely familiar to Continental mathematicians, and was more the *oeuvre* of English mathematicians - Cayley, Sylvester, Hammond, Whitworth and Dickson were mentioned by name. However, although Cahen acknowledged that MacMahon was the master of the subject, ("on ne saurait trouver dans l'étude de ces questions un guide plus qualifié"), he noted that the subject was not fashionable and many would regret the effort expended on pointless problems. By way of defence, he went on to say[333], "Mais pourquoi sont-elles futiles ? Est-ce parce qu'on y parle d'échiquiers et de jeux des cartes? Ce ne sont la que des représentations. Si l'on y parlait d'arrangements de molécules, la théorie prendrait toute de suite un air plus sérieux, tout en restant la meme au fond. Si elle n'a pas d'applications maintenant, elle en aura peut-etre plus tard." As an analogy he pointed out that the Greeks were also engaged in solving pointless problems that two thousand years later were invaluable in establishing modern astronomical theories and in building cars and aeroplanes. This is evidence that MacMahon was considered to be working in a field whose time had yet to come.

[326] [*Nature*, vol. 96, 30.12.1915, p. 478].
[327] [Vol. 10, 1916, pp. 601 - 606]. MacMahon himself was very pleased with this review, and wondered who C was (Source: Letter dated 24 April 1916 to Ronald Ross from LSHTM Archive).
[328] George Ballard Mathews, 1861 - 1922, was senior Wrangler in 1883 (St John's College), and taught at Bangor and Cambridge; he wrote, among others, *Theory of Numbers*, 1892.
[329] [*Mathematical Gazette*, 1915].
[330] [Mathews, 1917].
[331] [Lovitt, 1917].

G. Scorza[334] reviewed the book for *Scientia*[335]. He traced the origins of combinatory analysis to early work by Paciolo, Cardano and Tartaglia, later extended by Fermat, Pascal and Huygens, but noted that it had few adherents in 1915. He ascribed this state of affairs to the lack of a coherent theory to weld the many disparate problems together. MacMahon's work was precisely the solution to this difficulty, well thought out and deep.

Volume 2 was also reviewed in *Nature*[336], the anonymous reviewer concluded with the observation that, "by collecting these researches, which are very much his own, from their hiding places in scientific memoirs into these two volumes, the author has done much towards the promotion of a more

[332] [Cahen, 1916].
[333] This translates as: "But why are they pointless? Is it because they speak of games of chess and cards? These are nothing but representations. If one speaks of arrangements of molecules, the theory immediately takes on a more serious character although it is basically the same. Even if there are no applications now, there will be in the future."
[334] Gaetano Scorza, 1876 - 1939, taught in secondary schools; he was a member of the editorial board of *Rendiconti del Circolo Matematico di Palermo*.
[335] [Scorza, 1918].
[336] [*Nature*, Vol. 99, 29 March 1917, p. 82].

general outlook on the whole range of analytical work usually classed somewhat vaguely as 'algebra'."

The reviewer of Volume 2 in *Science Progress*[337] was again the enthusiastic C, who welcomed the publication of Volume 2 so quickly after Volume 1 with the explanation, "it has been an unamiable weakness of English mathematicians often to interrupt their tasks for so long that the first parts of their works have become tired of waiting for the companion volume, and settle down into a state of unblessed singleness upon our shelves." C went on to say how the second volume met the "high expectations based upon the part first published, and suggested that "the two volumes are so much one book that they might have been bound together within a single cover."

The most effusive praise was in the *Journal of the Royal Society of Arts*[338]: "Major MacMahon has proved himself easily first amongst his contemporary combinarians; so far as we know, he is only equalled (in this field) by Fermat, Pascal, Euler, E. Lucas, and Sylvester. And he is certainly the first to *methodise* the theory. All who are interested will eagerly look forward to the completion of this highly delectable and artistic work."

An indication of the importance of the work is the fact that a reprint of both volumes bound as one was published by Chelsea in 1960. This edition was reviewed in 1963 by T. H. Southard for the American Mathematical Society, who commented, "Because of the revived interest in combinatorial analysis, including such items as orthogonal Latin squares, finite projective geometries, and other discrete-variable problems which are being studied both theoretically and experimentally (with the use of high-speed digital computing machines), this work will prove to be a worthwhile addendum to the modern literature."

Southard's prediction was correct; the two volumes are MacMahon's most cited work. A search of

337 [*Science Progress*, Vol. 11, 1917, p. 695 -697].
338 [Anonymous, 1916]. This review was also reproduced on the dust jacket of the 1930 reprint of *New Mathematical Pastimes*.

the *Science Citation Index (Expanded Version)* in August 2003 showed about 200 citations of the original volumes, and 95 references to the 1960 Chelsea reprint[339], compared with just 79 citations of all MacMahon's other works. Many of these references are in papers in journals of number theory, combinatorial analysis and discrete mathematics, as might be expected, but there are also citations in papers from journals on probability, operational research and statistical subjects (which would no doubt meet with the approval of Darmois' reviewer cited above). To put this in some context, Arthur Cayley's only book, *Elementary Elliptic Functions*[340] is cited only 7 times, and J. J. Sylvester, who wrote no books on mathematics, has just 49 citations of his papers.

The only other contemporary book on the subject was Eugen Netto's[341] *Lehrbuch der Combinatorik*, first published in 1901. Netto made some references to MacMahon's work (and two chapters added to the book in 1927 by Skolem & Brun made reference to *Combinatory Analysis*), but Netto was writing a summary of the work done in the field up to that point, not trying to extend and generalise the topic. MacMahon, in his introduction to Volume 2, recommended Netto's book as an adjunct to his own work, for the reader who "desires a comprehensive account." However, he pointed out that *Combinatory Analysis* has a differerent purpose, the "presentation of processes of great generality, and of new ideas."

It is not proposed to describe in detail the content of the books, but rather to outline why they were so important to the development of combinatory analysis. It is important to note that the books were not designed as traditional textbooks to be used in teaching, in the manner of, say, G. H. Hardy's *A Course of Pure Mathematics*; there were very few examples and no exercises. The books are important because of the generalisations and extensions of the idea of a partition they contain. MacMahon did not concern himself with, for example, the work of Cayley and Sylvester in finding expressions for the

339 The accuracy of the references is not guaranteed; many authors refer to both volumes but attribute only one date (1915 or 1916).
340 [Cayley, 1876].
341 Eugen Netto, 1848 - 1919, was a pupil of Weierstrass, Kronecker and Kummer, and a Professor at Strasbourg, Berlin and Giessen.

coefficients of the generating functions, since this was more arithmetical than algebraic.

MacMahon set out a reason for writing the books quite clearly in the introduction to Volume 1. He felt that insufficient effort had been made to co-ordinate the work that had been done in the field: "Little attempt has been hitherto made either to make a general attack on the territory to be won or to co-ordinate or arrange the ground that has been already gained[342]."

In private, MacMahon gave an different reason for writing the book. During the course of correspondence with Ronald Ross[343] about Ross's own mathematical work, MacMahon wrote: "The important work in Combinatory Analysis that I did 29 years ago and at which I have been working at intervals ever since has received no notice and was not even mentioned in Netto's *Combinatorik* published 16 years ago - The result has been tho' that I have collared the 'stuff' and others haven't and it was in order to rub this in that I have written my two volumes." This comment shows that despite the public acclaim accorded to MacMahon, he was acutely aware that he was working in an unpopular field.

The achievement that MacMahon wanted to celebrate was the marrying together of algebra, or the mathematics of the continuous, and the 'higher arithmetic' or theory of numbers, essentially the mathematics of the discontinuous. It should be noted here that in the 19th Century, the modern distinction between algebra and analysis had not yet been fully worked out, and the mathematics of the continuous referred to by MacMahon is what we would now call analysis. The question of the misclassification of combinatory analysis described earlier was touched upon again in the introduction to Volume 1, but only as a way of explaining that the theories contained in the book were essentially algebraic in character, and that the ideas of number theory came into play only during the actual

342 Introduction to Volume 1, second paragraph. His use of military language is interesting, and suggests that MacMahon may have seen mathematics as a struggle to conquer, much as he may have struggled to help conquer India during his time there with the Royal Artillery.

343 Ronald Ross, 1857 - 1932, discovered the connection between mosquitos and malaria in 1897; he was awarded the Nobel prize for medicine in 1902 and was knighted in 1911.

evaluation of a particular coefficient.

Laplace was credited with the idea of creating an enumerating generating function, to assist him with his researches in probability theory. Here again, MacMahon showed his knowledge of the history of mathematics and his appreciation of the importance of a grasp of that history as an aid to understanding the subject matter. Of course, it was Euler who had first drawn attention to the use of generating functions in partition theory, but it was Laplace who took the idea and developed it more fully for use in his researches into probability theory[344].

MacMahon reiterated the importance of the differential calculus to combinatorial analysis, by virtue of the ability to design differential operators to achieve a particular transformation. A simple example is the operator $\left(\frac{d}{dx}\right)^n \equiv \delta_x^n$ to enumerate permutations that he mentioned in the *Encyclopaedia Britannica* article described earlier. The use of the Master Theorem to solve the generalised 'problème des rencontres', and the ease with which Simon Newcomb's problem may be solved by means of differential operators, were mentioned specifically. The latter problem was chosen because "from the time of Euler to that of Cayley inclusive, its solution was regarded as being beyond the powers of mathematical analysis." MacMahon was clearly proud of his achievements in this field, and made a point of saying that the work was original, whilst acknowledging the foundations laid by Cayley, Sylvester and Hammond.

In the introduction to Volume 2, MacMahon again refers to the history of the subject, recommending a study of Netto's *Combinatorik*[345]. But he was at pains to point out that the work he was presenting was of greater generality than any previously published, and that it contained new ideas. MacMahon hoped that these new ideas would provide "starting points for further investigations in an exceedingly

[344] [Laplace, 1779].
[345] [Netto, 1927].

interesting field of further mathematics."

This introduction also corrects an historical error made in the first volume where, as a result of a statement in *Grundzüge der Antiken und Modernen Algebra* by Ludwig Matthiessen (published in Leipzig in 1878), MacMahon had been led to ascribe Waring's formula to Albert Girard. MacMahon's attention had been drawn to this matter by Dr R. F. Muirhead.[346] By consulting the original sources MacMahon was able to show that Matthiessen had been wrong on two counts: Girard had *not* had a general formula, and Waring *had* given a proof. This was further evidence for the importance that MacMahon accorded to the history of the subject.

H. B. C. Darling[347], a colleague at the Board of Trade with whom MacMahon published two joint papers[348], was acknowledged as the proof reader for both volumes.

*Introduction to Combinatory Analysis* (1920)

This short book was published four years after Volume 2 of *Combinatory Analysis*. Its purpose was to overcome the objections of critics who had found the algebra of symmetric functions "difficult or troublesome reading." MacMahon apologised in the introduction that "much must be written that will appear to the reader to be self evident and unworthy of statement"[349]. The book comprises six short chapters, totalling 71 pages:

I. Elementary theory of symmetric functions
II. Opening theory of distributions
III. Distribution into different boxes
IV. Distribution when objects and boxes are equal in number
V. Distributions of given specification
VI. The most general case of distribution

346 Robert F. Muirhead, 1861 - 1941, taught at Göttingen, Glasgow and Cambridge. He was Ferguson scholar in 1882, and second Smith's prize winner in 1886.
347 Horatio B. C. Darling, b. 1866, also published a number of papers on Ramanujan's work.
348 [MacMahon 1918, [90:17]] and [MacMahon, 1919, [91:17]].
349 MacMahon was assisted in the composition of the book by Professor John E. A. Steggall, 1855 - 1935, who taught at Clifton College, Manchester and Dundee, and was first Smith's prize winner in 1878.

The chapters take the reader slowly and clearly through the elementary theory of symmetric functions and the recasting of partition theory as distribution theory. Although MacMahon claimed that he was explaining a "complicated if not difficult matter to an untrained mind", he nevertheless assumed in the reader a high degree of confidence with algebraic notation and manipulation.

A review was printed in *Nature* in May 1921.[350] In it, Dorothy Wrinch[351] noted that "this theory must enter into such a question as the formation of a muddy liquid ... and might very well prove another example of the extraordinary way in which abstract mathematics leads the way in applied science." The book became a recommended text of the American Statistical Society in 1924[352] .

*Concluding remarks about the books*

MacMahon's achievement in effectively creating the modern discipline of combinatory analysis cannot be understated. The above reviews show that this was recognised by his contemporaries.

As described by Norman Biggs in *Roots of Combinatorics*[353], the components of combinatory analysis (counting, permutations, combinations, and partitions) have an ancient history going back some 4000 years[354]. The Chinese were concerned with counting arrangements of symbols in the *I Ching* of the 7th century BC. There is evidence that similar problems were studied by Hindu mathematicians as early as 628 AD. Magic squares were written about by Chinese and Arab mathematicians from as early as the 1st century AD, and interest in partitions can be traced back to medieval gamblers. Many famous mathematicians have been involved with combinatorial problems: e.g. Bhaskara, Cardano, Galileo, Descartes, Fermat and Euler. But it was MacMahon who brought all the disparate parts together in the early 20th century, provided a method for solving problems, and suggested interesting avenues for

350 [Wrinch, 1921] This review was reproduced on the dustjacket of the second (1930) printing of *New Mahtematical Pastimes*.

351 Dorothy Wrinch, 1894 - 1976, a Wrangler in 1916, explained protein structure using mathematical principles. She was a fellow of Lady Margaret Hall, Oxford.

352 *Journal of the American Statistical Society*, Volume 19, number 147, p. 401, September 1924. It cost $2.60 in America, 7/6 in the UK .

353 [Biggs, 1979].

354 For example, Problem 79 of the Rhind Papyrus appears to be a counting problem.

further exploration.

## Recreational Mathematics

### *New Mathematical Pastimes (NMP)*

This book was first published in 1921, at a cost of 12 shillings, with a second edition[355] in 1930, the year after MacMahon's death. Albert A. Bennett reviewed it for the *American Mathematical Monthly* in September 1922[356], describing it as "an interesting little volume filled with strange and bizarre figures, and punctuated with quaint quotations in verse." This is an accurate summary: there are 114 pages, only 13 of which have no illustration, and there are 57 quotations scattered throughout.

A review also appeared in *Nature*[357] anonymously. MacMahon was described as "... the author of the well-known 'Combinatory Analysis' ..." Although there is praise for the book, "Everything [MacMahon] writes is carefully finished, and recreations invented by him are sure to be worth attention on their merits...", there is a caveat: "The questions considered in the second and third parts are of a more technical character, and are likely to appeal to the specialist rather than to the general reader: to the former they will open new and interesting lines of development."

The text is divided into three parts. The first part describes the triangle puzzles which were the subject of the patent number 3927 discussed in Chapter 4, and the extension of the idea to squares, right-angled triangles, pentagons and hexagons. Also described is the cube puzzle which was the subject of the 1893 paper, *On the thirty cubes constructed with six differently coloured squares*[358]. What makes this different from other puzzle books is that it is not a collection of puzzles to be solved; it is instead a detailed explanation of how the puzzles may be constructed and solved, giving insight into how the art of mathematical thinking can be applied to a recreation.

355 Strictly speaking, the second edition was simply a reprint, since it was published posthumously with no changes from the orginal 1921 impression. Indeed, the dustjacket was headed 'second impression'.
356 [Bennett, 1922].
357 [*Nature*, 1922, vol 109, pp200 - 201].
358 This paper is also described in Chapter 4.

An essential feature of the method is the subdivision of the shape in view (equilateral triangle, square, etc.) into (congruent) compartments, and then the matching of the compartments according to a rule, called a *contact system*. For example, the square can be divided into four compartments, numbered 1, 2, 3, 4, and then such squares can be arranged in a number of different ways, by matching the numbers 1 - 1, 2 - 2, 3 - 3, 4 - 4, or 1 - 2, 3 - 3, 4 - 4, or 1 - 2, 3 - 4 or 1 - 4, 2 - 3.

The second part of the book shows how the idea of edge-matching using colours can be replaced by altering the profile of adjoining edges to force the desired contact system. This method was used to produce some exotic jigsaw puzzles.

The final part of the book deals with the creation of repeating patterns, that would nowadays be described as tilings, beginning with the three regular tilings of the plane into equilateral triangles, squares or regular hexagons. From these *bases*, which MacMahon observed occur "constantly before our eyes", he described how the methods developed in the first two parts of the book can be used to create an infinite variety of ways of subdividing the shapes and altering their boundaries according to certain rules. He emphasised the role that symmetry has to play, and briefly mentioned the extension of the methods to three dimensions, where he explained how a space-filling array of cubes might be transformed into an array of rhombic dodecahedra. MacMahon also drew attention to the importance of such matters in crystallography.

Many of the planar designs in this part of the book are reminiscent of the designs of M. C. Escher[359], although they predate Escher's work by several decades. For example, Doris Schattschneider in *Tiling the Plane with Congruent Pentagons*[360] described the pattern that appears as (b) on page 101 of *New Mathematical Pastimes* as having "special aesthetic appeal". A favourite of Escher's, it is said to appear as street paving in Cairo and is the cover illustration for the first edition of Coxeter's *Regular Complex Polytopes*. Escher found the inspiration for it in a 1923 work of F. Haag[361] in the *Zeitschrift für Kristallographie*, and it appears in his notebooks from the 1930s, described by Schattschneider in *Visions of Symmetry*[362] on page 28. The detail and symmetries of the construction of such patterns were taken up in greater depth in the papers described in the next section.

[359] Maurits Escher, 1898 - 1972, was a Dutch artist, famous for his drawings of regular divisions of the plane and his use of optical illusions.
[360] [Schattschneider, 1978].

The book concluded with a short set of recommendations concerning the construction of pieces for the puzzles, including suggested materials and colours.

***Repeating patterns***

This work comprises a suite of three papers which expand upon themes from *NMP*. All were written or read in 1922, the year after *NMP* was published.

MacMahon's 1922 paper, *The design of repeating patterns - Part I*[363], was coauthored[364] with his nephew, William Percy Dartray MacMahon[365], the son of his youngest brother Reginald. W. P. D. MacMahon also published a paper in his own right, printed in the *Proceedings of the London Mathematical Society* in 1925[366], which is described in this section. These two papers, and the paper *The design of repeating patterns for decorative work* which was read[367] before the Royal Society of Arts (RSA) on 10 May 1922[368], represent a unique body of work somewhat ahead of its time. According to Lockwood and MacMillan[369], although much work on space filling had been done by several people, from Kepler's work on sphere packing in 1611[370] to Fedorov and Schoenflies' enumeration of the space groups in 1890[371], little had been done on two-dimensional repeating patterns. Indeed, the concentration on three-dimensional crystallography meant that the seven 'frieze groups' and seventeen 'wallpaper groups' were not fully described until 1924 (by Pólya[372]). So MacMahon's interest in the subject at this time was unique.

---

361 [Haag, 1923].
362 [Schattschneider, 1999].
363 [MacMahon, 1922, [98;15]].
364 The two MacMahons presented their work to the Royal Society on 17 November 1921, as reported in *Nature* on 24 November 1921 (Vol CVIII p. 421).
365 W. P. D. MacMahon, 1886 - 1949, served as a Chaplain in the First World War from 1916 to 1917 with the rank of 2nd lieutenant.
366 [MacMahon, 1925].
367 [MacMahon, 1922, [98;15]].
368 Published in the *Journal of the RSA* on 30 June 1922, and reviewed in *Nature* on 29 July 1922 (Vol CX, p162).
369 [Lockwood et al., 1978].
370 J. Kepler, *The six cornered snowflake*, 1611 (translated by L. L. Whyte in 1966, OUP).

The purpose of the paper was "to establish a simple method for the design of repeats" and to introduce a calculus of symmetry. In the paper read for the RSA[373], this purpose is further expanded as "necessary to explore this field of thought as a preliminary to the study of the corresponding division of space of three dimensions which is required in the treatment of crystallography, crystal structure and the structure of the atom." Indeed, Andrews has observed that MacMahon's goal may been related to Hilbert's 18th problem, which concerned the existence of *fundamental regions* of *n*-dimensional Euclidean space, congruent copies of which might fill the space.

There is evidence that other mathematicians also saw MacMahon's work in the context of Hilbert's problem. At a meeting of the Philadelphia Section of the Mathematical Association of America in February 1928, Professor A. H. Wilson of Haverford College presented a paper entitled *Space Filling Polyhedra*, which began with an illustrated account of MacMahon's work on repeating polygons.

MacMahon's work on the subject stopped after the publication of the two papers described below, so that the 'Part II' implied by the title *The design of repeating patterns - Part I* was never written, for reasons mentioned below[374]. On 5 February 1922 MacMahon wrote to H. F. Baker describing a space-filling tetrahedron he had discovered. On the letter, held in St. John's College, Cambridge, there

371 [Schoenflies, 1891] and [Fedorov, 1971]. Schoenflies, 1853 -1928,and Fedorov, 1853 - 1919, were both 19th century crystallographers. Fedorov worked in St Petersburg, publishing his work the year before Schoenflies, who was working in Alsace and Göttingen.

372 [Pólya, 1924]. In his article Pólya adopts the now traditional approach of describing the possible symmetries that a tiling of the plane might possess (translation, the cyclic groups $C^2$, $C^3$, $C^4$ and $C^6$, and the dihedral groups $D^1$, $D^2$, $D^3$, $D^4$ and $D^6$) and then trying to create patterns of a suitable nature. On page 281 there are illustrations of patterns, many of which are the same as those exhibited by MacMahon.

373 MacMahon had been a member of the council of the RSA since July 1919.

374 [Andrews, 1986, Chapter 15].

is a pencilled note, presumably by Baker, "?Is this one of Schoenflies Crystallosysteme p299." Andrews felt that this had dampened MacMahon's enthusiasm for the subject. Baker was correct and on page 299 of Schoenflies' *Krystallsysteme und Krystallstructur*[375] there is a description of a construction[376] that produces the same tetrahedron as MacMahon's. It is pity that MacMahon was apparently dissuaded from further investigation of this object, since it has been 'rediscovered' by puzzle enthusiasts and used to create, for example, rotating rings of joined tetrahedra, which would have delighted MacMahon. In his letter to Baker, MacMahon wrote: "It is most fruitful in further results. I have made a rough model and am having some made by a pattern maker." This suggests that MacMahon was intending to pursue the possibilities afforded by this polyhedron, and may have done so without getting around to publishing his results.

In *The design of repeating patterns - Part I*, MacMahon presented a series of definitions and results which provided a language for talking about repeating patterns, since, as was made clear in the RSA paper, extensive investigations into both artistic and mathematical literature had revealed a lack of development in the field. MacMahon thus clearly felt, quite justifiably, that he was pioneering a new area of work[377].

Modern treatments of the subject of tiling are essentially descriptive, whereas MacMahon's approach was constructive. By using the subdivision of a polygon and the rules of contact (edge matching) as developed from the original patent of 1892 in *NMP*, MacMahon was able to show how tilings of the plane with specific symmetries could be created from triangles, squares and hexagons. The pattern from page 101 of *NMP* mentioned above is derived from a square with contact system 1 - 4, 2 - 3.

---

375 [Schoenflies, 1891].

376 Schoenflies started with a square based prism with a length twice as long as the side of the square base. He then joined the midpoints of three adjacent faces to one another and to the vertex of the prism at the intersection of the faces. The result was a tetrahedron with four identical isosceles faces, where the ratio of the long side to the short side was $1:\frac{1}{2}\sqrt{3}$. MacMahon did not indicate how he had (re)discovered it.

377 Letter 5 to D'Arcy Thompson in appendix 6 shows that he still held this view in 1923.

MacMahon was also prepared to consider designs in which the repeating unit contained holes ('stencil repeats') or comprised disconnected pieces ('archipelago repeats'). The nomenclature is ascribed to his Cambridge friend, the geometer G. T. Bennett[378], so it seems clear that MacMahon felt that the topic was important enough to discuss with eminent colleagues.

The calculus of symmetry that MacMahon developed was based on the contact systems between the polygons, and the rules of symmetry which must apply to the edges in contact in order to force the desired contact system. MacMahon identified two types of symmetry for the edge contacts: mirror symmetry in a line at right angles to the untransformed edge, and point symmetry about the centre point of the untransformed edge, illustrated below (the diagrams are adapted from MacMahon's own).

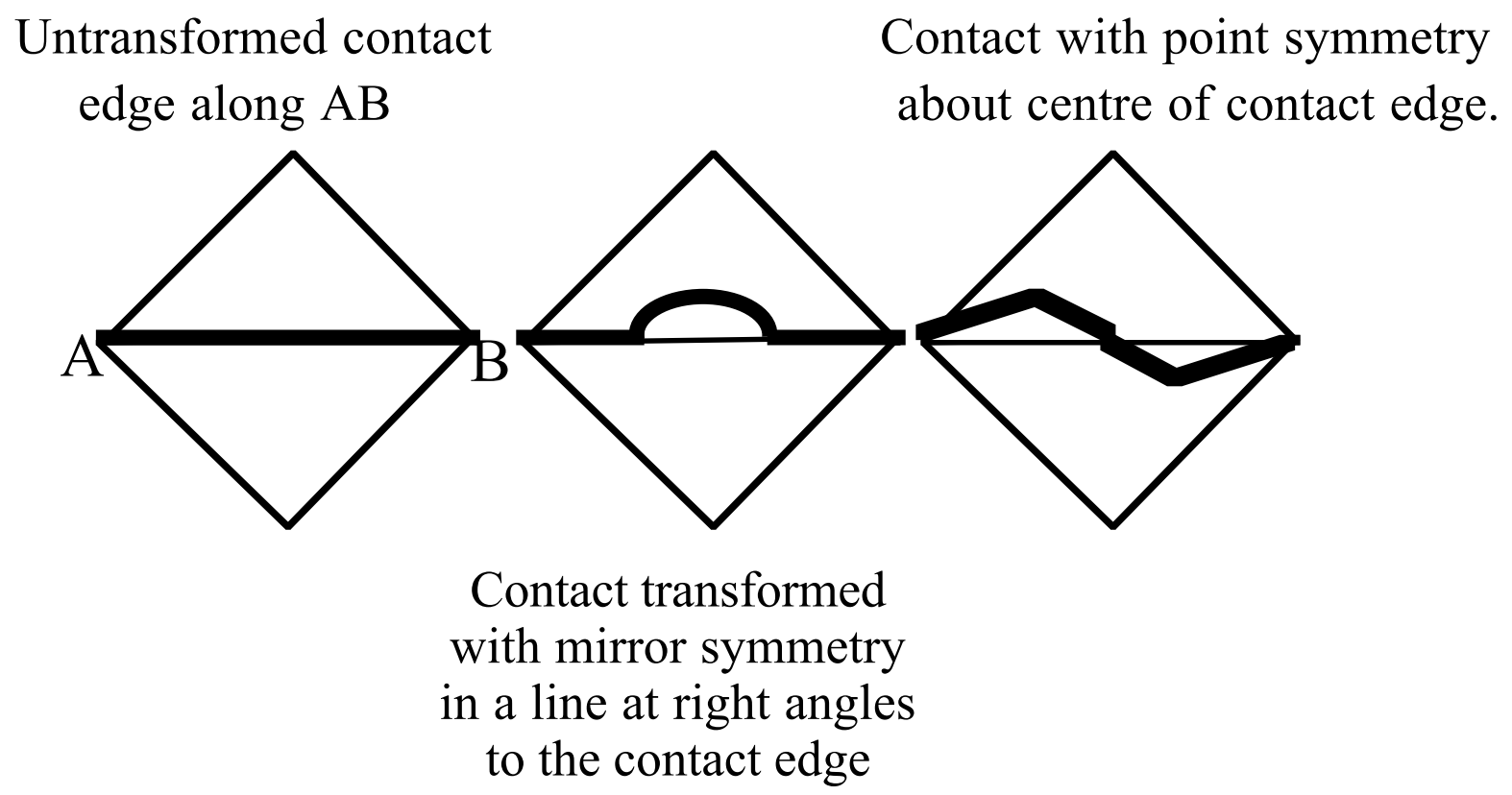

Fig 15: Contact symmetries 1

The symbol **I** was used to indicate mirror symmetry, **P** to indicate point symmetry. The suffices 1 and 2 were used to indicate the nature of the contact: 1 for the same colours touching and 2 otherwise. The diagram below illustrates the five classes of edge transformations identified in the paper. The notations S, U and V denote, respectively, transformations that are symmetrical about the line at right angles to the original contact edge, those with point symmetry about the midpoint of original contact edge, and those which have neither.

[378] See Appendix 6 containing the correspondence with D'Arcy Thompson for more about Bennett.

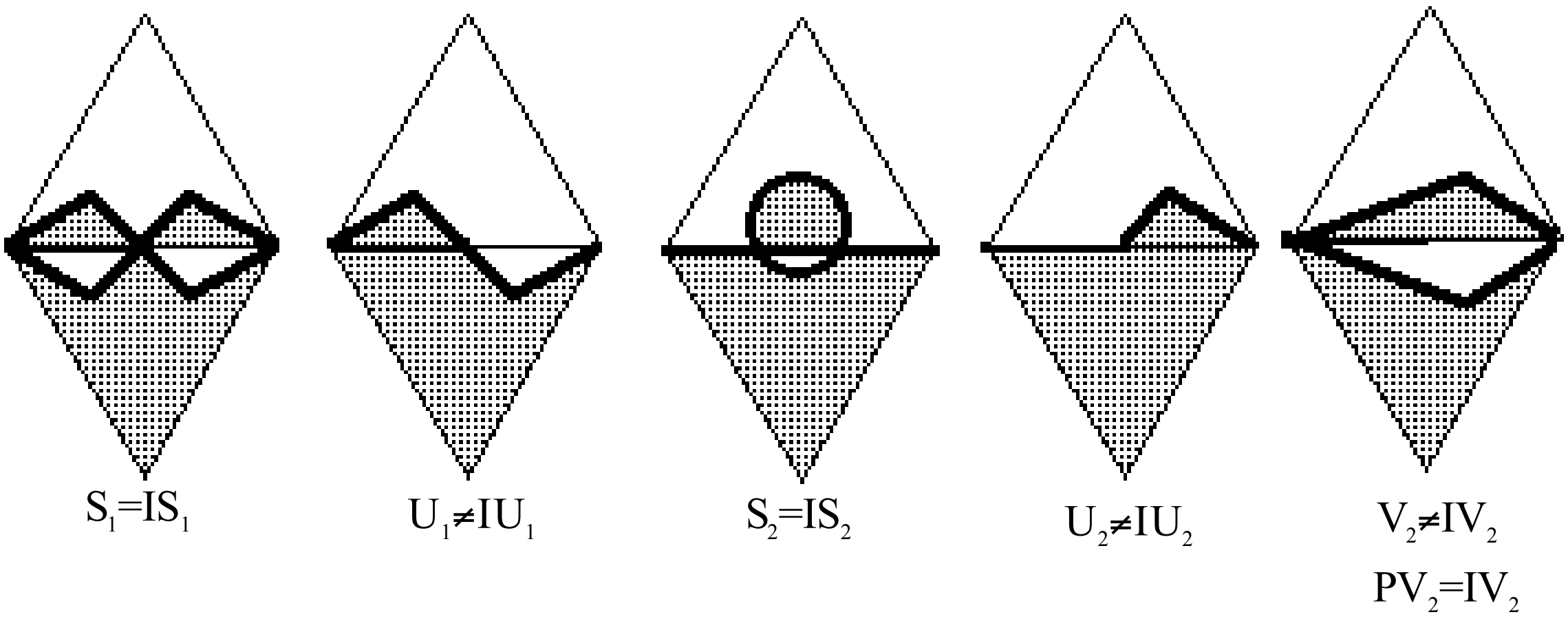

Fig 16: Contact symmetries 2

The equations under the diagrams show that, for example, an S transformation is the same as its mirror image: S = IS; or for a V transformation, the point symmetry image is the same as the mirror symmetry image (PV = IV). Assemblages of the shapes according to the contact systems creates repeat units which are translated to fill the plane. MacMahon identified the possible symmetries of repeat units, but did not seek to enumerate all the possible symmetries derivable from those symmetries; it may be that this was an idea he intended to pursue in the unwritten part II.

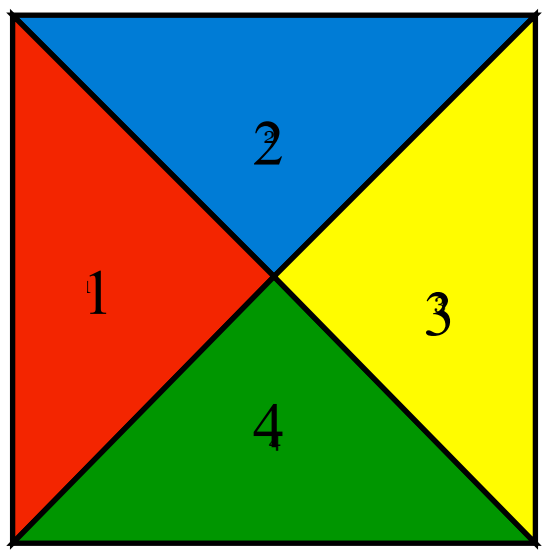

Fig 17: Four colour square

MacMahon does not seem to have considered in detail the enumeration of contact systems. For example, in the case of the square above, divided into four compartments numbered 1, 2, 3, 4 in clockwise order, there are ten different contact systems available; the ones MacMahon had illustrated in *NMP* are indicated with the page reference and figure number:

{11,24,33}, (page 96 fig 116) {11,22,34},

{12,34},
{13,24}, (page 97 fig 120a)
{11,22,33,44} (page 92 fig 110b),
{12,33,44},
{13, 22,44},
{14,22,33},
{11, 23,44},
{14, 23}.

Of these, only those printed in black are achievable without changing the order of the numbers within the square. MacMahon did not provide examples of all the possible contact systems.

With five numbered compartments in a pentagon, there are 26 potential contact systems, and six compartments in a hexagon provide 76 possibilities. The sequence 1, 2, 4, 10, 26, 76, ... ,[379] generated by the enumeration of possible contact systems is the same as the number of self-inverse permutations on *n* letters.

The direction in which MacMahon was going is clear: from the consideration of edge-matching the coloured tiles, through the replacement of the colours by transformed contact edges to the creation of repeat units, he created a fundamental region in the Hilbert sense, which generates the tiling from its symmetries.

The paper stopped rather abruptly, and the baton was passed to W. P. D. MacMahon[380], who explained how he had been involved with his uncle in "an investigation of the design of repeating patterns in two and three dimensions, a subject the idea of which was first broached by Major MacMahon in a little book entitled *New Mathematical Pastimes*, recently published." He went on to say, "the author believes that many of the results in the present paper are new", evidence that both MacMahons believed they were breaking new ground.

The paper started with a definition of a polygon or 'repeat', set out the labelling conventions for the diagrams and defined the idea of a 'system of contact' that formed the backbone of *New Mathematical Pastimes*. These polygons could be either convex or concave (that is, one or more of the angles could be re-entrant) - for the quadrilateral this was a new result. From these fundamental ideas, he derived the 'Law of Angle Distribution': *a closed polygon may be a repeat if, and only if, its angles are such that they may be associated together in sets of two, three, or four differently lettered angles, the sum of the angles in each set being equal to $\pi$ or $2\pi$.*

[379] This sequence was well known to MacMahon. He had written a paper on the subject of self-conjugate permutations [MacMahon, 1895, [48;4]].
[380] [MacMahon, 1925].

From this, W. P. D. MacMahon produced a table showing how the associations referred to in the theorem must be distributed. All triangles and quadrilaterals are repeats, but for a pentagon to be a repeat, it must have a set of three angles that sum to $\pi$ or $2\pi$. The possible kinds of contact system possible for repeat pentagons and hexagons were considered in some detail.

The final paper in this series was read to the Royal Society of Arts on 10 May 1922, and published in the *Journal of the Royal Society of Arts* on 30 June 1922. This must have been a very entertaining talk, illustrated with many slides and pictures 'thrown upon' a screen, probably using an epidiascope. The main body of the talk outlined the content of *New Mathematical Pastimes* and the papers referred to above, but what was remarkable about the talk was its audience, which included many artists. MacMahon's opening remarks provided a justification for the study of two-dimensional repeating patterns as a precursor to the study of the three dimensional case, which was of great importance to crystallography and the study of the structure of the atom. He described how he had conducted a search of geometrical and art literature before commencing his investigation, and then stated his purpose in making a presentation to the Royal Society of Arts "to bring to the notice of those who are concerned with the practical use of patterns the results that have been arrived at for scientific purposes."

His main claim was that the method of edge transformation and the calculus of symmetry he had developed were sufficient to create and classify all possible repeating patterns. In the report of the discussion afterward, it was clear that some members of the audience were not entirely convinced of this, and also cast some doubt on the aesthetic qualities of some of the patterns shown by MacMahon.

Apart from a letter in *Nature* in 1922[381] and a note in the *Proceedings of the London Mathematical Society* in 1924[382], mentioned but not reproduced in Andrews, there is no evidence that MacMahon continued to pursue this line of work, which had occupied him on and off for 30 years, possibly for the reason already mentioned.

***Later developments.***

Many people have developed the puzzles invented by MacMahon. A paper by Farrell on the *Mayblox* puzzle based on the 30 coloured cubes that can be made with six differently coloured squares has already been described in Chapter 4. The mathematical puzzle author Martin Gardner wrote three articles on them in his *Scientific American* column *Mathematical Pastimes* which were subsequently included in his books. In the discussion which follows, references are to these books rather than to the original magazine articles[383].

The first article[384] described the 24 three-colour squares, and the cubes. In it, Gardner told how the squares problem led to an article in *New Scientist* in 1961[385], an enumeration by Federico Fink[386] in Buenos Aires of 12 224 different possible patterns with the squares, and an exhaustive computer search in 1964 which found 12 261 patterns.

[381] Andrews reference [100;19]. In this letter, MacMahon described how a reinterpretation of the diagram of a well-known proof of Pythagoras's Theorem in terms of edge transformations could yield a tile that filled the plane, in three different ways. J. R. Cotter of Trinity College replied to this [*Nature*, 1922, Vol 109, p. 579] with a description of a proof of Pythagoras's Theorem based upon a chessboard. It is not the case, as claimed by G. N. Frederickson [Frederickson, 2002], that MacMahon was 'superposing tessellations to derive an unhingeable arrangement of the dissection.' There is no evidence that MacMahon had anything in mind other than repeating patterns.
[382] In *Note on Euler's Theorem in regard to certain polyhedra S-E+F=2*, MacMahon claimed that it would be more sensible to consider the dihedral angle between the faces of a polyhedron as a solid angle. Expressed in spherical measure, the "hemispherical deficiency" of a solid angle $\Omega$ is $2\pi - \Omega$, and the sum of the hemispherical deficiencies at the vertices is equal to the sum of the deficiencies at the edges in a polyhedron which satisfies Euler's theorem. [MacMahon, 1924]
[383] Descriptions of 'Instant Insanity' type puzzles, in which a column or tower of cubes must be constructed according to certain rules, have been ignored here since there is no evidence that MacMahon was ever concerned with such puzzles.
[384] [Gardner, 1966].

The second article[387] dealt with the 24 four colour triangles, and mentioned the work of Wade Philpott (1918 - 1985), an engineer from Ohio. Philpott became interested in MacMahon's puzzles during a period of hospitalisation following a shooting accident in 1947; his archive is now stored in the University of Calgary (as part of the Eugene Strens Recreational Mathematics Collection). Philpott created versions of the square and triangle puzzles that he called 'Multimatch[388]'.

The final article[389] dealt exclusively with the 30 coloured cubes. This is the only one of the three articles to include a bibliography.

In 1956, Paul B. Johnson[390] wrote an article about the cube puzzle in which he discussed the number of ways of selecting eight cubes to form a duplicate of a randomly selected 'key cube' from the full set of thirty possible cubes, and the number of ways in which possible sets of cubes 'stack' to form the required duplicate. He ignored the requirement that internal faces must also match. This article also has a bibliography.

In 1971, Norman T. Gridgeman[391], in *The 23 colored cubes*, took a wider view of the question of colouring cubes by allowing the cubes to be coloured in a single colour, two colours, three colours, etc., up to six different colours. These he grouped into 23 species (so in that sense the title may be misleading - it would more accurately be called the 23 *species* of coloured cubes - but only to avoid confusion with MacMahon's *30 coloured cubes*). This article also includes a bibliography.

---

385 [O'Beirne, 1961].
386 In 1972, Fink wrote to George Andrews describing his efforts to locate 'MacMahon's literary executrix' in Los Angeles in 1968. This would have been MacMahon's daughter Florence, who had died seven years previously. Fink's emissaries made contact with the sons of the woman who had run the nursing home in which Florence spent her last years, but found no *nachlass* from Florence.
387 [Gardner, 1985].
388 'Multimatch' puzzles are made and sold today by Kathy Jones of Kadon Enterprises; there is correspondence between Philpott and Kathy Jones from the period 1981 to 1984 in the University of Calgary archive.
389 [Gardner, 1992].
390 [Johnson, 1956].
391 [Gridgeman, 1971].

Two articles in the early 1970s[392] described a variant of the *Mayblox* puzzle, invented by Eric Cross, which had recently been put on the market. The eight blocks had to be assembled to form a larger cube after the manner of the *Mayblox*, but without the requirement that internal faces match.

**Some famous colleagues**

This chapter concludes with a discussion of some of the interactions between MacMahon and other famous mathematicians and scientists of the time. MacMahon's Fellowship of the Royal Society and other learned societies brought him into contact with most of the leading scientific figures of the period, but few records of any meetings survive. The choice of material below has been dictated by what documentation has survived, but it gives a good idea of the exalted circles into which the son of a soldier had moved from his scientifically humble beginnings.

***Ramanujan (1887 - 1920)***

In *The Man who Knew Infinity*, Robert Kanigel says[393]: "MacMahon was a whirlwind of a calculator. Sometimes, in fact, he would take on Ramanujan in friendly bouts of mental calculation - and regularly thrash him." A similar story is reported by G.-C. Rota in the preface he wrote for Andrews' *Percy Alexander MacMahon: Collected Papers*: "It would have been fascinating to be present at one of the battles of arithmetical wits at Trinity College, when MacMahon would regularly trounce Ramanujan by the display of superior ability for fast mental computation." This was heard from G. H. Hardy by D. C. Spencer, accorong to Rota. Hardy's own comment, at the end of his obituary for Ramanujan, is less dramatic. Referring to the calculation by both MacMahon and Ramanujan of a table of the number

[392] Steven J. Kahan, *Eight blocks to madness- a logical solution*, *Mathematics Magazine*, Vol. 45, 1972, pp. 57 - 65 and Andrew Sobczyk , *More progress to madness via eight blocks, Mathematics Magazine*, Vol. 47, 1974, pp. 115 - 124

[393] p. 250.

of free partitions $p(n)$ for all integers $n$ up to 200, Hardy says, "Major MacMahon was, in general, slightly the quicker and more accurate of the two."[394] These skills as a calculator were used by Hardy and Ramanujan to check their asymptotic formula for the value of $p(200)$ in 1916[395]. A later study of MacMahon's table of $p(n)$ led Ramanujan to notice the 'Ramanujan Congruences', $p(5n+4) \equiv 0 \pmod{5}$ and $p(7n+5) \equiv 0 \pmod{7}$.

It was not just as a computer that MacMahon was involved with Ramanujan. After Ramanujan's election to a Fellowship of Trinity College, Cambridge, in October 1918, Ramanujan wrote to Hardy[396], asking him to thank both "Mr. Littlewood and Major MacMahon" for "their encouragement." Earlier, MacMahon supported Ramanujan's election to the Cambridge Philosophical Society in February 1918, and in 1917 was the only person who had not been a Cambridge Wrangler to recommend Ramanujan to the Royal Society. After Ramanujan's death, MacMahon spoke about his work to a meeting of the London Mathematical Society on 10 June 1920.

***D'Arcy Wentworth Thompson (1860 - 1948)***

D'Arcy Thompson is best known as the author of *On Growth and Form*[397], first published in 1917, republished in 1942, and then in abridged form regularly since 1959. In this book, Thompson explored the connection between the growth patterns of animals and plants and their final adult forms.

The first firm evidence that MacMahon and Thompson were acquainted is in the minutes of the BAAS meeting held in Bradford on 7 September 1900 which record that D'Arcy Thompson was elected to be one of the Vice-Presidents for the Glasgow meeting held in 1901.

There are nine letters between MacMahon and Thompson spanning the years 1917 to 1926. They are

394 [Hardy, 1921].
395 [Hardy, 1918].
396 [Berndt & Rankin, 1995, p. 192].
397 [Thompson, 1997].

transcribed in full in Appendix 6.

The letters show that Thompson held MacMahon in high regard, both as a friend and a mathematician. They are about combinatorial and geometric matters, reflecting Thompson's interest in morphology in nature. In the early letters, Thompson was trying to understand how to enumerate combinations of three or more things taken 1, 2, 3, etc., at a time, using visual methods. MacMahon referred him to the liberal use of diagrammatic methods in *Combinatory Analysis*. MacMahon was particularly proud of having introduced Thompson to the geometer G. T. Bennett. The most interesting letter is number 8 in the list, in which Thompson poses the problem of the number of ways *n* hexagons might be arranged; his interest in this problem arose from studies of cells in developing embryos. There is no record of MacMahon's reply, nor that he worked on the problem at all.

***Sir Arthur Eddington[398] (1882 - 1944)***

In 1909 MacMahon wrote a paper concerning the possible use of lunar occultations to determine the diameter of distant stars[399]. The technique he proposed was to use photographic methods - either a moving film to give a trail or the "separate exposures of the kinematograph" - to time the occultation of a star by the moon. He had consulted Professor Dyson[400] on the matter, and had been assured that for a bright enough star it would be technically feasible.

The proposal was heavily criticised by Eddington. His objection was that diffraction of light by the limb of the moon would create a series of light and dark fringes during occultation whose separation would be beyond the resolving power of the telescope and thus overwhelm the tiny effect being measured.

[398] Arthur Eddington, 1882 - 1944, was senior Wrangler in 1904, chief assistant to the Astronomer Royal 1906 - 1913 and Plumian professor of astronomy 1913 - 1944. He was knighted in 1930.
[399] [MacMahon, 1909, [72;17]].
[400] Frank Watson Dyson, 1868 - 1939, was second Wrangler in 1889 and Astronomer Royal from 1910 - 1933. MacMahon quoted a letter from Dyson dated 1 November 1908, at which time Dyson was Astronomer Royal of Scotland.

In a private letter to George Andrews dated 9th February 1972, Professor R. Edward Nather[401] commented: "[MacMahon] was, in the light of long hindsight, very nearly correct even in detail. Certainly his suggestion was good one and the objections that Eddington raised proved, finally, to be incomplete and less well informed than one might expect from so famous a personage. I think it's fair to say that Eddington's comments held back the application of MacMahon's proposed technique for quite a long time."

**Personal life**

At the age of 52, MacMahon married for the second time. His wedding to Grace Elizabeth Howard took place on 9 February 1907 at the Church of St. John the Evangelist in Kilburn. MacMahon described himself as 'single' on the marriage certificate for this union, making no reference to his previous marriage. It is possible that after the break-up of the first marriage, MacMahon had not mentioned the episode to anyone. Certainly none of his three obituarists mentioned the first marriage; whether this was through ignorance or out of respect we cannot know. Neither can we know whether Grace was aware of the previous union.

However, divorced persons are not allowed to marry in church under the rules of the Church of England during the lifetime of the divorced partner. So, either MacMahon did not tell the truth to the Church (which seems unlikely) or he must have known - or believed - that Aimee Rose Leese was dead[402]. Supporting this latter position is the family history of the Leese family, which claims that Aimee Rose met her death in the San Francisco earthquake of 1906.[403] MacMahon's daughter Florence is also reputed to have died in this catastrophe, but this is not true. Florence was also MacMahon's literary executrix after his death, and continued to receive royalty cheques from the Cambridge University Press until her death in 1961.

[401] R.E. Nather is Professor Emeritus of Astronomy at the University of Texas at Austin and Director Emeritus of the Whole Earth Telescope project.
[402] MacMahon described himself as a widower in the census of 1901, which may mean that Aimee had actually died, or that MacMahon had claimed that this was the case.
[403] Private correspondence with Linda Cates, great-granddaughter of Percy Henry Leese.

At the time of his marriage to Grace Elizabeth, MacMahon was living at 13 Cambridge Gardens in West London. In 1908, the couple moved to 27 Evelyn Mansions in Carlisle Place, close to Victoria Station, where they lived until they moved to Cambridge in 1922.

Some time between 1906 and 1911, the impresario Luther Munday[404] made a bronze bust of MacMahon. This may have been commissioned by one of the newly-weds as a gift for the other. Luther Munday's book *A Chronicle of Friendships*[405] does not mention why the bust was made.

[404] Luther Munday, 1857 - 1922, was well known in London as a *bon viveur*, theatre agent and manager of the Lyric Club.
[405] [Munday, 1912].

## Chapter 7 Cambridge 1922 - 1929

In 1922 MacMahon, already retired from the Board of Trade, decided to leave London. He and his wife moved to 31 Hertford Street, Cambridge, from where MacMahon continued to work and publish mathematical papers. The year after his move to Cambridge was MacMahon's most prolific year in terms of the number of papers published, with 10 printed, against an annual average of 3.

### Cambridge

MacMahon continued his association with St John's College, Cambridge. According to the *Cambridge University Reporter*, MacMahon offered a short lecture course at Cambridge University under the auspices of the Special Board of Mathematics. It was entitled *Some processes in combinatory analysis*, was first advertised in October 1923, and continued to be offered until January 1925. From that date until January 1928, MacMahon's name appears in the list but without the title of the lecture series. No records of the content of, or attendance at, these lectures exists, although Philip Hall[406] is reported to have attended in 1925[407].

The statistician Maurice Kendall[408] also remembered seeing MacMahon at the College in the late 1920s, and described him as a "old man with a loud booming voice."[409] Kendall later went on to use MacMahon's work on symmetric functions in his work on Fisher's $k$-statistics; a $k$-statistic is an unbiased estimator of a distribution expressed in terms of the sums of $r$th powers of the data points.

[406] Philip Hall, 1904 - 1982, was Sadleirian Professor of Mathematics from 1953 - 1967, elected FRS in 1942, President of the LMS from 1955 - 1957, and was awarded the Royal Society's Sylvester Medal in 1961.
[407] [Green, et al., 1984].
[408] The statistician Maurice George Kendall, 1907 - 1983, was knighted in 1974, and was a Guy medallist. The Guy medal is awarded by the Royal Statistical Society every three years for innovative contributions to statistical theory.
[409] Private letter from Kendall to George Andrews, 3 January 1974.

**R. A. Fisher (1890 - 1962)**

In July 1924, the statistician R. A. Fisher[410] consulted MacMahon about Latin squares[411] on the strength of MacMahon's writing on the subject in *Combinatory Analysis*[412]. However, Fisher used his own method of direct enumeration rather than MacMahon's algebraic method and discovered four pairs of *reduced* 5 × 5 Latin squares that MacMahon had missed. A reduced Latin square is one in which the first row and column are in alphabetical or numerical order; all other Latin squares can be obtained from the reduced squares by interchanging rows and columns. In *Combinatory Analysis*, MacMahon had given the number of reduced 5 × 5 squares as 52, calculated using the Hammond operator *D* on a suitably chosen symmetric function. Fisher wrote in 1934[413] ,

> The corrected number was communicated by one of the authors in 1924 to Professor MacMahon in time to be incorporated in the copies of *Combinatory Analysis* then unsold; commenting to his co-author that the reason for this self-evident statement was that MacMahon had insisted that the Cambridge University Press should reprint and replace the offending pages! The correct number, 56, had actually been given 'd'après un *dénombrement exact*' by Euler in 1782, but MacMahon, in spite of extensive quotations from this paper (still in 1934 marked in pencil in the Royal Society copy), had overlooked this.

Both Fisher and MacMahon attended the International Congress of Mathematics in Toronto in August 1924, and in a letter to Fisher dated 19 September 1924 MacMahon wrote, "I have lost my Toronto note. Is the correct number of Latin squares (reduced) now 56 instead of my number 52?" so it is likely that the communication referred to by Fisher occurred at that meeting. This episode shows that even meticulous mathematicians are prone to the occasional oversight or miscalculation. MacMahon had noted that for Latin squares of order higher than 4 the calculations were impracticable[414], so it is not surprising that he may have made an error.

[410] Ronald Aylmer Fisher, 1890 - 1962, worked at the Rothamsted Agricultural Experimental Station and introduced the concept of *likelihood*. He was elected FRS in 1929, and received the Royal Society's Royal Medal in 1938, the Darwin Medal in 1948 and the Copley Medal in 1955. The correspondence between MacMahon and Fisher is in the archives of the University of Adelaide, Australia.
[411] [Box, 1991, p. 155] and [Bennett, 1990, p. 273]. A Latin square is an *n* by *n* array of *n* distinct symbols arranged so that each symbol occurs exactly once in each row and column.
[412] [MacMahon, 1915].
[413] [Fisher & Yates, 1934].
[414] [MacMahon, 1984, p. 251].

### The International Mathematical Union (IMU)[415] / International Congress of Mathematics (ICM).

MacMahon attended only two International Congresses. The first was the Fifth International Congress of Mathematics held in Cambridge, in August 1912, where he was chosen as a Vice President.[416] No papers or lectures were presented by MacMahon at this meeting, nor is he recorded as having taken a prominent part in any of the discussions.

MacMahon and his wife attended the 7th ICM in Toronto in 1924; this was also the second General Assembly of the IMU. He is shown as representing the BAAS, the Royal Society of Arts, the University of Cambridge (along with Sir Arthur Eddington), and as the Royal Society's representative on the British National Committee of the IMU. In addition to presenting a paper, *The expansion of determinants and permanents in terms of symmetric functions*, MacMahon also spoke at the closing session on Saturday 16 August. The entry in the Proceedings of the Congress[417] reads:

> Professor S Pincherle, President of the International Mathematical Union, and Major P. A. MacMahon, then spoke.

One interpretation of this is that MacMahon was speaking on behalf of the British National Committee of the IMU, but Dr Lehto, the present archivist and historian of the IMU, disagrees[418].

### De Morgan medal

The De Morgan medal is the highest accolade of the London Mathematical Society and is awarded every three years. It was first awarded to Cayley in 1884. MacMahon was proposed for the award on five occasions: in 1884, by Cayley; in 1890, and again in 1893[419] (on this occasion MacMahon requested that his name be withdrawn, and proposed Klein instead); and in 1920 Hardy put

[415] The IMU was founded in 1920, an occurrence supported by MacMahon but not by Hardy, as described earlier in the section on the LMS. A history of the IMU may be found in [Lehto, 1998].
[416] The Honorary President was Lord Rayleigh, the President was Sir G. H. Darwin, and the other Vice-Presidents were W. von Dyck, L. Fejér, R. Fujisawa, J. Hadamard, J. L. W. V. Jensen, G. Mittag-Leffler, E. H. Moore, F. Rudio, P. H. Schoute, M. S. Smoluchowski, V. A. Steklov and V. Volterra. The Secretaries were E. W. Hobson and A. E. H. Love.
[417] [Fields, 1928].
[418] Private e-mail; see also [Lehto, 1998].
[419] The dates refer to the LMS Council minutes when the proposals were recorded.

MacMahon's name forward. This may have been some sort of recompense for the dispute described earlier over the formation of the IMU. MacMahon himself proposed Forsyth and Hilbert for the De Morgan Medal in 1896, Burnside in 1899, Baker in 1905, and W. H. Young in 1911.

MacMahon was finally awarded the De Morgan Medal in 1923. The award was announced by the President, W. H. Young, at the Council meeting on 14 June 1923, and the medal was presented to him at the Annual General Meeting on 8 November 1923. There is no record of any citation to accompany the award. After MacMahon's death, his medal was returned to the LMS in 1956 by Colonel McLaverty, the third husband of MacMahon's niece Louise, the daughter of MacMahon's elder brother George.

## Mathematics

### *Partition theory*

MacMahon's 1921 paper on partition theory, *Note on the parity of the number which enumerates the partition of a number*[420], is a short description of the solution to a problem posed by Ramanujan to MacMahon in 1919 in a private letter. Using Euler's pentagonal number formula, MacMahon had already provided Hardy and Ramanujan with values for $p(n)$ up to $p(200)$, which they used to test their asymptotic formula for $p(n)$. The following picture shows part of the table of the enumeration of unrestricted partitions of all numbers from 1 to 200 calculated by MacMahon, and published in the Proceedings of the London Mathematical Society in 1918.

[420] [MacMahon, 1921, [97;9]].

| | | | |
|---|---|---|---|
| 27 ... 3010 | 65 ... 2012558 | 103 ... 271248950 | 141 ... 16670689208 |
| 28 ... 3718 | 66 ... 2323520 | 104 ... 304801365 | 142 ... 18440293320 |
| 29 ... 4565 | 67 ... 2679689 | 105 ... 342325709 | 143 ... 20390982757 |
| 30 ... 5604 | 68 ... 3087735 | 106 ... 384276336 | 144 ... 22540654445 |
| 31 ... 6842 | 69 ... 3554345 | 107 ... 431149389 | 145 ... 24908858009 |
| 32 ... 8349 | 70 ... 4087968 | 108 ... 483502844 | 146 ... 27517052599 |
| 33 ... 10143 | 71 ... 4697205 | 109 ... 541946240 | 147 ... 30388671978 |
| 34 ... 12310 | 72 ... 5392783 | 110 ... 607163746 | 148 ... 33549419497 |
| 35 ... 14883 | 73 ... 6185689 | 111 ... 679903203 | 149 ... 37027355200 |
| 36 ... 17977 | 74 ... 7089500 | 112 ... 761002156 | 150 ... 40853235313 |
| 37 ... 21637 | 75 ... 8118264 | 113 ... 851376628 | 151 ... 45060624582 |
| 38 ... 26015 | 76 ... 9289091 | 114 ... 952050665 | 152 ... 49686288421 |

* The numbers in this table were calculated by Major MacMahon, by means of the recurrence formulae obtained by equating coefficients in the identity

$$(1-x-x^2+x^5+x^7-x^{12}-x^{15}+\dots)\sum_0^\infty p(n)\,x^n = 1.$$

We have verified the table by direct calculation up to $n = 156$. Our calculation of $p(200)$ from the asymptotic formula then seemed to render further verification unnecessary.

Fig 18. Part of the table of the enumeration of unrestricted partitions

Ramanujan subsequently wanted to calculate $p(1000)$ and asked MacMahon whether he could determine whether the result would be even or odd. In the paper MacMahon describes the congruence relations with which he was able to establish "in about five minutes work, that $p(1000)$ is an uneven number."

The next short paper followed two years later, in 1923. In *The connexion between the sums of the squares of the divisors and the number of partitions of a given number*[421] MacMahon showed that by the application of a differential operator to the logarithm of the enumerating function for plane partitions of $n$, the coefficients of the result are the sums of the squares of the divisors of $n$.

Four further papers on aspects of partition theory followed in quick succession in the same year. *The partitions of infinity with some arithmetic and algebraic consequences*[422] developed some details from two earlier papers[423]. The idea of a perfect partition was extended in that paper to partitions of an infinite integer, and was developed further.

[421] [Macmahon, 1923, [105;11]].
[422] [Macmahon, 1923, [107;6]].
[423] [MacMahon, 1886, [20;6]] and [MacMahon, 1891, [38;6]].

The general expression of a partition of infinity is:

$$1^{\alpha_1}\left(1+\alpha_1\right)^{\alpha_2}\left\{\left(1+\alpha_1\right)\left(1+\alpha_2\right)\right\}^{\alpha_3}\left\{\left(1+\alpha_1\right)\left(1+\alpha_2\right)\left(1+\alpha_3\right)\right\}^{\alpha_4}\ldots,$$

where the $\alpha_i$ are positive integers, and (when used as exponents) represent repetitions; there may also be an integer $\alpha_s$ such that every $\alpha_i$ is 0 for $i > s$. MacMahon drew attention to some special cases (e.g. $\alpha_1 = \infty$, so that the sequence is an infinite succession of units) and noted that every partition of infinity corresponds to a scale of numeration - for example, the binary scale corresponds to $\alpha_s = 1$ for all $s$). He also gave some applications to the functions that enumerate partitions and symmetric functions.

The idea is further developed in *The prime numbers of measurement on a scale*[424], which might also be classified as a recreational paper. Here, MacMahon discussed the problem of dividing a rod into segments so that it might be used as a measuring rod, but without the usual redundancy in the markings. By considering a scale of infinite length, MacMahon avoided the difficulties caused by the boundaries, but he added the condition that each length from 1 to infinity must be measurable by a single measurement, since otherwise the problem is trivial (because it would then be possible simply to use a unit measurement repeatedly). The result is that segments of successive lengths 1, 2, 4, 5, 8, 10, 14, ... are required, with divisions at 0, 1, 3, 7, 12, 20, 30, 44, .... MacMahon had originally intended to list the first 347 terms of the first series, but printing costs required that this was reduced to 42.

The next paper, *The theory of modular partitions*[425], dealt with a generalisation of the Ferrers graph, where the dots are replaced by units, to give a 'modulo 1 partition' of a number, and then successively collapsed using 2s, 3s, and so on. An example will illustrate the concept; using the partition (8521) of

[424] [MacMahon, 1923, [109;6]].
[425] [MacMahon, 1923, [110;9]].

16, we have:

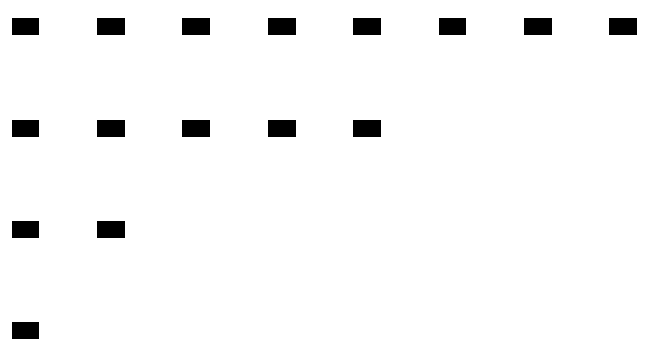

Fig 19: Ferrers' graph of (8521)

The modular partitions are then:

| | | | | | | | |
|---|---|---|---|---|---|---|---|
| Mod 1: | 11111111 | Mod 2: | 2222 | Mod 3: | 332 | Mod 4: | 44 |
| | 11111 | | 221 | | 32 | | 41 |
| | 11 | | 2 | | 2 | | 2 |
| | 1 | | 1 | | 1 | | 1 |
| Mod 5: | 53 | Mod 6: | 62 | Mod 7: | 71 | Mod 8: | 8 |
| | 5 | | 5 | | 5 | | 5 |
| | 2 | | 2 | | 2 | | 2 |
| | 1 | | 1 | | 1 | | 1 |

MacMahon used his considerable computational abilities to determine generating functions to enumerate these partitions.

The final 1923 paper was *Dirichlet series and the theory of partitions*[426] in which MacMahon defined the terms *potency* and *multiplicity* of an integer. If an integer $n$ has the prime factorisation $p_1^{\pi_1} p_2^{\pi_2} p_3^{\pi_3} \ldots p_s^{\pi_s}$ then the potency $pt(n)$ is $\sum_{i=1}^{s} \pi_i p_i$ and the multiplicity $mt(n)$ is $\sum_{i=1}^{s} \pi_i$. In the theory of partitions, $p_1^{\pi_1} p_2^{\pi_2} p_3^{\pi_3} \ldots p_s^{\pi_s}$ denotes a partition of the number $\sum_{i=1}^{s} \pi_i p_i$, which is the potency. Associated with the potency is the multipartite number $\overline{\pi_1 \pi_2 \pi_3 \cdots \pi_s}$ formed by the powers of the prime factors. From these correspondences, MacMahon was able to recast the famous Goldbach Conjecture (every even number greater than 2 is the sum of two primes) into partition theory form[427]: every even potency greater than 4 is compatible with association with a multipartite number $\overline{11}$. For example, if $n = 33$, then the prime factorisation is $3 \times 11$, the potency is $(1 \times 3) + (1 \times 11) = 14$, the multiplicity is $1 + 1 = 2$, and the associated multipartite number is $\overline{11}$. The number of integers $N_v$ with potency $v$ is given by the infinite product $\prod \frac{1}{1-b^p} = \sum_v N_v b^v$.

[426] [MacMahon, 1924, [111;8]].
[427] Veteran number theorists are certain that this reformulation is of no practical use in any attempt to prove Goldbach's conjecture.

MacMahon then showed that the Dirichlet series, $f(s) = \sum_{n=1}^{\infty} a_n n^{-s}$ is equal to the infinite product $\prod_{n=2}^{\infty} \frac{1}{1-n^{-s}}$ which enumerates the partitions of the multipartite number $\overline{\pi_1 \pi_2 \pi_3 \cdots \pi_s}$ He drew attention to the fact that this infinite product may be factorised, and the first factor is the Riemann zeta function, and observed that "the functions herein seem to be worthy of further study". He extended this work in 1925 in the paper *The enumeration of the partitions of multipartite numbers*[428].

MacMahon was now in his seventies, but he produced three more papers on partition theory before his death, two in 1926 and one in 1927.

The first of these was *Euler's $\phi$–function and its connexion with multipartite numbers*[429] in which MacMahon generalised Euler's totient function (the number of integers less than $n$ that are relatively prime to $n$), by recasting the definition to be: $\phi(n)$ is the number of bipartite numbers $\overline{a, n-a}$ which are such that $n$ and $n - a$ are relatively prime and both greater than or equal to 1. This then lends itself to extension to other multipartite numbers in a simple way. MacMahon also mentioned alternative generalisations of the totient function investigated by Schemmel[430], C. Jordan[431] and W. E. Story[432], showing that he was still reading widely.

[428] [MacMahon, 1925, [114;8]].
[429] [MacMahon, 1926, [117;8]].
[430] V. Schemmel - a Continental mathematician about whom little is known.
[431] M. E. C. (Camille) Jordan, 1838 - 1922, was elected FRS in 1919.
[432] William Edward Story, 1850 - 1930, awarded a PhD in Leipzig in 1875, was an American mathematician at Johns Hopkins University.

His second 1926 paper[433], *The parity of p(n), the number of partitions of n, when* $n \le 1000$, is essentially a list of all the values of $n$ for which $p(n)$ is even, following directly from his 1921 *Note on the parity of the number which enumerates the partitions of a number*[434]. MacMahon draws no consequences nor makes any conjectures based on this list. A later consequence of this work is the little known *MacMahon Number*, 1.74264258..., named by Parkin and Shanks[435]. The binary expansion of this number is 1.101111100001110111 01..., where the $k$th digit to the right of the point is 0 or 1, according as $p(k)$ is even or odd.

MacMahon's final paper on partition theory was *The elliptic products of Jacobi and the theory of linear congruences*[436]. He considered some properties of elliptic products which depend upon partitioning the product according to certain congruences, and gave some extensive tables of results. The final sentence of the paper was, "I hope to resume consideration of this subject on a future occasion", but it was not to be.

MacMahon's work on partition theory is characterised by two main features: his use of visual imagery to explain and motivate algebraic developments, sometimes in a recreational format (card games, weighing problems and measuring problems), and the generalisation of the relationship between partitions and existing ideas (divisors, compositions, Euler's theorem, Ferrers graphs and the totient function). Along the way he invented several new ideas (perfect partition, potency of an integer, multiplicity of an integer), and produced tables of results that we would nowadays obtain only with the aid of a computer[437] .

---

433 [MacMahon, 1926, 118;9]].
434 [MacMahon, 1921, [97;9]].
435 [Parkin & Shanks, 1967].
436 [MacMahon, 1927, [119;9]].
437 Making the calculations with a computer is also not easy: see [Leijenhorst, 2000, p. 16].

## MacMahon's final year

In 1928 MacMahon's health had deteriorated and he moved to Bognor Regis in Sussex[438] on the advice of his doctors. He died on Christmas Day 1929, aged 75, at Springfield, Normanton Avenue.[439] A funeral service was held at St Wilfred's in Bognor on Tuesday 31 December at 2:30pm, then part of the parish of St John's and very close to Normanton Avenue. MacMahon was buried in the town cemetery in grave number 3506. There is a six-foot stone cross on a plinth to mark the grave, which bears the inscription:

In loving memory of
Major Percy Alexander
MacMahon
Royal Artillery, F.R.S., D.Sc.
died on Christmas Day 1929
aged 75 years
Also of
Grace Elizabeth
widow of the above
who died March 20th 1935
Peace Perfect Peace

There is a small stone urn just in front of the plinth.

MacMahon's death was reported to the LMS at the Council meeting on 16 January 1930. On 6 February 1930, Dr F. S. Macaulay and Mr. Russell were charged with writing an obituary, but the published obituary was actually written by H. F. Baker.

His estate, amounting to £2,753 16s 5d, passed to his wife and a nephew, Lt. Col. George Dudley Ruadh MacMahon, the son of MacMahon's elder brother George. Six years later, his wife's estate, £10,686 15s, passed to G. D. R. MacMahon and another nephew, Lt. Commander Brian Patrick MacMahon, the son of MacMahon's younger brother, Ernest.

438 The obituarist A. R. Forsyth mistakenly placed him in Bournemouth, Dorset, at his death.
439 Mrs MacMahon is not listed as resident of the house after 1930. At her death on 20 March 1935 she was living at Harcourt, Nyewood Lane, Bognor Regis. Normanton Avenue still exists, but the house was demolished between 1959 and 1964 as part of a redevelopment programme.

**MacMahon's legacy**

We have seen how MacMahon's mathematical work began. He came from a family of soldiers but, unlike his brothers, was forced to abandon hope of a glorious career filled with military honours and decorations. Instead, he created an alternative career for himself in mathematics, filled with academic honours and decorations. He used his computational skills to help George Greenhill, and his work on ballistics led to the discovery of the correspondence between seminvariants and non-unitary symmetric functions, which marked his entry into the world of British mathematics. Cayley and Sylvester, the two most prominent figures in mathematics at the time, became friends of MacMahon and introduced him into the inner circle of Victorian mathematicians. His unusual abilities and outgoing personality enabled him to befriend many of the leading mathematicians and scientists of his day. His early election to the Royal Society in 1890, just seven years after his first mathematical breakthrough, made this process much easier.

MacMahon made his name in a demanding field (mathematics) populated by the intellectual elite of Victorian England, despite coming from a non-academic background. He did this without the education and social contacts that came from attendance at a traditional university. But it appears as if, once he had realised that he could become successful in the field, he made sure he was in the right places at the right time. He joined the appropriate organisations and volunteered for service on councils and committees, ensuring that as many people as possible knew who he was. Whether this was indeed a conscious effort on MacMahon's part[440], or simply that he was caught up in the excitement of London's intellectual life cannot be known. Certainly his military training and rank would have accustomed him to expect attention when he demanded it, and his physical appearance would have helped him to be noticed wherever he was. The evidence of his letter to Larmor in 1897 (see page 48) suggests that MacMahon was nevertheless surprised by his success.

[440] See Parshall and Seneta [Parshall, et al., 1997] for an example of a lifetime of deliberate action to establish a reputation in mathematics - the life of J. J. Sylvester.

MacMahon's mathematical work was special because he combined a strong computational ability with great insight, assisted by the use of visual imagery. He published new work more or less continuously for 47 years, and his style was clear and concise. However, having made such a name for himself, why is MacMahon less well known than many of his contemporaries ?

The nature of mathematics changed in the early twentieth century. The British obsession with exhibiting "the things themselves" gave way to a more symbolic way of working. The driving force behind the changes was David Hilbert, whose work on the finite basis theorem effectively destroyed general interest in the British style of invariant theory[441]. Hilbert's 1900 speech at the Paris ICM set the scene for future developments, and it did not include combinatory analysis. The emphasis was on "number theory and higher and abstract algebra, most of real and complex variable analysis and the still emerging branch of topology."[442] In Britain, the influence of Hardy shifted the area of interest toward analysis, until then a Continental strength with few devotees in the United Kingdom[443]. The new theories of relativity and quantum mechanics in mathematical physics were also attracting interest away from less obviously useful studies, such as partition theory.

MacMahon's work in the design of repeating patterns was too far from mainstream mathematics, and partition theory and combinatory analysis had yet to achieve a sufficient degree of importance. In a sense, then, MacMahon was disadvantaged on these fronts: he was working in already neglected or still unfashionable areas of mathematics, and the focus of mathematical research had shifted, both at home and in Europe (where MacMahon's work had in any case made little impact - see his comment to Ronald Ross on page 117).

[441] See [Parshall, 1990] for a detailed analysis of the decline of invariant theory, as practised by British mathematicians towards the end of the 19th century.
[442] [Grattan-Guinness, 2000].
[443] [Rice et al., 2003].

It is also significant that despite his renown, MacMahon never held a full academic post in a university. A consequence of this is that there is no convenient archive of MacMahon's writing in a university library to provide an accessible source of research material, nor any students who might have kept his memory alive in a popular sense. The move to Bognor in the year before his death meant that all his papers were left in the care of his wife, and then lost when she died; this loss may have included MacMahon's attempts to reconstruct Sylvester's results in partition theory and the second part of the work on repeating patterns, and would certainly have included his own sets of coloured triangles and squares. All these factors combined to rob him of wider historical recognition.

Was he an amateur or a professional mathematician ? There is an entire debate, which is beyond the scope of this thesis, to be had about the nature of professionalism and whether or not mathematics was professionalised in the Victorian era. John Heard[444] has characterised professionalisation in this era as "closely allied to the rise of the middle classes and their quest for gentlemanly status, the essential elements of a profession being ethics, public service, maintenance of standards, and official recognition." Heard went on to argue that these elements were not present in the activities of Victorian mathematicians in Oxford, Cambridge or the London Mathematical Society, and as a result it is "difficult to justify the claim that pure mathematics 'professionalised'." However, since MacMahon was paid to use mathematics for military purposes, and then to teach mathematics, whilst in the Royal Artillery, subsequently worked at the Board of Trade in a post which required a high level of mathematical ability and was consulted on mathematical problems (electoral reform, for example), his mathematical practice certainly included public service and official recognition. His arguments about the status of combinatory analysis in his 1896 LMS presidential speech and his remarks about education at the 1901 BAAS meeting show that he was concerned with maintenance of standards. So MacMahon's mathematical work certainly met three of Heard's four criteria; in this sense, he was a professional. MacMahon himself said[445]: "We are all of us amateurs who do not hold mathematical

[444] John Heard is a Senior Research Fellow at Imperial College, London. This comment quoted is from the abstract of a talk he gave to the Danish Society for the History of Science in 2002.
[445] Remark made in a letter dated 24 April 1916 to Ronald Ross, transcribed as Letter 4 in Appendix 9.

chairs", reflecting the contemporary view of what it meant to be a professional mathematician.

MacMahon's legacy to mathematics is predominantly in two areas: partition theory and recreational mathematics. In the latter field, there are three important contributions: the cube puzzle *Mayblox*, the edge-matching puzzles, and the work on the design of repeating patterns, all summarised in his popular book *New Mathematical Pastimes*.

The legacy from partition theory is the work on combinatory analysis set out in the eponymous book, and the conjectures and unsolved problems therein, related to two and three dimensional partitions, as described in earlier chapters. Many of the results in *Combinatory Analysis* still reappear in modern work, as mentioned in previous chapters. The book organised and extended an existing set of disparate results, providing mathematicians with a clear reference work and pointers for further development. He left a useful and interesting body of work which is still providing a source of results and inspiration for mathematicians in the 21st century.

# Appendices

## The life and work of Major Percy Alexander MacMahon
## PhD Thesis by Dr Paul Garcia

## Appendix 1

The following table shows a typical weekly timetable of the Royal Military Academy from around 1892, to illustrate the way in which MacMahon would have spent his time as a gentleman cadet.

*Weekly Detail of Studies and Drills,* 1892.

| Day | Time | 1st CLASS. Artillery Divn. | 1st CLASS. Engineer Divn. | 2nd CLASS. Artillery Division. | 2nd CLASS. Engineer Division. |
|---|---|---|---|---|---|
| MONDAY | 8-15 to 9-45 a.m. | Artillery | Mathematics | Military Topography | Military Topography |
| | 10-0 to 11-0 a.m. | | Drawing | | |
| | 11-15 to 1 p.m. | Riding Drill | Riding Drill | | |
| | 4-0 to 5-0 p.m. | | Fortification | Drawing | |
| | 5-15 to 7-15 p.m. | Artillery | | Fortification | Drawing |
| TUESDAY | 8-15 to 9-45 a.m. | Military Topography | Military Topography | Riding Drill | Riding Drill |
| | 10-0 a.m. to 3-0 p.m. | | | 10-0 a.m. to 1-0 p.m ½ Divn. at Arsenal ½ Divn. Artillery | 10-0 a.m. to 1-0 p.m ½ Divn. at Arsenal ½ Divn. Artillery |
| | 4-0 to 5 p.m. | Drawing | | | |
| | 5-15 to 7.15 p.m. | Tactics, Military Admin. & Law. | Tactics, Military Admin. & Law. | Artillery | Artillery |
| WEDNESDAY | 8-15 to 9-45 a.m. | Riding Drill | Riding Drill | Fortification | Fortification |
| | 9-45 a.m. to 1-0 p.m. | ½ Divn. at Arsenal | ½ Divn. at Arsenal | | |
| | 10-30 to 11-45 a.m. | ½ Divn. Artillery | ½ Divn. Artillery ½ Divn. Mathematics | Artillery Exercises | Artillery Exercises |
| | 12-0 noon to 1-0 p.m. | | | | |
| | 5-15 to 7-15 p.m. | Fortification | Mathematics | Artillery | Fortification |
| THURSDAY | 8-15 to 9-45 a.m. | Fortification | Artillery | Artillery Exercises | Mathematics |
| | 10-0 to 11-0 a.m. | | | | |
| | 11-15 a.m. to 1-0 p.m. | Riding Drill (2 Sections) | Riding Drill (2 Sections) | Battalion Drill | Battalion Drill |
| | 12-0 noon to 1-0 p.m. | Battalion Drill (1 Section) | Battalion Drill (1 Section) | | |
| | 2-30 to 3-30 p.m. | | | Riding Drill | Riding Drill |
| | 4 to 5 p.m. | | | Drawing | |
| | 5-15 to 7-15 p.m. | Chemistry and Physics | Chemistry and Physics | Tactics, Military Admin. & Law | Tactics, Military Admin. & Law |
| FRIDAY | 8-15 to 9-45 a.m. | Artillery Exercises. | Fortification | Tactics, Military Admin. & Law. | Tactics, Military Admin. & Law. |
| | 10-0 to 11-45 a.m. | | | Artillery Exercises | Mathematics |
| | 12-0 noon to 1-0 p.m. | Battalion Drill | Battalion Drill | Battalion Drill | Battalion Drill |
| | 2-10 to 3-45 p.m. | Riding Drill | Riding Drill | | |
| | 4-0 to 5-0 p.m. | Drawing | | | |
| | 5-15 to 7-15 p.m. | Tactics, Military Admin. & Law. | Tactics Military Admin. & Law. | Chemistry and Physics | Chemistry and Physics |
| SAT. | 8-15 to 9-45 a.m. | Tactics, Military Admin. & Law. | Tactics, Military Admin. & Law. | Artillery | Fortification |
| | 10-0 to 11-30 a.m. | Riding Drill | Riding Drill | | |

| Day | Time | 3rd CLASS. | 4th CLASS. |
|---|---|---|---|
| MONDAY | 8-15 to 9-45 a.m. | Fortification | Mathematics |
| | 10-0 to 11-45 a.m. | | |
| | 12-0 noon to 1-0 p.m. | Gun Drill | Squad Drill or Gymnastics |
| | 2-10 to 3-45 p.m. | Riding Drill | |
| | 4-0 to 5-0 p.m. | Drawing | French or German |
| | 5-15 to 7-15 p.m. | French or German | Chemistry and [Physics |
| TUESDAY | 8-15 to 9-45 a.m. | Mathematics | Fortification |
| | 10-0 to 11-45 a.m. | | |
| | 12-0 noon to 1-0 p.m. | Gun Drill | Squad Drill or Gymnastics |
| | 2-10 to 3-45 p.m. | Riding Drill | |
| | 4-0 to 5-0 p.m. | French or German | Drawing |
| | 5-15 to 7-15 p.m. | Chemistry and Physics | French or German |
| WEDNESDAY | 8-15 to 9-45 a.m. | Military Topography | Mathematics |
| | 10-0 to 11-45 a.m. | | |
| | 12-0 noon to 1-0 p.m. | Gun Drill | Squad Drill or Gymnastics |
| | 2-10 to 3-10 p.m. | | |
| | 5-15 to 7-15 p.m. | Mathematics | Chemistry and [Physics |
| THURSDAY | 8-15 to 9-45 a.m. | Mathematics | Military Topography |
| | 10-0 to 11-45 a.m. | | |
| | 12-0 noon to 1-0 p.m. | Battalion Drill | Battalion Drill |
| | 2-10 to 3-10 p.m. | | Squad Drill |
| | 2-30 to 3-30 p.m. | Gymnastics | |
| | 4-0 to 5-0 p.m. | Drawing | French or German |
| | 5-15 to 7-15 p.m. | French or German | Model Drawing |
| FRIDAY | 8-15 to 9-45 a.m. | Chemistry and Physics. | Fortification |
| | 10-0 to 11-45 a.m. | | |
| | 12 noon to 1-0 p.m. | Battalion Drill | Battalion Drill or Gymnastics |
| | 2-30 to 3-30 p.m. | Gymnastics | Squad Drill |
| | 4-0 to 5-0 p.m. | French or German | Drawing |
| | 5-15 to 7-15 p.m. | Model Drawing | French or German |
| SAT. | 8-15 to 9-45 a.m. | Fortification | Military Topography |
| | 10-0 to 11-30 a.m. | | |

**Appendix 2**

The following map shows the places in India where MacMahon is known to have been stationed.

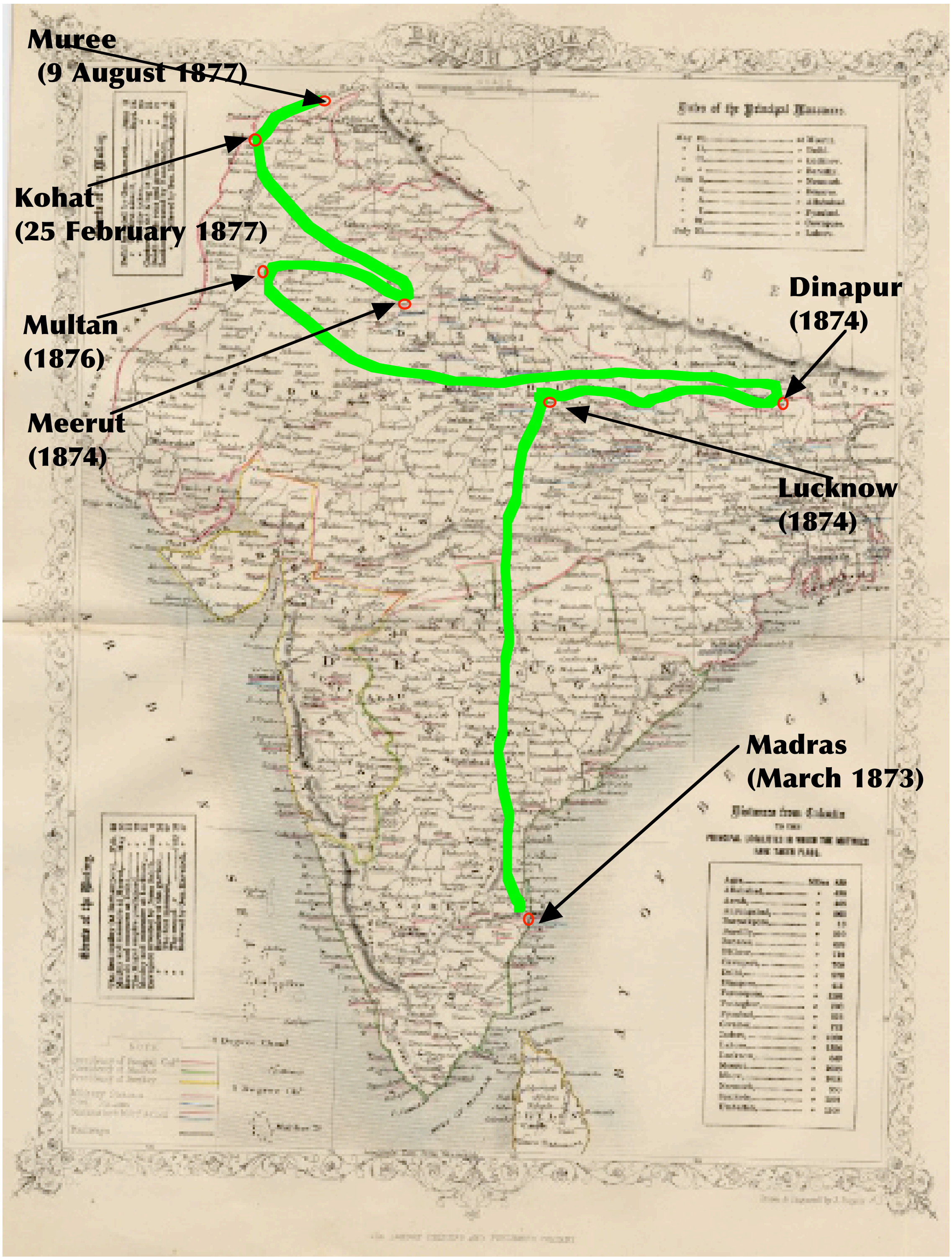

**Appendix 3**

The following is a family tree of the MacMahon family from MacMahon's father through to the present day descendants of his brothers.

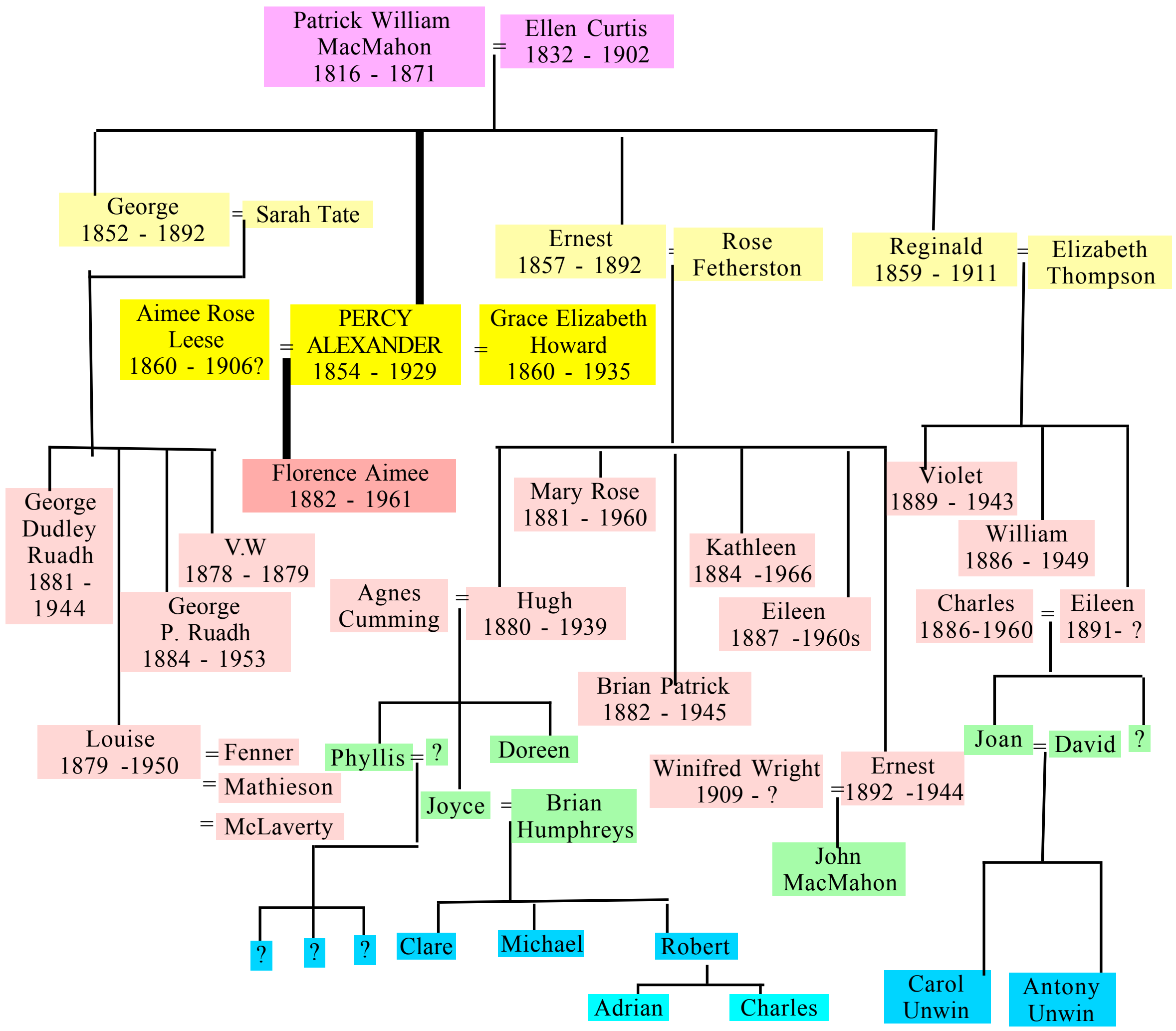

## Appendix 4
## From ballistics to invariant theory

This appendix outlines the main steps taken by MacMahon in moving from the military study of ballistics to mathematical work in invariant theory.

### *Ballistics*

In the *Proceedings of the Royal Artillery Institution* Greenhill wrote a paper[446] in which he began with the equation:

$$v\frac{d\psi}{dt} = -g\cos\psi,$$

where $v$ is the velocity of a projectile at a point on the trajectory where the tangent is inclined at an angle $\psi$ to the horizon, and the resistance of the air varies as some power of the velocity. Bashforth had determined that the resistance varied as the cube of the velocity, so by using this in his general solution, Greenhill obtained three equations for the *t, x* and *y* coordinates of the projectile in terms of the tangent $p$ of $\psi$ :

$$t = \frac{w}{g}\int_p^a \frac{dp}{\left(3a + a^3 - 3p - p^3\right)^{\frac{1}{3}}}$$

$$x = \frac{w^2}{g}\int_p^a \frac{dp}{\left(3a + a^3 - 3p - p^3\right)^{\frac{2}{3}}}$$

$$y = \frac{w^2}{g}\int_p^a \frac{pdp}{\left(3a + a^3 - 3p - p^3\right)^{\frac{2}{3}}}$$

where $w$ is the terminal velocity of the projectile in air, and $a$ is the point on the trajectory where the velocity would be infinite (and then using this point as the origin of the coordinate system); that is, the actual trajectory is part of a larger curve with vertical asymptotes at $\pm a$:

[446] [Greenhill, 1884].

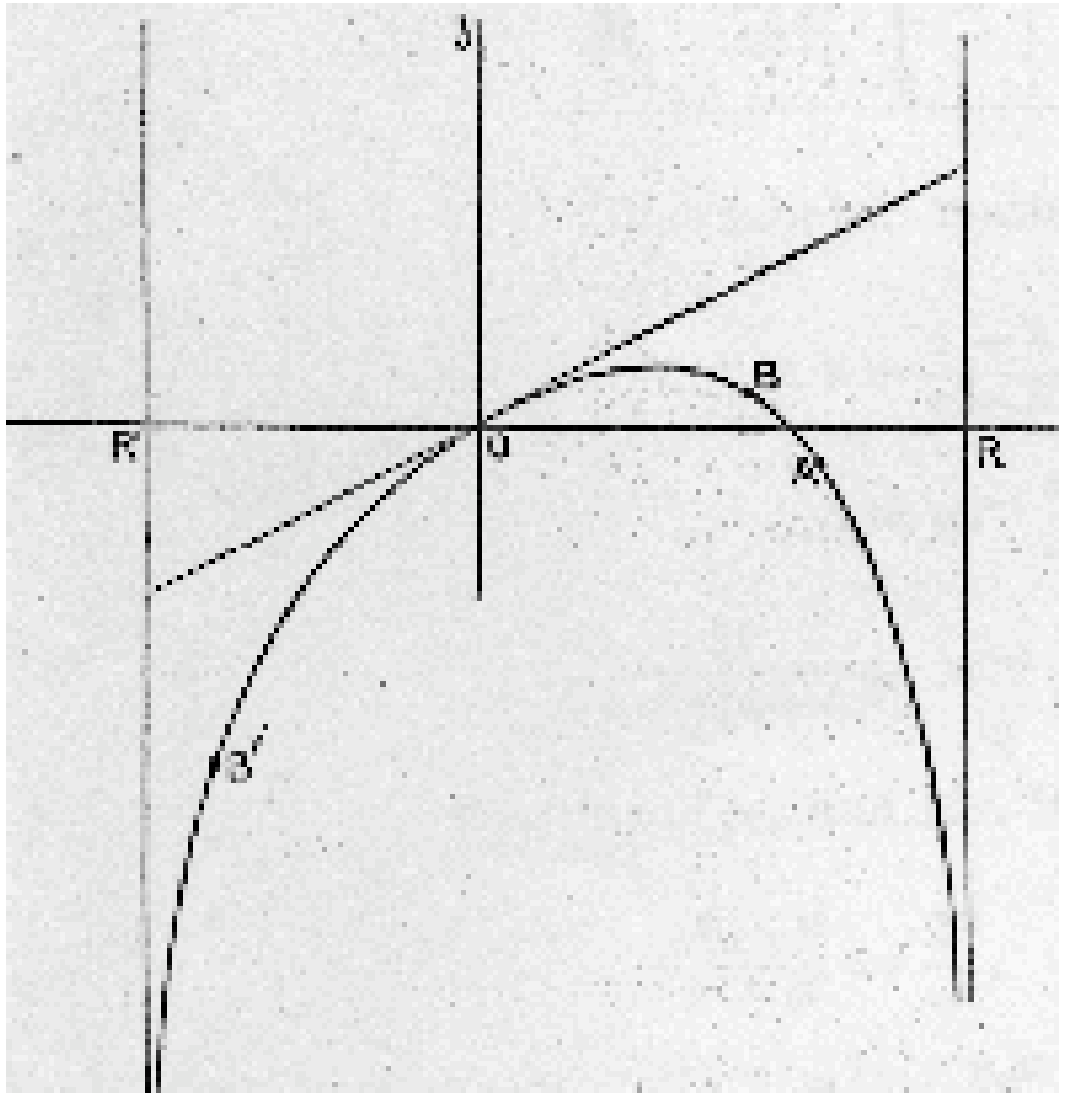

Fig 20: Greenhill's diagram for trajectory calculations

Greenhill solved these integrals using substitutions that convert them into elliptic integrals, and provided a complete solution to the motion when the resistance varies as the cube of the velocity. As a method for calculating the actual curves, however, he remarked: "The extreme complication of the formulae ... render this practically useless." Appended to the paper containing this mammoth calculation are tables, calculated by MacMahon from Bashforth's results of 1880, giving the reduced terminal velocity of a projectile and values of various quantities needed to plot the trajectories.

MacMahon also published his own analysis in 1884[447] of the integral formulas for the co-ordinates of the projectile, starting with that for the abscissa:

$$x = \frac{w^2}{g}\int_p^a \frac{dp}{\left(3a + a^3 - 3p - p^3\right)^{\frac{2}{3}}}$$

He had already published two papers, one in the *Comptes Rendus*[448] and one in the *Quarterly Journal of Mathematics*[449] on a problem he stated had been posed by M. Allégret[450] concerning the differential equation:

[447] [MacMahon, 1884].
[448] [MacMahon, 1881, [2;17]].
[449] [MacMahon 1883, [9;17]].
[450] Alexandre Allégret, a French mathematician about whom little is known.

$$\frac{dx}{\left(A+3Bx+3Cx^2+Dx^3\right)^{\frac{2}{3}}} + \frac{dy}{\left(A+3By+3Cy^2+Dy^3\right)^{\frac{2}{3}}} = 0$$

and a natural extension of the problem arising out of the observation that the solution is symmetrical:

$$X^{-\frac{2}{3}}dx + Y^{-\frac{2}{3}}dy + Z^{-\frac{2}{3}}dz = 0.$$

This extension to three terms has the same general form as Greenhill's trajectory equations, and MacMahon tackled it as a problem in symmetric functions. Using Euler's step-by-step method to solve the integration problem, MacMahon found that for the symmetric solution $z_2^2 z_3^2 + z_3^2 z_1^2 + z_1^2 z_2^2 = 2\left(z_1 z_2 z_3 - 2b^3\right)\left(z_1 + z_2 + z_3\right)$ we have $x_1 + x_2 + x_3 = 0$. From these relationships MacMahon was able to deduce other relationships which were important in tracing the curve.

MacMahon's energies were thus started on a road which led from the very practical considerations of gunnery to symmetric functions, and via partitions to combinatorial analysis of a very abstract nature.

***Invariant theory in the late 19th century***

The roots of MacMahon's work on symmetric functions lie in the researches into invariant theory started in the 1840s by Boole, Salmon, Cayley and Sylvester.[451] Below is a brief outline of where the study had reached when MacMahon became involved.

The example of a homogeneous binary quartic with binomial coefficients is used to illustrate the methods developed by Cayley, and then Sylvester. Two forms were commonly used:

(form 1) $\left(a,b,c,d,e\right)\left(x,y\right)^4 = ax^4y^0 + 4bx^3y^1 + 6cx^2y^2 + 4dx^1y^3 + ex^0y^4$

or (form 2) $\left(a_0,a_1,a_2,a_3,a_4\right)\left(x,y\right)^4 = a_0x^4y^0 + 4a_1x^3y^1 + 6a_2x^2y^2 + 4a_3x^1y^3 + a_4x^0y^4$.

[451] Descriptions of these researches may be found in K. H. Parshall's book, *James Joseph Sylvester: Life and Work in Letters* [Parshall, 1998], in a paper by the same author, *Toward a history of nineteenth century invariant theory* [Parshall, 1989], and also in [Crilly, 1986], [Crilly, 1988] and [Crilly, 2006].

This quartic has *order* 4 (the highest power of the variables $x$ and $y$), and *degree* 1 (the highest power of the coefficients $a$, $b$, $c$, $d$ and $e$), and is written as degree-order (1.4); each term has *weight* 4. It has two independent invariants; the first is $I \equiv ae - 4bd + 3c^2$, with degree 2 and weight 4; this is a different definition of weight from that used for the quantic and is the sum of the suffices of the coefficients. To see this it is easier to express the invariant using the notation of form 2 above: $I \equiv a_0a_4 - 4a_1a_3 + 3a_2^2$. The second invariant, the *catalecticant* of Sylvester[452], has degree 3 and weight 6: $J \equiv ace + 2bcd - ad^2 - b^2e - c^3$,

and can be written as the determinant:

$$\begin{vmatrix} a & b & c \\ b & c & d \\ c & d & e \end{vmatrix}$$

Any determinant of the coefficients written in this cyclic manner is the *catalecticant invariant.* $I^3$ and $J^2$ are thus both of degree 6 and weight 12, and the product $I^3J^{-2}$ is known as an absolute invariant of the quartic, where the power $s$ of $M$ in the definition on page 24 is 0; that is, the linear transformation of the variables leaves the invariant absolutely unchanged.

Another, non-independent, invariant of the quartic is the discriminant, $I^3 - 27J^2$. Such a relationship between invariants was named a *syzygy* by Sylvester; it means 'conjunction' or 'opposition' and is usually associated with planetary alignments.

The quartic also has two covariants. One is called the *Hessian* (see below for a definition):

$$H \equiv \left(ac - b^2\right)x^4 + 2\left(ad - bc\right)x^3y + \left(ae + 2bd - 3c^2\right)x^2y^2 + 2\left(be - cd\right)xy^3 + \left(ce - d^2\right)y^4.$$

The other is a sextic covariant, or cubicovariant:

[452] Sylvester was renowned for inventing terms to describe mathematical objects. "Catalecticant" was coined in 1851 - full details of the derivation and use of the word are at *members.aol.com/jeff570/c.html.*

$$\begin{aligned}G \equiv & \left(a^2d - 3abc + 2b^3\right)x^6 + \left(a^2e + 2abd - 9ac^2 + 6b^2c\right)x^5y \\ & +\left(5abe - 15acd + 10b^2d\right)x^4y^2 + \left(10b^2e - 10ad^2\right)x^3y^3 \\ & +\left(15bce - 5ade - 10bd^2\right)x^2y^4 + \left(9c^2e - ae^2 - 2bde - 6cd^2\right)xy^5 \\ & +\left(3cde - be^2 - 2d^3\right)y^6\end{aligned}$$

Cayley felt that the enumeration and tabulation of invariants and covariants was an important task. In his 1896 LMS Presidential address, MacMahon said:

> In no case did Cayley's interest seem to cease with solving the question of enumeration. He desired, up to a certain point, to see the things themselves, and, moreover, considered it very important to attain this end.

The first problem to be tackled was how to calculate invariants for a given quantic. Cayley's first method was the *hyperdeterminant calculus*. He noticed that by writing a quantic twice using different variables (called *cogredient variables*, and subject to the same linear transformation), $(x_1, y_1), (x_2, y_2)$, and then using a differential operator gives an invariant or covariant of the system of two quantics. By taking the two cogredient quantics as identical, an invariant or covariant of a single function is then obtained.

The differential operator he used was

$$\begin{vmatrix} \dfrac{\partial}{\partial x_1} & \dfrac{\partial}{\partial x_2} \\ \dfrac{\partial}{\partial y_1} & \dfrac{\partial}{\partial y_2} \end{vmatrix} = \frac{\partial}{\partial x_1}\frac{\partial}{\partial y_2} - \frac{\partial}{\partial x_2}\frac{\partial}{\partial y_1}.$$

This operator was written symbolically as $\overline{12}$, where the 1 and 2 refer to the suffices of the variables; any power of the operator produces a *concomitant* (the collective term used for invariants and covariants). The first power is called the *Jacobian*[453]; the second power is called the *Hessian*[454] .

[453] This term is used to denote the matrix of partial derivatives of *n* equations in *n* variables.
[454] This term is used to denote the matrix of second partial derivatives of *n* equations in *n* variables.

The invariant $I \equiv ae - 4bd + 3c^2$ of the quartic is obtained by applying the fourth power of the hyperdeterminant,

$$\overline{12}^4 = \begin{vmatrix} \frac{\partial}{\partial x_1} & \frac{\partial}{\partial x_2} \\ \frac{\partial}{\partial y_1} & \frac{\partial}{\partial y_2} \end{vmatrix}^4 = \left( \frac{\partial}{\partial x_1}\frac{\partial}{\partial y_2} - \frac{\partial}{\partial x_2}\frac{\partial}{\partial y_1} \right)^4 .$$

The calculation proceeds as follows. First two cogredient binary quartics are defined:

$$U \equiv ax_1^4 + 4bx_1^3 y_1 + 6cx_1^2 y_1^2 + 4dx_1 y_1^3 + ey_1^4$$
$$V \equiv a'x_2^4 + 4b'x_2^3 y_2 + 6c'x_2^2 y_2^2 + 4d'x_2 y_2^3 + e'y_2^4$$

The hyperdeterminant is expanded:

$$\overline{12}^4 = \begin{vmatrix} \frac{\partial}{\partial x_1} & \frac{\partial}{\partial x_2} \\ \frac{\partial}{\partial y_1} & \frac{\partial}{\partial y_2} \end{vmatrix}^4 = \left( \frac{\partial}{\partial x_1}\frac{\partial}{\partial y_2} - \frac{\partial}{\partial x_2}\frac{\partial}{\partial y_1} \right)^4$$
$$= \frac{\partial^4}{\partial x_1^4}\frac{\partial^4}{\partial y_2^4} - 4\frac{\partial^4}{\partial x_1^3 \partial y_1}\frac{\partial^4}{\partial x_2 \partial y_2^3} + 6\frac{\partial^4}{\partial x_1^2 \partial y_1^2}\frac{\partial^4}{\partial x_2^2 \partial y_2^2} - 4\frac{\partial^4}{\partial x_1 \partial y_1^3}\frac{\partial^4}{\partial x_2^3 \partial y_2} + \frac{\partial^4}{\partial x_2^4}\frac{\partial^4}{\partial y_1^4}$$

Applying this operator to the product $UV$ gives:

$$24a \cdot 24e' - 4(24b \cdot 24d') + 6(24c \cdot 24c') - 4(24b' \cdot 24d) + 24a' \cdot 24e$$

Ignoring the common factor of $2 \times 24^2$ gives:

$$ae' - 4(bd') + 6(cc') - 4(b'd) + a'e$$

The coefficients are then put equal to one another, and the result follows (again ignoring the common factor 2): $ae - 4bd + 3c^2$.

The second invariant of the above quartic, $J \equiv ace + 2bcd - ad^2 - b^2e - c^3$, is of degree 3 in the coefficients, and is generated by extending the method to a system of three quantics and using the operator $\overline{12}^2\,\overline{23}^2\,\overline{31}^2$.

Salmon[455] gave the following example to show how covariants are used to solve equations. If the cubic $(a,b,c,d)(x,y)^3 = ax^3 + 3bx^2y + 3cxy^2 + dy^3$ can be linearly transformed to the correct form, it can immediately be factorised. So it is supposed that the linear transformation is such that $ax^3 + 3bx^2y + 3cxy^2 + dy^3$ maps to $AX^3 + 3BX^2Y + 3CXY^2 + DY^3$. Then the Hessians of both equations must be equal:

$$(ac - b^2)x^2 + (ad - bc)xy + (bd - c^2)y^2 = (AC - B^2)X^2 + (AD - BC)XY + (BD - C^2)Y^2$$

Suppose now that the transformation is such that *B* and *C* vanish. Then the Hessian reduces to *ADXY,* so that *X* and *Y* can be taken as the factorisation of the quadratic Hessian and *A* and *D* can be obtained by comparing coefficients.

Salmon's numerical example was this: the Hessian of the cubic $4x^3 + 9x^2y + 18xy^2 + 17y^3$ is $180(x + 3y)(3x + y)$. So on taking $y = 1$ the cubic equation $4x^3 + 9x^2 + 18x + 17 = 0$ can be transformed into $A(x + 3)^3 + D(3x + 1)^3$. Expanded, this gives:

$$(A + 27D)x^3 + (9A + 27D)x^2 + (27A + 9D)x + 27A + D.$$

So *A* and *D* can be found by solving the simultaneous equations formed by comparing, say, the first and last coefficients:

$$\left.\begin{aligned} A + 27D = 4 \\ 27A + D = 17 \end{aligned}\right\} \Rightarrow A = \tfrac{5}{8}, \quad D = \tfrac{1}{8}$$

The transformed equation $\frac{5}{8}(x + 3)^3 + \frac{1}{8}(3x + 1)^3 = 0$ may now be readily solved.

During the 1850s Cayley stopped using the hyperdeterminant calculus, which suffers from the defect

455 [Salmon, 1859, p. 59].

that although all hyperdeterminants correspond to concomitants, the converse is not true. He returned to an earlier theorem which asserts that any invariant must satisfy two differential equations. Two operators are defined:

$$\Omega \equiv a_0 \frac{\partial}{\partial a_1} + 2a_1 \frac{\partial}{\partial a_2} + 3a_2 \frac{\partial}{\partial a_3} + \cdots + pa_{p-1} \frac{\partial}{\partial a_p}$$
$$\mathrm{O} \equiv pa_1 \frac{\partial}{\partial a_0} + (p-1)a_2 \frac{\partial}{\partial a_1} + (p-2)a_3 \frac{\partial}{\partial a_2} + \cdots + a_p \frac{\partial}{\partial a_{p-1}}$$

Any invariant $I$ must then satisfy $\Omega I = \mathrm{O}I = 0$. A proof of this may be found in Elliott[456].

$\Omega$ may then be used to generate invariants. The process is illustrated by calculating the invariant $J$ of degree 3 of the above quartic. The first step in the process is to note that, in an invariant of degree $i$ = 3 of a quartic where the order $p = 4$, the weight $w$ of each term must satisfy[457] $w = \frac{1}{2}ip$, so the weight of each term in the invariant sought is 6. Thus the form of the invariant must be: $a_0a_2a_4 + \lambda a_0a_3^2 + \rho a_2^3 + \upsilon a_1a_2a_3 + \mu a_1^2 a_4$. Applying the operator gives:

$$\Omega\left(a_0a_2a_4 + \lambda a_0a_3^2 + \rho a_2^3 + \upsilon a_1a_2a_3 + \mu a_1^2a_4\right) = 0$$
$$a_0\left(\upsilon a_2a_3 + 2\mu a_1a_4\right) + 2a_1\left(a_0a_4 + 3\rho a_2^2 + \upsilon a_1a_3\right) + 3a_2\left(2\lambda a_0a_3 + \upsilon a_1a_3\right) + 4a_3\left(a_0a_2 + \mu a_1^2\right) = 0$$
$$\left(\upsilon + 6\lambda + 4\right)a_0a_2a_3 + \left(2\mu + 2\right)a_0a_1a_4 + \left(6\rho + 3\upsilon\right)a_1a_2^2 + \left(2\upsilon + 4\mu\right)a_1^2a_3 = 0$$

Equating the coefficients to 0 then yields the values for the invariant:

$$J \equiv a_0a_2a_4 - a_0a_3^2 - a_2^3 + 2a_1a_2a_3 - a_1^2a_4$$

This illustrates the connection with partitions, since for each term in the invariant to have weight 6, the suffixes must sum to 6 with no part greater than 4 (since the coefficient with the highest suffix is $a_4$). Similarly, for the invariant of degree 2, the weight $I$ of each term must be 4, and so the suffices can only be (0,4), (1,3) and (2,2).

For an invariant of degree 4, the weight of each term must be 8 and the process yields the invariant

[456] [Elliott, 1895, p. 13].
[457] [Elliott, 1895, pp. 31 - 39].

$a^2e^2 - 8abde + 6ac^2e + 16b^2d^2 - 24bc^2d + 9c^4$. However, this is simply $I^2$, and so is not an independent invariant.

Similarly, the invariant of degree 5 is simply *IJ*:

$$a^2ce^2 - ab^2e^2 - a^2d^2e - 2abcde + 4b^3de + 2ac^3e - 3b^2c^2e$$
$$+ 4abd^3 - 3ac^2d^2 - 8b^2cd^2 + 10bc^3d - 3c^5$$

Several questions then arise: are there any other independent invariants ? how many independent, or *asyzygetic*, invariants are there ? can the syzygies be detected easily ?

On the Continent, Paul Gordan[458], encouraged by Clebsch[459], had provided a constructive proof in 1868 that the invariants and covariants of a binary quantic of any degree could be formed as linear combinations of a finite basis set; this corrected Cayley's 1856 error for binary quintics for which he asserted that the basis set was infinite. The method used was similar in spirit to Cayley's original hyperdeterminant methods, but more abstract. It became known as a 'symbolic' method, in contrast to the more direct 'algebraic' methods of Cayley and Sylvester.

Despite twenty years of work, Gordan was unable to extend the proof to quantics with more than two variables. In 1888, Hilbert[460] provided a proof that a finite basis exists for higher degree quantics, but gave no clue as to the construction of such a basis, which Gordan found deeply unsatisfactory; he is reputed to have described Hilbert's proof as theology, rather than mathematics[461].

MacMahon's work in invariant theory became important in 1884, when showed that there was correspondence between semivariants and non-unitary partitions, as described in Chapter 3.

[458] Paul Gordan,1837 - 1912.
[459] R. F. A. Clebsch, 1833 - 1872.
[460] David Hilbert, 1862 - 1943.
[461] Quoted in Max Noether's obituary for Gordan, in the *Mathematische Annalen* **75** (1914), 1 - 41 - from an e-mail to the Historia Mathematica electronic discussion list dated 1 September 2001 by Colin McLarty.

## Appendix 5
## Partitions - early English work

This appendix describes early English work on partitions during the middle of the 19th century. using recurrence relations. It is this tradition of algebra that typifies the algebraic environment in which MacMahon worked. He sketched the history of the subject[462] in his *Encyclopaedia Britannica* article on Combinatorial Analysis. MacMahon was always conscious of the need to know and understand the history of a topic, so this detailed description has been included to set the scene for the style of partition theory in England at the time.

### *Nicholson and J. Hamilton*

The fundamental problem in partition theory is to enumerate the partitions of a given number. Euler had started a detailed analysis of this problem in *Introductio in analysi infinitorum*[463], and a great deal of development had been done by Continental mathematicians over the following half century. The idiosyncratic English work began in 1818, when Peter Nicholson[464] included an article on partitions in his *Essays on the Combinatorial Analysis*[465]. The article was not by Nicholson himself but by J. A. Hamilton, described as a Professor of Music. He referred to Euler's work and then described the use of a table to calculate the number of ways $\binom{n}{m}$ of partitioning $n$ into $m$ parts. Part of his table is reproduced below. Hamilton did not describe its construction, claiming that it would be obvious after some study. It was evidently based on a recurrence relation, rather than the generating functions of Euler and the Continental mathematicians, although Hamilton was not explicit about the detail. It was this use of recurrence relations and difference equations that characterised the English work on the subject.

[462] [MacMahon, 1910, [76;19]].
[463] [Euler, 1748] The question of partitions dates back to 1669 and a letter from Leibniz to Bernoulli - a detailed history of the subject prior to the 19th century can be found in [Dickson, 1992, Volume 2].
[464] Peter Nicholson, 1765 - 1844, architect and private teacher of mathematics [source: Lloyd, 2003].
[465] [Nicholson, 1818, p. 30Q].

| *n* | 1 | 2 | 3 | 4 | 5 | 6 | 7 | 8 | 9 | 10 | 11 | 12 | 13 | 14 | 15 | 16 | 17 | 18 |
|---|---|---|---|---|---|---|---|---|---|---|---|---|---|---|---|---|---|---|
| | 1 | 1 | 1 | 1 | 1 | 1 | 1 | 1 | 1 | 1 | 1 | 1 | 1 | 1 | 1 | 1 | 1 | 1 |
| | | 1 | 1 | 1 | 1 | 1 | 1 | 1 | 1 | 1 | 1 | 1 | 1 | 1 | 1 | 1 | 1 | 1 |
| | | | | 1 | 1 | 1 | 1 | 1 | 1 | 1 | 1 | 1 | 1 | 1 | 1 | 1 | 1 | 1 |
| | | | | | | 1 | 1 | 1 | 1 | 1 | 1 | 1 | 1 | 1 | 1 | 1 | 1 | 1 |
| | | | | | | | | 1 | 1 | 1 | 1 | 1 | 1 | 1 | 1 | 1 | 1 | 1 |
| | | | | | | | | | | 1 | 1 | 1 | 1 | 1 | 1 | 1 | 1 | 1 |
| | | | | | | | | | | | | 1 | 1 | 1 | 1 | 1 | 1 | 1 |
| | | | | | | | | | | | | | | 1 | 1 | 1 | 1 | 1 |
| | | | | | | | | | | | | | | | | 1 | 1 | 1 |
| | | | | | | | | | | | | | | | | | | 1 |
| | | 1 | 1 | 2 | 2 | 3 | 3 | 4 | 4 | 5 | 5 | 6 | 6 | 7 | 7 | 8 | 8 | 9 |
| | | | 1 | 1 | 2 | 2 | 3 | 3 | 4 | 4 | 5 | 5 | 6 | 6 | 7 | 7 | 8 | 8 |
| | | | | | | 1 | 1 | 2 | 2 | 3 | 3 | 4 | 4 | 5 | 5 | 6 | 6 | 7 |
| | | | | | | | | | 1 | 1 | 2 | 2 | 3 | 3 | 4 | 4 | 5 | 5 |
| | | | | | | | | | | | | 1 | 1 | 2 | 2 | 3 | 3 | 4 |
| | | | | | | | | | | | | | | | 1 | 1 | 2 | 2 |
| | | | | | | | | | | | | | | | | | | 1 |
| | | | 1 | 1 | 2 | 3 | 4 | 5 | 7 | 8 | 10 | 12 | 14 | 16 | 19 | 21 | 24 | 27 |
| | | | | 1 | 1 | 2 | 3 | 4 | 5 | 7 | 8 | 10 | 12 | 14 | 16 | 19 | 21 | 24 |
| | | | | | | | | 1 | 1 | 2 | 3 | 4 | 5 | 7 | 8 | 10 | 12 | 14 |
| | | | | | | | | | | | | 1 | 1 | 2 | 3 | 4 | 5 | 7 |
| | | | | | | | | | | | | | | | | 1 | 1 | 2 |
| | | | | 1 | 1 | 2 | 3 | 5 | 6 | 9 | 11 | 15 | 18 | 23 | 27 | 34 | 39 | 47 |

Table 5: Hamilton's table of partitions

To read the table, find *n* along the top, and look at the row labelled $\binom{n}{m}$; for example, the intersection of column 5 and row $\binom{n}{2}$ gives 2 as the number of ways that 5 can be partitioned into 2 parts. The rows intermediate to the labelled rows are for calculation purposes only, and are the previous labelled row simply shifted to the right by a number of spaces corresponding to the next value of *m*.

Nicholson's own contribution to the subject occurs in his *Essay 1*[466], in Problem V, ('to find the *m*th term of any order *n* of combinations with repetitions raised from a base *a, b, c*, &c in alphabetical order'), Example II lists all the compound divisors of 360 , and in Example III where all the possible combinations of 1s, 2s and 3s in parcels of 4 are listed. Similarly, Example IV lists "all the possible ways of adding 1, 3 and 7 in one's, two's and three's". Problem VI Example I shows how to "decompose 14 into parcels of 6 digits in every possible combination of the 9 digits", and Example II partitions 9 "in all possible ways". Problem 10 Example II is to "break 10 with the nine digits as constituents."

In each case, Nicholson gave a constructive algorithm for generating the desired partitions, but confessed that he had been unable to find a formula for calculating the number of partitions. He felt that such a formula would have been very useful as a check on his algorithm to ensure that no partitions had been missed.

After Hamilton's article in Nicholson's book, the earliest interest in partitions in England[467] was in a paper by A. De Morgan published in 1843 in the *Cambridge Mathematical Journal,* entitled *On a new species of equation of differences*[468]. This was not attributed to him at the time, but was referred to by Henry Warburton in a paper read on 1 March 1847 and printed in the *Transactions of the Cambridge Philosophical Society* in 1849[469]. MacMahon, in his article for the tenth edition of the *Encyclopaedia*

[466] [Nicholson, 1818].
[467] The source of this information is principally L. E. Dickson's work on the *History of the Theory of Numbers,* which provides detailed references for the Continental and English work, but does not mention the Hamilton article.
[468] [De Morgan, 1843].
[469] [Warburton, 1849].

*Britannica* published in 1902[470] , stated that interest in England in the subject was first aroused by a letter from De Morgan to Warburton in 1846. Warburton refers to this letter thus:

> In the Autumn of 1846, having communicated a Theorem on the Partition of Numbers to Professor A. De Morgan, I received from him an obliging reply, wherein he intimated a wish that I would turn my attention to Combinations, as a department in Mathematics, which, he thought, needed much cultivation.

This suggests that MacMahon was unaware of De Morgan's 1843 paper. MacMahon's article then went on to trace the history of the subject to a subsequent paper by Herschel in 1850[471] and the work of Cayley[472] and Sylvester in 1856. No further advances were published until 1882 when graphical methods were introduced by Sylvester ("Ferrers graphs[473] "), Franklin and Durfee. Although aware of Kirkman's work in 1855[474] and 1857[475], MacMahon did not feel it important enough to mention, either in the *Britannica* article, nor in the address he delivered as retiring President to the London Mathematical Society in 1896[476], despite listing Kirkman's papers in the bibliography. A description now follows of the unique English approach to partitions made by De Morgan, Warburton, Herschel and Kirkman.

### *De Morgan*

Continental work on the topic of partitions followed Euler's approach of using generating functions, but De Morgan seemed to be unaware of any of this previous work:

> It is well known that the number of ways in which the number $n$ may be put together is $2^{n-1}$ [...].But in this theorem orders are counted[477] .[...].The corresponding problem, namely to find in how many ways the number $n$ can be constructed, different orders not being counted as different ways, is of much greater difficulty and I am not aware of anything having been written upon it.

[470] [MacMahon, 1902].
[471] [Herschel, 1850].
[472] [Cayley, 1856].
[473] Norman MacLeod Ferrers, 1829 - 1903, vice chancellor of Cambridge University from 1884, was the first to point out the conjugacy of partitions represented by arrays or dots, but the use of Ferrers name to label such diagrams is due to Sylvester. See [Kimberling, 1999] for more details.
[474] [Kirkman, 1855].
[475] [Kirkman, 1857].
[476] [Andrews, 1986, Vol II pp. 839 -866].
[477] Such a partition is called a *composition* - see, for example, [Riordan, 1958].

His approach was to use a method of finite differences. He used the notation $u_{x,y}$ to represent the number of ways of making the number $x$ from numbers not exceeding $y$. It was then clear that the number of ways of making, say, 16 from numbers less than or equal to 6 is equal to the sum of the ways of making 10 from numbers less than or equal to 1, 2, ... 6, since you can simply add 6 to each of these ways. In symbols:

$$u_{16,6} = u_{10,1} + u_{10,2} + u_{10,3} + u_{10,4} + u_{10,5} + u_{10,6}$$

Generalising gives $u_{x,y} = u_{x-y,1} + u_{x-y,2} + u_{x-y,3} + \ldots + u_{x-y,y}$ (1). De Morgan then noted that $u_{x+1,y+1} = u_{x-y,1} + u_{x-y,2} + u_{x-y,3} + \ldots + u_{x-y,y} + u_{x-y,y+1}$ (2). Subtracting (1) from (2) gives $u_{x+1,y+1} - u_{x,y} = u_{x-y,y+1}$ or, subtracting 1 from each suffix, $u_{x,y} - u_{x-1,y-1} = u_{x-y,y}$ (3). By noting that $u_{x,1} = u_{x,x} = u_{x,x-1} = 1$ he constructed a table of values of $u_{x,y}$ , for values of $x$ up to 10 (reproduced below). Although his table is correct, his description of how it was constructed contained an error, probably the result of poor proof-reading.

| | **Values of y (greatest part)** | | | | | | | | | |
|---|---|---|---|---|---|---|---|---|---|---|
| **x** | 1 | 2 | 3 | 4 | 5 | 6 | 7 | 8 | 9 | 10 |
| 1 | 1 | | | | | | | | | |
| 2 | 1 | 1 | | | | | | | | |
| 3 | 1 | 1 | 1 | | | | | | | |
| 4 | 1 | 2 | 1 | 1 | | | | | | |
| 5 | 1 | 2 | 2 | 1 | 1 | | | | | |
| 6 | 1 | 3 | 3 | 2 | 1 | 1 | | | | |
| 7 | 1 | 3 | 4 | 3 | 2 | 1 | 1 | | | |
| 8 | 1 | 4 | 5 | 4 | 3 | 2 | 1 | 1 | | |
| 9 | 1 | 4 | 7 | 6 | 5 | 3 | 2 | 1 | 1 | |
| 10 | 1 | 5 | 8 | 9 | 7 | 5 | 3 | 2 | 1 | 1 |

Table 6: Part of De Morgan's table of partitions

De Morgan then derived a general formula for $u_{x,y}$ . He started with equation (3) for $y$ = 2, $u_{x,2} - u_{x-1,1} = u_{x-2,2}$ which, rearranged and noting that $u_{x,1} = 1$ for all values of $x$ (since there is only one way of making a number out of 1s), gives the recurrence relation $u_{x,2} - u_{x-2,2} = 1$.

De Morgan then stated "the complete integral of which, the constants being found from $u_{1,2}$ = 0 and $u_{2,2}$ = 1 is $u_{x,2} = \frac{1}{2}x - \frac{1}{4} + \frac{1}{4}(-1)^x$." The integration referred to is a process of summation which works as follows. Writing $u_{x,2} - u_{x-2,2} = A_x$, where $A_x$ is a function of $x$ gives:

$$\begin{aligned}
u_{3,2} - u_{1,2} &= A_3 \\
u_{4,2} - u_{2,2} &= A_4 \\
u_{5,2} - u_{3,2} &= A_5 \\
u_{6,2} - u_{4,2} &= A_6 \\
&\vdots \\
u_{x-2,2} - u_{x-4,2} &= A_{x-2} \\
u_{x-1,2} - u_{x-3,2} &= A_{x-1} \\
u_{x,2} - u_{x-2,2} &= A_x
\end{aligned}$$

In summing this sequence of equations, most of the terms on the left-hand sides cancel out[478], leaving $u_{x,2} + u_{x-1,2} - u_{2,2} - u_{1,2} = \sum_{i=3}^{x} A_i$, which gives, from $u_{1,2} = 0$, $u_{2,2} = 1$ and $A_x = 1$, $u_{x,2} + u_{x-1,2} - 1 = x - 2$.

Now it is necessary to note that $u_{x,2} = u_{x-1,2}$ if $x$ is odd and $u_{x,2} = u_{x-1,2} + 1$ if $x$ is even.

Thus, for $x$ odd:

$$\begin{aligned}
u_{x,2} + u_{x-1,2} - 1 &= x - 2 \\
2u_{x,2} - 1 &= x - 2 \\
u_{x,2} &= \frac{x-1}{2} = \frac{x}{2} - \frac{1}{2}
\end{aligned}$$

and for $x$ even:

$$\begin{aligned}
u_{x,2} + u_{x-1,2} - 1 &= x - 2 \\
2u_{x,2} + 1 - 1 &= x - 2 \\
u_{x,2} &= \frac{x}{2}
\end{aligned}$$

De Morgan combined these into the single formula $u_{x,2} = \frac{1}{2}x - \frac{1}{4} + \frac{1}{4}(-1)^x$. The next case was $y = 3$, which De Morgan disposed of in just two lines.

The next problem was to express $u_{x-1,3}$ and $u_{x-2,3}$ in terms of $u_{x,3}$. This resulted in a *circulating*

[478] A process now called *telescoping.*

function (that we would now call a periodic function) which depends on the value of $x$ (mod 6). De Morgan also obtained a single formula:

$$u_{x,3} = \frac{6x^2 - 7 - 9(-1)^x + 8(\beta^x + \gamma^x)}{72},$$

where β and γ are the complex cube roots of 1. The terms $-7 - 9(-1)^x + 8(\beta^x + \gamma^x)$ are De Morgan's solution to the problem of finding a periodic function which takes the values
0, -6, -24, 18, -24 or -6 according as $x$ is equivalent to 0, 1, 2, 3, 4 or 5 (mod 6).

It is not clear how he came across this solution to the problem, given his apparent unfamiliarity with the idea of generating functions. The use of generating functions, makes the formulas for $u_{x,2}$ and $u_{x,3}$ easy to calculate, as follows.[479]

We can describe a partition of $n$ into 1s and 2s by the equation $n = 1x + 2y$. If $r(n)$ is the number of integer solutions to this equation, then we can use the generating function:

$$f(x) = 1 + \sum_{n=1}^{\infty} r(n)x^n = 1 + r(1)x + r(2)x^2 + r(3)x^3 + r(4)x^4 \ldots$$
$$= (1 + x + x^2 + x^3 + \ldots)(1 + x^2 + x^4 + \ldots) = \frac{1}{(1-x)(1-x^2)} = \frac{1}{(1-x)^2(1+x)}$$

This final expression is decomposed into partial fractions and expanded:

[479] This explanation is based on remarks made by G. H. Hardy in his Inaugural Lecture as Savilian Professor of Geometry at Oxford in 1920.

$$\frac{1}{(1-x)^2(1+x)} = \frac{\frac{1}{4}}{(1+x)} + \frac{\frac{1}{4}}{(1-x)} + \frac{\frac{1}{2}}{(1-x)^2}$$

$$= \tfrac{1}{4}\left(1 - x + x^2 - x^3 + x^4 - \ldots + (-1)^n x^n \ldots\right)$$

$$+ \tfrac{1}{4}\left(1 + x + x^2 + x^3 + x^4 + x^5 \ldots + x^n \ldots\right)$$

$$+ \tfrac{1}{2}\left(1 + 2x + 3x^2 + 4x^3 \ldots + (n+1)x^n \ldots\right)$$

$$= 1 + x + 2x^2 + 2x^3 + \ldots + \left(\frac{n+1}{2} + \tfrac{1}{4} + \tfrac{1}{4}(-1)^n\right)x^n \ldots$$

From this it is clear that the general term $r(n) = \frac{n+1}{2} + \tfrac{1}{4} + \tfrac{1}{4}(-1)^n$ agrees with De Morgan's result $u_{x,2} = \frac{x}{2} - \tfrac{1}{4} + \tfrac{1}{4}(-1)^x$ when $x = n+1$. If De Morgan had been through this process, he chose not to mention it, or thought it so obvious as to be unnecessary to explain.

### *Henry Warburton*[480]

Warburton's approach was similar, in the sense that he used difference equations rather than generating functions. He began by describing how he was led to investigate the topic by considering an application of partitions made by Waring to the development of a power series. He apparently communicated some of his results to De Morgan, who drew his attention to the anonymous 1843 paper described above.

Warburton's first attempt was to consider the number $N$ to be partitioned in the form $mt + r$, so that the number of partitions is a function of $t$. In the case of partitions[481] into two parts (bi-partitions), the modulus $m$ is 2; for partitions into three parts (tri-partitions), the modulus is 6, for four parts (quadri-partitions) it is 12.

[480] Henry Warburton, 1784 - 1858, was educated at Eton and Trinity College, Cambridge, and elected FRS in 1809. He was an MP from 1826 - 1847.

[481] [Warburton, 1849, p. 472]. The fact that this results in two formulas for bi-partitions, six for tri-partitions, twelve for quadri-partitions, and so on, led Warburton to abandon this line of enquiry, since in his view, "the functions obtained would be so many, and of such complexity as to be of little or no practical utility."

His attack on the problem then followed a similar line to De Morgan's, except that instead of considering partitions of $N$ into parts not exceeding $y$ he started by looking at partitions into $p$ parts where each part is not less than some number *h.* Warburton introduced the notation $\left[N, p_{\eta}\right]$ to describe this situation, where $N$ is the number to be partitioned, $p$ is the number of parts and $h$ is the *smallest* part. De Morgan took no account of the number of parts, so it is not possible to get a precise correspondence between the two notations.

He noted that partitions with each part not less than $h$ (the *lower limit*) necessarily include all the partitions with each part not less than $h$ + 1; in more modern notation, $\left[N, p_{\eta}\right] \supset \left[N, p_{\eta+1}\right]$, $\left[N, p_{\eta+1}\right] \supset \left[N, p_{\eta+2}\right]$ etc.

He then sought to derive a formula, using similar techniques to De Morgan as follows. First, if in any partition of $N$ a fixed quantity $q$ is added to, or subtracted from, each part, then $N$ is increased or diminished by $pq$ and the lower limit is increased or diminished by $q$, so giving the identity

$$\left[N, p_{\eta}\right] = \left[N \pm p\theta, p_{\eta\pm\theta}\right] \qquad (1)$$

The first difference equation arises from noting that if take all the partitions of N - $h$ with $p$ - 1 parts and add $h$ as an extra part, then we have partitions of $N$ with $p$ parts. Moreover, each of these partitions has at least one part equal to *h,* thus excluding all the partitions of $N$ with $h$ + 1 as the lower limit. This gives the difference equation

$$\left[N, p_{\eta}\right] - \left[N, p_{\eta+1}\right] = \left[N - \eta, \left(p-1\right)_{\eta}\right] \qquad (2)^{482}$$

This can be seen from the following table of the 14 tri-partitions of 13 with lower limit 1:

[482] Here Warburton used a slightly different form of the suffix notation. Where there was a single letter to be suffixed, he wrote: $p_{\eta}$ but where a suffix was necessary on a compound term, he did not use brackets and a suffix, but instead wrote: $\underset{\eta}{p-1}$ with no brackets. Later on, he also omitted the suffix if it was unity.

| | | | |
|---|---|---|---|
| 11+1+1 | 9+2+2 | 7+3+3 | 5+4+4 |
| 10+2+1 | 8+3+2 | 6+4+3 | |
| 9+3+1 | 7+4+2 | 5+5+3 | |
| 8+4+1 | 6+5+2 | | |
| 7+5+1 | | | |
| 6+6+1 | | | |

Table 7: Tripartitions of 13, after Warburton

So, for example, $\left[13,3_1\right]-\left[13,3_2\right]=\left[12,2_1\right]$ is 14 - 8 = 6 and $\left[13,3_3\right]-\left[13,3_4\right]=\left[10,2_3\right]$ is 4 - 1 = 3.

The 6 bi-partitions of 12 with lower limit 1 are obtained by deleting the 1 from the six tri-partitions of 13 with lower limit 1, and the 3 bi-partitions of 10 with lower limit 3 are similarly obtained from the 3 tri-partitions of 13 with lower limit 3 by deleting the 3.

From (2) it is easy to get $\left[N+\eta,(p+1)_\eta\right]-\left[N+\eta,(p+1)_{\eta+1}\right]=\left[N,p_\eta\right]$. He also derived a further sequence of identities: (3) If $N < p\eta$ then $\left[N,p_\eta\right]=0$, (3a) $\left[0,p_\eta\right]=0$, (4) $\left[N,0_\eta\right]=0$, (5) $\left[0,0_\eta\right]=1$ and (6) if $N = p\eta$ then $\left[p\eta,p_\eta\right]=\left[0,p_0\right]=\left[p,p_1\right]=1$.

Warburton then wrote:

> Professor De Morgan (as he informs me) has, in the *Cambridge Mathematical Journal*, traced Equation (2) to its consequences, in the case where the number of parts is preserved constant, and the variation is thrown on the number to be partitioned.

This seems to imply that De Morgan knew of Euler's result, that the number of partitions of *N* into parts not exceeding *m* is the same as number of partitions of *N* into at most *m* parts, because in De Morgan's work the greatest part is kept constant:

$$u_{x,3}+u_{x-1,3}+u_{x-2,3}-u_{3,3}-u_{2,3}-u_{1,3}=\sum_{i=4}^{x}A_i\,.$$

So perhaps De Morgan knew more than he was revealing.

Warburton then used equation (2) to create a series of equations which could then be telescoped to give

$$\left[N, p_{\eta}\right] - \left[N, p_{\eta+\theta+1}\right] = \sum_{z=0}^{\theta}\left[N - \eta - zp, \left(p-1\right)_{\eta}\right]$$ [483].

Warburton's problem was how to remove the second term on the left-hand side and at the same time transform the right-hand side to a calculable form. He explained how he tackled the former, but felt the latter was clearly trivial enough for the reader to see immediately. A possible explanation for both techniques is outlined below.

The second term in the equation effectively excludes all the partitions of $N$ where every part is greater than $\eta + \theta$. Writing $Y'$ for $\eta + \theta$ and letting $N = pY + r$, where $r < p$ gives $\left[N, p_{\eta+\theta+1}\right] = \left[pY + r, p_{Y'+1}\right]$. Then letting $Y = Y'$ so that $r = 0$ gives:

$$\begin{aligned}\left[pY + r, p_{Y'+1}\right] &= \left[p(\eta + \theta), p_{\eta+\theta+1}\right] \\ &= \left[p\eta, p_{\eta+1}\right] = 0\end{aligned}$$

This uses identity (1) above to eliminate $pq$ and $q$ and identity (6) with the fact that $\left[N, p_{\eta}\right] \supset \left[N, p_{\eta+1}\right]$. So now we have, when $Y = Y'$, $\left[N, p_{\eta}\right] = \sum_{z=0}^{Y-\eta}\left[N - \eta - zp, \left(p-1\right)_{\eta}\right]$.

Now we can deal with the RHS, using identity (1) with $\theta = \eta - 1$ and $p$ - 1 in place of $p$:

This gives the final formula $\left[N, p_{\eta}\right] = \sum_{z=0}^{y-\eta}\left[N - \left(1 + p(\eta - 1)\right) - zp, \left(p-1\right)_{1}\right]$.

This is a recurrence formula, not a single formula like De Morgan's. Warburton illustrated its use with the example:

$$\begin{aligned}\left[31, 5_3\right] &= \sum_{z=0}^{6-3}\left[31 - (1 + 5\times 2) - 5z, 4_1\right] = \left[20, 4_1\right] + \left[15, 4_1\right] + \left[10, 4_1\right] + \left[5, 4_1\right] \\ &= 64 + 27 + 9 + 1 = 101\end{aligned}$$

[483] Warburton did not use the Greek sigma notation for summation, but a simple capital S.

Warburton did not explain where the results for the partitions into four parts with lowest part 1 came from. From these results it is clear that the variation in the summation terms falls upon the number to be partitioned (the 31 is reduced to 20, 15, 10 and 5, whilst the number of parts is constant at 4 and the smallest part constant at 1). In his final section on this topic, he took equation (2), $\left[N, p_\eta\right] - \left[N, p_{\eta+1}\right] = \left[N - \eta, (p-1)_\eta\right]$, and transformed it differently, arriving at a formula in which the number to be partitioned remains constant, but the number of parts varies. He did this by using $N + h$ in place of $N$ and $q+1$ in place of $p$, and using equation (1) to get $\left[N + \eta, (q+1)_\eta\right] - \left[N, q_\eta\right] = \left[N - q\eta, (q+1)_1\right]$.

This then formed the basis for the recurrence relation and the telescoping, by writing equation (2) with $N + 2h$ in place of $N$ and $q+2$ in place of $p$, $N + 3h$ in place of $N$ and $q+2$ in place of $p$, etc., to finish up with $\left[N + p\eta, (q+p)_\eta\right] - \left[N, q_\eta\right] = \sum_{z=1}^{p}\left[N - q\eta, (q+z)_1\right]$.

If $q = 0$ and $h = 1$ the second term on the left disappears, the final result is:

$$\left[N, p_\eta\right] = \sum_{z=0}^{p}\left[N - p\eta, z_1\right]$$

Using the same numerical example as above gives:

$$\begin{aligned}\left[31, 5_3\right] &= \sum_{z=0}^{5}\left[31 - (5 \times 3), z_1\right] = \left[16, 0_1\right] + \left[16, 1_1\right] + \left[16, 2_1\right] + \left[16, 3_1\right] + \left[16, 4_1\right] + \left[16, 5_1\right] \\ &= 0 + 1 + 8 + 21 + 34 + 37 = 101\end{aligned}$$

Having derived this formula, Warburton then proceeded to give an elementary proof of it, based on the observation that it exhausts all the ways of constructing $N - ph$ from $p$ or fewer parts. This is very similar to De Morgan's first recurrence relation which was used to construct Table 6. Warburton then went on to construct a similar table showing the number of partitions $p$ for each value of $N$ up to 12. This table is reproduced below (Table 8) but not in the same form as Warburton. He listed values of $p$ vertically down the left, and values of $N$ horizontally across the top, and provided a sum of each column for the total number of partitions of each *N* (which De Morgan did not give). Below is the transpose of Warburton's table so that it may be more readily compared with De Morgan's table[484].

| | Values of p (number of partitions) | | | | | | | | | | | | |
|---|---|---|---|---|---|---|---|---|---|---|---|---|---|
| N | 1 | 2 | 3 | 4 | 5 | 6 | 7 | 8 | 9 | 10 | 11 | 12 | sum |
| 1 | 1 | | | | | | | | | | | | 1 |
| 2 | 1 | 1 | | | | | | | | | | | 2 |
| 3 | 1 | 1 | 1 | | | | | | | | | | 3 |
| 4 | 1 | 2 | 1 | 1 | | | | | | | | | 5 |
| 5 | 1 | 2 | 2 | 1 | 1 | | | | | | | | 7 |
| 6 | 1 | 3 | 3 | 2 | 1 | 1 | | | | | | | 11 |
| 7 | 1 | 3 | 4 | 3 | 2 | 1 | 1 | | | | | | 15 |
| 8 | 1 | 4 | 5 | 4 | 3 | 2 | 1 | 1 | | | | | 21 |
| 9 | 1 | 4 | 7 | 6 | 5 | 3 | 2 | 1 | 1 | | | | 30 |
| 10 | 1 | 5 | 8 | 9 | 7 | 5 | 3 | 2 | 1 | 1 | | | 42 |
| 11 | 1 | 5 | 10 | 11 | 10 | 7 | 5 | 3 | 2 | 1 | 1 | | 56 |
| 12 | 1 | 6 | 12 | 15 | 13 | 11 | 7 | 5 | 3 | 2 | 1 | 1 | 77 |

Table 8: (after Warburton)

This table is identical to De Morgan's, because the number of ways of making *N* from parts less than or equal to *y* is the same as the number of ways of making *N* out of *y* parts (Euler's Theorem), but Warburton made no mention of this.

### *Herschel*

The work of Warburton inspired J. F. W. Herschel to write a paper in 1850, *On the algebraic expression of the number partitions of which a given number is susceptible*[485]. In this paper Herschel first developed the "general properties of the differences of the powers of the natural numbers" based on previous work done by him in 1815[486] and the Rev. John Brinkley[487] on finite differences. Also mentioned was a paper on circulating functions[488]. It is not until paragraph 26, after ten pages of preliminary work, that Herschel got to the main point of the paper, namely the expression of the number of ways of partitioning $x$ into a given number of parts.

[484] The row and column for $N = 0$ and $p = 0$ have also been omitted.
[485] [Herschel, 1850]. John Frederick William Herschel, 1792 - 1871, is best known as the astronomer who discovered Uranus. He was elected FRS in 1813, and won the Royal Society's Copley Medal twice.
[486] [Herschel, 1815].
[487] [Brinkley, 1807]. John Brinkley, 1766 - 1855, was senior Wrangler in 1788, and Astronomer Royal from 1792. He became Bishop of Cloyne in 1826.
[488] [Herschel, 1818].

He first described his notation, which is different from both De Morgan's and Warburton's:

$$^{s}\prod(x)$$

Here $x$ is the number to be partitioned and $s$ is the number of parts. Later in the paper (paragraph 41) Herschel referred to De Morgan's paper[489] and compared his notation to that of De Morgan, noting that $u_{x,y} \equiv {}^{y}\prod(x)$. This is only true if Euler's Theorem is brought into play, since Herschel's $s$ is the number of parts and De Morgan's $y$ is the largest part.

Herschel's first observation was that ${}^{1}\prod(x) = 1$. The development then proceeded by considering all the possible bi-partitions of $x$ (that is, the case when $s = 2$):

$$1, x-1; \quad 2, x-2; \quad 3, x-3; \ldots$$

The enumeration of these was given as $\frac{x}{\overline{2}}$. This notation, using the double separator between numerator and denominator, means the integer part of the quotient. Paragraph 12 of the paper develops this idea and uses the notation for a circulating function (from paragraph 6) $s_x$ to denote the sum of the $x$th powers of the $s$th roots of unity. In the case in question here we have:

$${}^{2}\prod(x) = \frac{x}{\overline{2}} = \frac{x}{2} - \tfrac{1}{2}2_{x-1}$$

Now $2_{x-1}$ is 0 if $x$ is even and 1 if $x$ is odd. This produces the same result as De Morgan's formula for $u_{x,2}$. However, the derivation of this result is more complex, and proceeds from the general solution to the particular cases, rather than the other way round.

Herschel observed that the general formula for $\prod^{s}(x)$ (note the slight change in the form of the notation) must consist of two portions: a non-periodical function $f(x)$ of $x$, and a periodical or circulating function $Q_x$.

[489] [De Morgan, 1843].

From his general recurrence relation $\prod^{s}(x) = \prod^{s-1}(x-1) + \prod^{s-1}(x-s-1) + \prod^{s-1}(x-2s-1) + \cdots$

he derived the relation $\prod^{s}(x) = \phi(x-1) + \phi(x-s-1) + \cdots + Q_{x-1} + Q_{x-s-1} + \cdots$.

Four pages of analysis of the circulating functions followed using the machinery developed in the first part of the paper, culminating in the result ${}^{2}\prod(x) = \frac{1}{2}(x - 2_{x-1})$. From this it is then relatively straightforward to derive the non-periodic part of ${}^{3}\prod(x)$ which is $\frac{1}{12}(x^2 + x)$, and the circulating parts which are:

$$Y = -\tfrac{1}{12}\left\{2.6_{x-1} + 6.6_{x-2} + 2.6_{x-4} + 6.6_{x-5}\right\}$$
$$Z = -\tfrac{x}{12} + \tfrac{1}{12}\left\{6_{x-1} + 2.6_{x-2} + 3.6_{x-3} - 2.6_{x-4} + 5.6_{x-5}\right\}$$

Putting these parts together gives: ${}^{3}\prod(x) = \frac{1}{12}\left\{x^2 - 6_{x-1} - 4.6_{x-2} + 3.6_{x-3} - 4.6_{x-4} - 6_{x-5}\right\}$ which agrees with De Morgan's result for $u_{x,3}$.

Herschel continued by calculating results for ${}^{s}\prod(x)$ for $s = 4$ and $s = 5$ explicitly, but he did not show the calculation for $s = 6$, claiming simply that it would not be materially more difficult than the earlier cases. He concluded by noting that the effect of the circulating part was "to adjust the final value to *the nearest integer* of the rational fraction."

***Kirkman***

Thomas Kirkman presented a paper on 10 January 1854[490] entitled *On the k-partitions of N*, in which he reproduced the results of Herschel (and De Morgan and Warburton) but using what he described as "more elementary methods." Herschel's work was described as "somewhat formidable analysis". His only knowledge of the work of De Morgan and Warburton came from the description in Herschel's paper, for he lamented:

> "I have not seen, nor have any means of seeing, so far as I know, within 150 miles of Manchester, what Prof. De Morgan has written on this subject[491] ."

[490] [Kirkman, 1855].

Kirkman used a different method from that employed by the earlier authors to calculate the algebraic and periodic parts of his formulas, although his starting point for the analysis was similar to De Morgan's. He began by observing that any number $N$ can be written as the sum $\left(a+x_1\right)+\left(a+x_1+x_2\right)+\cdots+\left(a+x_1+\cdots x_k\right)$, where there are $k$ numbers, none less than $a$, and the $x_i$ "may be anything, positive or zero". By taking $a = 1$ and using the notation ${}_kP_x$ to denote the number of $k$ partitions of $x$, a process similar to De Morgan's leads to the recurrence relation ${}_kP_x - {}_kP_{x-k} = {}_{k-1}P_{x-1}$. He claimed it is then "easily seen" that ${}_2P_x = \frac{1}{2}\left(x - {}^*2_{x-1}\right)$.

The notation ${}^*2_{x-1}$ is "that well known function of the square roots of unity which is = 1 for $x$ odd, and = 0 for $x$ even." These functions are called *circulators*, and are defined thus:

$$ {}^*s_e = \begin{cases} 1, \text{ if } \dfrac{e}{s} \text{ is a positive integer} \\ 0, \text{ otherwise} \end{cases} $$

Note that Kirkman wrote the quotient $e/s$ as *e:s*. The notation ${}^*s_e$ is similar to that used by Herschel.

From the result for ${}_2P_x$ the result for ${}_3P_x$ is derived in two stages. First the non-periodic (algebraic) part is calculated by noting that the recurrence relation ${}_3P_x - {}_3P_{x-3} = {}_2P_{x-1}$ must hold purely for the algebraic parts. The subtraction on the left eliminates at most one power of $x$, and since the right hand side

[491] [Kirkman, 1855, p. 145].

contains the first power, the left hand parts can be at most quadratic expressions. So Kirkman could write

$$\left({}_3P_x\right) = Ax^2 + Bx + C \text{ and } \left({}_3P_{x-3}\right) = A\left(x-3\right)^2 + B\left(x-3\right) + C,$$

where the brackets around ${}_3P_x$ indicate that only the algebraic part is being considered. By summing these two expressions and setting them equal to the known result for ${}_2P_{x-1}$, Kirkman was able to calculate values for $A$, $B$ and $C$. This gave him the result $\left({}_3P_x\right) = \frac{1}{12}\left(x^2 + x - 2\cdot{}^*3_{x-1}\right)$. The calculation of the circulating part then follows by summing the $e + 1$ terms in $-\frac{1}{2}\left({}^*2_{x-2} + {}^*2_{x-5} + \cdots + {}^*2_{c+1} + {}^*2_{c-2}\right)$.

Kirkman did this by separately considering the cases when $e$ is odd or $x$ is 0, 1 or 2 (mod 3) and $e$ is even, or $x$ is 0, 1, 2, 3, 4 or 5 (mod 6) and he arrived at the same result as Herschel. He then proceeded to calculate explicit expressions for $k$ = 4, 5 and 6.

In a subsequent paper[492], *On the 7-partitions of X,* he exhibited an explicit formula for the case $k$ = 7. In this paper he referred to Euler's work on the subject, and stated that, "This view of the matter, however, led to no practical results, until, by the profound researches of Mr Cayley, ... , and the more decisive discoveries of Professor Sylvester, who has, if I mistake not, thoroughly worked out the idea of Euler, and given the results of it in a calculable form, the problem of partitions has been directly solved, without the aid from any induction."

***Cayley and Sylvester***

The work of Cayley and Sylvester referred to by Kirkman consisted of two papers by Cayley, *Researches on the Partition of Numbers*[493] and *Apropos of Partitions*[494], and a paper by Sylvester, *On a discovery in the partition of numbers*[495].

492 [Kirkman, 1857].
493 [Cayley, 1857].
494 [Cayley, 1857].
495 [Sylvester, 1857].

Cayley started with Euler's result:

$$P\left(a,b,c,\ldots\right)q = \text{the coefficient of } x^q \text{ in } \frac{1}{\left(1-x^a\right)\left(1-x^b\right)\left(1-x^c\right)\cdots},$$

where $P(a,b,c,...)q$ is the number of ways that $q$ can be made up of the elements $a, b, c, \ldots$ He then proceeded to develop a theory based on algebraical fractions and partial fractions, rather than difference equations. The results are similar to those already discussed, in so far as they have a non-periodic part and a periodic part.

Cayley introduced the notion of a *prime circulator*, pcr $a_q$, defined by:

$$pcr\ a_q = \begin{cases} 1 \text{ if } q = 0 (\operatorname{mod} a) \\ 0 \text{ otherwise} \end{cases}$$

to deal with the circulating parts of the formulas. He gave the results of several calculations; for example, the number of ways of making $q$ out of 1s and 2s is:

$$P\left(1,2\right)q = \tfrac{1}{4}\left\{2q+3+\left(1,-1\right)pcr2_q\right\}$$

where

$$\left(1,-1\right)pcr2_q = 1\cdot 2_q + -1\cdot 2_{q-1}$$

So for $q = 5$, the calculation is:

$$P\left(1,2\right)5 = \tfrac{1}{4}\left\{2\cdot 5+3+\left(1,-1\right)pcr2_5\right\} = \tfrac{1}{4}\left\{13+1\cdot 2_5 + -1\cdot 2_4\right\} = 3$$

The three partitions are $1^5$, $1^3\,2$ and $1\ 2^2$. Cayley's formula differed from those of De Morgan, Warburton, Herschel and Kirkman in that Cayley specified the actual parts to be included in the partition, whereas the others specified either the least part (Warburton and Kirkman), the greatest part (De Morgan) or the number of parts (Warburton, Herschel and Kirkman). So, for example, $P(1,2)5 = u_{5,1} + u_{5,2}$ in De Morgan's notation, because in each partition De Morgan included the greatest part.

Sylvester's approach, typically, introduced a lot of new vocabulary. It was similar to Cayley's, in that he specified the elements to be included in the partition of *n*, and defined the number of ways of composing *n* from these elements to be the *quotity* of *n*. The formula for the quotity comprises two parts: a non-periodic part he called the *quot-additant*, and a periodic or circulating part he called the *quot-undulant*. Again, the analysis was not based on difference equations but on expanding rational functions and involved Bernoulli numbers, logarithms and integrals.

In 1859, Sylvester gave a course of seven lectures[496] at King's College, London, on partitions. The records of these lectures consist only of his notes, which are very detailed in places, and merely sketches in others. In the first lecture, when introducing the subject, he made reference to the work of De Morgan, Warburton, Herschel and Kirkman, but did not detail what it was he actually intended to say. He defined the general problem of simple partitions as finding the number of ways of composing *n* from elements *a, b, c, ..., k*. The problem was then reformulated as a problem of solving systems of equations and he developed the solution using matrices, geometry, groups, symmetric functions and other deep ideas which are beyond the scope of this thesis.

The particularly English contribution to the topic was the use of difference equations and recurrence relations by De Morgan, Warburton, Herschel and Kirkman, as an almost self-contained theory which did not require the more detailed and difficult mathematical concepts used by Cayley, Sylvester, and those who followed.

***Table to compare notations***

To summarise, the different notations described in this appendix are collected together below, along with examples of the formulas developed.

---

[496] [Sylvester, 1897]

| **Date** | **Name** | **Notation** | **Example Formula** |
|---|---|---|---|
| 1818 | Hamilton: | $\binom{n}{m}$ | none |
| 1843 | De Morgan | $u_{x,y}$ | $u_{x,2} = \frac{x}{2} - \frac{1}{4} + \frac{1}{4}(-1)^x$ |
| 1847 | Warburton | $[N, p_\eta]$ | $[N, p_\eta] = \sum_{z=0}^{y-\eta}\left[N - (1 + p(\eta - 1)) - zp, (p-1)_1\right]$ |
| 1850 | Herschel | ${}^s\prod(x)$ | ${}^2\prod(x) = \frac{1}{2}(x - 2_{x-1})$ |
| 1857 | Cayley | $P(a,b,c,...)q$ | $P(1,2)q = \frac{1}{4}\{2q + 3 + (1,-1)pcr2_q\}$ |

### List of MacMahon's partition theory papers

The list below of all of MacMahon's papers relating to partition theory in chronological order illustrates the way in which the work was spread over his entire career. Each line is in the order: *Date, Andrews Reference, Title.*

1886 20;6 Certain special partitions of numbers
1888 27;18 The expression of syzygies among perpetuants by means of partitions
1889 33;10 On play "à outrance'
1890 37;19 Weighing by a series of weights
1891 38;6 The theory of perfect partitions and the compositions of multipartite numbers
1894 46;5 Memoir on the theory of the compositions of numbers
1897 51;9 Memoir on the theory of the partition of numbers Part I
1899 54;10 Memoir on the theory of the partition of numbers Part II
1899 56;12 Partitions of numbers whose graphs possess symmetry
1900 58;10 Partition analysis and any system of consecutive integers
1904 62;10 The Diophantine inequality $\lambda x \geq \mu y$
1905 65;10 Note on the Diophantine inequality $\lambda x \geq \mu y$
1905 66;10 On a deficient multinomial expansion
1906 67;10 Memoir on the theory of the partition of numbers Part III
1907 68;10 The Diophantine equation $x^n - Ny^n = z$
1908 71;5 Second memoir on the composition of numbers
1909 73;10 Memoir on the theory of the partition of numbers Part IV
1912 77;11 Memoir on the theory of the partition of numbers Part V
1912 78;12 Memoir on the theory of the partition of numbers Part VI
1912 79;8 On compound denumeration
1917 87;8 Memoir on the theory of the partition of numbers Part VII
1920 93;16 Divisors of numbers and their continuations in the theory of partitions
1920 94;9 On partitions into unequal and into uneven parts
1921 97;9 Note on the parity of the number which enumerates the partitions of a number
1923 105;11 The connexion between the sum of the squares of the divisors and the number of partitions of a given number
1923 107;6 The partitions of infinity with some arithmetic and algebraic consequences
1923 109;6 The prime numbers of measurement on a scale
1923 110;9 The theory of modular partitions
1924 111;8 Dirichlet series and the theory of partitions
1925 114;8 The enumeration of the partitions of multipartite numbers
1926 117;8 Euler's phi function and its connection with multipartite numbers
1926 118;9 The parity of $p(n)$, the number of partitions of $n$, when $n$ is less than 1000
1927 119;9 The elliptic products of Jacobi and the theory of linear congruences

**Appendix 6**
**The correspondence between MacMahon and D'Arcy Wentworth Thompson.**

On 14 September 1911 MacMahon was awarded an LLD *honoris causa* by St Andrews University, Fife. It is likely that he was recommended for this honour by D'Arcy Wentworth-Thompson[497] ,Professor of Zoology at the University. Thompson is best known as the author of *On Growth and Form*[498], first published in 1917, republished in 1942 and then in abridged form regularly since 1959.

The University has correspondence between MacMahon and Thompson comprising the following nine letters:

1. Thompson to MacMahon, Dundee, 2 July 1917
2. MacMahon to Thompson, London, 8 July 1917
3. MacMahon to Thompson, London, 15 January 1922
4. MacMahon to Thompson, Cambridge, 22 August 1923
5. MacMahon to Thompson, Antwerp, 23 August 1923
6. Thompson to MacMahon, St Fillans, 27 August 1923
7. MacMahon to Thompson, Cambridge, 7 December 1923
8. Thompson to MacMahon, St Andrews, 23 May 1924
9. MacMahon to Thompson, Cambridge, 20 January 1926

A description and the text of each letter follows, with notes appended as footnotes.

**Letter 1, Thompson to MacMahon**
*(Double spaced typewritten letter over two pages, with some pencil diagrams attached.)*

Dundee. 2/7/17

My dear MacMahon,

Coming home from Edinburgh on Saturday, I found your new Phil. Trans. paper[499] on my table; and, for a moment, I though you had paid me the delicate, though exaggerated, compliment of sending me a copy. I soon discovered, of course, that I had flattered myself unduly. However, it reminded me that I had never sent you a copy of the little paper of which we talked when it was being written, and this, belated as it is, I send you now. It is all printed over again, slightly expanded, as a chapter of my recent book[500], - which (by the way) I hope you will have a look at. There was a copy on the table in the Athenaeum when I was last there.

Also you[sic] own work is sky high above me, and I can only see very dimly, here and there, what you are driving at. But in a kind of way it interests me all the same.

I can't help thinking (and if I talk a little nonsense you must forgive me) that the hole thing might be tuned into, and would look much simpler to my eyes, as a spatial problem.

[497] Sir D'Arcy Wentworth-Thompson, 1860 -1948, was knighted in 1937.
[498] [Thompson, 1997].
[499] This paper was the *Seventh memoir on the partition of numbers: a detailed study of the enumeration of multi-partite numbers*, published in 1917 [MacMahon, 1917, [87;8]].
[500] This is the original 1917 edition of *On Growth and Form*, although it is not clear which chapter is being referred to. The reply from MacMahon (letter 2) suggests that it is Chapter 9, 'On the Theory of Transformations, or the Comparison of Related Forms'.

Take for instance, the Partition of Multipartite Numbers (p.83).
I was lying in bed this morning thinking of this, and the elementary cases seemed to me ever so much simpler in a diagrammatic form, like this:
Here are your three 'factors Δ[501] ; and here, obviously are your three ways of combining them, two and one,
and the fourth way, i.e. all three together.[502]
Here again are four factors. And it is equally obvious that you have six ways of (..).., three ways of (..)(..) , four ways of (...). and one way (....)[503]
e.g. 1,2,3,4,5 1,3,2,4,5 etc., might it not be amusing to treat this a sort of problem in 'beknottedness', imagining (as it were) one's numbers to be written on a long ribbon, to be wound or kinked in various ways. You would get, as it seems to me, a series of symmetrical patterns or knots, whose classification would be that of your numerical permutations.
Alas, I am now to old to learn; but I wish to goodness I had even a ghost of a knowledge of mathematics.

Yours ever faithfully, D'Arcy W Thompson

*There are three pages of pencil sketches apparently attached to this letter, comprising hexagons divided into six triangles with different arrangements of the numbers 1, 2, 4 and 5 in them. Their relevance is to the context of this letter is not clear, and it may be that they do not belong to it but have become misplaced during storage over the years.*

**Letter 2: MacMahon to Thompson**
*(Three page handwritten letter with three sides of sketches attached. The letter is on Board of Trade paper, so presumably MacMahon was writing whilst at work).*

8/7/17
My dear Thompson,
In my work on Combinatory Analysis two vols. lately published by the Cambridge University Press[504] you will see how largely graphical methods are employed in every part of the subject.
They may of course be simply used so as to be visual representations but they are mostly of service for obtaining new theorems and correspondences by graphical transformation.
I will look out a selection of my original papers which have special reference to these methods and will transmit them to you in a day or two.
My extension of line-combinations etc. to two-dimensional spaces and some work on those in three-dimensions arose entirely from my having the space representation constantly in mind.
Sylvester pointed the way in the first instance.
I am much obliged for your morphology and mathematics which is evidently the result[?] of deep thinking. I expect that if you took rubber sheets and drew your fishes etc.[505] the straining of the sheets in various ways would give you the transformations but I have no doubt you that you have thought of this.
I have not yet read your paper very carefully but I intend to do so to see if there are any analogies which one is constantly on the look out for.

[501] The triangle is hand-drawn in pencil, with clear dots at each vertex.
[502] There is a handwritten amendment dated 5.1.39 at this point: [and a fifth, all separate]. This note was written in pencil and subsequently inked over. Since it was written 10 years after MacMahon's death, it suggests that Thompson kept a copy and referred to it prior to the reprint of his book in 1942.
[503] There is another later addition here, also dated 5.1.39: " and one way (.)(.)(.)(.)."
[504] [MacMahon, 1984], originally published in 1915/6.
[505] This is the clue that the chapter of *On Growth and Form* which Thompson had sent to MacMahon was Chapter 9.

With kind regards, yrs sincerely, Percy A MacMahon

*The three pages of sketches do not appear to have any bearing on the letter and are in Thompson's rather than MacMahon's hand. They seem to be attempts to create all the combinations of two or three things taken three at a time with repetitions allowed. It is possible that they belong to Letter 1 above.*

**Letter 3: MacMahon to Thompson.**
*(A single sheet, handwritten, covering one and a half sides of the paper, written from MacMahon's home address at 27 Evelyn Mansions, Carlisle Place, SW1.*)

15 Jan 22
Dear Thompson,
It is very kind of you to send me a model of your interesting polyhedron[506] . I am at present researching in the subject of repeating configurations in space and have [illegible word] the values of various polyhedral angles in spherical measure. There are restrictions upon the angles possible in such forms, but the subject is difficult and will proceed I think best by analogy with two dimensions which has been, quite recently, fully worked out.
With kind regards and hoping that you will continue to think about these subjects as suggested by natural objects.
sincerely yrs, Percy A MacMahon

**Letter 4: MacMahon to Thompson**
*(A single side, written from MacMahon's address in Cambridge)*

22 August 1923
Dear Thompson,

Thanks for yours.
Your 1 nz[?]: DWT$^1$ is a circular percession of PAM$^1$ which I regard as not differing in anyway to your 2 nz. DWT$^2$ is a circular percession of my PAM$^3$ and practically identical with it.
Neither of them are in relation to PAM$^2$.
What I mean is that I take

```
   4          5
3     5    4     6

2     6    3     1
   1          2
```

to be identical hexagons. Your letter gave me rather a flutter because I have in vain sought for a fourth contact system. My address for ten days is Hotel des Flandres Antwerp and you will hear from me thence.

[506] There is no clue as to the nature of this model.

Yrs, PAM

**Letter 5: MacMahon to Thompson**

*(A single handwritten side of hotel notepaper, headed Hotel Restaurant, Cafe Suisse, Place Verte 2-3 Anvers.)*

Anvers, le 23 August

Dear Thompson,

I hope that you received the card about the hexagons that I sent you before leaving Cambridge. I am very grateful to you for the attention that your giving to the little book and I gladly take up the faults that you may point out that I may avoid them in a new edition which may shortly be called for. I look upon part III as a really important part of the book. It is strange that this space filling geometry has been ignored by geometricians. It has since led on to quite a revolution in Decorative Geometry which I am accentuating[?] in a new work. Geometricians in Cambridge admit but cannot explain the subject. That reminds me that I have recently had your book on Natural History Geometry in my hands and sincerely congratulate you upon it. It may interest you to hear that G. T. Bennett[507] of Emmanuel, the leading geometrician in this country has gone carefully thro [sic] your book and has made many observations and notes while reading it. He had found some slips and inaccuracies which I am sure he will put at your disposal for a new edition. They are few in number but their correction would increase the value of your book. You might look him up next time you are in Cambridge and have an interesting hour.

G.T.B:W.D.T::W.D.T:P.A.M and you are square as the result! I should add that Part 3 of Pastimes looks fair to have an influence on the general theory of periodic functions. This was pointed out by H. F. Baker and I have since realised that it must be so - seeing what wonderful results Klein obtained from an assemblage of parallelograms.

yrs Sincerely P A MacMahon

**Letter 6: Thompson to MacMahon**

*(A three page handwritten letter from St Fillans, with a single diagram attached, which is a correction to page 15 of 'New Mathematical Pastimes'[508].)*

St Fillans 27 VIII 1923

My dear Macmahon[sic],

Very many thanks for your Antwerp letter. I shall write at once to G.T.B. I am a little surprised that he should have read my work - but not at all surprised that he should have found

[507] Geoffrey Thomas Bennett ,1868 - 1943, senior (male)Wrangler 1890 (Phillipa Fawcett, 1868 - 1948, was placed first), FRS 1914. He was an interesting character with a fondness for billiards and cycling as well as geometry. He used a racing cycle around Cambridge and rode in a prone position with a mirror fitted to the cycle so that he didn't have to raise his head to see where he was going (the despatch note from the maker in London is preserved in the Emmanuel College archives). He kept detailed records of the lengths of his many rides. He also kept records of his walks around the area, writing down the readings from his pedometer on little scraps of paper. Bennett's *nachlass* contains no reference to MacMahon or to D'Arcy Thompson (although it does contain catalogues of billiard tables and samples of the cloth used to cover them and many diagrams of the paths of balls). However, he was a member of the Council of the London Mathematical Society from 1908 to 1917, and Honorary Secretary to the Society from 1915 to 1917. He attended many meetings of the council with MacMahon. MacMahon consulted him about his work on repeating patterns.

[508] [MacMahon, 1921].

mistakes in it. I have a pretty long list of them already, detected some by myself and some by my friends. One of the worst - a real howler - I detected in the proof and then forgot to correct. The edition is all but sold out, but the Press want to reprint it from the plates. That is to say, they want merely to correct verbal mistakes, but not to expand or alter materially. They say they can't afford to re-set the book, and so they have the whip hand of me.
As to your own book, the only thing I have noted is the omission of B5,0,0,4 and B1,0,0,8[509] in the table on p. 15. The former I have not yet tried; the latter comes out easily.
On p.42 §34 it might be worth pointing out (obvious as it is) how enormously the number of cases increases when we come to 12- and 20-hedra. In other words, you might expand, in a single sentence, qv[?] 'and so on for any regular solid.'
The method of notating one figure or another, on a help to inspection, of which spoke in my last letter, reminds me of an old idea of mine in another field. Here is the old Eulerian Diagram,

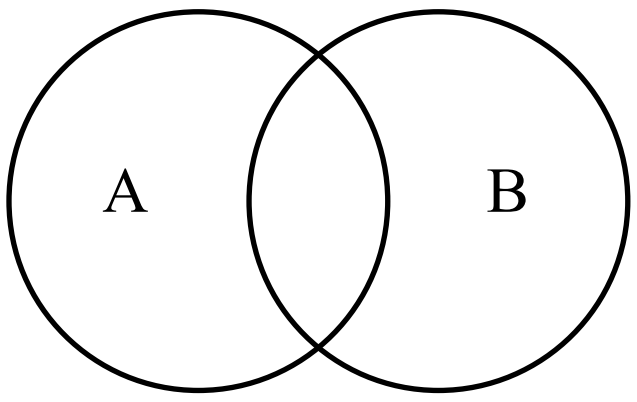

of elementary Logic, Some A's and B's, i.e. certain A's are B's. But let $A \cap B$ revolve and you get a very different meaning of "some", - e.g. Some of the crew are always on deck.
But that is off the point, meanwhile,

Yrs faithfully, D'Arcy W Thompson

I was interested to come across Col. Jocelyn's name in yr. book. I met him once years ago at a dinner of the Savages - where I rather think I was your guest. I met him again in a railway carriage the next day; and have never seen him since! - to my great regret.

**Letter 7: MacMahon to Thompson.**
*(A handwritten letter over two sides of paper, written from MacMahon's Cambridge address.)*

7 Dec. 1923
Dear Thompson,
I had your letter several weeks ago and went round to see Bennett and found he had also heard [from Thompson]. We discussed the matter and he shewed to me a German book with diagrams which seemed to be very like what you had been doing and he (Bennett) said he would write to you on behalf of me as well as himself. However he postponed doing this and I fear may not have even by this time have replied to you. He will no doubt do so almost at once as I spoke again about a day or two ago. He has all the information that will be useful to you - I unfortunately do not.

[509] This is the subject of the diagram attached to this letter.

I congratulate myself too in bringing you and Bennett together and I know you will be great friends.

Yrs sincerely, P A MacMahon"

**Letter 8: Thompson to MacMahon**

*(A three page typewritten letter from St Andrews, with a number of hand drawn diagrams on the second and third pages.)*

23rd May 1924
"My dear Macmahon[sic],
Here is a little problem that might perhaps interest you; and if, or when, you have time to think of it you might help me in doing so.
The essence of it is, In how many different patterns can you arrange *n* Hexagons? The practical application and interest of it is, to find in how many different ways the the partition-walls of *n* soap bubbles, or of *n* cells in a developing embryo can be arranged in a plane configuration.
As you know, in the soap bubble or the living cell, the partition-walls meet (in surface view) three in a point, at angles of $120^0$.
When there are three of them, there is obviously only one possible arrangement; there is likewise one only when you have four, or when you have five. In the case of six, you have three arrangements, when you have seven, there are five,

Now with eight cells, Plateau says there are thirteen possible arrangements, and so I said in my book; but now, on going into the matter again, I am pretty certain there are only twelve.
But I cannot find any rule or formula for the increasing number; in other words I don't know how many possibilities exist in the case of nine cells, or ten.
The problem is the more interesting to me because one of my students has been examining large numbers of developing frogs' eggs; and he has found, in the eight-celled stage, every one of the twelve possible configurations.

I go to work by the following method, - perhaps a clumsy one. In the first place, it is plain that one need only study the variations in pattern of these partition-walls which do not extend to the boundary of the system: i.e. such walls as are thickened in the following figure.

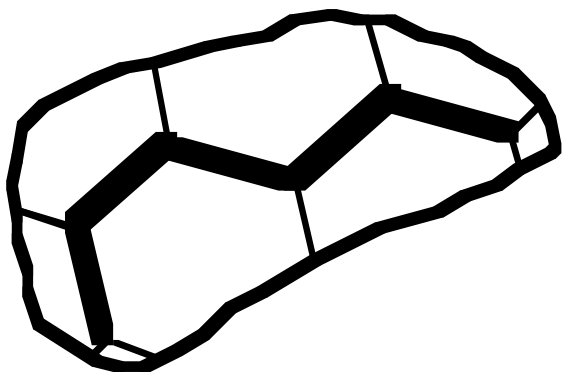

Given such a simplified diagram, one knows that one has to add on a line at each angle, and a fork at each free end.
Next, in such a diagram we may call an angle (or flexure) a when in one (positive) direction, ^, and b when in the negative. V

Then we obviously have, to begin with, the following six arrangements:

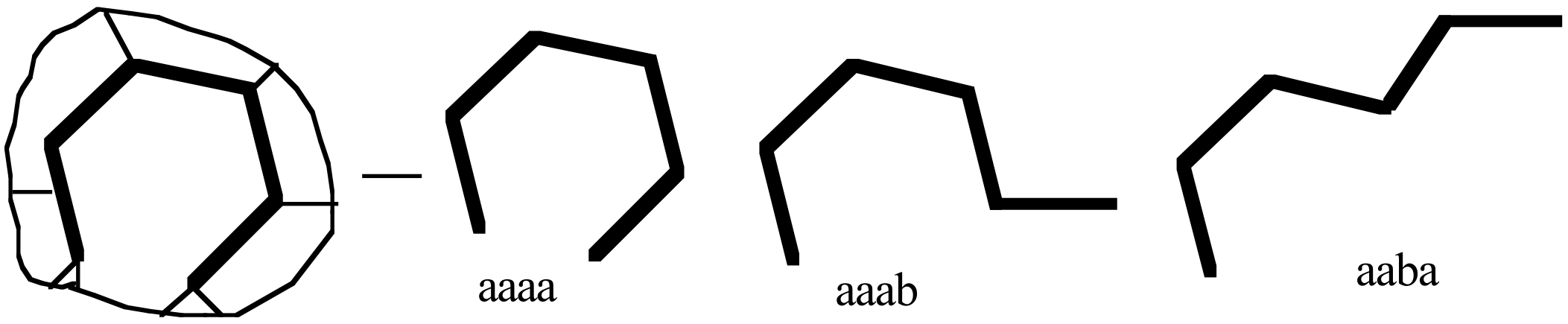

etc aaaa aaab aaba

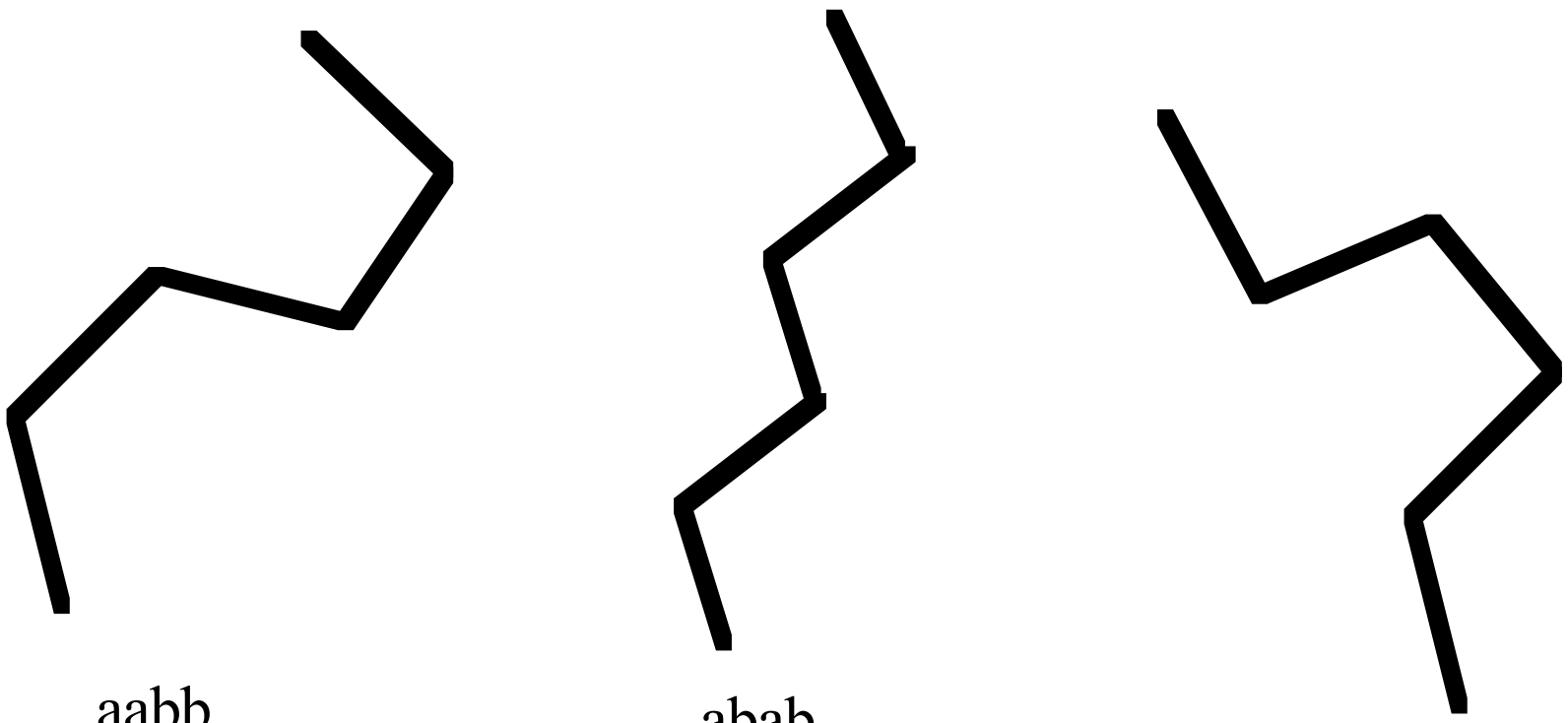

aabb abab abba

Next, instead of these four successive angles, we may take three of them and then apply to one of the the three another partition-wall, whose free end will give rise to a fork (like the two free ends already existing). Instead of the above six cases there will now be only three, viz.

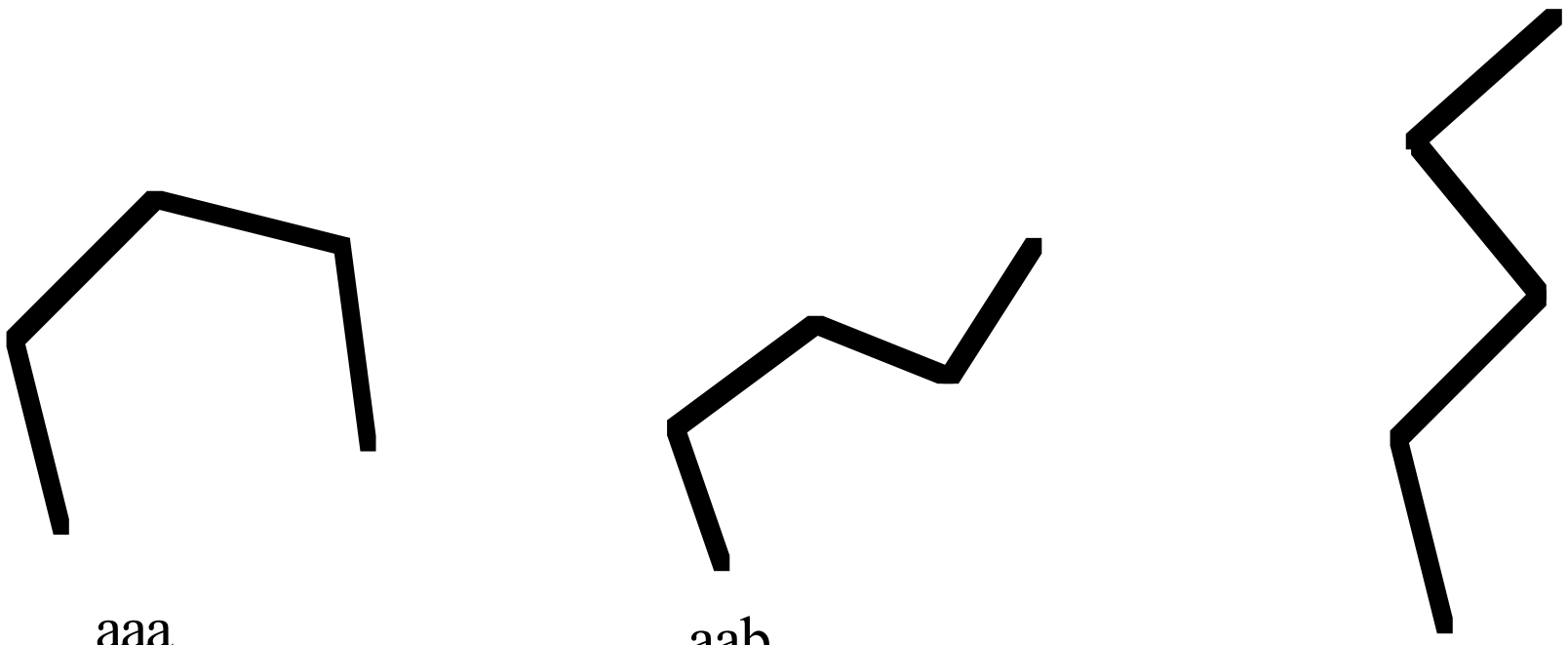

aaa aab aba

But when we add on the new partition (which we may call c), we may in each case add it in one or other of two places, and get what we may designate as

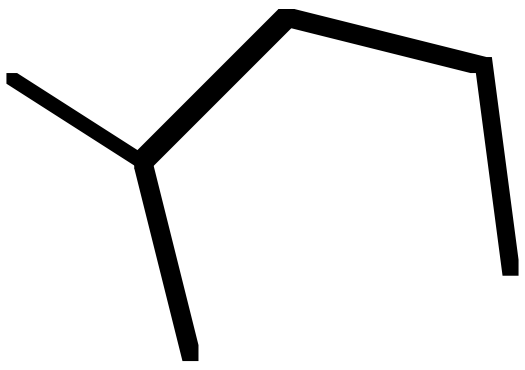

$\overset{c}{aaa}$
*

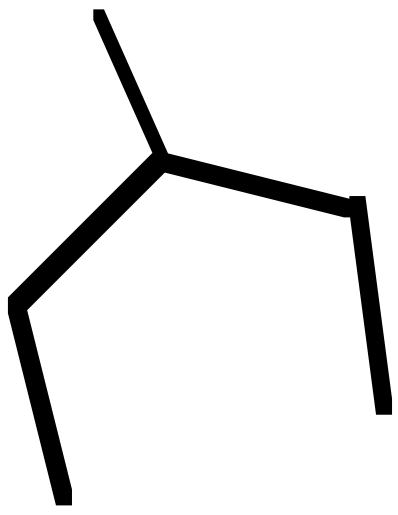

$\overset{c}{aaa}$

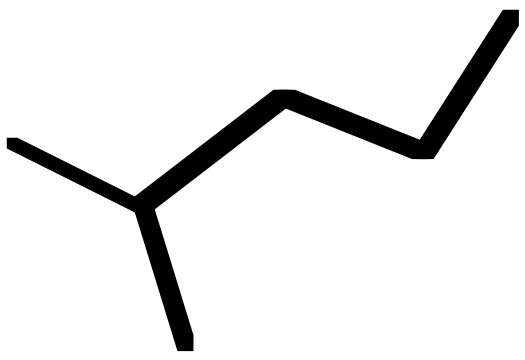

$\overset{c}{aab}$
╫

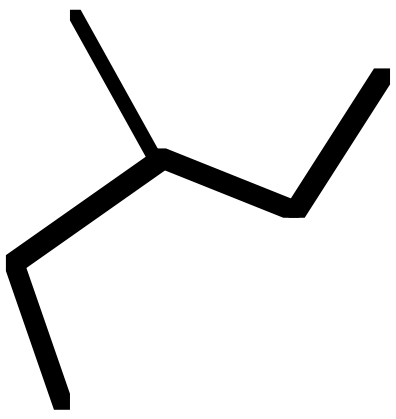

$\overset{c}{aab}$

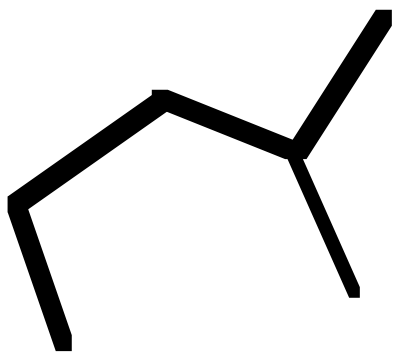

$\overset{c}{aab}$
*

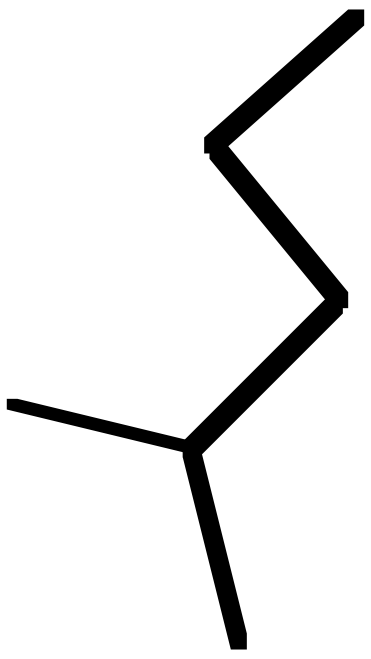

$\overset{c}{aba}$
╫

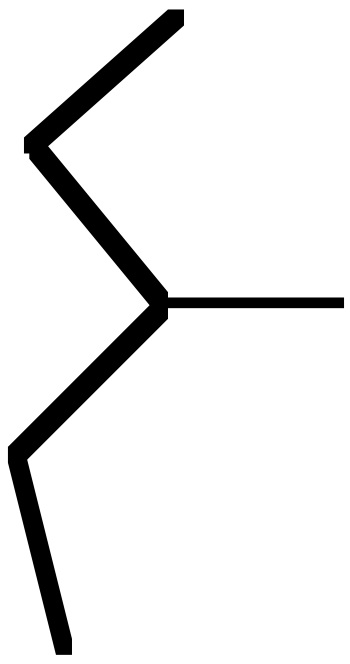

$\overset{c}{aba}$

Here, however, we see at once that $\overset{c}{aaa}$ is identical with $\overset{c}{aab}$ ; and $\overset{c}{aab}$ identical with $\overset{c}{aba}$ . You have therefore five variants in this second series.

Lastly, you may take only two successive flexures, viz. 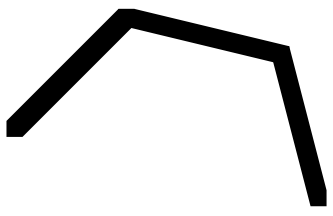 aa or ab 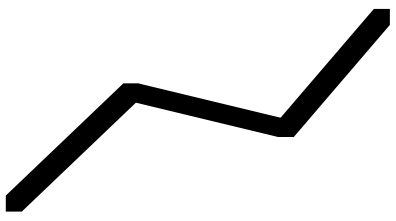 . You must now add on two c's ; but, on doing so you find that the result is the same, - $\overset{cc}{aa}$ is identical with $\overset{cc}{ab}$

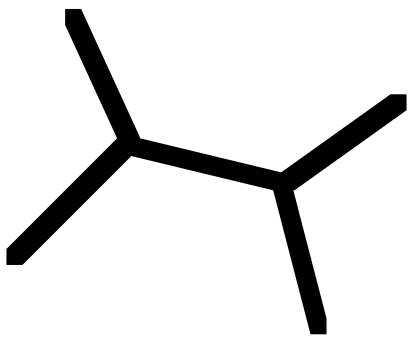 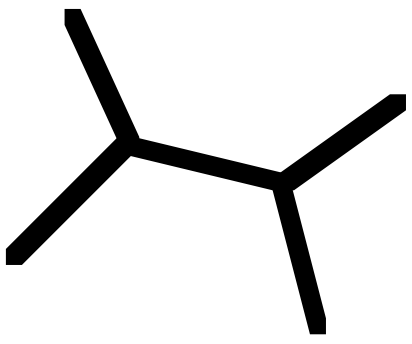

.

These then are the twelve, and only twelve, possibilities which I can discover.

Plateau gives thirteen, on the authority of Van Rees. But I have looked up Van Rees in the R. Soc. Catalogue, and cannot find that he published anything on the subject. Plateau must have got it from him in conversation or by letter.
Don't bother about this unless, perchance, it happens to amuse you.
I am going to Toronto, - though I don't know where the money is coming from, as I have just had a long and very unwonted holiday in Portugal. I hope you are coming too.
Yours faithfully, D'Arcy W Thompson

**Letter 9: MacMahon to Thompson**
*(A single handwritten page from MacMahon's Cambridge address.)*

20 Jan 26
My dear Thompson,
It will truly give me pleasure to meet Somerville when he comes to Cambridge. I have long admired his work and in particular his mastery of questions in solid geometry in which I am interested. I hope that he will acquaint in advance of the day or days that he can be here and of the address that will reach him.
yrs sincerely, Percy A MacMahon"

## Appendix 7
## Correspondence with James Parker Smith

This Appendix contains the survivng correspondence between James Parker Smith MP and MacMahon concerning evidence to be put before a Royal Commission on electoral reform in 1909. The letters are transcribed as accurately as possible. There is also a document written by MacMahon with some detailed calculations, the purpose of which is described on page 100.

### Letter 1 James Parker Smith to MacMahon

11th Jan 09
26 Murrayfield Road
Edinburgh

My dear MacMahon,

All good wishes to you for a prosperous new year.
You will have seen that Govt have appointed a Royal Commission to go into the question of Proportional Repn of which Ld Rich. Cavendish is chairman.
I wrote to him saying that there were various problems connected with the question which needed to be investigated by a mathematician of high rank, and have a sympathetic reply[510] asking me to see him when I am next in London - which will be for the 26th.
Are you inclined to take up the subject ? I think the chance problems very interesting and am sure no one would do them better justice than yourself.
Of course the problem is to find some system of selection by which the composition of the House of Commons corresponds as nearly as may be with the strength of the parties in the country.
Everyone knows that this is not the case now and that the tremendous changes at successive elections do not really correspond to similar changes in the country.
The general idea however is that if you had the country divided up into single member constituencies of equal size throughout - you would get fair representation of everyone.
In reality that system would give most exaggerated power to the majority and gives 'the balancing elector' i.e. the small proportion of electors whose views change from time to time or who have special fads, quite disproportionate power.
This is true whether constituencies are of uniform or variable size - the error tends to be cured if the size is very small.
You can state the problem this way.
A bin full of balls of two colours white and red in known proportions, say 53% white and 47% red. A set of boxes of uniform or varying sizes are filed with a shovel from this bin, and then the balls in each box are counted; in how many boxes will white balls be in the majority ?

[510] In a letter dated 10 January 1909, Lord Richard Cavendish said: "Many thanks for your letter re Proportional Representation. Your suggestion is I think a very good one and I should be grateful to you if, when you are in London I could have a talk with you on the subject. I must say that some of the problems connected with the subject are fairly complicated and as you say require an expert mathematician to deal with them. I shall be in London next Wednesday and Thursday. I shall be there more or less permanently after the 19th, if you happen to be going up to town I hope you will let me know."

The ordinary man who has not considered the matter will say, in the same proportion as the whole number of balls. A moments thought however will shew that the chance is that each box will contain balls of the two colours in nearly the same [proportions as the bin and therefore that white balls will be the majority in nearly all the boxes - say in 90 out of every hundred. In other words the white balls will have an enormous majority in H of C and the red party will be nearly wiped out at the polls.
That is the case with Unionists in Ireland, Wales and Scotland just now. In the last Parlt. it was the same with the Liberals in London and the South of England.
As a matter of fact the actual results are not so bad because the balls are not uniformly mixed throughout the bin; in some parts one colour in other parts the other predominates.
If you selected your hundred balls by drawing them singly at random out of the bin, you would of course get nearly the same proportion as the whole no. of balls i.e. 53:47. that was more or less the case in the old pre reform Parlt. and explains why the H of C in old days represented the country at least as accurately as it does now.
But nowadays when every one wants a vote you can't go back to any arrangement of that sort, and the only way of getting anything like fair representation is to take some system of plural numbered constituencies with cumulative or transferable votes. That is to say, let a big city e.g. Glasgow which returns 7 members vote as a unit under a system which lets a majority return 4 members and a minority 3 - or majority 5 and minority 2 according to the proportions. that can be done either by the cumulative or the transferable vote. If you divide Glasgow in 7 separate single member constituencies you will often have seven members returned by the majority and none by the minority - so that the minority though large gets no representation at all.
As you know I worked at this problem 20 years ago[511] , and sent you my results before the last election,[512] but my mathl. equipment is not sufficient to deal with the question adequately, and moreover is now quite rusty.
Do you feel disposed to go with these problems with a view to giving evidence before the Commn. - and if you do not can you suggest any other mathematician whose tastes lie in this direction, and who would be interested in working at the problem.

I am yrs sincerely, Parker Smith,

**Letter 2 MacMahon to Parker Smith**

27 Evelyn Mansions
Carlisle Place
SW

[511] Under proposals made in the early in 1880s by the Proportional Representation Society, voters would be able to use a single transferable vote. Once a candidate had reached the minimum number of votes required for election, known as the *quota*, then the remaining votes for than candidate would be allocated to the second choice of each voter. There were concerns that this system would be open to abuse, since the order in which the transferable votes were selected could affect which candidate was the second to be elected. Parker Smith had calculated that a random choice of such transferable votes would result in a gain of or loss of more than 60 votes only once in every 2000 elections. This was reported in 1884 by Sir John Lubbock, of the Proportional Representation Society, in a letter to the *Daily News*.
[512] In 1906, a large Liberal victory.

1909 Jan 13

My dear Parker Smith,

As you know the questions connected with probabilities at elections may be formulated in such a manner as to involve very long and intricate calculations. I hesitate to take up the subject if the work involved is to be additional to and not partly in substitution of the work I already do for Government in the Board of Trade.
I dare say that I can think on the matter pending your interview with Lord Richard Cavendish in the 26th inst.
I will consider particularly if I ~~think~~ can be of real service and turn over in my mind the names of mathematicians whom it may be well to interest in the matter.

Sincerely yrs
P A MacMahon

**Document 1 - 15 handwritten pages apparently by MacMahon (see p. 102 for analysis of the purpose of this document.**

*[A handwritten note attached to this document reads: "MacMahon. As you know I worked this problem out 20 years ago and sent you my results before the last election. Smith's 'cubic' rule ?".*

*The first sentence is a quote from Parker Smith's letter above. The handwriting is neither Parker Smith's nor MacMahon's].*

At an election there are two candidates $A$, $B$ who finally have $a$, $b$ votes respectively where $a > b$.

The probability that at no time during the voting A has a minority of votes is

[crossed out formula]

$$1 - \frac{b}{a+1}$$

I have just generalised the above (which is the well known 'Problème du Scrutin') so as to refer to any number of candidates $A$, $B$, $C$....[513]

There are $A$ white and black balls in a bag, $b$ white and $c$ black

I draw out $n$ of them What is the chance that ~~there are $p$ white~~ the ~~balls~~ set drawn are composed of $p$ and $q$ black (all the nos being large)

The chance of drawing a white ball to begin with is $\frac{b}{A}$ of drawing another is $\frac{b-1}{A-1}$

[513] This refers to the fourth *Memoir on the partitions of numbers* [MacMahon, 1909,[73;10]] read in 1908 and published in 1909.

$\therefore$ Chance of drawing $p$ white and $q$ black in one particular order is[514]

$$\frac{b\cdots(b-p+1).c\cdots(c-q+1)}{A\cdots(A-(p+q)+1)}$$

But the no. of possible orders is

$$\frac{(p+q)!}{p!\,q!}$$

$\therefore$ whole chance is

$$C_{pq}=\frac{b!.c!.(A-(p+q))!}{(b-p)!.(c-q)!.A!}.\frac{(p+q)!}{p!\,q!} \qquad (1)$$

--------------------------------------------------------------------------------

*[The following two lines are crossed through]*

$$\frac{C_{p+1,q-1}}{C_{pq}}=\frac{b-p}{c-q+1}.\frac{q}{p+1}$$

$$\frac{C_{p+2,q-2}}{C_{pq}}=\frac{(b-p)(b-p-1)}{(c-q+1)(c-q+2)}.\frac{q(q-1)}{(p+1)(p+2)}$$

--------------------------------------------------------------------------------

Now we know that

$$s!=s^s e^{-s}\sqrt{2\pi s}\left(1+\tfrac{1}{12}s^{-1}+\cdots\right) \qquad (2)$$ [515]

All the numbers being large. We will neglect neg. powers of $s$

$$\therefore\ C_{pq}=\frac{b^b.c^c.(A-n)^{(A-n)}.n^n}{(b-p)^{(b-p)}.(c-q)^{(c-q)}.A^A.p^p.q^q}\sqrt{\frac{b.c.(A-n).n}{2\pi(b-p).(c-q).A.p.q}}$$

(the exponentials disappearing)[516] .

Put

$$c=\lambda b$$
$$q=\lambda p \qquad \therefore A=(1+\lambda)b$$
$$n=(1+\lambda)p$$

---

[514] Brackets have been insertedwhere necessary and modern factorial notation used for clarity.

[515] A pencilled note referenced to this formula reads: “This is obtained by means of the formula $s!=\int_0^\infty x^s e^{-x}dx$ by a method of Laplace (??) for expanding expanding integrals (MS notes of lectures on the [unreadable word] of [unreadable word] by Glaisher”

[516] This has been checked by the author and it is correct.

Substituting these values in the last expression for $C_{pq}$ it will be found that the powers all disappear[517] and we are left with

$$C_{p,\lambda p} = \sqrt{\frac{(1+\lambda).b}{2\pi\lambda.p.(b-p)}} \qquad (3)$$

Call this $C_0$ and call $C_{p+r,\,\lambda p-r}$ $C_r$. From (1) we see that

$$\frac{C_1}{C_0} = \frac{b-p}{c-q+1}.\frac{q}{p+1}$$

Similarly $\frac{C_2}{C_0} = \frac{(b-p)(b-p-1)}{(c-q)+\cdots(c-q+2)}.\frac{q(q-1)}{(p+1)(p+2)}$

and $\frac{C_r}{C_0} = \frac{(b-p)!}{(b-p)!-r}.\frac{(c-q)!}{(c-q)!+r}.\frac{q!}{q!-r}.\frac{b!}{b!+r}$ (4)

~~Let~~ Suppose $r$ small compared with the other numbers, so that all the factorials are large.

Thus (2) may be used in the form

$$s! = s^{s+\frac{1}{2}}e^{-s}\sqrt{2\pi}$$

By means of this, and writing

$$b-p=d$$

$$\frac{C_r}{C_0} = \frac{d^{d+\frac{1}{2}}.\left(\lambda d^{\lambda d+\frac{1}{2}}\right)(\lambda p)^{\lambda p+\frac{1}{2}}p^{p+\frac{1}{2}}}{(d-r)^{d-r+\frac{1}{2}}(\lambda d+r)^{\lambda d+r+\frac{1}{2}}(\lambda p-r)^{\lambda p-r+\frac{1}{2}}(p+r)^{p+r+\frac{1}{2}}} \qquad (5)$$

This may be written

$$\frac{C_r}{C_0} = \left(1-\frac{r}{d}\right)^{-\left(d-r+\frac{1}{2}\right)}\left(1+\frac{r}{\lambda d}\right)^{-\left(\lambda d+r+\frac{1}{2}\right)} \times \left(1-\frac{r}{\lambda p}\right)^{-\left(\lambda p-r+\frac{1}{2}\right)}\left(1+\frac{r}{p}\right)^{-\left(p+r+\frac{1}{2}\right)} \qquad (6)$$

Consider these factors separately

Let $y = \left(1-\frac{r}{d}\right)^{-\left(d-r+\frac{1}{2}\right)}$

[517] The author has been unable to verify this. Calculations, assisted by *Mathematica*, reduce the powers to $\frac{(b-p)^p b^{b-1} p^{\lambda}}{(b-p)^b p^{p-1} b^{\lambda}}$.

$$\therefore \log y = -\left(d - r + \tfrac{1}{2}\right)\log\left(1 - \frac{r}{d}\right)$$

$$= \left(d - r + \tfrac{1}{2}\right)\left(\frac{r}{d} + \frac{1}{2}\frac{r^2}{d^2} + \cdots\right)$$

$$= r + \frac{1}{2}\frac{r^2}{d} + \frac{1}{3}\frac{r^3}{d^2} + \frac{1}{4}\frac{r^4}{d^3}$$

$$- \frac{r^2}{d} - \frac{1}{2}\frac{r^3}{d^2} - \frac{1}{3}\frac{r^4}{d^3}$$

$$+ \frac{1}{2}\frac{r}{d} + \frac{1}{4}\frac{r^2}{d^2}$$

$$= r - \frac{1}{2}\frac{r^2}{d} + \log 1 + \frac{1}{2}\frac{r}{d} - \frac{1}{6}\frac{r^3}{d^2}$$

$$\therefore y = e^{r - \frac{1}{2}\frac{r^2}{d}}\left(1 + \frac{1}{2}\frac{r}{d} - \frac{1}{6}\frac{r^3}{d^2}\right)$$

neglecting[518] $\frac{r^2}{d^2}$ and $\frac{r^4}{d^3}$

The second factor can be got from this by writing *-r* for *r* and $\lambda d$ for *d* - and the factors involving *p* in a similar way multiplying together and as before neglecting terms of the orders $\frac{r^2}{d^2}$ and $\frac{r^4}{d^3}$. We get

$$\frac{C_r}{C_0} = e^{-\frac{1}{2}r^2\left\{\frac{1}{d} + \frac{1}{\lambda d} + \frac{1}{\lambda p} + \frac{1}{p}\right\}} \bullet \left\{1 + \frac{r}{2}\left(1 - \frac{1}{\lambda}\right)\left(\frac{1}{d} - \frac{1}{p}\right) - \frac{r^3}{6}\left(1 - \frac{1}{\lambda^2}\right)\left(\frac{1}{d^2} - \frac{1}{p^2}\right)\right\}$$

$$\frac{C_r + C_{-r}}{C_0} = e^{-\frac{1}{2}r^2\left\{\frac{1}{d} + \frac{1}{\lambda d} + \frac{1}{\lambda p} + \frac{1}{p}\right\}} \bullet \left\{\ 2\ \right\}$$

or writing

$$\mu^2 = \frac{1}{2}\left(1 + \frac{1}{\lambda}\right)\left(\frac{1}{d} + \frac{1}{p}\right) \qquad (7)$$

we get

[518] A marginal note by this calculation explains "This ~~follows~~ is Laplace's method given by Todh. Hist. of prob. p. 549". Todhunter's paragraphs 993 to 994, on pages 548 to 552 of *A history of the mathematical theory of probability from the time of Pascal to that of Laplace* reprise Laplace's calculation of the probability that, given an event with probability *p*, in a specified number of trials the probability of the event occuring lies between certain limits. The result is, in $\mu = m + n$ trials, where the desired event occurs *m* times, the probability of *m* being between $\mu p + z \pm \frac{\tau\sqrt{2mn}}{\sqrt{\mu}}$, where $m = \mu p + z$ and $\tau = \frac{r\sqrt{\mu}}{\sqrt{2mn}}$, is $\frac{2}{\sqrt{\pi}}\int_0^\tau e^{-t^2}dt + \frac{\sqrt{\mu}}{\sqrt{2mn}}e^{-r^2}$

$$C_r + C_{-r} = 2C_0 e^{-\mu^2 r^2} \qquad (8)$$

or from (3) we have

$$C_0 = \frac{\mu}{\sqrt{\pi}} \qquad (9)$$

Let $S_r = C_r + C_{r-1} + \cdots C_0 + \cdots C_{-r}$ (10)

$$\therefore S_r = 2C_0 e^{-\mu^2 r^2} + 2C_0 \sum_{x=0}^{x=r-1} e^{-\mu^2 x^2} - C_0$$

But by Euler's theorem (Todh. Hist. of prob. p. 551)

$$\sum_0^{r-1} y = \int y dr + \frac{1}{2} y - \frac{1}{2} y \qquad (11)$$

that is in our case

$$2C_0 \sum_{x=0}^{x=r-1} e^{-\mu^2 x^2} = 2C_0 \int_0^r e^{-\mu^2 x^2} dx + C_0 - C_0 e^{-\mu^2 r^2}$$

$$\therefore S_r = 2C_0 \int_0^r e^{-\mu^2 x^2} dx + C_0 e^{-\mu^2 r^2}$$

or writing $\mu x = z$ and using (9)

$$\therefore S_r = \frac{2}{\sqrt{\pi}} \int_0^{\mu r} e^{-z^2} dz + \frac{\mu}{\sqrt{\pi}} e^{-\mu^2 r^2} \qquad (12)$$

The integral $\int_0^\tau e^{-z^2} dz$ is tabulated by De Morgan for hundredths up to 2.00 in his article on Probabilities in Encyclop. Metropol. and by Galloway for hundredths up to 3.00 in his article on Prob. in Encyclop. Britan.

In obtaining (12) the highest powers neglected have been $\frac{r^2}{d^2}, \frac{r^2}{p^2}, \frac{r^4}{d^3}, \frac{r^4}{p^3}$ which it will be seen in the cases that concern us are small. In the present case this theorem may be written

$$\sum_0^{r-1} e^{-\mu^2 x^2} = \int_0^r e^{-\mu^2 r^2} dr$$

Let us consider the meaning of our results.
We have

$$C_0 = \frac{1}{\sqrt{2\pi}}\sqrt{\frac{(1+\lambda).b}{\lambda p.(b-p)}}$$
$$= \frac{1}{\sqrt{2\pi}}\sqrt{\frac{A}{q.(b-p)}}$$
$$= \frac{1}{\sqrt{2\pi}}\sqrt{\frac{A}{\frac{c}{b}.p.(b-p)}}$$

Let us suppose $\lambda$ to remain constant then $C_0 \propto \frac{1}{\sqrt{A}}$ that is to say $C_0$ diminishes as $A$ increases, but at a much slower rate. and the chance that the numbers drawn be nearly exact proportionality than any assigned fraction - say $\frac{1}{\mu}$th - of A

$$= S_{\frac{A}{\mu}} = \left(2\frac{A}{\mu}+1\right)C_0\{1-\cdots terms\}$$
$$\propto AC_0$$
$$\propto \sqrt{A}$$

(when $\mu$ is small but the series in brackets does not converge when $A$ increases - which leaves the proof imperfect).

And if

$$b = c = \frac{A}{2}$$
$$p = q = \frac{1}{2}b = \frac{A}{4}$$
$$\lambda = 1$$
$$C_0 = \frac{2\sqrt{2}}{\sqrt{\pi A}}$$

We have from (4)

$$\frac{C_r}{C_0} = \frac{(p!)^4}{((p-r)!)^2((p+r)!)^2}$$

or using the formula

$$s! = s^{s+\frac{1}{2}}e^{-s}.\sqrt{2\pi}$$

$$\frac{C_r}{C_0} = \frac{p^{4p+2}}{(p-r)^{2(p-r)+1}(p+r)^{2(p+r)+1}}$$

$$\therefore \log C_r = \log C_0 + (4p+2)\log p$$
$$-(2p+1)\{\log(p+r) + \log(p-r)\}$$
$$-2r\{\log(p+r) - \log(p-r)\}$$

which is the formula from which I have obtained my numerical results. [The first column of the table which follows is unlabeled, but is $r$].

$A = 10000$
$b = 5000$
$p = 2500$
results given by taking the means of the given values of C in the 2nd column

| | C | S |
|---|---|---|
| 0 | .01596 | .01596 |
| 5 | .01564 | .1739 |
| 10 | .01476 | .3259 |
| 20 | .01165 | .5900 |
| 30 | .00779 | .7844 |
| 40 | .00444 | .9067 |
| 50 | .00215 | |
| 60 | .00100 | |
| 70 | .00032 | |

$A = 10000$
$b = 5000$
$p = 2500$
$$\mu = \frac{\sqrt{2}}{50}$$

| r | $\mu r$ | $\frac{2}{\sqrt{\pi}}\int_0^{\mu r} e^{-z^2}\,dz$ | $\frac{\mu}{\sqrt{\pi}} e^{-\mu^2 r^2}$ |
|---|---|---|---|
| 0 | 0 | 0 | .01595 |
| 5 | .1414 | .1585 | .0156 |
| 10 | .2828 | .3108 | .0147 |
| 16 | .4525 | .4778 | .0129 |
| 17 | .4808 | .5035 | .0125 |
| 20 | .5656 | .5762 | .0116 |
| 25 | .7070 | .6826 | .0096 |
| 30 | .8484 | .7698 | .0076 |
| 40 | 1.1312 | .8903 | .0044 |
| 50 | 1.4140 | .9544 | .0021 |
| 60 | 1.6968 | .9836 | .0009 |
| 70 | 1.9796 | .9948 | .0003 |
| 80 | | | |

*End of document*

## Appendix 8
## Correspondence from MacMahon listed in the Royal Society archives

| Royal Society reference | Description |
|---|---|
| 06412 | 9 October 1906: An official Board of Trade letter requesting that the supply of Series B of the Transactions be stopped |
| 07215 | (missing) |
| 07222 | (missing) |
| 08281 | 24 February 1908: A letter on Board of Trade paper in reply to someone about a request for information. |
| 10156 | (missing) |
| 11150 | 7 September 1911: An official Board of Trade letter requesting that the standard yard and pound be made available for testing. |
| 12241 | 15 February 1912: A letter on Board of Trade paper returning a photo and description of the Balta Sound tablet-stone. |
| 16110 | 20 November 1916: A personal letter complaining about the misspelling of his surname in official publications. |

## List of applications signed by MacMahon, with date of election.

1891 Elliott, Edwin Bailey
1893 Burnside, William
1894 Hill, Micaiah John Muller
1894 Love, Augustus Edward Hough
1895 Holden, Sir Henry Capel Lofft
1896 Clark, George Sydenham, Baron Sydenham of Combe
1896 Pearson, Karl
1897 Wislicenus, Johannes
1897 Mathews, George Ballard
1898 Lindley, Nathaniel, Baron Lindley
1898 Baker, Henry Frederick
1899 Tanner, Henry William Lloyd
1900 Muir, Sir Thomas
1901 Dyson, Sir Frank Watson
1901 Macdonald, Hector Munro
1904 Joly, Charles Jasper
1905 Campbell, John Edward
1905 Whittaker, Sir Edmund Taylor
1906 Cowell, Philip Herbert
1906 Lyons, Sir Henry George
1907 Young, William Henry
1908 Grace, John Hilton
1909 Barnes, Ernest William
1910 Hardy, Godfrey Harold
1910 Bornet, Jean Baptiste Edouard
1910 Ehrlich, Paul
1910 Volterra, Sir Vito
1910 Weismann, August Friedrich Leopold
1911 Backlund, Jons Oskar
1911 Groth, Paul Heinrich Ritter von
1911 Kayser, Heinrich Johannes Gustav
1911 Bel, Joseph Achille Le
1911 Timiriazeff, Clement Arkadevich
1912 Dixon, Arthur Lee
1913 Fields, John Charles
1914 Bennett, Geoffrey Thomas
1914 Eddington, Sir Arthur Stanley
1915 Evershed, John
1916 Bousfield, William Robert
1918 Ramanujan, Srinivasa Aaiyangar
1919 Arden-Close, Sir Charles Frederick
1919 Watson, George Neville
1920 Lindemann, Frederick Alexander
1924 Rogers, Leonard James
1925 Fowler, Sir Ralph Howard
1928 Macaulay, Francis Sowerby
1930 Spencer Jones, Sir Harold

**Appendix 9**
**Letters to Ronald Ross**

This appendix contains the correspondence with Sir Ronald Ross, 1857 - 1932, held in the archive of the London School of Hygiene and Tropical Medicine (LSHTM). Ronald Ross was best known as the man who discovered the connection between mosquitos and malaria. He was also an amateur mathematician, and indeed in his memoirs[519] Ross claims that, along with literature, mathematics was one of his first loves. A detailed biography of Ross may be found in [Nye et al., 1997].

**Letter 1**
29 July 1905

Dear Major Ross,
I have enjoyed the study of your verb functions. The method seems analogous to certain methods already known in differential equations. It seems to me quite new and worth prolonged study with the particular object of determining the conditions under which the operations are legitimate. I am of [the] opinion that questions of convergency enter and that you should determine what restrictions there must be on operands to render the series convergent. I hope you will persevere and produce an extensive paper which is bound to enlarge our views on algebraic equations.
I congratulate you on having got so far in the development of the subject and must apologise for a somewhat long delay in answering your letter.
Believe me, sincerely yours, Percy A. MacMahon

**Letter 2**
13 October 1908

Have received your interesting note and am making a few inquiries and will write to you in a few days.
P. A. MacMahon

**Letter 3**
18 October 1908

Dear Ross,
The enclosed from H. F. Baker may interest you - you have already already decided to go to *Nature* so I shall look forward to reading your account.
Yrs sincerely, P. A. MacMahon

**Letter 4**
24 April 1916

My dear Ross,
I am much indebted to you for the three numbers of *Science Progress* - I am much pleased with the review of my book and am speculating as to the identity of C.

[519] [Ross, 1923].

I am also glad to seize the opportunity of becoming acquainted with your work on Operative Division. You have really great mathematical intuition and power of developing an idea - you need not stigmatise yourself as a mathematical amateur if there be any stigma upon an amateur. We are all amateurs who do not hold mathematical chairs and i am sure that you are quite able to hold your own although you have in the first place devoted yourself to medicine. I feel strongly that your work should have been noticed. Your notation is somewhat difficult to get accustomed to and I hope that your new important results you will express in language and notation that can be understood by the mathematical multitude. I find the same thing happens to myself - others are so immersed in their own alphabets and in forming words from them that they have not the patience to learn another man's alphabet and the words he forms. The important work in Combinatory Analysis that I did 29 years ago and at which I have been working at intervals ever since has received no notice and was not even mentioned in Netto's *Combinatorik* published 16 years ago - the result has been tho' that I have collared the 'stuff' and others haven't and it was in order to rub this in that I have written my two volumes.
I think that i shall get some useful ideas for my own work out of your operative division.
With kindest regards and thanks, sincerely yrs, Percy A. MacMahon.

**Letter 5**
3 June 1916

I am sorry that my travelling about delayed receipt of your first letter.
I shall be very pleased to read your paper on 'Iteration' and will expect to receive it from you shortly.
Yrs sincerely, P. A. MacMahon.

**Letter 6**
20 June 1916

My dear Ross,
I think that your paper on iteration is really a branch of the Calculus of Functions. Do you know Boole's Finite Differences and his discussion of certain cases of iteration by means of Finite Difference Equations ? His discussion however leads up to periodical functions and not to the solution of algebraic equations. The subject is not quite my line and I am sure that I could not be one of the referees selected for the paper as many men can give a better opinion. The paper is very interesting and I would be glad to see it published. It is impossible to predict what the referees report will be. I think however that you should connect the subject up with the Calculus of Functions, a well known branch of mathematics, and Boole's book.
I return your paper.
Sincerely yrs, Percy A. MacMahon.

**Letter 7**
30 April 1917

My dear Ross,
Excuse the delay in replying to your letter as I have to think over the reply. The fact is that the London Mathematical Society has more papers than it will be able to print for a long time, so that a paper communicated now and passed by the referees would have some time to wait for publication.

But having said this I am glad to add that I will willingly communicate it if you so desire - I think your operational work is good.
Yrs sincerely, Percy A. MacMahon.

**Letter 8**
20 July 1921

Dear Ross,
Thanks much for *Science Progress* just received. Mathews is generous about my work - I shall try to think that I deserve the nice things he says. Your work upon operations is very interesting. I am enjoying the reading of it.
Yrs sincerely, Percy A. MacMahon.

## Appendix 10

## The Master Theorem

In Chapter 4, the genesis of MacMahon's so-called 'Master Theorem' was described, and its use in calculating the coefficients of enumerating generating functions was alluded to. In this Appendix, an example of the use of the theorem, to calculate the solution to the Problème des Rencontres in the case $n = 4$, is given.

The problem in this instance is to find the number of rearrangements of four symbols, say *1234*, so that no symbol is in its original position; e.g. *3142*.

The solution can easily be found by multiplying out the expression $(x_2+x_3+x_4)(x_1+x_3+x_4)(x_1+x_2+x_4)(x_1+x_2+x_3)$ and noting the coefficient of the term $x_1x_2x_3x_4$. This is the *redundant generating function*, since it provides more information than is required by the problem, and for larger values of $n$ the algebra is very hard to do accurately by hand. The Master Theorem says that when we need to find a coefficient in the expansion of $\frac{1}{(1-X_1)(1-X_2)(1-X_3)(1-X_4)}$, where the $X_n$ are linear functions of the form $a_{11}x_1 + a_{12}x_2 + a_{13}x_3 + a_{14}x_4$, and the coefficient required is of the term $x_1x_2x_3x_4$, then it can be found in the expansion of $V_4^{-1}$, where $V_4$ is the determinant $(-1)^4 x_1x_2x_3x_4 \begin{vmatrix} a_{11}-\frac{1}{x_1} & a_{12} & a_{13} & a_{14} \\ a_{21} & a_{22}-\frac{1}{x_2} & a_{23} & a_{24} \\ a_{31} & a_{32} & a_{33}-\frac{1}{x_n} & a_{34} \\ a_{41} & a_{42} & a_{43} & a_{44}-\frac{1}{x_4} \end{vmatrix}$.

In this example, we need all the $a_{nn} = 0$, so we have to calculate $(-1)^4 x_1x_2x_3x_4 \begin{vmatrix} -\frac{1}{x_1} & a_{12} & a_{13} & a_{14} \\ a_{21} & -\frac{1}{x_2} & a_{23} & a_{24} \\ a_{31} & a_{32} & -\frac{1}{x_n} & a_{34} \\ a_{41} & a_{42} & a_{43} & -\frac{1}{x_4} \end{vmatrix}$, which may be written $1 - p_2 - 2p_3 - 3p_4$, where $p_2 = x_1x_2 + x_1x_3 + x_1x_4 + x_2x_3 + x_2x_4 + x_3x_4 = \sum x_1x_2 = (1^2)$, $p_3 = \sum x_1x_2x_3 = (1^3)$ and $p_4 = \sum x_1x_2x_3x_4 = (1^4)$. The *condensed* generating function is thus $\frac{1}{1 - p_2 - 2p_3 - 3p_4}$.

MacMahon then had to deal with the problem of expanding this fraction, which he did by discovering a recurrence formula for the coefficients $P_n$ in the expansion: $\frac{1}{1 - p_2 - 2p_3 - 3p_4} = 1 + P_1p_1 + P_2p_2 + P_3p_3 + P_4p_4 + \cdots$.

The formula is $P_n = \binom{n}{2}P_{n-2} + 2\binom{n}{3}P_{n-3} + 3\binom{n}{4}P_{n-4} + \cdots (n-1)\binom{n}{n}$, which in this example becomes $P_4 = \binom{4}{2}P_2 + 2\binom{4}{3}P_1 + 3$. Here $P_1$ is number of ways of permuting 1 symbol so no symbol is in its original place, so clearly 0, and $P_2$ is the number of ways of permuting 2 symbols, so 1. Hence $P_4$ = 6 + 3 = 9.

MacMahon noted that this method for determining the coefficients was limited in its applicability, and went on to develop more general methods in *Combinatory Analysis*, Volume 1, Section 3, Chapter 3.

MacMahon gave an alternative method for calculating the value of $P_4$, using the redundant generating function. Take $p_1 = x_1 + x_2 + x_3 + x_4 = \sum x_1 = (1)$, then the generating function can be written

$(p_1 - x_1)(p_1 - x_2)(p_1 - x_3)(p_1 - x_4)$. Expanded, this becomes

$$p_1^4 - p_1^3\sum x_1 + p_1^2\sum x_1x_2 - p_1\sum x_1x_2x_3 + \sum x_1x_2x_3x_4 = p_1^4 - p_1^3 p_1 + p_1^2 p_2 - p_1 p_3 + p_4$$
$$= p_1^2 p_2 - p_1 p_3 + p_4 .$$

In general, in the expansion of $p_1^{n-s} p_s$, the coefficient of $p_n$ is $\frac{n!}{s!}$. So in the expression above the value is $\frac{4!}{2!} - \frac{4!}{3!} + \frac{4!}{4!} = 12 - 4 + 1 = 9$, as before.

MacMahon also used the notation $\{0;1^n\}$ to represent the coefficient that enumerates the number of ways that $n$ different symbols may be permuted so that none are in their original positions. This can then be extended to $\{m;\xi_1\xi_2\ldots\xi_n\}$ to represent the number of ways in which exactly $m$ of the $\xi_1 + \xi_2 + \ldots + \xi_n$ symbols in $x_1^{\xi_1}x_2^{\xi_2}\ldots x_n^{\xi_n}$ remain in their original places.

# Bibliographies

## Manuscript sources

This is a list of the manuscript sources consulted during the course of research for this thesis. Where appropriate they are also listed in the main bibliography.

**The National Archive** (formerly **The Public Record Office**), Kew, London.
War Office papers W032/6157, WO33/42
Board of Trade papers BT101/626 to 845

**St John's College, Cambridge**
Papers of Joseph Larmor: eleven letters (details at http://janus.lib.cam.ac.uk)

**University College, London**
Karl Pearson archive: three letters reference 756/7 described in Chapter 4.

**London School of Hygiene and Tropical Medicine**
Ronald Ross papers: eight letters (not indexed as at March 2006).

**St Andrews University, Fife**
D.W. Thompson papers: nine letters as described in Appendix 6.

**Glasgow City Library, Glasgow**
Parker Smith archive: two letters and one document as described in Appendix 7.

**Royal Society, London**
MacMahon correspondence: eight letters as described in Appendix 8.

**The Bodleian Library, Oxford**
British Association for the Advancement of Science archive.

**The Sandhurst Collection**
The Royal Military Academy Sandhurst, Camberley, Surrey.

**Chronological list of MacMahon's papers**

Below is a bibliography of the Mathematical Papers of Percy Alexander MacMahon adapted from George Andrews' *Percy Alexander MacMahon Collected Papers,* grouped according to the chapter headings in this thesis. Papers marked with an asterisk (*) are also referenced in the main bibliography, so there is a degree of repetition.

**<u>Chapter 4: 1890 - 1898</u>**

**<u>Chapter 5: 1899 - 1906</u>**

**Chapter 6: 1907- 1922**

## Chapter 7: 1923 - 1929

**Main Bibliography**